\documentclass[12 pt]{amsart}
\title{Twistor Spaces for Supersingular K3 Surfaces}
\author{Daniel Bragg and Max Lieblich}
\usepackage{macros}

\begin{document}

\maketitle

\begin{abstract}
  We develop a theory of twistor spaces for supersingular K3 surfaces,
  extending the analogy between supersingular K3 surfaces and complex
  analytic K3 surfaces. Our twistor spaces are obtained as relative
  moduli spaces of twisted sheaves on universal gerbes associated to
  the Brauer groups of supersingular K3 surfaces.  In rank 0,
  this is a geometric incarnation of the Artin-Tate isomorphism.
  Twistor spaces give rise to curves in moduli spaces of twisted
  supersingular K3 surfaces, analogous to the analytic moduli space of
  marked K3 surfaces. We describe a theory of crystals for twisted
  supersingular K3 surfaces and a twisted period morphism from the
  moduli space of twisted supersingular K3 surfaces to this space of
  crystals.
  
  As applications of this theory, we give a new proof of the Ogus-Torelli
  theorem modeled on Verbitsky's proof in the complex analytic setting and a
  new proof of the result of Rudakov-Shafarevich that supersingular K3 surfaces
  have potentially good reduction. These proofs work in
  characteristic 3, filling in the last remaining gaps in the theory. As a further
  application, we show that each component of the supersingular locus in each 
  moduli space of polarized K3 surfaces is unirational.
\end{abstract}

\tableofcontents

\section{Introduction}
A K3 surface over an algebraically closed field is said to be
\textit{supersingular} if the rank of its Picard group is 22. This is
the maximal value possible among K3 surfaces, and only occurs 
in positive characteristic. The first examples of such surfaces were
found by Tate 
\cite{MR0225778} and Shioda \cite{MR0429918, MR0435084}. The
systematic study of supersingular K3 surfaces was initiated by Artin
\cite{Artin74}.\footnote{Artin noticed that Picard rank 22 has a
  cohomological implication, namely, that the formal Brauer group
  $\widehat{\Br}(X)$ is isomorphic to the formal additive group
  $\widehat{\G}_a$. He was able to show that any elliptic K3
  satisfying this condition must also have Picard number 22, and
  conjecture that this should be true in general. The Tate conjecture
  for supersingular K3 surfaces is equivalent to Artin's
  conjecture. Since this has been proven for K3 surfaces in all
  characteristics \cite{Charles13,Charles14,KMP16, MP15,Mau14}, we
  will not make a distinction between the cohomological condition
  (unipotence of the formal Brauer group) and the algebraic one
  (maximality of the Picard number).} Among other things, he described
a series of analogies between supersingular K3 surfaces over a fixed
algebraically closed field of positive characteristic and analytic K3
surfaces over the complex numbers.

In this paper we add to this series of analogies by
developing a positive characteristic analog of \emph{twistor families\/}. 
Recall that in the complex setting, a twistor family associated to a K3 surface
$X$ is given by fixing a Kähler class $\alpha\in\H^2(X,\R)$ and varying the
complex structure of $X$ using a sphere spanned by the hyper-Kähler 
structure $I, J, K\in\H^2(X,\R)$ associated to 
$\alpha$. There results a non-algebraic family $X(\alpha)\to\P^1$ of analytic K3
surfaces. If we take the 3-plane $W$ in $\H^2(X,\R)$ spanned by $I, J,
K$, we can describe the periods of the line $X(\alpha)$ by intersecting the
complexification $W_{\C}$ with the period domain $D$, viewed as a subvariety of
$\P(\H^2(X,\C))$ (see, for example, Section 6.1 and Section 7.3 of \cite{Huy16}).

In place of the Kähler class in $\H^2(X,\R)$, our supersingular
twistor spaces are determined by a choice of isotropic vector $v$ in
the quadratic $\F_p$-space $p\NS(X)^\ast/p\NS(X)$. In
place of the sphere of complex structures in $\H^2(X,\R)$, we use
continuous families of classes in the cohomology group
$\H^2(X,\m_p)$. In place of a fixed differentiable surface carrying a
varying complex structure, we end up with a fixed algebraic surface
carrying a varying Brauer class (which is an algebraic avatar of
transcendental structure). The group $\H^2(X,\m_p)$ may be naturally
viewed as a subvariety of $\H^2_{dR}(X/k)$, and just as in the
classical setting, this family of cohomology classes traces out a line
in an appropriate period domain classifying supersingular K3 crystals.

Here is a more precise description of supersingular
twistor families. If $X$ is a supersingular K3 surface, then by a
result of Artin 
(proved here as Theorem \ref{thm:relativeartin}) the flat cohomology
group $\H^2(X,\m_p)$ can be naturally viewed as the $k$-points of a
group scheme whose connected component is isomorphic to
$\bA^1$. This 
results in a \emph{universal twistor family\/} $\sX\to\bA^1$ of
twisted supersingular 
K3 surfaces, whose central fiber is the trivial $\m_p$-gerbe over
$X$. Although the underlying K3 surface of each fiber is constant, the
gerbe structure is varying. This situation is only
possible in positive characteristic: if $k$ has characteristic zero,
then the groups $\H^2(X,\m_n)$ 
are all finite, and so there do not exist such continuously varying
families of $\m_n$-gerbes. By taking relative moduli spaces of twisted
sheaves on $\sX\to\bA^1$ with suitably chosen topological invariants,
we can produce new families $\sX(v)\to\bA^1$ of twisted supersingular
K3 surfaces. We call such families \textit{(supersingular) twistor
  families}. (Note that the original universal family of gerbes $\sX\to\bA^1$
is the stack of $\sX$-twisted sheaves of length $1$, so it is itself a
family of moduli spaces of $\sX$-twisted sheaves.) While the
universal twistor family has constant coarse moduli space, the coarse
moduli family of a general twistor family will typically be non-isotrivial.

Certain special cases of our twistor families have been studied before
by various authors, and their connection to the Brauer group was
noticed already by Artin (see the footnote on page 552 of
\cite{Artin74}). The novelty in our approach lies in realizing that
they are most naturally viewed as families of $\m_p$-gerbes over
supersingular K3 surfaces (that is, as families of twisted K3
surfaces). To develop this viewpoint, we take as our basic object of
study a $\m_p$-gerbe $\sX\to X$ over a supersingular K3 surface, and
seek to extend the material of \cite{Ogus78,Ogus83} to the twisted
setting. We introduce an appropriate moduli space of such objects,
define their crystalline periods, and study the corresponding period
morphism. By working throughout with twisted K3 surfaces, we are able
to treat the universal twistor families $\sX\to\bA^1$ on the same
footing as any other family of supersingular K3 surfaces. This allows
us to obtain precise control over the periods of twistor families.

The moduli space of supersingular K3 surfaces is stratified by the
\textit{Artin invariant}, which takes integral values
$1\leq\sigma_0\leq 10$. The dimension of the locus of Artin
invariant at most $\sigma_0$ is $\sigma_0-1$. In particular, the moduli space
of supersingular K3 surfaces has dimension 9. We extend the Artin
invariant to a twisted supersingular K3 surface, and find that
$\sigma_0$ now goes up to 11. Accordingly, the moduli space of twisted
supersingular K3 surfaces has dimension 10, and the moduli space of
non-twisted surfaces is naturally embedded in it as a divisor. This
bears some resemblance to the situation over the complex numbers,
where loci of algebraic K3 surfaces are embedded in the space of
analytic K3 surfaces as divisors. In some sense, there is only one
such locus for twisted K3 surfaces because all of the cohomology is
already algebraic up to $p$-torsion.

Here is a rough dictionary translating between the classical analytic
situation and the supersingular one, extending Artin's Table 4.12
\cite{Artin74}. For a K3 surface over the complex numbers, we let
$\rho$ be the Picard number and $\rho_0$ be the rank of the
transcendental lattice. For a supersingular K3 surface with Artin
invariant $\sigma_0$ we will write $\sigma=11-\sigma_0$.  \vskip
1\baselineskip

\newcommand\Tstrut{\rule{0pt}{2.6ex}}         
\newcommand\Bstrut{\rule[-0.9ex]{0pt}{0pt}}   

\begin{center}
\begin{tabular}{ c|c } 
  \bf{Complex} & \bf{Supersingular}\Bstrut\\
  \hline
  Complex structure, class in $\H^2(X,\R)$ & Gerbe, class in $\H^2(X,\m_p)$\Tstrut \\ 
  Twistor line $\bP^1$                     & Connected component $\bA^1\subset\R^2\pi_\ast\m_p$ \\
  Algebraic moduli space                   & Moduli space of supersingular K3 surfaces \\
  Analytic moduli space                    & Moduli space of twisted supersingular K3 surfaces \\
  Period domain                            & Space of supersingular K3 crystals \\
  $\rho_0\geq 2$                           & $\sigma_0\geq 1$   \\ 
  $\rho+\rho_0=22$                         & $\sigma+\sigma_0=11$  \\
  Generic analytic K3 has $\rho=0$         & Generic twisted supersingular K3 has $\sigma=0$
\end{tabular}
\end{center}

\vskip 1\baselineskip

The most interesting aspects of this dictionary are the places where it breaks
down. Notably, while a twistor line gives a $\P^1$ in the analytic moduli space, our
construction only yields an $\A^1$ in the moduli space of twisted supersingular
K3 surfaces. We know in retrospect from the crystalline Torelli theorem that the underlying family of surfaces of a twistor family extends to a family over $\bP^1$ (after possibly normalizing the base). However, the theory developed in this paper only yields a direct description of such a family over the open locus $\bA^1\subset\bP^1$. It would be very interesting to find a natural interpretation of the fiber over infinity. While the complex analytic twistor spaces are decidedly non-algebraic, ours are
sufficiently algebraic to give rational curves in moduli spaces of algebraic
(supersingular) K3 surfaces.

It is our hope that similar constructions will be useful in the study of other
classes of varieties in positive characteristic. In particular, it seems
reasonable that some the constructions in this paper might extend to
supersingular holomorphic symplectic varieties. We will explore this in future
work.

\subsection{A brief outline of this paper}
\label{sec:outline}

In Section \ref{sec:cohomology} we discuss some preliminaries on the
cohomology of supersingular K3 surfaces, focusing on the cohomology of
$\m_p$ (Section \ref{sec:flat-coho-mu-p}) and the Brauer group
(Section \ref{sec:BrauerGroup}). We introduce the relative étale site
(Section \ref{sec:relativeetalesite}), which is a technical tool that
will be useful for calculations later in the paper.

The goal of Section \ref{sec:periods} is to define and study the
periods of twisted supersingular K3 surfaces, building on work of Ogus
\cite{Ogus78, Ogus83}. In Section \ref{sec:perioddomain} we recall
Ogus's moduli space of characteristic subspaces, which is a
crystalline analog of the classical period domain, and study the
relationship between the spaces classifying K3 crystals of different
Artin invariants. Twistor lines in this period domain are defined in
Section \ref{sec:twistorlines}. In Section \ref{sec:twistedk3crystals}
we develop a crystalline analog of the classical B-fields associated
to a Brauer class, which we use to construct a K3 crystal associated
to a twisted supersingular K3 surface. In Section
\ref{sec:twistedperiodmorphism} we study the moduli space of marked
twisted supersingular K3 surfaces. Using the twisted K3 crystals of
the previous section, we define a period morphism on this space, and
compute its differential.

In Section \ref{sec:moduli} we construct supersingular twistor
spaces. We begin by discussing moduli spaces of twisted sheaves, which
we then use to construct our twistor families. These naturally fall
into two types, which we describe separately. The main result of
Section \ref{sec:constructionoftwistorlines} is a construction of
twistor families of positive rank in the moduli space of twisted
supersingular K3 surfaces. In Section \ref{sec:artintate}, we
construct Artin-Tate twistor families, which is somewhat more involved
than the positive rank case. We show how in the presence of an
elliptic fibration the twistor line construction yields a geometric
form of the Artin-Tate isomophism (Proposition \ref{prop:Br-iso-Sha}).

In Section \ref{sec:applications} we give applications of
theory. Section \ref{sec:torellitheorem} is devoted to a supersingular
form of Verbitsky's proof of the Torelli theorem: the existence of the
twistor lines gives a relatively easy proof that the twisted period
morphism is an isomorphism (after adding ample cones), which is a
twisted version of Ogus's crystalline Torelli theorem. This implies
immediately that Ogus's period morphism is an isomorphism, giving a
new proof of the main result of \cite{Ogus83}. Our proof works
uniformly in all odd characteristics, and so we are able to extend
Ogus's result to characteristic $p=3$, where it was not previously
known. In Section \ref{sec:misc-applications}, we deduce some
consequences for the moduli of supersingular K3 surfaces. As a
corollary of Theorem \ref{thm:twistedtorelli}, we obtain a new proof
of a result of Rudakov and Shafarevich on degenerations of
supersingular K3 surfaces (Theorem 3 of \cite{RTS83}). This result is
also new in characteristic 3. We also prove that the supersingular
locus in each moduli spaces of polarized K3 surfaces is
unirational.

\begin{Remark} We fix throughout a prime number $p$. Unless noted
  otherwise, we will assume that $p$ is odd. The main reason for this
  restriction will be our use of bilinear forms. We think it is likely
  that our techniques should yield similar results in characteristic
  $p=2$ as well, with suitable modifications to certain definitions,
  but we will not pursue this here.
\end{Remark}

\subsection{Acknowledgments}

François Charles first pointed out to the second-named author that the
curves he found in \cite{Lie15} should be thought of as positive
characteristic twistor families.  Much of the work described here on
twistor lines, the twisted period morphism, and the new proofs of the
Ogus-Torelli theorem and the Rudakov-Shafarevich theorem was done by
the first-named author as part of his graduate work. The present paper
also includes improvements to some of the arguments of \cite{Lie15}, many due to the first author, and therefore supersedes this
preprint. We thank Josh Swanson for the proof of Lemma
\ref{BigTangentSpace}. We had helpful conversations with Bhargav Bhatt, François
Charles, Brendan Hassett, Daniel Huybrechts, Aise Johan de Jong,
Sándor Kovács, David Benjamin Lim, Keerthi Madapusi Pera, Davesh Maulik, and
Martin Olsson during the 
course of this work. We received comments from Bhargav Bhatt, François Charles,
and Christian Liedtke on an earlier version of this manuscript, as well as especially detailed comments from Yuya Matsumoto.

\section{Cohomology}\label{sec:cohomology}
\subsection{Flat cohomology of $\mu_p$}\label{sec:flat-coho-mu-p}
Let $X$ be a scheme. Fix a prime number $p$, and consider the group scheme $\m_p$ over $X$. This group is related to the sheaf of units $\cO_X^\times=\bG_{m,X}$ by the Kummer sequence
\[
  1\to\m_p\to\cO_X^\times\xrightarrow{x\mapsto x^p}\cO_X^\times\to 1
\]
which is exact in the flat topology (if $X$ has characteristic not equal to $p$, then it is also exact in the \'{e}tale topology). By general results of Giraud \cite{MR0344253}, if $A$ is any sheaf of abelian groups on $X$, then $\H^2(X,A)$ is in bijection with isomorphism classes of $A$-gerbes on $X$. We will denote the cohomology class of an $A$-gerbe $\sX\to X$ by $[\sX]\in\H^2(X,A)$. We begin this section by collecting some results related to gerbes banded by the sheaves $\m_p$ and $\cO_X^\times$. We refer to Section 2 of \cite{Lie04} for background, and for the definition of twisted and $n$-twisted sheaves. If $\sX\to X$ is a $\m_p$ or $\cO_X^\times$-gerbe and $n$ is an integer, we will write $\Qcoh^{(n)}(\sX)$ for the category of $n$-twisted quasicoherent sheaves on $\sX$, and $\Coh^{(n)}(\sX)$ for the category of $n$-twisted coherent sheaves on $\sX$. 

Consider the long exact sequence
\[
  \ldots\to \H^1(X,\cO_X^\times)\too{\cdot p}\H^1(X,\cO_X^\times)\too{\delta} \H^2(X_{\fl},\m_p)\to \H^2(X,\cO_X^\times)\too{\cdot p} \H^2(X,\cO_X^\times)\to\ldots
\]
induced by the Kummer sequence.
\begin{Definition}
  If $\sL\in\Pic (X)$ is a line bundle, then the \textit{gerbe of p-th
    roots} of $\sL$ is the stack $\{\sL^{1/p}\}$ over $X$ whose objects
  over an $X$-scheme $T\to X$ are line bundles $\sM$ on $T$ equipped
  with an isomorphism $\sM^{\otimes p}\iso\sL_T$.
\end{Definition}
The stack $\{\sL^{1/p}\}$ is canonically banded by $\m_p$, and $\{\sL^{1/p}\}\to X$ is a $\m_p$-gerbe. Its cohomology class is the image of $\sL$ under
the boundary map $\delta$.
\begin{Lemma}\label{lem:firstlemma}
  Let $\sX\to X$ be a $\m_p$-gerbe with cohomology class
  $\alpha\in \H^2(X_{\fl},\m_p)$. The following are equivalent.
  \begin{enumerate}
    \item There exists an invertible twisted sheaf on $\sX$,
    \item there is an isomorphism $\sX\cong\{\sL^{1/p}\}$ of
      $\m_p$-gerbes for some line bundle $\sL$ on $X$,
    \item the cohomology class $\alpha$ is in the image of the
      boundary map $\delta$, and
    \item the associated $\cO_X^\times$-gerbe is trivial.
  \end{enumerate}
\end{Lemma}
\begin{proof}
  We will show $(1)\Leftrightarrow(2)$. Suppose that there exists an
  invertible twisted sheaf $\sM$ on $\sX$. The $\m_p$-action on
  $\sM^{\otimes p}$ is trivial, so the adjunction map
  $p^*p_*\sM^{\otimes p}\iso\sM^{\otimes p}$ is an isomorphism. We
  obtain by descent an isomorphism
  $\sX\cong \{(p_*\sM^{\otimes p})^{1/p}\}$ of
  $\m_p$-gerbes. Conversely, on $\{\sL^{1/p}\}$ there is a universal
  invertible twisted sheaf $\sM$. The image of a line bundle $\sL$
  under the boundary map $\delta$ in non-abelian cohomology is
  $\{\sL^{1/p}\}$, so $(2)\Leftrightarrow(3)$. Finally,
  $(3)\Leftrightarrow (4)$ by the exactness of the long exact sequence
  on cohomology.
\end{proof}
We say that a $\m_p$-gerbe satisfying the conditions of Lemma
\ref{lem:firstlemma} is \textit{essentially trivial}. We will also
need the following slightly more general construction. Let $\sX\to X$
be a $\m_p$-gerbe and $\sL\in\Pic(X)$ a line bundle.
\begin{Definition}\label{def:relativeroot}
  Let $\sX\{\sL^{1/p}\}$ be the stack on $X$ whose objects over an
  $X$-scheme $T\to X$ are invertible sheaves $\sM$ on $\sX_T$ that are $(-1)$-twisted, together with an isomorphism
  $\sM^{\otimes p}\iso\sL_{\sX_T}$.
\end{Definition}
The stack $\sX\{\sL^{1/p}\}\to X$ is canonically banded by $\m_p$,
and is again a $\m_p$-gerbe over $X$. In the special case where $\sX=\mathbf{B}\m_p$ is the trivial
gerbe, there is a canonical isomorphism
$\sX\{\sL^{1/p}\}\cong \sL^{1/p}$ of $\m_p$-gerbes over $X$. Furthermore, there is a universal invertible twisted sheaf
\[
  \sM\in\Pic^{(\text{-}1,1)}(\sX\times_X\sX\{\sL^{1/p}\}).
\]

\begin{Lemma}
  The cohomology class $[\sX\{\sL^{1/p}\}]\in\H^2(X,\m_p)$ is equal to
  $[\sX]+[\{\sL^{1/p}\}]$. 
\end{Lemma}
\begin{proof}
  The group structure on $\H^2(X,\m_p)$ is induced by the contracted
  product of $\m_p$-gerbes along the anti-diagonal
  $\m_p\to\m_p\times\m_p$. That is, given two $\m_p$-gerbes $\ms A\to
  X$ and $\ms B\to X$, a $\m_p$-gerbe $\ms C\to X$ represents the
  class $[\ms A] + [\ms B]$ if and only if there is a morphism of
  $X$-stacks $\ms A\times_X\ms B\to\ms C$ such that the induced
  morphism of bands $\m_p\times\m_p\to\m_p$ is the
  multiplication map. (The reader is referred to Section IV.3.4 and
  the associated internal references of \cite{MR0344253} for details.)

  Consider the $\m_p$-gerbe $\sY\to X$ whose sections over $T\to X$
  are given by the groupoid of pairs $(M,\tau)$ where $M$ is an
  invertible $\sX_T$-twisted sheaf and $\tau:M^{\tensor p}\simto\ms O_{\sX_T}$
  is a trivialization; these are ``$\sX_T$-twisted
  $\m_p$-torsors''. Such a pair $(M,\tau)$ corresponds precisely to a morphism
  $\sX_T\to\B\m_{p,T}$ that induces the identity map 
  $\m_p\to\m_p$ on bands. There results a morphism $\sX\times\sY\to\B\m_p$ that
  induces the multiplication map $\m_p\times\m_p\to\m_p$ on
  bands. Thus, $\sY$ represents $-[\sX]\in\H^2(X,\m_p)$. Similarly,
  $\sX$ is isomorphic to the stack of $\sY$-twisted $\m_p$-torsors. Since
  $(-1)$-twisted sheaves on $\sX$ correspond to $\sY$-twisted sheaves,
  we conclude that $\sX$ is itself isomorphic to the $X$-stack of pairs
  $(N,\gamma)$ with $N$ an invertible $(-1)$-twisted sheaf on $\sX$
  and $\gamma:N^{\tensor p}\simto\ms O$ a trivialization of the $p$th
  tensor power. We will use this identification in what follows.

  Now consider the morphism $\sX\times_X\{\sL^{1/p}\}\to\sX\{\sL^{1/p}\}$ defined
  as follows. As above, an object of $\sX$ over $T\to X$ can be
  thought of as a pair $(\sN,\gamma)$ with $\sN$ a $(-1)$-twisted
  invertible sheaf on $\sX_T$ and $\gamma:\sN^{\tensor
    p}\simto\ms O_{\sX_T}$. An object of $\{\sL^{1/p}\}$ over $T$ is a pair
  $(\sM,\chi)$ with $\sM$ an invertible sheaf on $T$ and
  $\chi:\sM^{\tensor p}\simto\sL$ an isomorphism. We get an object of
  $\sX\{\sL^{1/p}\}$ by sending $(\sN,\gamma),(\sM,\chi)$ to
  $(\sN\tensor\sM, \gamma\tensor\chi)$. This induces the multiplication
  map on bands, giving the desired result.
\end{proof}

Let $\pi\colon X\to S$ be a morphism. Consider the group scheme $\m_p$
as a sheaf of abelian groups on the big fppf site $X_{\fl}$ of
$X$. The morphism $\pi$ induces a map
$\pi^{\fl}\colon X_{\fl}\to S_{\fl}$ of sites. We will study the sheaf
$\R^2\pi^{\fl}_*\m_p$ on $S_{\fl}$. This sheaf is of formation
compatible with arbitrary base change, in the sense that for any
$S$-scheme $T\to S$ there is a canonical isomorphism
\[
  \R^2\pi^{\fl}_*\m_p(T)=\H^0(T,\R^2\pi_{T*}^{\fl}\m_p)
\]
The sheaf $\R^2\pi^{\fl}_*\m_p$ may be described as the sheafification in the flat topology of the functor 
\[
  T\mapsto \H^2(X_{T\fl},\m_p)
\]
In this paper, we will only consider the cohomology of $\m_p$ with respect to the flat topology, so we will when convenient omit the subscript indicating the topology.

\begin{Lemma}\label{lem:mup1}
  Suppose that $\pi\colon X\to S$ is a morphism such that $\pi^{\fl}_*\m_p=\m_p$ and $\R^1\pi^{\fl}_*\m_p=0$ (for instance, a relative K3 surface). For any $T\to S$, there is an exact sequence
  \[
    0\to \H^2(T,\m_p)\to \H^2(X_T,\m_p)\to \H^0(T,\R^2\pi_{T*}^{\fl}\m_p)\to \H^3(T,\m_p)
  \]
\end{Lemma}
\begin{proof}
  This follows from the Leray spectral sequence for the sheaf $\m_p$ on $X_T\to T$.
\end{proof}
In particular, if $S=\Spec k$ is the spectrum of an algebraically closed field, then there is a canonical identification
\[
  \R^2\pi^{\fl}_*\m_p(k)=\H^2(X_{\fl},\m_p)
\]

\begin{Theorem}\label{thm:relativeartin}
  If $\pi\colon X\to S$ is a relative K3 surface, then the sheaf $\R^2\pi^{\fl}_*\m_p$ on the big flat site of $S$ is representable by a group algebraic space locally of finite presentation over $S$.
\end{Theorem}

This is a relative form of a result of Artin, published as Theorem 3.1 of \cite{Artin74} with the caveat that the proof would appear elsewhere. We have not found the proof in the literature, so we provide a moduli-theoretic proof here. Note that we make no assumptions as to the characteristic of $S$, although we will be most interested in the case when $S$ has characteristic $p$. We will write $\sS=\R^2\pi^{\fl}_*\m_p$ in the remainder of this section.

\begin{Proposition}
  \label{P:diag}
  The diagonal
  $$\sS\to\sS\times\sS$$
  is representable by closed immersions of finite presentation.
\end{Proposition}
\begin{proof}
  Let $a,b\colon T\to\sS$ be two maps corresponding to classes
  $$a,b\in\H^0(T,\R^2\pi^{\fl}_{T*}\m_p).$$
  To prove that the locus where $a=b$ is represented by a closed
  subscheme of $T$, it suffices by translation to assume that $b=0$
  and prove that the functor
  $Z(a)$ sending a $T$-scheme $T'\to T$ to $\emptyset$ if $a_{T'}\neq 0$
  and $\{\emptyset\}$ otherwise is representable by a closed subscheme
  $Z_a\subset T$.

  First, suppose $a$ is the image of a class $\alpha\in\H^2(X_T,\m_p)$
  corresponding to a $\m_p$-gerbe $\ms X\to X_T$. Since
  $\R^1\pi^{\fl}_\ast\m_p=0$, we see that the functor $Z(a)$ parametrizes
  schemes $T'\to T$ such that there is an $\alpha'\in\H^2(T',\m_p)$
  with
  \[
    \alpha'|_{X_{T'}}=\alpha|_{X_{T'}}.
  \]

  Let $\ms P\to T$ be the stack whose objects over $T'\to T$ are
  families of $\ms X_{T'}$-twisted
  invertible sheaves $\ms L$ together with isomorphisms $\ms
  L^{\tensor p}\simto\ms O_{\ms X_{T'}}$. The stack $\ms P$ is a $\m_p$-gerbe over a
  quasi-separated algebraic space $P\to T$ that is locally of finite
  presentation (Proposition 2.3.1.1 of \cite{Lie04} and Section C.23 of \cite{MR2786662}). Moreover, since $\Pic_X$ is torsion free, the
  natural map $P\to T$ is a monomorphism.
  Note that if we change $\alpha$ by the preimage
  of a class $\alpha'\in\H^2(T,\m_p)$ we do not change $P$ 
  (but we do change the class of the gerbe $\ms P\to P$ by $\alpha'$).

  \begin{Lemma}
    The algebraic space $P\to T$ is a closed immersion of finite presentation.
  \end{Lemma}
  \begin{proof}
    First, let us show that $P$ is of finite presentation. It
    suffices to show that $P$ is quasi-compact under the assumption
    that $T$ is affine. Moreover, since $\sS$ is locally of finite
    presentation, we may assume that $T$ is Noetherian. By Gabber's
    Theorem (Theorem 1.1 of \cite{dejong-gabber}) there is a Brauer-Severi scheme $V\to X$ such that
    $\alpha|_V$ has trivial Brauer class, i.e., so that there is an
    invertible sheaf $L\in\Pic(V)$ satisfying
    $$\alpha|_V=\delta(\sL)\in\H^2(V,\m_p),$$
    Writing $$\ms V=\ms X\times_{X_T} V,$$ we know
    from the isomorphism $\ms V\cong \{L^{1/p}\}$ that there is an
    invertible $\ms V$-twisted sheaf $\ms L$ such that
    $$\ms L^{\tensor p}\cong L.$$

    Let $W$ denote the algebraic space parametrizing invertible $\ms
    V$-twisted sheaves whose $p$-th tensor powers are trivial. By
    the argument in the preceding paragraph, tensoring with $\ms
    L^{\vee}$ defines an isomorphism between $W$ and the fiber of the
    $p$-th power map
    $$\Pic_{V/T}\to\Pic_{V/T}$$
    over $[L^\vee]$. Since the $p$-th power map is a closed immersion
    ($X$ being K3), we
    see that $W$ is of finite type. The following lemma then applies
    to show that $P$ is of finite type.

    \begin{Lemma}
      The pullback map $\Pic^{(1)}_{\ms X/T}\to\Pic^{(1)}_{\ms V/T}$
      is of finite type.
    \end{Lemma}
    \begin{proof}
      It suffices to prove the corresponding results for the stacks of
      invertible twisted sheaves. Since $V\times_X V$ and $V\times_X
      V\times_X V$ are proper over $T$ this follows from descent
      theory: the category of invertible $\ms X$-twisted sheaves is
      equivalent to the category of invertibe $\ms V$-twisted sheaves
      with a descent datum on $\ms V\times_{\ms X}\ms V$. Thus, the
      fiber over $L$ on $\ms V$ is a locally closed subspace of
      $$\Hom_{\ms V\times_{\ms X}\ms V}(\pr_1^\ast L,\pr_2^\ast L).$$
      Since the latter is of finite type (in fact, a cone in an affine
      bundle), the result follows.
    \end{proof}

    We claim that $P$ is proper over $T$. To see this, we may use the
    fact that it is of finite presentation (and everything is of
    formation compatible with base change) to reduce to the case in
    which $T$ is Noetherian, and then we need only check the valuative
    criterion over DVR's. Thus, suppose $E$ is a DVR with fraction
    field $F$ and $p\colon \spec F\to P$ is a point. Replacing $E$ by a finite
    extension, we may assume that $p$ comes from an invertible $\ms
    X_F$-twisted sheaf $\ms L$. Taking a reflexive extension and using
    the fact that $\ms X_E$ is locally factorial, we see that $\ms L$
    extends to an invertible $\ms X_E$-twisted sheaf $\ms L_E$. Since
    $\Pic_X$ is separated, it follows that $\ms L_E$ induces the
    unique point of $P$ over $E$ inducing $p$.

    Since a proper monomorphism is a closed immersion, we are done.
  \end{proof}

We now claim that the locus $Z(a)\subset T$ is represented by the
closed immersion $P\to T$. Consider the $\m_p$-gerbe $\ms P\to P$. Subtracting the pullback from $\alpha$ yields a class such that the associated Picard stack
$\ms P'\to P$ is trivial, whence there is an invertible twisted sheaf
with trivial $p$-th power. In other words,
$$a|_P=0\in\sS(P).$$
On the other hand, if $a|_{T'}=0$ then up to changing $a$ by the pullback
of a class from $T$, there is an invertible twisted with trivial $p$-th
power. But this says precisely that $T'$ factors through the moduli
space $P$. This completes the proof of Proposition \ref{P:diag}.
\end{proof}

Let $\Az$ be the stack over $S$ whose objects over a scheme $T\to S$ are Azumaya
algebras $\ms A$ of degree $p$ on $X_T$ such that for every geometric
point $t\to T$, we have
$$\ker(\Tr:\H^2(X,\ms A)\to\H^2(X,\ms O))=0.$$
It is well known that $\Az$ is an Artin stack locally of finite type over $S$
(see, for example, Lemma 3.3.1 of \cite{MR2579390}, and note that the trace condition is open). 
There is a morphism of stacks
$$\chi:\Az\to\R^2\pi^{\fl}_\ast\m_p$$
(with the latter viewed as a stack with no non-trivial automorphisms
in fiber categories)
given as follows. Any family $\ms A\in\Az_T$ has a corresponding class
$$[\ms A]\in\H^1(X_T,\PGL_p).$$
The non-Abelian coboundary map yields a class in $\H^2(X_T,\m_p)$
which has a canonical image
$$\chi(\ms A)\in\H^0(T,\R^2\pi^{\fl}_*\m_p).$$

\begin{Proposition}
  \label{P:surj-by-az}
  The morphism $\chi$ described above is representable by smooth Artin stacks.
\end{Proposition}
\begin{proof}
  Since $\Az$ is locally of finite type over $S$ and the diagonal of
  $\sS$ is representable by closed immersions of finite
  presentation, it suffices to show that $\chi$ is formally smooth.
  Suppose $A\to A'$ is a square-zero extension of Noetherian rings and
  consider a diagram of solid arrows
  \[
    \begin{tikzcd}
      \Spec A'\arrow{r}\arrow{d} & \Az\arrow{d}\\
      \Spec A\arrow{r}\arrow[dashed]{ur} & \sS
    \end{tikzcd}
  \]
  We wish to show that we can produce the dashed diagonal
  arrow. Define a stack $\ms T$ on $\Spec A$ whose objects over an
  $A$-scheme $U\to\spec A$ are dashed arrows in the restricted diagram
    \[
      \begin{tikzcd}
      U'\arrow{r}\arrow{d} & \Az\arrow{d}\\
      U\arrow{r}\arrow[dashed]{ur} & \sS
      \end{tikzcd}
    \]
    where $U'=U\tensor_A A'$, and whose morphisms are isomorphisms between the objects of $\Az$ over $U$ restricting to the identity on the restrictions to $U'$.

    \begin{Claim}
      The stack $\ms T\to\Spec A$ is an fppf gerbe with coherent band.
    \end{Claim}
    \begin{proof}
      First, it is clear that $\ms T$ is locally of finite
      presentation. Suppose $U$ is the spectrum of a complete local
      Noetherian ring with algebraically closed residue field. Then
      $\H^2(U,\m_p)=0$ and the section $U\to\ms S$ is equivalent to a
      class $$\alpha\in\H^2(X_U,\m_p).$$ Let $\ms X\to X_U$ be a
      $\m_p$-gerbe representing $\alpha$ and write $\ms X'=\ms
      X\tensor_A A'$ the restriction of $\ms X$ to $U'$. An object of
      $\Az_{U'}$ is then identified with $\send(V)$ where $V$ is a
      locally free $\ms X'$-twisted sheaf of rank $p$ with trivial
      determinant. The obstruction to deforming such a sheaf lies in
      $$\ker(\H^2(X_{U'},\send(V)\tensor I)\to\H^2(X_{U'},\ms O\tensor
      I)),$$
      and deformations are a pseudo-torsor under
      $$\ker(\H^1(X_{U'},\send(V)\tensor I)\to\H^1(X_{U'},\ms O\tensor
      I))=\H^1(X_{U'},\send(V)\tensor I).$$
      Standard arguments starting from the assumption on the geometric
      points of $\Az$ show that the obstruction group is trivial,
      while the band is the coherent sheaf $\R^1\pi^{\fl}_*\send(V)$, as desired.
    \end{proof}
Since any gerbe with coherent band over an affine scheme is neutral,
we conclude that $\ms T$ has a section. In other words, a dashed arrow
exists, as desired.
\end{proof}

\begin{proof}[Proof of Theorem \ref{thm:relativeartin}]
  Using Proposition \ref{P:surj-by-az}, a smooth cover $B\to\Az$ gives
  rise a to a smooth cover $B\to\sS$. Thus, by Proposition
  \ref{P:diag}, $\sS$ is a group algebraic space locally
  of finite presentation over $S$.
\end{proof}

\subsection{The Brauer group of a supersingular K3 surface}\label{sec:BrauerGroup}

Fix an algebraically closed field $k$ of characteristic $p$, and let $X$ be a supersingular K3 surface over $k$ with structure morphism $\pi\colon X\to\Spec k$. In this section we will record some facts about the Brauer group of $X$ and the flat cohomology group $\sS=\R^2\pi^{\fl}_*\m_p$. We recall the following definition.
 
\begin{Definition}\label{def:completion}
  Suppose that $E$ is a sheaf of abelian groups on the big fppf site of $\Spec k$. We define a functor on the category of local Artinian $k$-algebras with residue field $k$ by
  \[
    \widehat{E}(A)=\ker\left(E(\Spec A)\to E(\Spec k)\right)
  \]
  which we call the \textit{completion of $E$ at the identity section}.
\end{Definition}

In particular, we let $\widehat{\Br(X)}$ denote the completion at the identity section of the functor $\R^2\pi^{\Etale}_*\cO_X^\times$. This is called the \textit{formal Brauer group} of $X$. The association $E\mapsto \widehat{E}$ gives a functor from the category of sheaves of abelian groups on $(\Spec k)_{\fl}$ to the category of presheaves on the opposite of the category of local Artinian $k$-algebras. It is immediate that this functor is left exact.

\begin{Lemma}\label{lem:completionisexact}
  Let $0\to H\to G\to K\to 0$ be a short exact sequence of sheaves of abelian groups on $(\Spec k)_{\fl}$. If $H$ is representable by a smooth group scheme over $\Spec k$, then the induced sequence
  \[
    0\to\widehat{H}\to\widehat{G}\to\widehat{K}\to 0
  \]
  of presheaves is exact.
\end{Lemma}
\begin{proof}
  Let $A$ be a local Artinian $k$-algebra with residue field $k$. Consider the diagram
  \[
    \begin{tikzcd}
      0\arrow{r}& H(A)\arrow{d}\arrow{r} & G(A) \arrow{d}\arrow{r} & K(A) \arrow{d}\arrow{r} &0\\
      0\arrow{r}& H(k)\arrow{r} &G(k) \arrow{r} & K(k) \arrow{r}& 0
    \end{tikzcd}
  \]
  We claim that the rows are exact. It will suffice to show that the flat cohomology groups $\H^1(\Spec A_{\fl},H|_{\Spec A})$ and $\H^1(\Spec k_{\fl},H)$ vanish. For the latter, note that by the Nullstellensatz every cover in $\Spec k_{\fl}$ admits a section, and hence the functor of global sections on the category of sheaves of abelian groups is exact. For the former, our assumption that $H$ is smooth implies by a theorem of Grothendieck (Theorem \ref{thm:grothendieck}) that $\H^1(\Spec A_{\fl},H|_{\Spec A})=\H^1(\Spec A_{\etale},H|_{\Spec A})$. As $A$ is a strictly Henselian local ring, the latter group vanishes.
  
  Finally, as $H$ is smooth, the map $H(A)\to H(k)$ is surjective. The result therefore follows from the snake lemma.
\end{proof}

\begin{Lemma}\label{lem:isooncompletions}
  If $\pi\colon X\to\Spec k$ is a supersingular K3 surface, then the natural map $\R^2\pi^{\fl}_*\m_p\to\R^2\pi^{\Etale}_*\cO_X^\times$ induces an isomorphism
  \[
    \widehat{\R^2\pi^{\fl}_*\m_p}\iso\widehat{\Br(X)}
  \]
  on completions at the identity section.
\end{Lemma}
\begin{proof}
  The flat Kummer sequence induces a short exact sequence
  \[
    0\to\Pic_X/p\Pic_X\to\R^2\pi^{\fl}_*\m_p\to\R^2\pi^{\Etale}_*\cO_X^\times[p]\to 0
  \]
  The quotient $\Pic_X/p\Pic_X$ is smooth, so by Lemma \ref{lem:completionisexact} the induced sequence on completions is exact. Moreover, it is discrete, and hence has trivial completion, so obtain an isomorphism
  \[
    \widehat{\R^2\pi^{\fl}_*\m_p}\iso\widehat{\R^2\pi^{\Etale}_*\cO_X^\times[p]}
  \]
  (with no assumptions on the height of $X$). We also have an exact sequence
  \[
    0\to \R^2\pi^{\Etale}_*\cO_X^\times[p]\to \R^2\pi^{\Etale}_*\cO_X^\times\too{\cdot p}\R^2\pi^{\Etale}_*\cO_X^\times
  \]
  Because $X$ is supersingular, multiplication by $p$ induces the zero map on the formal Brauer group. Taking completions, we get an isomorphism
  \[
    \widehat{\R^2\pi^{\Etale}_*\cO_X^\times[p]}\iso\widehat{\Br(X)}
  \]
  This completes the proof.
\end{proof}

\begin{Proposition}\label{prop:relativeartin1}
  If $\pi\colon X\to S$ is a supersingular K3 surface, then $\sS=\R^2\pi^{\fl}_*\m_p$ is a smooth group scheme over $\Spec k$ of dimension 1 with connected component isomorphic to $\G_a$.
\end{Proposition}
\begin{proof}
  By Lemma \ref{lem:isooncompletions}, the completion of $\sS$ at the identity section is isomorphic to $\widehat{\G_a}$, which is formally smooth and $p$-torsion. It follows that $\sS$ is smooth over $k$ with 1-dimensional connected component. The only 1-dimensional connected group scheme over $\Spec k$ that is smooth and $p$-torsion is $\G_a$, so the identity component of $\sS$ is isomorphic to $\G_a$.
\end{proof}

\begin{Lemma}\label{lem:flatcohomologyofaffine}
  If $S=\Spec A$ is an affine scheme of characteristic $p$, then $\H^i(S,\m_p)=0$ for all $i\geq 3$.
\end{Lemma}
\begin{proof}
  This follows from the exact sequence of Lemma \ref{lem:4term} and Milne's isomorphisms 
  \[
    \H^i(S_{\etale},\nu(1))\iso \H^{i+1}(S_{\fl},\m_p)
  \]
  (see the proof of Corollary 1.10 of \cite{Milne76}), together with the vanishing of the higher cohomology of quasi-coherent sheaves on affine schemes.
\end{proof}

As a scheme, $\sS=\R^2\pi^{\fl}_*\m_p$ is a disjoint union of finitely many copies of $\bA^1$. Let $\bA^1\subset\sS$ be a connected component. By Tsen's theorem, $\H^2(\bA^1,\m_p)=0$, and by Lemma \ref{lem:flatcohomologyofaffine}, $\H^3(\bA^1,\m_p)=0$. So, the exact sequence of Lemma \ref{lem:mup1} gives an isomorphism
\[
  \H^2(X_{\sS},\m_p)\iso \H^0(\sS,\R^2\pi_{\sS*}^{\fl}\m_p)
\]
In particular, the universal cohomology class $\alpha\in \H^0(\sS,\R^2\pi_{\sS*}^{\fl}\m_p)$ is the image of a unique element of $\H^2(X_{\sS},\m_p)$. This element corresponds to an isomorphism class of $\m_p$-gerbes over $X\times\sS$. Let $\sX\to X\times\sS$ be such a $\m_p$-gerbe. We will think of $\sX$ as a family of twisted K3 surfaces over $\sS$ via the composition
\[
  \sX\to X\times\sS\to\sS
\]
of the coarse space morphism with the projection, and we will refer to such a gerbe $\sX$ as a \textit{tautological} or \textit{universal} family of $\m_p$-gerbes on $X$. Restricting $\sX\to\sS$ to a connected component $\bA^1\subset\sS$, we find a family
\[
  \sX'\to\bA^1
\]
of twisted supersingular K3 surfaces. This is the basic example of a twistor family.

To obtain more precise information regarding the groups $H^2(X,\m_p)$ and $\Br(X)$, we recall some consequences of flat duality for supersingular K3 surfaces. By Lemma \ref{lem:mup1}, there is a canonical identification
\begin{equation}\label{eq:canonicalidentification}
  \R^2\pi^{\fl}_*\m_p(k)=\H^2(X_{\fl},\m_p)
\end{equation}
In \cite{Artin74}, Artin identifies a certain subgroup $\U^2(X,\m_p)\subset \H^2(X,\m_p)$ and a short exact sequence
\begin{equation}\label{eq:artinsshortexactsequence}
  0\to \U^2(X,\m_p)\to \H^2(X,\m_p)\to \D^2(X,\m_p)\to 0
\end{equation}
(see also Milne \cite{Milne76}). It is an immediate consequence of the definitions that under the identification (\ref{eq:canonicalidentification}) the subgroup $\U^2(X,\m_p)\subset \H^2(X,\m_p)$ is equal to the subgroup of $\sS(k)$ consisting of the $k$-points of the connected component of the identity. Thus, as a group $\U^2(X,\m_p)$ is isomorphic to the underlying additive group of the field $k$. We will use the following result.
\begin{Theorem}\label{thm:Artinsflatduality}
  Write $\Lambda=\Pic(X)$. The Kummer sequence induces a diagram
    \begin{equation}\label{eq:artinsflatduality}
      \begin{tikzcd}
              &0\arrow{d}                       &0\arrow{d}                    &               & \\
      0\arrow{r}&\dfrac{p\Lambda^*}{p\Lambda}\arrow{r}\arrow{d}&\dfrac{\Lambda}{p\Lambda}\arrow{r}\arrow{d}&\dfrac{\Lambda}{p\Lambda^*}\arrow{r}\isor{d}{}&0\\
      0\arrow{r}&\U^2(X,\m_p)\arrow{r}\arrow{d}   &\H^2(X,\m_p)\arrow{r}\arrow{d}&\D^2(X,\m_p)\arrow{r}&0\\
      &\Br(X)\arrow[equals]{r}\arrow{d} &\Br(X)\arrow{d} &               & \\
              &0                                &0                             &               &
    \end{tikzcd}
  \end{equation}
  In particular, the natural map $\U^2(X,\m_p)\to\Br(X)$ is surjective, and the Brauer group $\Br(X)$ is $p$-torsion.
\end{Theorem}
\begin{proof}
  In Theorem 4.2 of \cite{Artin74}, this result is shown under the assumptions that $X$ admits an elliptic fibration, which is now known to always be the case by the Tate conjecture, and under the then-conjectural existence of a certain flat duality theory, which was subsequently developed by Milne \cite{Milne76}.
\end{proof}

\subsection{The relative \'{e}tale site}\label{sec:relativeetalesite}


If $X$ is a scheme, we write $X_{\Etale}$ for the big \'{e}tale site of $X$ and $X_{\etale}$ for the small \'{e}tale site of $X$. Given a morphism $\pi\colon X\to S$ of schemes (or algebraic spaces), we define a site $(X/S)_{\etale}$, which is halfway between $X_{\etale}$ and $S_{\Etale}$. A closely related construction is briefly discussed in Section 3 of \cite{MR0419450}.

\begin{Definition}
  The \textit{relative \'{e}tale site} of $X\to S$ is the category $(X/S)_{\etale}$ whose objects are pairs $(U,T)$ where $T$ is a scheme over $S$, and $U$ is a scheme \'{e}tale over $X\times_ST$. An arrow $(U',T')\to (U,T)$ consists of a morphism $U'\to U$ and a morphism $T'\to T$ over $S$ such that the diagram
  \[
    \begin{tikzcd}
      U'\arrow{r}\arrow{d}&U\arrow{d}\\
      T'\arrow{r}&T
    \end{tikzcd}
  \]
  commutes. Note that this data is equivalent to giving a morphism $T'\to T$ over $S$ and a morphism $U'\to U\times_{X_{T}}X_{T'}$ over $X_{T'}$. A family of morphisms $\left\{(U_i,T_i)\to (U,T)\right\}_{i\in I}$ is a covering if both $\left\{U_i\to U\right\}_{i\in I}$ and $\left\{T_i\to T\right\}_{i\in I}$ are \'{e}tale covers.
\end{Definition}

The small, relative, and big \'{e}tale sites are related by various obvious functors, which induce a diagram of sites
\[
  \begin{tikzcd}
    X_{\Etale}\arrow{r}\arrow[bend left=30]{rr}{\alpha_X}\arrow{dr}[swap]{\pi^{\Etale}}&(X/S)_{\etale}\arrow{r}\arrow{d}{\pi^{\etale}}&X_{\etale}\arrow{d}\\
    &S_{\Etale}\arrow{r}{\alpha_S}&S_{\etale}
  \end{tikzcd}
\]
\begin{Remark}
  The relative \'{e}tale site is well suited to studying structures that are compatible with arbitrary base change on $S$ and \'{e}tale base change on $X$, such as sheaves of relative differentials and the relative Frobenius. In this work, we will be concerned with various representable functors on $S_{\Etale}$ obtained as higher pushforwards of sheaves along the map $\pi^{\etale}\colon (X/S)_{\etale}\to S_{\Etale}$. The relative \'{e}tale site will be useful to us because it is more constrained than $X_{\Etale}$, while still being large enough to admit a map to the big \'{e}tale site of $S$.
\end{Remark}
 
If $\pi'\colon X'\to S'$ is a morphism, and $g\colon S'\to S$ and $f\colon X'\to X$ are morphisms such that $\pi\circ f=g\circ \pi'$, then there is an induced map of sites
\[
  (f,g)\colon (X'/S')_{\etale}\to (X/S)_{\etale}
\]
On underlying categories, this map sends a pair $(U,T)$ to $(U\times_{X}X',T\times_SS')$. We obtain a corresponding pushforward functor $(f,g)_*$, which has an exact left adjoint $(f,g)^{-1}$, and together give a morphism of topoi.
\begin{Lemma}\label{lem:heythatscool}
  If $f$ and $g$ are universal homeomorphisms, then the pushforward
  \[
    (f,g)_*\colon \Sh_{(X'/S')_{\etale}}\to\Sh_{(X/S)_{\etale}}
  \]
  between the respective categories of sheaves of abelian groups on the relative \'{e}tale sites is exact.
\end{Lemma}
\begin{proof}
  Being the right adjoint to the pullback, $(f,g)_*$ is left exact. That it is also right exact follows from the topological invariance of the small \'{e}tale site under universal homeomorphisms \cite[05ZG]{stacks-project}. Let $\sE\to\sF$ be a surjective map of sheaves of abelian groups on $(X'/S')_{\etale}$. We will show that
  \[
    (f,g)_*\sE\to (f,g)_*\sF
  \]
  remains surjective. That is, we will show that for any object $(U,T)\in (X/S)_{\etale}$, any section $x\in (f,g)_*\sF(U,T)$ is in the image after taking a cover of $(U,T)$. By a change of notation, we reduce to the case when $T=S$. Let $U_{S'}=U\times_{S}S'$, and choose a cover $\left\{(U_i',T_i')\to (U_{S'},S')\right\}_{i\in I}$ such that the restriction of $x\in (f,g)_*\sF(U,T)=\sF(U_{S'},S')$ is in the image of each $\sF(U_i',T_i')\to\sE(U_i',T_i')$. By the topological invariance of the \'{e}tale site applied to $g\colon S'\to S$, we may find morphisms $T_i\to S$ such that $T_i\times_SS'=T_i'$. Applying the same to the universal homeomorphisms $X'_{T_i'}\to X_{T_i}$, we find morphisms $U_i\to X_{T_i}$ such that $U_i\times_{S}S'=U_i'$. The restriction of $x$ to the cover $\left\{(U_i,T_i)\to (U,S)\right\}_{i\in I}$ is then in the image by construction.
\end{proof}

There is a sheaf of rings $\cO_X$ on $(X/S)_{\etale}$ given by
\[
  \cO_X(U,T)=\Gamma(U,\cO_U)
\]
The pair $((X/S)_{\etale},\cO_X)$ is a ringed site, in the sense of \cite[Tag 03AD]{stacks-project}.
\begin{Definition}\label{def:quasi-coherent}
  A sheaf $\sE$ of $\cO_X$-modules on $(X/S)_{\etale}$ is \textit{quasi-coherent} if
  \begin{enumerate}
      \item for every object $(U,T)$ of $(X/S)_{\etale}$ there exists an open cover $U_i\to U$ such that for each $i$ the restriction $\sE|_{U_{i\zariski}}$ of $\sE$ to the small Zariski site of $U_i$ is a quasi-coherent sheaf (in the usual sense), and
      \item if $(f,g)\colon (U',T')\to (U,T)$ is a morphism in $(X/S)_{\etale}$ then the natural comparison map
      \[
        f^*(\sE|_{U_{\zariski}})\to \sE|_{U'_{\zariski}}
      \]
      is an isomorphism.
  \end{enumerate}
  The sheaf $\sE$ is \textit{locally free of finite rank} if it is quasi-coherent, and for each object $(U,T)$ the restriction of $\sE$ to the small Zariski site of $U$ is locally free of finite rank in the usual sense.
\end{Definition}
\begin{Remark}
  This definition is similar to \cite[Tag 06WK]{stacks-project}. In \cite[Tag 03DL]{stacks-project} there is given a general definition of quasi-coherent sheaf and locally free sheaf on a ringed site. One can show that these agree with the above definitions in our special case. 
\end{Remark}

Our main interest in this paper is when $X$ and $S$ are schemes over $\bF_p$. In this case, we have a diagram
\[
  \begin{tikzcd}
    X\arrow{r}{F_{X/S}}\arrow{dr}[swap]{\pi}\arrow[bend left=45]{rr}{F_X}&X^{(p/S)}\arrow{r}{W_{X/S}}\arrow{d}{\pi^{(p/S)}}&X\arrow{d}{\pi}\\
    &S\arrow{r}{F_S}&S
  \end{tikzcd}
\]
where $F_X$ and $F_S$ are the absolute Frobenius morphisms, and the square is Cartesian. If there is no risk of confusion, we may write $F=F_X$, $X^{(p/S)}=X^{(p)}$, $\pi^{(p/S)}=\pi^{(p)}$, and $W_{X/S}=W$. We will denote the various induced maps on the relative \'{e}tale sites by $F_X=(F_X,F_S)$, $F_{X/S}=(F_{X/S},\id_S)$, and $W_{X/S}=(W_{X/S},F_S)$.

\begin{Lemma}\label{lem:relativefrobenius}
  If $\pi$ is smooth, then $F_{X/S}$ is finite and flat. If $\pi$ is locally of finite presentation, then the relative Frobenius $F_{X/S}$ is an isomorphism if and only if $\pi$ is \'{e}tale. 
\end{Lemma}
\begin{proof}
  See Proposition 3.2 of \cite{Ill96}. The second claim is Proposition 2 of Expos\'{e} XV 1 in \cite{MR0491704}.
\end{proof}

\begin{Lemma}\label{lem:relativefrobeniusbasechange}
  The relative Frobenius is compatible with arbitrary base change on $S$ in the following sense. For any morphism $T\to S$, there is a canonical identification 
    \[
      (X^{(p/S)})\times_ST=(X\times_ST)^{(p/T)}
    \]
  and a commutative diagram
  \[
    \begin{tikzcd}
      X_T\arrow[bend left=30]{rr}{\pi_T}\arrow{r}[swap]{F_{X_T/T}}\arrow{d}[swap]{\pi_X}&(X_T)^{(p/T)}\arrow{d}\arrow{r}[swap]{\pi_T^{(p/T)}}&T\arrow{d}\\
      X\arrow[bend right=30]{rr}[swap]{\pi}\arrow{r}{F_{X/S}}&X^{(p/S)}\arrow{r}{\pi^{(p/S)}}&S
    \end{tikzcd}
  \]
  where both squares are Cartesian.
\end{Lemma}

\begin{Lemma}\label{lem:relativeetalesiteisgood}
  The relative Frobenius is compatible with \'{e}tale base change on $X^{(p)}$, in the following sense. Suppose that $U\to X^{(p)}$ is \'{e}tale, and define $U_X=X\times_{X^{(p)}}U$ by the Cartesian diagram
  \[
    \begin{tikzcd}
      U_X\arrow{r}\arrow{d}&U\arrow{d}\\
      X\arrow{r}{F_{X/S}}&X^{(p)}
    \end{tikzcd}
  \]
  The relative Frobenius $F_{U/X^{(p)}}\colon U\to U^{(p/X^{(p)})}=(U_X)^{(p/S)}$ is an isomorphism, and there is a diagram
  \[
    \begin{tikzcd}
      U_X\arrow{dr}\arrow{r}{F_{U_X/S}}&(U_X)^{(p/S)}\\
     &U\isor{u}{F_{U/X^{(p)}}}
    \end{tikzcd}
  \]
\end{Lemma}
\begin{proof}
  Consider the diagram
  \[
    \begin{tikzcd}
      U_X\arrow{r}\isor{d}{F_{U_X/X}}&U\isol{d}{F_{U/X^{(p)}}}&&&\\
      (U_X)^{(p/X)}\arrow{r}\arrow{d}&(U_X)^{(p/S)}\arrow{r}\arrow{d}&U_X\arrow{d}\arrow{r}&U\arrow{d}\\
      X\arrow{r}{F_{X/S}}\arrow{dr}&X^{(p)}\arrow{r}{W_{X/S}}\arrow{d}&X\arrow{d}\arrow{r}{F_{X/S}}&X^{(p)}\\
      &S\arrow{r}{F_S}&S&
    \end{tikzcd}
  \]
  where the squares are Cartesian. By Lemma \ref{lem:relativefrobenius}, $F_{U/X^{(p)}}$ is an isomorphism, and the composition $U_X\to (U_X)^{(p/X)}\to(U_X)^{(p/S)}$ is equal to $F_{U_X/S}$.
\end{proof}

Consider the sheaf of units $\cO_X^\times$ on $X_{\Etale}$ given by $\cO_X^\times(T)=\Gamma(T,\cO_T)^\times$. Its pushforward to $(X/S)_{\etale}$, which we will also denote by $\cO_X^\times$, is given by
\[
  \cO_X^\times(U,T)=\Gamma(U,\cO_U)^\times
\]
Given a morphism $T\to X^{(p)}$, we consider the diagram
\[
  \begin{tikzcd}
    T\arrow{r}{F_{T/X^{(p)}}}\arrow{dr}&T^{(p/X^{(p)})}\arrow{r}\arrow{d}&T'\arrow{r}\arrow{d}&T\arrow{d}\\
    &X^{(p)}\arrow{r}{W_{X/S}}&X\arrow{r}{F_{X/S}}&X^{(p)}
  \end{tikzcd}
\]
where the squares are Cartesian. Pulling back along $T\to T'$ and $T'\to T$ induces natural maps
\[
  F_{X/S*}\cO_X^\times\to\cO_{X^{(p)}}^\times\xrightarrow{F_{X/S}^*}F_{X/S*}\cO_X^\times
\]
of big \'{e}tale sheaves whose composition is given by the $p$-th power map $x\mapsto x^p$. We therefore have a commuting square
\begin{equation}\label{eq:itsasquare}
  \begin{tikzcd}
    F_{X/S*}\cO_X^\times\arrow{r}{x\mapsto x^p}\arrow{d}&F_{X/S*}\cO_X^\times\arrow[equal]{d}\\
    \cO_{X^{(p)}}^\times\arrow{r}{F_{X/S}^*}&F_{X/S*}\cO_X^\times
  \end{tikzcd}
\end{equation}
of sheaves on $X^{(p)}_{\Etale}$, or by restriction on $(X^{(p)}/S)_{\etale}$. It is shown in Lemma 2.1.18 of \cite{Ill79} that the restriction of the map $F_{X/S}^*$ to the small \'{e}tale site is injective. By Lemma \ref{lem:relativefrobeniusbasechange} and Lemma \ref{lem:relativeetalesiteisgood}, it is also injective on the relative \'{e}tale site.
\begin{Definition}
  We define a sheaf $\nu(1)$ on $(X^{(p)}/S)_{\etale}$ by the cokernel of the pullback map $F_{X/S}^*\colon \cO_{X^{(p)}}^\times\to F_{X/S*}\cO_X^\times$ If $X\to S$ is smooth, then this map is injective, and we have a short exact sequence
  \begin{equation}\label{eq:defNu}
    1\to\cO_{X^{(p)}}^\times\to F_{X/S*}\cO_X^\times\to\nu(1)\to 1
  \end{equation}
\end{Definition}
If $X$ is a scheme, we let $\varepsilon_X\colon X_{\fl}\to X_{\Etale}$ and $\varepsilon_{X/S}\colon X_{\fl}\to (X/S)_{\etale}$ denote the natural maps of sites. These are related by a diagram
\begin{equation}\label{eq:changeoftopology}
  \begin{tikzcd}
    X_{\fl}\arrow{r}{\varepsilon_{X/S}}\arrow{d}[swap]{\pi^{\fl}}\arrow{dr}{\Theta}&(X/S)_{\etale}\arrow{d}{\pi^{\etale}}\\
    S_{\fl}\arrow{r}{\varepsilon_S}&S_{\Etale}
  \end{tikzcd}
\end{equation}

We recall the following theorem of Grothendieck.

\begin{Theorem}[\cite{MR0244271}, Th\'{e}or\`{e}me 11.7]\label{thm:grothendieck}
 If $X$ is a scheme and $A$ is a smooth group scheme over $X$, then $\R^p\varepsilon_{X*}A=0$ for $p>0$.
\end{Theorem}

We consider the long exact sequence induced by applying $\varepsilon_{X/S*}$ to the short exact sequence
\begin{equation}\label{eq:kummermup}
  1\to\m_p\to\cO_X^{\times}\xrightarrow{x\mapsto x^p}\cO_X^\times\to 1
\end{equation}
of sheaves on $X_{\fl}$. Because pushforward along the map $X_{\Etale}\to (X/S)_{\etale}$ is exact, Grothendieck's Theorem implies that $\R^p\varepsilon_{X/S*}\cO_X^\times=0$ for $p>0$. Thus, we obtain an exact sequence
\[
  1\to \m_p\to \cO_X^\times\xrightarrow{x\mapsto x^p}\cO_X^\times\to \R^1\varepsilon_{X/S*}\m_p\to 1
\]
of sheaves on $(X/S)_{\etale}$. Let us write
\[
  \R^1\varepsilon_{X/S*}\m_p=\cO_X^\times/\cO_X^{\times p}
\]
where it is understood that the quotient is taken in the \'{e}tale topology. By Lemma \ref{lem:heythatscool}, this sequence remains exact after applying $F_{X/S*}$, and so the square (\ref{eq:itsasquare}) induces a morphism
\begin{equation}\label{eq:themorphism}
  F_{X/S*}(\cO_X^\times/\cO_X^{\times p})\to\nu(1)
\end{equation}
of sheaves on $(X^{(p)}/S)_{\etale}$.

\begin{Proposition}\label{prop:flattoetalecohomology}
  If $\pi\colon X\to S$ is a smooth morphism of $\bF_p$-schemes such that
  \begin{enumerate}
    \item The adjunction $\cO_S\to\pi_*\cO_X$ is universally an isomorphism, and
    \item $\R^1\pi^{\fl}_*\m_p=0$,
  \end{enumerate}
  then there is a natural map 
  \[
    \Upsilon\colon \varepsilon_{S*}\R^2\pi^{\fl}_*\m_p\to \R^1\pi^{(p)\etale}_*\nu(1)
  \]
  of big \'{e}tale sheaves. For any perfect scheme $T$ over $S$, this map is a bijection on $T$-points.
\end{Proposition}
\begin{proof}
  First, note that by Lemma \ref{lem:heythatscool}
  \[
    \R^1\pi^{(p)\etale}_*F_{X/S*}(\cO_X^\times/\cO_X^{\times p})=\R^1\pi^{\etale}_*(\cO_X^\times/\cO_X^{\times p})
  \]
  Thus, the morphism (\ref{eq:themorphism}) induces a map
  \begin{equation}\label{eq:morphism11}
    \R^1\pi^{\etale}_*(\cO_X^\times/\cO_X^{\times p})\to \R^1\pi^{(p)\etale}_*\nu(1)
  \end{equation}
  Next, we consider the Grothendieck spectral sequences
  \[
    E^{p,q}_{2}=\R^p\pi^{\etale}_*(\R^q\varepsilon_{X/S*}\m_p)\Longrightarrow \R^{p+q}\Theta_*\m_p
  \]
  \[
    \overline{E}^{p,q}_{2}=\R^p\varepsilon_{S*}(\R^q\pi^{\fl}_*\m_p)\Longrightarrow \R^{p+q}\Theta_*\m_p
  \]
  induced by the square (\ref{eq:changeoftopology}). The first spectral sequence gives an exact sequence
  \[
    E^{2,0}_2\to\ker(E^2_{\infty}\to E^{0,2}_2)\to E^{1,1}_2\to E^{3,0}_2
  \]
  Applying Theorem \ref{thm:grothendieck} to $\bG_m$, the Kummer sequence shows that $\R^p\varepsilon_{X/S*}\m_p=0$ for $p\geq 2$
  Therefore, $E^{0,2}_{2}=0$, so we get a map 
  \[
    \R^2\Theta_*\m_p\to \R^1\pi^{\etale}_*\R^1\varepsilon_{X/S*}\m_p=\R^1\pi^{\etale}_*(\cO_X^\times/\cO_X^{\times p})
  \]
  The second spectral sequence gives a map $\R^2\Theta_*\m_p\to\varepsilon_{S*}\R^2\pi^{\fl}_*\m_p$. Hence, we have a diagram
  \[
    \begin{tikzcd}
      &\R^2\Theta_*\m_p\arrow{dl}[swap]{f}\arrow{dr}{g}&&\\
      \varepsilon_{S*}\R^2\pi^{\fl}_*\m_p&&\R^1\pi^{\etale}_*(\cO_X^\times/\cO_X^{\times p})\arrow{r}&\R^1\pi^{(p)\etale}_*\nu(1)
    \end{tikzcd}
  \]
  It follows from assumptions (1) and (2) that $f$ is an isomorphism. We therefore obtain the desired map
  \[
    \Upsilon\colon \varepsilon_{S*}\R^2\pi^{\fl}_*\m_p\to \R^1\pi^{(p)\etale}_*\nu(1)
  \]
  Suppose that $T$ is a perfect $S$-scheme. For any \'{e}tale $U\to X^{(p)}_T$, the map $F_{X/S*}\cO_X^\times(U,T)\to\cO_{X^{(p)}}^\times(U,T)$ is given by pullback along a base change of $W_{X_T/T}$ followed by pullback along the relative Frobenius $F_{U/X_T^{(p)}}$, both of which are isomorphisms. Therefore, by restricting (\ref{eq:itsasquare}) to the small \'{e}tale site (which is an exact functor), we get an isomorphism
  \[
    \alpha_{T*}\R^1\pi^{\etale}_{T*}\R^1\varepsilon_{X_T/T*}\m_p\iso \alpha_{T*}\R^1\pi^{(p)\etale}_{T*}\nu(1)
  \]
  of sheaves on $S_{\etale}$. We have an exact sequence
  \[
    \R^2\pi_{T*}^{\etale}\varepsilon_{X_T/T*}\m_p\to\varepsilon_{T*}\R^2\pi_{T*}^{\fl}\m_p\to \R^1\pi_{T*}^{\etale}\R^1\varepsilon_{X_T/T}\m_p\to \R^3\pi_{T*}^{\etale}\varepsilon_{X_T/T*}\m_p
  \]
   Because $X_T$ is smooth over a perfect scheme, $\alpha_{X_T*}\varepsilon_{X_T/T*}\m_p=0$, and therefore the middle arrow becomes an isomorphism. We have shown that the induced map
   \[
     \alpha_{T*}\varepsilon_{T*}\R^2\pi_{T*}^{\fl}\m_p\iso\alpha_{T*}\R^1\pi^{\etale}_{T*}\R^1\varepsilon_{X_T/T*}\m_p\iso \alpha_{T*}\R^1\pi^{(p)\etale}_{T*}\nu(1)
   \]
   of small \'{e}tale sheaves is an isomorphism. In particular, it is a bijection on global sections, which gives the result.
\end{proof}

In this paper we are interested in the case when $\pi\colon X\to S$ is a relative supersingular K3 surface. If this is so, then the conditions of Proposition \ref{prop:flattoetalecohomology} are satisfied. Moreover, by Theorem \ref{thm:relativeartin} the sheaf $\R^2\pi^{\fl}_*\m_p$ is representable, and we will see that $\R^1\pi^{(p)\etale}_*\nu(1)$ is as well in Proposition \ref{prop:nuRepresentable}. The following result therefore applies to show that $\Upsilon$ is a universal homeomorphism.

\begin{Proposition}\label{prop:homeomorphism}
  Let $S$ be a scheme over $\bF_p$ and $f\colon X\to Y$ a morphism of $S$-schemes. If $f$ induces a bijection on $T$-points for every perfect $S$-scheme $T$, then $f$ is a universal homeomorphism.
\end{Proposition}
\begin{proof}
  Recall that the \textit{perfection} $X^{\pf}\to X$ of an $\bF_p$-scheme $X$ is characterized by the property that any map $T\to X$ from a perfect scheme $T$ factors uniquely through $X^{\pf}\to X$. There is therefore a commutative diagram
  \[
    \begin{tikzcd}
      X^{\pf}\arrow{r}{f^{\pf}}\arrow{d}&Y^{\pf}\arrow{d}\\
      X\arrow{r}{f}&Y
    \end{tikzcd}
  \]
  Applying our assumption on $f$ to the $S$-scheme $Y^{\pf}\to Y\to S$, we get a map $Y^{\pf}\to X$, and hence a map $Y^{\pf}\to X^{\pf}$. By the universal property of the perfection, this map is an inverse to $f^{\pf}$, so $f^{\pf}\colon X^{\pf}\to Y^{\pf}$ is an isomorphism. By Lemma 3.8 of \cite{BS15}, this is equivalent to $f$ being a universal homeomorphism.
\end{proof}

Finally, we will compute the tangent spaces of $\R^2\pi^{\fl}_*\m_p$ and $\R^1\pi^{(p)\etale}_*\nu(1)$. The former is slightly subtle, as the pushforward under a universal homeomorphism (such as the projection $\pi_X\colon X[\varepsilon]\to X$) is not necessarily exact in the flat topology. We record the following lemma.
\begin{Lemma}\label{lem:kindoflame}
  If $X$ is a scheme and $\pi_X\colon X[\varepsilon]\to X$ is the natural projection, then $\R^1\pi^{\fl}_{X*}\cO_{X[\varepsilon]}^\times=0$.
\end{Lemma}
\begin{proof}
  Consider the diagram
  \begin{equation}\label{eq:diagram1aaaa}
    \begin{tikzcd}
      X[\varepsilon]_{\fl}\arrow{d}[swap]{\varepsilon_{X[\varepsilon]}}\arrow{r}{\pi_X^{\fl}}\arrow{dr}{\Theta}&X_{\fl}\arrow{d}{\varepsilon_X}\\
      X[\varepsilon]_{\Etale}\arrow{r}[swap]{\pi_X^{\Etale}}&X_{\Etale}
    \end{tikzcd}
  \end{equation}
  and the two induced Grothendieck spectral sequences. On the one hand, we have by Theorem \ref{thm:grothendieck} that $\R^p\varepsilon_{X[\varepsilon]*}\cO_{X[\varepsilon]}^\times=0$ for $p\geq 1$. As in Lemma \ref{lem:heythatscool}, the pushforward $\pi_{X*}^{\Etale}$ is exact, and we conclude that $\R^p\Theta_*\cO_{X[\varepsilon]}^\times=0$ for $p\geq 1$. We now consider the other spectral sequence. The exact sequence of low degree terms gives an isomorphism
  \[
    \varepsilon_{X*}\R^1\pi_{X*}^{\fl}\cO_{X[\varepsilon]}^\times\iso\R^2\varepsilon_{X*}\pi_{X*}^{\fl}\cO_{X[\varepsilon]}^\times
  \]
  We will show that the right hand side vanishes. Consider the standard short exact sequence
  \begin{equation}\label{eq:morphism1a}
    0\to\cO_X\to\pi^{\fl}_{X*}\cO^\times_{X\left[\varepsilon\right]}\to\cO^\times_X\to 1
  \end{equation}
  of sheaves on the big flat site of $X$, where the first map is given by $f\mapsto 1+f\varepsilon$ and the second by $g+f\varepsilon\mapsto g$. Note that this sequence is split. Therefore, $\pi^{\fl}_{X*}\cO^\times_{X\left[\varepsilon\right]}$ is represented by a smooth group scheme, and hence by Theorem \ref{thm:grothendieck} we have $\R^p\varepsilon_{X*}\pi_{X*}^{\fl}\cO_{X[\varepsilon]}^\times=0$ for all $p\geq 1$. This completes the proof.
\end{proof}

\begin{Lemma}\label{lem:tangentspacetomodulispace}
  Suppose that $\pi\colon X\to S$ is a smooth morphism of $\bF_p$-schemes such that
  \begin{enumerate}
      \item The sheaf $\R^2\pi^{\fl}_*\m_p$ is representable by an algebraic space,
      \item The adjunction $\cO_S\to\pi_*\cO_X$ is universally an isomorphism, and
      \item $\R^1\pi_*\cO_X=0$.
  \end{enumerate}
  Set $\sS=\R^2\pi^{\fl}_*\m_p$ and let $p\colon \sS\to S$ be the forgetful morphism. There is a natural identification
  \[
    T^1_{\sS/S}\iso p^*\R^2\pi_*\cO_{X}
  \]
\end{Lemma}
\begin{proof}
  Let $\sigma_e\colon S\to\sS$ be the identity section. We will construct an isomorphism
  \[
    \R^2\pi_*\cO_X\iso\sigma_e^*T^1_{\sS/S}
  \]
  Consider the Cartesian diagram
  \begin{equation}\label{eq:diagram1a}
    \begin{tikzcd}
      X[\varepsilon]\arrow{d}[swap]{\pi_{[\varepsilon]}}\arrow{r}{\pi_X}&X\arrow{d}{\pi}\\
      S[\varepsilon]\arrow{r}{\pi_S}&S
    \end{tikzcd}
  \end{equation}
  The group of sections $\sigma_e^*T^1_{\sS/S}(U)$ over an open set $U\subset S$ is in natural bijection with the collection of morphisms $t\colon U[\varepsilon]\to\sS$ over $S$ such that the diagram
  \[
    \begin{tikzcd}
    &&\sS\arrow{d}\\
    U\arrow{r}\arrow[bend left=10]{urr}{\sigma_e}&U[\varepsilon]\arrow{r}\arrow{ur}{t}&S
    \end{tikzcd}
  \]
  commutes. Because $\R^2\pi^{\fl}_*\m_p$ is compatible with base change, this is the the same as the kernel of the natural map
  \begin{equation}\label{eq:thetangentspace}
    \pi_{S*}\R^2\pi^{\fl}_{[\varepsilon]*}\m_p\to \R^2\pi^{\fl}_*\m_p
  \end{equation}
  of small Zariski sheaves. Combining (\ref{eq:morphism1a}) with the Kummer sequence, and using the vanishing result of Lemma \ref{lem:kindoflame}, we obtain a diagram
  \[
    \begin{tikzcd}
              &0\arrow{d}&1\arrow{d}&1\arrow{d}&\\
    0\arrow{r}&\cO_X\arrow{r}\isor{d}{}&\pi^{\fl}_{X*}\m_p\arrow{r}\arrow{d}&\m_p\arrow{r}\arrow{d}&1\\
    0\arrow{r}&\cO_X\arrow{r}\arrow{d}{0}&\pi^{\fl}_{X*}\cO^\times_{X[\varepsilon]}\arrow{r}\arrow{d}{\cdot p}&\cO^\times_X\arrow{r}\arrow{d}{\cdot p}&1\\
    0\arrow{r}&\cO_X\arrow{r}\isor{d}{}&\pi^{\fl}_{X*}\cO^\times_{X[\varepsilon]}\arrow{r}\arrow{d}&\cO^\times_X\arrow{d}\arrow{r}&1\\
    0\arrow{r}&\cO_X\arrow{r}{\sim}\arrow{d}&\R^1\pi^{\fl}_{X*}\m_p\arrow{r}\arrow{d}&1&\\
    &0&1&&
    \end{tikzcd}
  \]
  of sheaves on the big flat site of $X$ with exact rows and columns. Taking cohomology of the split exact sequence
  \[
    0\to\cO_X\to\pi^{\fl}_{X*}\m_{p}\to\m_{p}\to 1
  \]
    we get an exact sequence
  \[
    0\to \R^2\pi^{\fl}_{*}\cO_X\to \R^2\pi^{\fl}_{*}(\pi^{\fl}_{X*}\m_{p})\to \R^2\pi^{\fl}_{*}\m_p\to 0
  \]
  of sheaves on the big flat site of $S$. To compare this to the kernel of the morphism (\ref{eq:thetangentspace}), we consider the spectral sequences
  \[
    E^{p,q}_{2}=\R^p\pi^{\fl}_*(\R^q\pi^{\fl}_{X*}\m_p)\Longrightarrow \R^{p+q}\Theta^{\fl}_*\m_p
  \]
  \[
    \overline{E}^{p,q}_{2}=\R^p\pi^{\fl}_{S*}(\R^q\pi^{\fl}_{[\varepsilon]*}\m_p)\Longrightarrow \R^{p+q}\Theta^{\fl}_*\m_p
  \]
  induced by the commuting square (\ref{eq:diagram1a}), where $\Theta=\pi\circ\pi_X=\pi_S\circ\pi_{[\varepsilon]}$. The first spectral sequence gives an exact sequence
  \[
    0\to E^{1,0}_2\to E^1_{\infty}\to E^{0,1}_2\to E^{2,0}_2\to\ker(E^2_{\infty}\to E^{0,2}_2)\to E^{1,1}_2\to E^{3,0}_2
  \]
  By the Kummer sequence and Grothendieck's Theorem, we get that $\R^q\pi^{\fl}_{X*}\m_p=0$ for $q\geq 2$, and in particular, $E^{0,2}_2=0$. The isomorphism 
  \begin{equation}\label{eq:morphism3}
    \cO_X\iso\R^1\pi^{\fl}_{X*}\m_p
  \end{equation}
  and condition (3) imply that $E^{1,1}_2=0$. Thus, we have an exact sequence
  \begin{equation}\label{eq:anexactsequence}
    0\to\R^1\pi^{\fl}_*\pi^{\fl}_{X*}\m_p\to\R^1\Theta^{\fl}_*\m_p\to\pi^{\fl}_*\R^1\pi^{\fl}_{X*}\m_p\to\R^2\pi^{\fl}_*\pi^{\fl}_{X*}\m_p\to\R^2\Theta^{\fl}_*\m_p\to 0
  \end{equation}
  We next examine the second spectral sequence. We have
  \begin{align*}
    \overline{E}^{p,0}_2&=\R^p\pi^{\fl}_{S*}(\pi^{\fl}_{[\varepsilon]*}\m_p)=\R^p\pi^{\fl}_{S*}\m_p=0\mbox{ for }p\geq2,\mbox{ and }\\
    \overline{E}^{p,1}_2&=\R^p\pi^{\fl}_{S*}(\R^1\pi^{\fl}_{[\varepsilon]*}\m_p)=0\mbox{ for }p\geq 0.
  \end{align*}
  Thus, we have natural isomorphisms
  \begin{align*}
      \R^1\pi^{\fl}_{S*}(\pi^{\fl}_{[\varepsilon]*}\m_p)&\iso\R^1\Theta^{\fl}_*\m_p\\
      \R^2\Theta^{\fl}_*\m_p&\iso\pi^{\fl}_{S*}\R^2\pi^{\fl}_{[\varepsilon]*}\m_p
  \end{align*}
  Comparing with the exact sequence (\ref{eq:anexactsequence}), we find a surjection
  \begin{equation}\label{eq:morphism1}
    \R^2\pi^{\fl}_*\pi^{\fl}_{X*}\m_p\to\R^2\Theta^{\fl}_*\m_p\iso\pi^{\fl}_{S*}\R^2\pi^{\fl}_{[\varepsilon]*}\m_p
  \end{equation}
  whose kernel is given by the cokernel of the map
  \begin{equation}\label{eq:morphism2}
    \R^1\pi^{\fl}_{S*}(\pi^{\fl}_{[\varepsilon]*}\m_p)\iso\R^1\Theta^{\fl}_*\m_p\to\pi^{\fl}_*\R^1\pi^{\fl}_{X*}\m_p
  \end{equation}
  We claim that (\ref{eq:morphism2}) is an isomorphism. Indeed, using the isomorphism (\ref{eq:morphism3}) and condition (2), we find a diagram
  \[
    \begin{tikzcd}
      \R^1\pi^{\fl}_{S*}(\pi^{\fl}_{[\varepsilon]*}\m_p)\arrow{r}\isor{d}{}&\pi^{\fl}_*\R^1\pi^{\fl}_{X*}\m_p\isor{d}{}\\
      \cO_S\arrow{r}&\cO_S
    \end{tikzcd}
  \]
  where the induced map $\cO_S\to\cO_S$ is the identity. Thus, (\ref{eq:morphism2}) and hence (\ref{eq:morphism1}) are isomorphisms, and restricting to the Zariski site we obtain a map of short exact sequences
  \[
    \begin{tikzcd}
      0\arrow{r}&\sigma_e^*T^1_{\sS/S}\arrow{r}\isor{d}{}&\pi_{S*}\R^2\pi^{\fl}_{[\varepsilon]*}\m_p\arrow{r}\isor{d}{}&\R^2\pi^{\fl}_*\m_p\arrow{r}\isor{d}{}&0\\
      0\arrow{r}&\R^2\pi^{\fl}_*\cO_X\arrow{r}&\R^2\pi^{\fl}_*(\pi^{\fl}_{X*}\m_p)\arrow{r}&\R^2\pi^{\fl}_*\m_p\arrow{r}&0
    \end{tikzcd}
  \]
  For any group space $p\colon G\to S$, there is a canonical isomorphism
  \[
    T^1_{G/S}\iso p^*\sigma_e^*T^1_{G/S}
  \]
  The composition
  \[
    T^1_{\sS/S}\iso p^*\sigma_e^*T^1_{\sS/S}\iso p^*\R^2\pi_*\cO_X
  \]
  gives the desired isomorphism.
\end{proof}

\begin{Definition}
  We define a quasicoherent sheaf $B^1_{X/S}$ on $X^{(p)}$ by the short exact sequence
  \[
    0\to\cO_{X^{(p)}}\to F_{X/S*}\cO_X\to B^1_{X/S}\to 0
  \]
\end{Definition}

\begin{Lemma}\label{lem:tangentspacetonu}
  Suppose that $\pi\colon X\to S$ is a smooth morphism of $\bF_p$-schemes such that the sheaf $\R^1\pi^{(p)\etale}_*\nu(1)$ is representable by an algebraic space. Set $\sS_\nu=\R^1\pi^{(p)\etale}_*\nu(1)$ and let $p_\nu\colon \sS_\nu\to S$ be the forgetful morphism. There is a natural identification
  \[
    T^1_{\sS_\nu/S}\iso p_\nu^*\R^1\pi^{(p)}_*B^1_{X/S}
  \]
\end{Lemma}
\begin{proof}
  Let $\sigma_e\colon S\to\sS_\nu$ be the identity section. We will construct an isomorphism
  \[
    \R^1\pi^{(p)}_*B^1_{X/S}\iso\sigma^*_eT^1_{\sS_\nu/S}
  \]
  As in Lemma \ref{lem:tangentspacetomodulispace}, we identify the sheaf $\sigma^*_eT^1_{\sS_\nu/S}$ with the kernel of the natural map
  \begin{equation}\label{eq:thetangentspacenu}
    \pi_{S*}\R^1\pi^{(p)\etale}_{[\varepsilon]*}\nu(1)\to\R^1\pi^{(p)\etale}_*\nu(1)
  \end{equation}
  where $\pi^{(p)\etale}_{[\varepsilon]}=(\pi_{X^{(p)}},\pi_S)$ is the map of relative \'{e}tale sites defined as in (\ref{eq:diagram1a}). Restricting (\ref{eq:morphism1a}) to the relative \'{e}tale site, we get a short exact sequence
  \[
    0\to\cO_{X}\to\pi^{\etale}_{X*}\cO_{X[\varepsilon]}^\times\to\cO_X^\times\to 1
  \]
  of sheaves on the relative \'{e}tale site $(X/S)_{\etale}$. We have a diagram
  \[
    \begin{tikzcd}
                    &                   0  \arrow{d}&1                                                   \arrow{d}&1\arrow{d}                            & \\
        0\arrow{r}  &\cO_{X^{(p)}}\arrow{r}\arrow{d}&\pi^{\etale}_{X^{(p)}*}\cO_{X^{(p)}[\varepsilon]}   \arrow{r}\arrow{d}&\cO_{X^{(p)}}^\times\arrow{r}\arrow{d}&1\\
        0\arrow{r}  &F_{X/S*}\cO_X\arrow{r}\arrow{d}&F_{X/S*}\pi^{\etale}_{X*}\cO_{X[\varepsilon]}^\times\arrow{r}\arrow{d}&F_{X/S*}\cO_X^\times\arrow{r}\arrow{d}&1\\
        0\arrow{r}  &B^1_{X/S}          \arrow{r}\arrow{d}&\pi^{\etale}_{X^{(p)}*}\nu(1)                       \arrow{r}\arrow{d}&\nu(1)              \arrow{r}\arrow{d}&1\\
                    &0                              &1                                                            &1                                     &
    \end{tikzcd}
  \]
  Taking cohomology of the split exact sequence
  \[
    0\to B^1_{X/S}\to\pi^{\etale}_{X^{(p)}*}\nu(1)\to\nu(1)\to 1
  \]
  we get an exact sequence
  \[
    0\to\R^1\pi^{(p)\etale}_*B^1_{X/S}\to\R^1\pi^{(p)\etale}_*\pi^{\etale}_{X^{(p)}*}\nu(1)\to\R^1\pi^{(p)\etale}_*\nu(1)\to 1
  \]
  By Lemma \ref{lem:heythatscool}, the functors $\pi^{\etale}_{X^{(p)}*}$ and $\pi_{S*}^{\etale}$ are exact. Thus, the kernel of the morphism (\ref{eq:thetangentspacenu}) is identified with $\R^1\pi^{(p)}_*B^1_{X/S}$, and we obtain an isomorphism
  \[
    \R^1\pi^{(p)}_*B^1_{X/S}\iso\sigma^*_eT^1_{\sS_\nu/S}
  \]
  As in Lemma \ref{lem:tangentspacetomodulispace}, this induces an isomorphism
  \[
    T^1_{\sS_\nu/S}\iso p_\nu^*\R^1\pi^{(p)}_*B^1_{X/S}
  \]
\end{proof}

\begin{Lemma}\label{lem:LemmaNumeroUno}
  If $X\to S$ is a morphism of $\bF_p$-schemes satisfying the assumptions of Proposition \ref{prop:flattoetalecohomology}, then the diagram
  \begin{equation}\label{eq:diagram111}
    \begin{tikzcd}
      \varepsilon_{S*}\R^2\pi^{\fl}_*\m_p\arrow{r}{\Upsilon}\arrow{d}&\R^1\pi^{(p)\etale}_*\nu(1)\arrow{d}\\
      \R^2\pi^{\etale}_*\cO_X^\times\arrow{r}{F}&\R^2\pi^{(p)\etale}_*\cO_{X^{(p)}}^\times
    \end{tikzcd}
  \end{equation}
  commutes, where the left vertical arrow is induced by the inclusion $\m_p\to\cO_X^\times$, the right vertical arrow is the boundary map induced by (\ref{eq:defNu}), and the lower horizontal arrow $F$ is induced by the map $F_{X/S*}\cO_X^\times\to\cO_{X^{(p)}}^\times$ of (\ref{eq:itsasquare}).
\end{Lemma}
\begin{proof}
  This follows from the construction of the exact sequence of low degree terms of a spectral sequence.
\end{proof}

\begin{Remark}\label{rem:RemarkNumeroUno}
  Let us interpret the infinitesimal information contained in the diagram (\ref{eq:diagram111}) in the case when $X\to S=\Spec k$ is a supersingular K3 surface. By Lemma \ref{lem:isooncompletions}, the left vertical arrow induces an isomorphism on completions at the identity. The right vertical arrow fits into the long exact sequence
  \[
    \ldots\to\R^1\pi^{\etale}_*\cO_X^\times\to\R^1\pi^{(p)\etale}_*\nu(1)\to \R^2\pi^{(p)\etale}_*\cO_{X^{(p)}}^\times\xrightarrow{V}\R^2\pi^{\etale}_*\cO_X^\times\to\ldots
  \]
  The formal Brauer group of $X$ has infinite height, which means that the map $V$ induces the zero map on completions at the identity. Using that the Picard group of $X$ is discrete, it follows as in Lemma \ref{lem:isooncompletions} that the right vertical arrow of (\ref{eq:diagram111}) induces an isomorphism on completions. We therefore have a diagram
  \begin{equation}\label{eq:diagram111b}
    \begin{tikzcd}
      \widehat{\R^2\pi^{\fl}_*\m_p}\arrow{r}{\widehat{\Upsilon}}\isor{d}{}&\widehat{\R^1\pi^{(p)\etale}_*\nu(1)}\isor{d}{}\\
      \widehat{\Br(X)}\arrow{r}{0}&\widehat{\Br(X^{(p)})}
    \end{tikzcd}
  \end{equation}
  In particular, $\widehat{\Upsilon}=0$. The corresponding diagram on the tangent spaces to the identity element is
  \begin{equation}\label{eq:diagram111a}
    \begin{tikzcd}
      \H^2(X,\cO_X)\arrow{r}{0}\isor{d}{}&\H^1(X^{(p)},B^1_{X/S})\isor{d}{}\\
      \H^2(X,\cO_X)\arrow{r}{0}&\H^2(X^{(p)},\cO_{X^{(p)}})
    \end{tikzcd}
  \end{equation}
\end{Remark}


\subsection{De Rham cohomology on the relative \'{e}tale site}

Let $\pi\colon X\to S$ be a smooth proper morphism of $\bF_p$-schemes. We will 
discuss the Hodge and conjugate filtrations on the de Rham cohomology of $\pi$ 
and their relationship under the Cartier operator. We will then relate the de 
Rham cohomology to the \'{e}tale cohomology of $\nu(1)$. The material in this section is essentially well 
known, although we have chosen to work throughout on the relative \'{e}tale 
site. This will allow us to cleanly obtain moduli theoretic results later in 
Section \ref{sec:sheaves on the moduli space}. For a thorough treatment of the 
special features of de Rham cohomology in positive characteristic, we refer the 
reader to \cite{Ill96} and \cite{MR0337959}. We remark that the 
results of this section can be seen as a direct generalization of Section 1 of 
\cite{Ogus78} to the case of a non-perfect base.
\begin{Notation}
If $\sE^\bullet$ is any complex of sheaves of abelian groups on a site, we define
\[
  Z^i\sE^\bullet=\ker(\sE^i\too{d}\sE^{i+1})\hspace{1cm}\mbox{ and }\hspace{1cm}B^i\sE^\bullet=\im(\sE^{i-1}\too{d}\sE^i)
\]
\end{Notation}

We consider the algebraic de Rham complex
\[
  \Omega^{\bullet}_{X/S}=\left[0\to\cO_X\too{d}\Omega^{1}_{X/S}\too{d}\Omega^{2}_{X/S}\too{d}\dots\right]
\]
where $\cO_X$ is placed in degree 0. The relation $\frac{d}{dx}(x^p)=0$ implies that the exterior derivative is $\cO_{X^{(p)}}$-linear. This means that $F_{X/S*}\Omega^{\bullet}_{X/S}$ is a complex of $\cO_{X^{(p)}}$-modules, and so the sheaves $F_{X/S*}\Omega^{i}_{X/S}$, $Z^i(F_{X/S*}\Omega_{X/S}^{\bullet})$, $B^i(F_{X/S*}\Omega_{X/S}^{\bullet})$, and $\sH^i(F_{X/S*}\Omega_{X/S}^{\bullet})$ on $X^{(p)}$ have natural $\cO_{X^{(p)}}$-module structures.
\begin{Lemma}\label{lem:illusie}
  The sheaves
  \[
    \Omega^i_{X/S}, F_{X/S*}\Omega^i_{X/S}, Z^i(F_{X/S*}\Omega_{X/S}^\bullet), B^i(F_{X/S*}\Omega_{X/S}^\bullet), \mbox{ and } \sH^i(F_{X/S*}\Omega_{X/S}^\bullet)
  \]
  of modules on the small Zariski sites of $X$ and $X^{(p)}$ are locally free of finite rank.
\end{Lemma}
\begin{proof}
  Because $\pi$ is smooth, $\Omega^i_{X/S}$ is locally free of finite rank. For the remainder, see Corollary 3.6 of \cite{Ill96}.
\end{proof}

We find a similar story on the relative \'{e}tale site. We define sheaves of $\cO_X$-modules on $(X/S)_{\etale}$ by the formula
\begin{equation}\label{eq:extendthesheaves}
  \Omega^{\etale\,i}_{X/S}(U,T)=\Gamma(U,\Omega^i_{U/T})
\end{equation}
Because of the compatibility of the differential with pullback, we obtain a complex
\[
  \Omega^{\etale\,\bullet}_{X/S}=\left[0\to\cO_X\too{d}\Omega^{\etale\,1}_{X/S}\too{d}\Omega^{\etale\,2}_{X/S}\too{d}\dots\right]
\]
of sheaves of abelian groups on $(X/S)_{\etale}$. We will show that the analog of Lemma \ref{lem:illusie} holds on the relative \'{e}tale site.
\begin{Proposition}\label{prop:compatiblewithbasechange}
  The sheaves 
  \[
    \Omega^{\etale\,i}_{X/S}, F_{X/S*}\Omega^{\etale\,i}_{X/S}, Z^i(F_{X/S*}\Omega_{X/S}^{\etale\,\bullet}), B^i(F_{X/S*}\Omega_{X/S}^{\etale\,\bullet}), \mbox{ and } \sH^i(F_{X/S*}\Omega_{X/S}^{\etale\,\bullet})
  \]
  are locally free of finite rank, in the sense of Definition \ref{def:quasi-coherent}.
\end{Proposition}
\begin{proof}
  The content of this claim is that the sheaves of Lemma \ref{lem:illusie} are of formation compatible with base change by \'{e}tale morphisms on the source and arbitrary morphisms on the base. Let us first show that these sheaves are quasi-coherent. If $(U,T)$ is an object of $(X/S)_{\etale}$, then there is a natural identification
  \[
    \Omega^{\etale\,i}_{X/S}|_{U_{\zariski}}=\Omega^i_{U/T}
  \]
  Thus, condition (1) of Definition \ref{def:quasi-coherent} holds. Suppose that $(f,g)\colon (U',T')\to (U,T)$ is a morphism in $(X/S)_{\etale}$. We have a diagram
  \begin{equation}\label{eqn:NiceLookinDiagram}
    \begin{tikzcd}
     U'\arrow{r}{a}\arrow[bend left=30]{rr}{f}\arrow{dr}&U_{T'}\arrow{r}{b}\arrow{d}&U\arrow{d}\\
      &X_{T'}\arrow{r}\arrow{d}&X_T\arrow{d}\\
      &T'\arrow{r}{g}&T
    \end{tikzcd}
  \end{equation}
  where the squares are Cartesian. Note that $a$ is \'{e}tale. Thus, the comparison map factors as a composition of isomorphisms
  \[
    \begin{tikzcd}
      f^*(\Omega^{\etale\,i}_{X/S}|_{U_{\zariski}})\arrow{rr}\arrow[equals]{d}&&\Omega^{\etale\,i}_{X/S}|_{U'_{\zariski}}\arrow[equals]{d}\\
      f^*\Omega^i_{U/T}\arrow{r}{\sim}&a^*\Omega^i_{U_{T'}/T'}\arrow{r}{\sim}&\Omega^i_{U'/T'}
    \end{tikzcd}
  \]
  Condition (2) follows, and therefore $\Omega^{\etale\,i}_{X/S}$ is quasi-coherent. Because $F_{X/S}$ is affine, the pushforward of a quasi-coherent sheaf under $F_{X/S}$ is again quasi-coherent, and therefore $F_{X/S*}\Omega^{\etale\,i}_{X/S}$ is quasi-coherent. Restriction to the small Zariski site is a left exact functor, so the natural map
  \[
    Z^i(F_{X/S*}\Omega^{\etale\,\bullet}_{X/S})|_{U_{\zariski}}\iso Z^i(F_{X/S*}\Omega^{\etale\,\bullet}_{X/S}|_{U_{\zariski}})
  \]
  is an isomorphism. Thus, $Z^i(F_{X/S*}\Omega^{\etale\,\bullet}_{X/S})$ satisfies condition (1). By Lemma \ref{lem:illusie}, the latter sheaf is locally free of finite rank. Because $F_{X/S*}\Omega^{\etale\,i}_{X/S}$ is quasi-coherent, $Z^i(F_{X/S*}\Omega^{\etale\,\bullet}_{X/S})$ therefore satisfies condition (2), and hence is quasi-coherent. The cokernel of a map of quasi-coherent sheaves is quasi-coherent, so $B^i(F_{X/S*}\Omega_{X/S}^{\etale\,\bullet})$ and $\sH^i(F_{X/S*}\Omega_{X/S}^{\etale\,\bullet})$ are quasi-coherent.
  
  By Lemma \ref{lem:illusie}, and the fact that the restriction to the small Zariski site is acyclic for quasi-coherent sheaves, it follows that these sheaves are locally free of finite rank.
\end{proof}

The \textit{Cartier operator} (as described in Section 3 of \cite{Ill96}, and Section 7 of \cite{MR0291177}) defines for each $i$ an isomorphism
\begin{equation}\label{eq:cartieroperator}
    C_{X/S}\colon \sH^i(F_{X/S*}\Omega^\bullet_{X/S})\iso\Omega^i_{X^{(p)}/S}
\end{equation}
of $\cO_{X^{(p)}}$-modules.
\begin{Lemma}\label{lem:cartieroperator}
  The Cartier operator extends to an isomorphism
  \[
    C_{X/S}^{\etale}\colon \sH^i(F_{X/S*}\Omega^{\etale\,\bullet}_{X/S})\iso\Omega^{\etale\,i}_{X^{(p)}/S}
  \]
  of sheaves of $\cO_{X^{(p)}}$-modules on $(X^{(p)}/S)_{\etale}$
\end{Lemma}
\begin{proof}
To define such a map, it will suffice to give morphisms
\[
  \sH^i(F_{X/S*}\Omega^{\etale\,\bullet}_{X/S})|_{U_{\zariski}}\to\Omega^{\etale\,i}_{X^{(p)}/S}|_{U_{\zariski}}\\
\]
for each object $(U,T)$ of $(X^{(p)}/S)_{\etale}$ that behave functorially under pullback. Define $U_{X_T}$ by the Cartesian diagram
\[
  \begin{tikzcd}
    U_{X_T}\arrow{r}\arrow{d}&U\arrow{d}\\
    X_T\arrow{r}{F_{X_T/T}}&(X_T)^{(p/T)}
  \end{tikzcd}
\]
where we are using the natural isomorphism $(X^{(p/S)})_T\iso (X_T)^{(p/T)}$ of Lemma \ref{lem:relativefrobeniusbasechange}. As in Lemma \ref{lem:relativeetalesiteisgood}, we have a diagram
\[
  \begin{tikzcd}
    U_{X_T}            \arrow{r}\arrow{dr}[outer sep=-4pt]{F_{U_{X_T}/T}}      \isor{d}{F_{U_{X_T}/T}} & U                \isol{d}{F_{U/X_T^{(p/T)}}} &\\
    (U_{X_T})^{(p/X_T)}\arrow{r}           \arrow{d}               & (U_{X_T})^{(p/T)}\arrow{d}& \\
    X_T                \arrow{r}{F_{X_T/T}}\arrow{d}               & (X_T)^{(p/T)} \arrow{d}\arrow{r}&T\arrow{d}           \\
        X\arrow{r}{F_{X/S}}&X^{(p/S)}\arrow{r}&S
  \end{tikzcd}
\]
with Cartesian squares. Consider the composite arrow $f\colon U\to X^{(p/S)}$. 
We define a map $C^{\etale}_{X/S}|_{U_{\zariski}}$ by the diagram
\[
  \begin{tikzcd}
    F_{U_{X_T}/T}^* \sH^i(F_{U_{X_T}/T*}\Omega^{\bullet}_{U_{X_T}/T})\isor{d}{}\arrow{r}{C_{U_{X_T}/T}}&F_{U_{X_T}/T}^*\Omega^i_{U_{X_T}^{(p/T)}/T}\isor{d}{}\\
    f^*(\sH^i(F_{X/S*}\Omega^{\bullet}_{X/S}))\isor{d}{}      &f^*(\Omega^{i}_{X^{(p)}/S})\isor{d}{}\\
    f^*(\sH^i(F_{X/S*}\Omega^{\etale\,\bullet}_{X/S})|_{X^{(p/S)}_{\zariski}})\isor{d}{}      &f^*(\Omega^{\etale\,i}_{X^{(p)}/S}|_{X^{(p/S)}_{\zariski}})\isor{d}{}\\
    \sH^i(F_{X/S*}\Omega^{\etale\,\bullet}_{X/S})|_{U_{\zariski}}\arrow{r}{C_{X/S}^{\etale}|_{U_{\zariski}}} &\Omega^{\etale\,i}_{X^{(p)}/S}|_{U_{\zariski}}
  \end{tikzcd}
\]
where the top horizontal arrow is the pullback of $C_{U_{X_T}/T}$ under the isomorphism $F_{U_{X_T}/T}$, and the two lower vertical arrows are the comparison morphisms, which are isomorphisms by Proposition \ref{prop:compatiblewithbasechange}. The explicit description of the Cartier operator in Section 7 of \cite{MR0291177} shows that the morphisms $C^{\etale}_{X/S}|_{U_{\zariski}}$ are compatible with pullback in the appropriate sense, and so we obtain the desired isomorphism $C_{X/S}$ of sheaves on the relative \'{e}tale site.

\end{proof}

We consider the map 
\begin{equation}
  d\log\colon \cO_X^\times\to\Omega^{\etale\,1}_{X/S}
\end{equation}
of sheaves on $(X/S)_{\etale}$ given by $f\mapsto df/f$. Note that the image of $d\log$ is contained in the subsheaf of closed forms. We also have a map 
\[
  Z^1(F_{X/S*}\Omega^{\etale\,\bullet}_{X/S})\too{W^*} F_{X/S*}W_*\Omega^{\etale\,1}_{X^{(p)}/S}=F_{X^{(p)}*}\Omega^{\etale\,1}_{X^{(p)}/S}\to\Omega^{\etale\,1}_{X^{(p)}/S}
\]
where the first map is the pushforward of the map $W^*$ under $F_{X/S}$, and the second is the canonical map.

\begin{Lemma}\label{lem:4term}
  The sequence
  \[
    0\to\cO_{X^{(p)}}^\times\to F_{X/S*}\cO_X^\times\too{d\log}Z^1(F_{X/S*}\Omega^{\etale\,\bullet}_{X/S})\xrightarrow{C_{X/S}^{\etale}-W^*}\Omega^{\etale\,1}_{X^{(p)}/S}\to 0
  \]
  of sheaves on $(X^{(p)}/S)_{\etale}$ is exact.
\end{Lemma}
\begin{proof}
  Corollaire 2.1.18 of \cite{Ill79} states that the restriction of this sequence to the small \'{e}tale site is exact. The result follows by definition of the sheaves involved.
\end{proof}

We next discuss the Hodge and conjugate spectral sequences in the setting of 
the relative \'{e}tale site, which converge to the relative de Rham cohomology
\[
  \H^m_{dR}(X/S)_{\etale}\defeq\R^m\pi^{\etale}_*\Omega^{\etale\,\bullet}_{X/S}
\]
of $X\to S$. The so-called na\"{i}ve filtration on $\Omega_{X/S}^{\etale\bullet}$ (see \cite{Ill96}, and Chapter 2 of \cite{Huy06}) induces the Hodge spectral sequence
\begin{equation}\label{eq:HodgeSS}
  E_{H1}^{p,q}=\R^q\pi^{\etale}_*\Omega^{\etale\,p}_{X/S}\Longrightarrow \R^{p+q}\pi^{\etale}_*\Omega^{\etale\,\bullet}_{X/S}
\end{equation}
On the other hand, we may compute the higher pushforwards of the de Rham 
complex using the Leray spectral sequence induced by the canonical 
factorization $\pi=\pi^{(p)}\circ F_{X/S}$. By Lemma \ref{lem:heythatscool} the 
functor $F_{X/S*}$ is exact on the relative \'{e}tale site, so this gives a spectral sequence
\begin{equation}\label{eq:ConjugateSS}
  E_{C2}^{p,q}=\R^p\pi^{(p)\etale}_*\sH^{q}(F_{X/S*}\Omega^{\etale\,\bullet}_{X/S})\Longrightarrow \R^{p+q}\pi^{\etale}_*\Omega^{\etale\,\bullet}_{X/S}
\end{equation}
called the \textit{conjugate spectral sequence}.
Fix an integer $m\geq 0$. The Hodge and conjugate spectral sequences induce filtrations
\begin{align}\label{eq:filtrations}
  0\subset F_H^{m,m}\subset F_H^{m-1,m} \subset\dots\subset  F_H^{i,m}\subset\dots\subset F_H^{0,m}=&\H^m_{dR}(X/S)_{\etale}\\
  0\subset F_C^{m,m}\subset F_C^{m-1,m} \subset\dots\subset  F_C^{i,m}\subset\dots\subset F_C^{0,m}=&\H^m_{dR}(X/S)_{\etale}
\end{align}
\begin{Remark}
  Working on the relative \'{e}tale site gives a functorial packaging of the usual Hodge and conjugate spectral sequences on the small Zariski site. If $T$ is a scheme over $S$, then the restriction of the spectral sequences (\ref{eq:HodgeSS},\ref{eq:ConjugateSS}) to the small Zariski site of $T$ is naturally isomorphic to the usual Hodge and conjugate spectral sequences induced by the complex $\Omega^{\bullet}_{X_T/T}$. In particular, there is a natural identification
  \[
    \R^q\pi^{\etale}_*\Omega^{\etale\,p}_{X/S}(T)= \Gamma(T,\R^q\pi_{T*}\Omega^p_{X_T/T})
  \]
  and similarly for $E^{p,q}_{C2}$, $\H^m_{dR}(X/S)_{\etale}$, and the filtrations $F^{i,m}_H$ and $F^{i,m}_C$.
\end{Remark}

We will make the following assumptions on $\pi$.
\begin{Definition}
  We say that the morphism $\pi$ satisfies $(*)$ if
  \begin{enumerate}
    \item the Hodge spectral sequence degenerates at $E_1$,
    \item the conjugate spectral sequence degenerates at $E_2$, and
    \item the $\cO_S$-modules $E^{p,q}_{H1}$, $E^{p,q}_{C2}$, $\H^m_{dR}(X/S)_{\etale}$, $F^{i,m}_H$ and $F^{i,m}_C$ are all locally free of finite rank, in the sense of Definition \ref{def:quasi-coherent}.
  \end{enumerate}
\end{Definition}
\begin{Remark}
  If $\pi\colon X\to S$ is a relative K3 surface, then $\pi$ satisfies $(*)$. This follows from the fact that the spectral sequences associated to the usual (small Zariski) De Rham cohomology degenerate, and are of formation compatible with base change.
\end{Remark}
Under these assumptions, the Hodge and conjugate spectral sequences are related by the Cartier operator.

\begin{Lemma}
  If $\pi\colon X\to S$ satisfies $(*)$, then there are natural isomorphisms
  \[
    F_H^{i,m}/F_H^{i+1,m}\iso E_{H1}^{i,m-i}=\R^{m-i}\pi^{\etale}_*(\Omega^{\etale\,i}_{X/S})
  \]
  \[
    F_C^{i,m}/F_C^{i+1,m}\iso E_{C2}^{i,m-i}=\R^{i}\pi^{(p)\etale}_*(\sH^{m-i}(F_{X/S*}\Omega^{\etale\,\bullet}_{X/S}))
  \]
  and the Cartier operator induces isomorphisms
  \[
    F_C^{m-i,m}/F_C^{m-i+1,m}\iso F_S^* (F_H^{i,m}/F_H^{i+1,m})
  \]
\end{Lemma}
\begin{proof}
  The first claim follows from the assumed degeneration of the Hodge and conjugate spectral sequences. For the second, we apply $\R\pi^{(p)\etale}_*$ to both sides of (\ref{eq:cartieroperator}) to get an isomorphism
  \[
    \R\pi^{(p)\etale}_*\sH^i(F_{X/S*}\Omega^{\etale\,\bullet}_{X/S})\iso\R\pi^{(p)\etale}_*\Omega^{\etale\,i}_{X^{(p)}/S}
  \]
  Applying the base change isomorphism in the derived category and using our assumption that the cohomology sheaves of $\R\pi^{\etale}_*\Omega^{\etale\,i}_{X/S}$ are locally free, we find isomorphisms
  \[
    \R\pi^{(p)\etale}_*\sH^i(F_{X/S*}\Omega^{\etale\,\bullet}_{X/S})\cong\R\pi^{(p)\etale}_*\Omega^{\etale\,i}_{X^{(p)}/S}\cong\R\pi^{(p)\etale}_*W^*\Omega^{\etale\,i}_{X/S}\cong F_S^*\R\pi^{\etale}_*\Omega^{\etale\,i}_{X/S}
  \]
  Applying the isomorphisms of the previous claim gives the result.
\end{proof}

For the remainder of this paper, we only consider the case $m=2$, and so we omit $m$ from the notation. We have sheaves $F^1_H$ and $F^1_C$ on the big \'{e}tale site of $S$ that are locally free of finite rank, and in particular quasi-coherent. By $F^1_H\cap F^1_C$ we will mean the fiber product of sheaves on the big \'{e}tale site of $S$. This will not in general be a quasi-coherent sheaf.
  
\begin{Lemma}\label{lem:closeddifferentialforms}
  If $X\to S$ satisfies $(*)$, then the natural map $\R^1\pi^{\etale}_*Z^1\Omega^\bullet_{X/S}\to \R^2\pi^{\etale}\Omega^\bullet_{X/S}$ induces an identification
  \[
    \R^1\pi^{\etale}_*Z^1\Omega^\bullet_{X/S}\iso F^1_H\cap F^1_C
  \]
\end{Lemma}
\begin{proof}
  We have truncations
  \[
    \Omega^{\etale\,\geq i}_{X/S}=\left[\dots \to 0\to\Omega^{\etale\,i}_{X/S}\to\Omega^{\etale\,i+1}_{X/S}\to\dots\right]
  \]
  where $\Omega^{\etale\,i}_{X/S}$ is placed in degree $i$. There is an obvious map of complexes $\Omega^{\etale\,\geq i}_{X/S}\to \Omega_{X/S}^{\etale\,\bullet}$, and the image of the induced map $\R^2\pi^{\etale}_*\Omega^{\etale\,\geq i}_{X/S}\to \R^2\pi^{\etale}_*\Omega^{\etale\,\bullet}_{X/S}$ is the $i$-th level of the Hodge filtration. There is an obvious map $Z^1\Omega^{\etale\,\bullet}_{X/S}\left[-1\right]\to\Omega^{\etale\,\geq 1}_{X/S}$, which induces an exact sequence
  \[
    0\to Z^1\Omega^{\etale\,\bullet}_{X/S}\left[-1\right]\to\Omega^{\etale\,\geq 1}_{X/S}\to Q^\bullet\to 0
  \]
  Therefore, we get a triangle
  \[
    \R\pi^{\etale}_*Z^1\Omega^{\etale\,\bullet}_{X/S}\left[-1\right]\to \R\pi^{\etale}_*\Omega^{\etale\,\geq 1}_{X/S}\to \R\pi^{\etale}_*Q^\bullet\to
  \]
  Note that $\sH^i(Q^\bullet)=0$ for $i\leq 1$, so by the spectral sequence associated to the canonical filtration on $Q^\bullet$ we get that $\R^1\pi^{\etale}_*Q^\bullet=0$ and $\R^2\pi^{\etale}_*Q^\bullet\iso \R^0\pi^{\etale}_*\sH^2(\Omega_{X/S}^{\etale\,\bullet})$. Therefore, we get an exact sequence
  \[
    0\to \R^1\pi^{\etale}_*Z^1\Omega^{\etale\,\bullet}_{X/S}\to F^1_H\to F^0_C/F^1_C
  \]
  This gives the result.
\end{proof}

\begin{Proposition}\label{prop:ogusexactsequence}
  If $X\to S$ satisfies $(*)$ and $\pi^{\etale}_*\Omega^{\etale\,1}_{X/S}=0$, then there is an exact sequence
  \[
    0\to \R^1\pi^{(p)\etale}_*\nu(1)\to F^1_H\cap F^1_C\xrightarrow{C\circ\pi_C-F^*_S\circ\pi_H}F_S^*(F^1_H/F^2_H)
  \]
  of sheaves of abelian groups on the big \'{e}tale site of $S$.
\end{Proposition}
\begin{proof}
  Here, the arrow $C\circ\pi_C$ is the projection $F^1_C\to F^1_C/F^2_C$ followed by the isomorphism $C^{\etale}_{X/S}\colon F^1_C/F^2_C\iso F_S^*(F^1_H/F^2_H)$, $\pi_H$ is the projection $F^1_H\to F^1_H/F^2_H$, and $F^*_S$ is the canonical $p$-linear map $F^1_H/F^2_H\to F_S^*(F^1_H/F^2_H)$ given by $s\mapsto s\otimes 1$.
  
  We apply $\R\pi^{(p)\etale}_*$ to the sequence
  \[
    0\to \nu(1)\to F_{X/S*}Z^1\Omega_{X/S}^{\etale\,\bullet}\to\Omega^{\etale\,1}_{X^{(p)}/S}\to 0
  \]
  of Lemma \ref{lem:4term} to get
  \[
    0\to \R^1\pi^{(p)\etale}_*\nu(1)\to \R^1\pi^{\etale}_*Z^1\Omega^{\etale\,\bullet}_{X/S}\to F^*_S(\R^1\pi^{\etale}_*\Omega^{\etale\,1}_{X/S})
  \]
  Applying Lemma \ref{lem:closeddifferentialforms}, we get a diagram
  \[
    \begin{tikzcd}
      0\arrow{r}& \R^1\pi^{(p)\etale}_*\nu(1)\arrow{r}\arrow[equal]{d}&\R^1\pi^{\etale}_*Z^1\Omega^{\etale\,\bullet}_{X/S}\arrow{r}{C^{\etale}_{X/S}-W^*}\isor{d}{}&F_S^*(\R^1\pi^{\etale}_*\Omega^{\etale\,1}_{X/S})\isor{d}{}\\
      0\arrow{r}&\R^1\pi^{(p)\etale}_*\nu(1)\arrow{r}&F^1_H\cap F^1_C\arrow{r}{C\circ\pi_C-F^*_S\circ\pi_H}&F_S^*(F^1_H/F^2_H)
    \end{tikzcd}
  \]
  which gives the result.
\end{proof}

\section{Periods}\label{sec:periods}

\subsection{Characteristic subspaces and the crystalline period domain}\label{sec:perioddomain}

In this section we consider the crystalline analog of the period domain, Ogus's moduli space of characteristic subspaces. We study certain rational fibrations relating the period spaces of different Artin invariants. We begin by recalling some definitions and results having to do with bilinear forms on vector spaces over $\bF_p$. Recall that we are assuming $p\geq 3$. 
\begin{Definition}\label{def:hyperbolicplane}
  Let $U_2\otimes\bF_p$ denote the \textit{hyperbolic plane over $\bF_p$}, which is the $\bF_p$-vector space generated by the two vectors $v,w$ which satisfy $v^2=w^2=0$ and $v.w=-1$.  
\end{Definition}
This notation is explained by Definition \ref{def:hyperbolicplane2}.
\begin{Proposition}\cite[Theorem 4.9]{MR1859189}\label{prop:classification}
  Let $V$ be a vector space over $\bF_p$ of even dimension $2\sigma_0$, equipped with a non-degenerate $\bF_p$-valued bilinear form. If there exists a totally isotropic subspace of $V$ of dimension $\sigma_0$, then $V$ is isometric to the orthogonal sum of $\sigma_0$ copies of $U_2\otimes\bF_p$. If there does not exist such a subspace, then $V$ is isometric to the orthogonal sum of $\sigma_0-1$ copies of $U_2\otimes\bF_p$ and one copy of $\bF_{p^2}$.
\end{Proposition}
Here, we view $\bF_{p^2}$ as a two dimensional vector space over $\bF_p$ equipped with the quadratic form coming from the trace. In the former case, we say that the form on $V$ is \textit{neutral}, and in the latter we say that it is \textit{non-neutral}.
\begin{Theorem}\cite[Theorem 5.2]{MR1859189}\label{thm:Witt}
  If $V$ is a vector space over a field of characteristic not equal to 2 equipped with a non-degenerate bilinear form, then any isometry $W\to W'$ between two subspaces of $V$ extends to an isometry $V\to V$.
\end{Theorem}
In particular, this implies that the group of isometries of $V$ acts transitively on the set of totally isotropic subspaces of any given dimension. 

Let $V$ be a vector space over $\bF_p$ equipped with a non-degenerate, non-neutral bilinear form. Let $S$ be a scheme over $\bF_p$. By viewing $V$ as a coherent sheaf on $\Spec\bF_p$, we see that there is a canonical isomorphism
\[
  F^*_S(V\otimes\cO_S)\iso V\otimes\cO_S
\]
Precomposing with the canonical $p$-linear map $V\otimes\cO_S\to F_S^*(V\otimes\cO_S)$, we obtain a map of sheaves
\[
  F^*_S\colon V\otimes\cO_S\to V\otimes\cO_S
\]
This map is given on sections by $v\otimes s\mapsto v\otimes s^p$. When $S=\Spec k$, we will denote this map by $\varphi$.
\begin{Definition}\cite[Section 4]{Ogus78}\label{def:charsub}
  A \textit{characteristic subspace} of $V\otimes\cO_S$ is a submodule $K\subset V\otimes\cO_S$ such that
  \begin{enumerate}
    \item both $K$ and $K+F^*_SK$ are locally free and locally direct summands of $V\otimes\cO_S$,
    \item $K$ is totally isotropic of rank $\sigma_0$, and
    \item $K+F_S^*K$ has rank $\sigma_0+1$.
  \end{enumerate}
  We define a functor $M_V$ on schemes over $\bF_p$ by letting $M_V(S)$ be the set of characteristic subspaces of $V\otimes\cO_S$.   
\end{Definition}

The scheme $M_V$ is studied in Section 4 of \cite{Ogus78}, where it is shown to be smooth and projective over $\bF_p$ of dimension $\sigma_0-1$. Moreover, $M_V$ is irreducible, and has no $\bF_p$-points. The base change $M_V\otimes_{\bF_p}\overline{\bF}_p$ has two irreducible connected components, which are defined over $\bF_{p^2}$ and are interchanged by the action of the Galois group.
\begin{Remark}\label{rem:smallArtininvariant}
  For small values of $\sigma_0$, the period space $M_V$ admits the following descriptions.
  \begin{itemize}
  \item If $\sigma_0=1$, $M_V\cong\Spec\bF_{p^2}$. Hence,
    $M_V\otimes_{\bF_p}\overline{\bF}_p$ is a disjoint union of two
    copies of $\Spec\overline{\bF}_p$.
  \item If $\sigma_0=2$, $M_V\cong\bP^1_{\bF_{p^2}}$. Hence,
    $M_V\otimes_{\bF_p}\overline{\bF}_p$ is a disjoint union of two
    copies of $\bP^1_{\overline{\bF}_p}$.
      \item If $\sigma_0=3$, $M_V$ is isomorphic to the Fermat surface
      \[
        V(x^{p+1}+y^{p+1}+z^{p+1}+w^{p+1})\subset\bP^3_{\bF_{p^2}}
      \]
      Hence, $M_V\otimes_{\bF_p}\overline{\bF}_p$ is a disjoint union
      of two copies of the Fermat surface of degree $p+1$.
  \end{itemize}
  The first two are described in \cite[Example 4.7]{Ogus78}, and the last case is
  proven in \cite[Proposition 5.14]{Ekedahl}. We note that
  \cite[Remark 4.4]{Liedtke15} erroneously quotes
  \cite[Example 4.7]{Ogus78} in claiming that
  $M_V\cong\P_{\F_{p^2}}^1\times\P_{\F_{p^2}}^1$ when
  $\sigma_0=3$. The reader is referred to Remark \ref{rem:!!}ff for
  further discussion of the global geometry of $M_V$.
\end{Remark}

By definition, $M_V$ comes with a tautological sub-bundle
$K_V\subset V\otimes\cO_{M_V}$, and the pairing induces a diagram
\begin{equation}\label{eq:duality}
  \begin{tikzcd}
    0\arrow{r}&K_V\arrow{r}\isor{d}{}&V\otimes\cO_{M_V}\arrow{r}\isor{d}{}&Q_V\arrow{r}\isor{d}{}&0\\
    0\arrow{r}&Q_V^*\arrow{r}&V^*\otimes\cO_{M_V}\arrow{r}&K_V^*\arrow{r}&0
  \end{tikzcd}
\end{equation}

\begin{Definition}\label{def:strictlycharacteristic}
  The \textit{Artin invariant} of a characteristic subspace $K\subset V\otimes k$ is 
  \[
    \sigma_0(K)=\sigma_0-\dim_{\bF_p}(K\cap V)
  \]
  The subspace $K$ is \textit{strictly characteristic} if $\sigma_0(K)=\sigma_0$.
\end{Definition}

\begin{Definition}\label{def:MVW}
  Let $W\subset V$ be a subspace. We consider the Cartesian diagram
  \[
    \begin{tikzcd}
      \underline{K}_V\cap \underline{W}\arrow{r}\arrow{d}&\underline{W}\arrow{d}\\
      \underline{K}_V\arrow{r}&\underline{V}\otimes\cO_{M_V}
    \end{tikzcd}
  \]
  of schemes over $M_V$, where $\underline{K}_V$ is the vector bundle on $M_V$ associated to $K_V$, and $\underline{W}$ and $\underline{V}$ are the constant group schemes associated to $W$ and $V$. The zero section $Z\subset \underline{K}_V\cap \underline{W}$ is open and closed, and so the morphism $\underline{K}_V\cap\underline{W}\setminus Z\to M_V$ is closed. We let $M^W_V\subset M_V$ be the open complement of the image of this morphism. A $k$-point $[K]\in M_V(k)$ lies in $M^W_V$ exactly when the intersection of $K$ and $W$ inside of $V\otimes k$ is trivial. In particular, we set $U_V=M^V_V$. This is the locus of points with maximal Artin invariant. 
  
  We let $M_{V,W}\subset M_V$ denote the closed complement of the union of the $M^{\langle w\rangle}_V$ as $w$ ranges over the non-zero vectors in $W$. We give $M_{V,W}$ the reduced induced closed subscheme structure. A $k$-point $[K]\in M_V(k)$ lies in $M_{V,W}$ exactly when $W\subset K$.
\end{Definition}

\begin{Notation}\label{notation:vectorspace} For the remainder of this section,
  the following notation will be consistently used.
  \begin{itemize}
      \item We let $\tV$ be a fixed vector space over $\bF_p$ of dimension $2\sigma_0+2$ equipped with a non-degenerate non-neutral bilinear form.
      \item We let $v\in\tV$ be a vector that is \textit{isotropic}, meaning that $v\neq 0$ and $v^2=0$.
      \item We set $V=v^\perp/v$. This is a vector space of dimension $2\sigma_0$, and is equipped with a natural bilinear form which is again non-degenerate and non-neutral.
      \item We let $\tV=V\oplus (U_2\otimes\bF_p)$ be an orthogonal decomposition, where $U_2\otimes\bF_p$ is defined in Definition \ref{def:hyperbolicplane} (such a decomposition exists for any $v$ by Theorem \ref{thm:Witt}).
  \end{itemize}
\end{Notation}

\begin{Lemma}\label{lem:definitionMorphism}
  Let $S$ be a scheme over $\bF_p$ and let $\tK\subset \tV\otimes\cO_S$ be a characteristic subspace. If for every geometric point $\Spec k\to S$ the fiber $\tK\otimes k\subset \tV\otimes k$ does not contain $v(=v\otimes 1)$, then the image $K$ of $\tK\cap (v^\perp\otimes\cO_{S})$ in $V\otimes\cO_{S}$ is a characteristic subspace, whose formation is compatible with arbitrary base change.
\end{Lemma}
\begin{proof}
  Every geometric fiber of $K$ is a maximal totally isotropic subspace of $\tV\otimes k$. Therefore by Nakayama's Lemma the map of sheaves $\tK\to\cO_{S}$ given by pairing with $v$ is surjective, and we have a short exact sequence
  \[
    0\to \tK\cap(v^\perp\otimes\cO_{S})\to\tK\to\cO_{S}\to 0
  \]
  This shows that $\tK\cap(v^\perp\otimes\cO_{S})$ is a locally free quasi-coherent sheaf of rank $\sigma_0$, whose formation is compatible with arbitrary base change. The quotient map $v^\perp\to v^\perp/v=V$ induces an isomorphism
  \[
    \tK\cap (v^\perp\otimes\cO_{S})\iso (\tK\cap v^\perp\otimes \cO_{S})/(v\otimes\cO_{S})=K
  \]
  Thus, $K$ is locally free of rank $\sigma_0$, is locally a direct summand, and is compatible with arbitrary base change. It is clear that $K$ is totally isotropic. The same reasoning applied to $\tK + F^*\tK$ shows that $(\tK + F^*\tK)\cap (v^\perp\otimes\cO_{S})$ and its image under the quotient map are direct summands of rank $\sigma_0+1$. 
\end{proof}

Note that this condition on the geometric fibers of $\tK$ is equivalent to the induced map $S\to M_{\tV}$ factoring through $M^{\langle v\rangle}_{\tV}$. Because of the compatibility with respect to base change, the association $\tK\mapsto K$ defines a morphism
\begin{equation}
  \pi_v\colon M^{\langle v\rangle}_{\tV}\to M_V.
\end{equation} 

\begin{Definition}
  We define a group scheme $\sU_V$ on $M_V$ by the kernel of the map
  \[
    \frac{\underline{V}\otimes\cO_{M_V}}{\underline{K}_V}\xrightarrow{1-F^*_{M_V}}\frac{\underline{V}\otimes\cO_{M_V}}{\underline{K}_V+F^*_{M_V}\underline{K}_V}
  \]
  of group schemes over $M_V$.
\end{Definition}

\begin{Proposition}\label{prop:isoofgroupschemes}
  There is an isomorphism $\sU_{V}\iso M^{\langle v\rangle}_{\tV}$ of schemes over $M_{V}$.
\end{Proposition}
\begin{proof}
    We will first define a morphism $f\colon \sU_{V}\to M^{\langle v\rangle}_{\tV}$. That is, we will give a characteristic subspace of $\tV\otimes\cO_{\sU_{\tV}}$ whose intersection with $\langle v\rangle\otimes\cO_{M^{\langle v\rangle}_{\tV}}$ is trivial. By definition, we have an exact sequence
  \[
    0\to\sU_V\to\frac{\underline{V}\otimes\cO_{M_V}}{\underline{K}_V}\xrightarrow{1-F^*_{M_V}}\frac{\underline{V}\otimes\cO_{M_V}}{\underline{K}_V+F^*_{M_V}\underline{K}_V}
  \]
  of group schemes over $M_V$. On $\sU_V$, there is a tautological section
  \[
    \sB_V\in\Gamma\left(\sU_V,\frac{V\otimes\cO_{M_V}}{K_V}\bigg\rvert_{\sU_V}\right)
  \]
  satisfying $(1-F^*)\sB_V\in K_V+F^*K_V$. Consider the sub-bundle $\tK\subset\tV\otimes\cO_{\sU_{V}}$ defined locally as
  \[
    \tK=\left\langle x_1+(x_1\innerproduct B)v,\dots,x_{\sigma_0}+(x_{\sigma_0}\innerproduct B)v, w+B+\frac{B^2}{2}v\right\rangle\subset\tV\otimes\cO_{\sU_{V}}
  \]
  where $x_1,\dots,x_{\sigma_0}$ is a local basis for $K_V|_{\sU_V}$, and $B$ denotes a local lift of the tautological section $\sB_V$ to $V\otimes\cO_{\sU_V}$. Note that as $K_V$ is totally isotropic, this does not depend on our choice of $B$. By construction, $\tK$ is locally free of rank $\sigma_0+1$, is locally a direct summand, and is totally isotropic. Using the defining property of $\sB_V$, and that $K_V$ is characteristic, it follows that $\tK$ is characteristic as well. Thus, we get a morphism $f\colon \sU_V\to M_{\tV}^{\langle v\rangle}$, which is compatible with the respective morphisms to $M_V$.
  
  Let us construct an inverse. We seek a global section of $\sU_V|_{M_{\tV}^{\langle v\rangle}}$. Pairing with -$v$ gives a short exact sequence
  \[
    0\to K_{\tV}|_{M_{\tV}^{\langle v\rangle}}\cap (v^\perp\otimes\cO_{M_{\tV}^{\langle v\rangle}})\to K_{\tV}|_{M_{\tV}^{\langle v \rangle}}\xrightarrow{\langle\_,\text{-}v\rangle}\cO_{M_{\tV}^{\langle v \rangle}}\to 0
  \]
  and the projection $v^\perp\to V$ induces an isomorphism
  \[
    K_{\tV}|_{M_{\tV}^{\langle v\rangle}}\cap (v^\perp\otimes\cO_{M_{\tV}^{\langle v\rangle}})\iso\pi_v^* K_V
  \]
  Thus, letting $\rho\colon \tV\to V$ denote the projection, we have maps
  \[
    \cO_{M_{\tV}^{\langle v \rangle}}\iso\frac{K_{\tV}|_{M_{\tV}^{\langle v \rangle}}}{K_{\tV}|_{M_{\tV}^{\langle v\rangle}}\cap (v^\perp\otimes\cO_{M_{\tV}^{\langle v\rangle}})}\too{\rho}\frac{\rho(K_{\tV}|_{M_{\tV}^{\langle v \rangle}})}{\pi_v^* K_V}\subset \frac{V\otimes\cO_{M_{\tV}^{\langle v\rangle}}}{\pi_v^* K_V}
  \]
  The image of $1\in\Gamma(M_{\tV}^{\langle v \rangle},\cO_{M_{\tV}^{\langle v \rangle}})$ under the left isomorphism can be lifted locally to a section of $K_{\tV}$ of the form $w+b+\frac{b^2}{2}v$ where $b$ is a section of $V\otimes\cO$, and the above composition sends $1$ to a global section which is locally given by the image of $b$ in the quotient. Notice that
  \[
    (1-F^*)\left(w+b+\frac{b^2}{2}v\right)\in (K_{\tV}+F^*K_{\tV})\cap (v^\perp\otimes\cO_{M_{\tV}})
  \]
  Because $v$ is fixed by the Frobenius, the projection $v^\perp\to V$ also induces an isomorphism 
  \[
    (K_{\tV}+F^*K_{\tV})\cap (v^\perp\otimes\cO_{M_{\tV}})\iso \pi_v^*(K_V+F^*K_V)
  \]
  Thus, the image of our global section under the map
  \[
    \frac{V\otimes\cO_{M_{\tV}^{\langle v\rangle}}}{\pi_v^* K_V}\xrightarrow{1-F^*}\frac{V\otimes\cO_{M_{\tV}^{\langle v\rangle}}}{\pi_v^*(K_V+F^*K_V)}
  \]
  is zero, and hence gives a global section of $\sU_V|_{M_{\tV}^{\langle v\rangle}}$. We have constructed a morphism $g\colon M_{\tV}^{\langle v\rangle}\to\sU_V$ over $M_V$. One checks that this map is an inverse to $f$, and this gives the result.
\end{proof}
\begin{Remark}
  Via crystalline cohomology, the period domain $M_{\tV}$ classifies supersingular K3 surfaces (in a sense that is made precise by Ogus's crystalline Torelli theorem). In the following sections, we will interpret $\sU_V$ as the connected component of the relative second cohomology of $\m_p$ on the universal K3 surface. That is, $\sU_V$ classifies (certain) twisted supersingular K3 surfaces. The isomorphism of Proposition \ref{prop:isoofgroupschemes} suggests that there should be a relationship between the collection of supersingular K3 surfaces of Artin invariant $\sigma_0+1$ and the collection of twisted supersingular K3 surfaces whose coarse space has Artin invariant $\sigma_0$. A central goal of this paper is to find a geometric interpretation of this isomorphism.
\end{Remark}
\begin{Lemma}\label{lem:smooth}
The morphisms $\sU_V\to M_V$ and $\pi_v\colon M^{\langle v\rangle}_{\tV}\to M_V$ are smooth.
\end{Lemma}
\begin{proof}
  We will show that $\sU_V\to M_V$ is smooth. Because the differential of the Frobenius vanishes, the morphism
  \[
    1-F^*_{M_V}\colon\frac{\underline{V}\otimes\cO_{M_V}}{\underline{K}_V}\to\frac{\underline{V}\otimes\cO_{M_V}}{\underline{K}_V+F^*_{M_V}\underline{K}_V}
  \]
  is smooth. But $\sU_V\to M_V$ is the pullback of this morphism along the zero section. By Proposition \ref{prop:isoofgroupschemes}, $\pi_v$ is smooth as well.
\end{proof}

By definition, $M_V$ comes equipped with two natural line bundles $K_V/K_V\cap F^*K_V$ and $F^*K_V/K_V\cap F^*K_V$. The pairing on $V$ induces an isomorphism
\[
  \dfrac{K_V}{K_V\cap F^*K_V}\otimes\dfrac{F^*K_V}{K_V\cap F^*K_V}\iso\cO_{M_V}
\]
and so these line bundles are naturally dual. Consider the open subset $U_V\subset M_V$ parameterizing strictly characteristic subspaces. Write $U=U_V$, $K_U=K_V|_{U}$, and $\sU_U=\sU_V|_U$. As discussed in Section 4 of \cite{Ogus78}, the subsheaf
\[
  \sL_U=K_U\cap F_{U}^*K_U\cap\ldots\cap F^{*\sigma_0-1}_UK_U
\]
is locally free of rank 1, and we have
\begin{align}\label{eq:miscequation2}
  \sL_U+F^*\sL_U+\cdots+F^{*\sigma_0-1}\sL_U&=F^{*\sigma_0-1}K_U\\
  \sL_U+F^*\sL_U+\cdots+F^{*2\sigma_0-1}\sL_U&=V\otimes\cO_M
\end{align}
Furthermore, the natural map
\[
  \sL_U\to F^{*\sigma_0-1}\left(\dfrac{K_V}{K_V\cap F^*K_V}\bigg\rvert_{U}\right)
\]
is an isomorphism, and the pairing on $V$ induces an isomorphism
\[
  \sL_U\otimes F^{*\sigma_0}\sL_U\iso\cO_U
\]
Consider the composition
\begin{equation}\label{eq:miscequation1}
  F^{*\sigma_0-1}\sU_U\subset\dfrac{\underline{V}\otimes\cO_U}{F^{*\sigma_0-1}\underline{K}_U}\iso (F^{*\sigma_0-1}\underline{K}_U)^\vee\to\underline{\sL}_U^\vee
\end{equation}
\begin{Lemma}\label{lem:bundlestructure}
  The map $F^{*\sigma_0-1}\sU_U\to\underline{\sL}_U^\vee$
  (\ref{eq:miscequation1}) is an isomorphism of group schemes.
\end{Lemma}
\begin{proof}
  We claim that it suffices to show that (\ref{eq:miscequation1}) is
  injective (as a map of group schemes). Indeed, we may check
  isomorphy on geometric fibers. By Lemma \ref{lem:smooth},
  $F^{*\sigma_0-1}\sU_U$ is a smooth, $p$-torsion group scheme of
  relative dimension 1. Hence, its geometric fibers are abstractly
  isomorphic to a disjoint union of copies of $\bA^1$. Of course, the
  geometric fibers of $\underline{\sL}_U^\vee$ are also isomorphic to
  $\bA^1$. The claim follows from the fact that any monomorphism
  $\bA^1\to\bA^1$ is an isomorphism.
  
  Suppose that $S$ is a scheme over $M_V$, and consider a section in
  $F^{*\sigma_0-1}\sU_U(S)$ whose image is equal to 0. Passing to an
  open cover of $S$, we may assume that $\sL_U$ is trivial, generated
  by some section $e\in\Gamma(S,\sL_U|_S)$, and that our section lifts
  to some element $b\in\Gamma(S,V\otimes\cO_S)$. We have that
  $b.e=0$. By assumption,
  \[
    b-F^*(b)\in \Gamma(S,F^{*\sigma_0-1}K+F^{*\sigma_0}K)
  \]
  and therefore $b-F^*(b)$ is perpendicular to
  $F^{*\sigma_0-1}K\cap F^{*\sigma_0}K$. It follows that
  \[
    b.F^*(e)=(b-F^*(b)).F^*(e)+F^*(b).F^*(e)=F^*(b).F^*(e)=F^*(e.b)=0
  \]
  Similarly, we find $b.F^{*i}(e)=0$ for $0\leq i\leq \sigma_0-1$. By
  (\ref{eq:miscequation2}), $b$ is orthogonal to $F^{\sigma_0-1}K$. As
  $F^{*\sigma_0-1}K$ is a maximal totally isotropic sub-bundle, it
  follows that $b\in\Gamma(S,F^{*\sigma_0-1}K)$. This means that our
  original section was equal to 0.
\end{proof}

\begin{Lemma}\label{lem:fiberdescription}
  If $K\subset V\otimes k$ is a characteristic subspace of Artin
  invariant $\sigma$, then the fiber of $\sU_V$ over $\left[K\right]$
  is isomorphic to $\bA^1\times(\bZ/p\bZ)^{\oplus
    \sigma_0-\sigma}$.
\end{Lemma}
\begin{proof}
  Consider $\cap_i\varphi^i(K)\subset V\otimes
  k$. This is fixed by $\varphi$, hence is equal to $W\otimes k$ for
  some totally isotropic subspace $W\subset V$, and has dimension
  $\sigma_0-\sigma$. Using Theorem \ref{thm:Witt}, we find a subspace
  $V_0\subset V$ such that the form on $V_0$ is neutral and
  nondegenerate and $W\subset V_0$ is a maximal totally isotropic
  subspace. There is a direct sum decomposition $V=V_0\oplus V/V_0$
  and an induced isomorphism
  \[
    \frac{V\otimes k}{K}\cong \frac{V_0\otimes k}{K\cap (V_0\otimes
      k)}\oplus\frac{V\otimes k/V_0\otimes k}{(K+V_0\otimes
      k)/V_0\otimes k}
  \]
  Because $W\otimes k$ is a maximal totally isotropic subspace of
  $V_0\otimes k$, and $K$ is totally isotropic,
  $W\otimes k=K\cap(V_0\otimes k)$. On the other hand,
  $K\cap V=K\cap V_0=W$. It follows that that
  $K\cap(V_0\otimes k)=(K\cap V_0)\otimes k$, so the first term is
  isomorphic to $(V_0/W)\otimes k$. The kernel of $1-F$ on this term
  is therefore $(\bZ/p\bZ)^{\oplus\sigma_0-\sigma}$. In the second
  term, note that the subspace
  $(K+V_0\otimes k)/V_0\otimes k\subset (V/V_0)\otimes k$ has trivial
  intersection with the $\bF_p$-subspace $V/V_0$ and therefore is
  strictly characteristic.
\end{proof}

\begin{Remark}\label{rem:!!}
  By Lemma \ref{lem:bundlestructure}, $\sU_U\to U$ is an fppf 
  $\bA^1$-bundle: its Frobenius pullback
  $F^{*\sigma_0-1}\sU_U$ is isomorphic to a line bundle, and hence
  locally trivial in the Zariski topology. (The Frobenius of $U$ is
  an fppf covering because $U$ is smooth over $\F_p$.) One might wonder
  if $\sU_U\to U$ is itself a Zariski $\A^1$-bundle. We do not believe
  this to be true for any $\sigma_0\geq 2$. Here is a proof for
  $\sigma_0=2$. Let $M_2$ be the period space of Artin invariant $2$
  and $M_3$ the period space of Artin invariant $3$. As noted in
  Remark \ref{rem:smallArtininvariant},
  $M_2\tensor_{\F_p}\overline{\F}_p$ is geometrically a disjoint union
  of two copies of $\bP^1$, and the 
  group scheme $\sU_U$ over $M_2$ is
  isomorphic to an open subset of $M_3$. Hence, if $\sU_U$ were
  Zariski-locally trivial, then each geometric component of $M_3$
  period space would be rational. On the other hand,
  as pointed out in  
  Remark \ref{rem:smallArtininvariant}, $M_3$ does not have rational
  geometric components for any $p\geq 3$. 
\end{Remark}

We conclude that the period space is \emph{not\/} in general an
iterated $\P^1$-bundle, contrary to what is claimed in
\cite[Theorem 4.3]{Liedtke15}.

We do have the following weaker consequence of Lemma
\ref{lem:bundlestructure}, which will suffice for our purposes.
\begin{Proposition}\label{prop:periodspaceisunirational}
  Each connected component of $M_V\otimes_{\bF_p}\bF_{p^2}$ is purely
  inseparably unirational over $\bF_{p^2}$.
\end{Proposition}
\begin{proof}
  We will induct on $\sigma_0$. The cases $\sigma_0=1$ and
  $\sigma_0=2$ are covered by Remark
  \ref{rem:smallArtininvariant}. Suppose that the result is true for
  some $\sigma_0$. By Proposition \ref{prop:isoofgroupschemes},
  $\sU_V$ is isomorphic to an open subset of $M_{\tV}$. Choosing an
  open subset of $M_V$ where $\sL_U$ is trivial, we find by Lemma
  \ref{lem:bundlestructure} a birational map
  \begin{equation}\label{eq:birational}
    M_{\tV}^{(p^{\sigma_0-1}/M)}\dashrightarrow M_V\times\bA^1
  \end{equation}
   where $M_{\tV}^{(p^{\sigma_0-1}/M)}$ is defined by the Cartesian diagram
   \[
     \begin{tikzcd}
       M_{\tV}^{(p^{\sigma_0-1}/M)}\arrow{r}{W^{\sigma_0-1}}\arrow{d}&M_V^{\langle v\rangle}\arrow{d}{\pi_v}\\
       M_V\arrow{r}{F_{M_V}^{\sigma_0-1}}&M_V
     \end{tikzcd}
   \]
   Composing the inverse of (\ref{eq:birational}) with $W^{\sigma_0-1}$, we find a dominant rational map
   \[
     M_V\times\bA^1\dashrightarrow M_{\tV}
   \]
   over $M_V$ which is generically finite and purely inseparable on fibers. The result follows by induction.
\end{proof}

\begin{Remark}
  In Proposition 10.3 of \cite{RS83} it is shown that $M_V$ is unirational. We have not checked whether their unirational parameterization is related to ours.
\end{Remark}

\begin{Remark}
  Here is an amusing observation. If $p=3$, then by Remark \ref{rem:smallArtininvariant} the $\sigma_0=3$ period space is geometrically a disjoint union of two copies of the Fermat quartic, which is known to be isomorphic to the unique supersingular K3 surface $X_1$ of Artin invariant $1$. The above gives an explicit description of a number of quasi-elliptic fibrations $X_1\to\P^1$.
\end{Remark}


\subsection{Twistor lines}\label{sec:twistorlines}
We will fix for the remainder of this paper an algebraically closed field $k$ of characteristic $p>0$. We will continue to assume that $p\geq3$, unless explicitly noted otherwise.
\begin{Notation}
  We write
  \[
    \oM_V=M_V\times_{\Spec\bF_p}\Spec k
  \]
  for the base change of $\oM_V$ to $\Spec k$. Equivalently, this is the functor on the category of schemes over $k$ whose value $\oM_V(S)$ on a $k$-scheme $S$ is the set of characteristic subspaces of $V\otimes\cO_S$. We let $\oM_V^W$ and $\oM_{V,W}$ be the base changes of the subschemes defined in Definition \ref{def:MVW}, and we will denote the base change of the universal bundle $K_V$ and the morphism $\pi_v$ by the same symbols.
\end{Notation}
It is shown in \cite{Ogus78} that $\oM_V$ has two connected components, each of which is defined over $\bF_{p^2}$. We maintain Notation \ref{notation:vectorspace}, so that $\tV$ is a vector space over $\bF_p$ of dimension $2\sigma_0+2$ equipped with a non-degenerate, non-neutral bilinear form.

\begin{Definition}
  A \textit{twistor line} in $\oM_{\tV}$ is a subvariety $L\subset\oM_{\tV}$ that is a connected component of a fiber of $\pi_v$ over a $k$-point $[K]\in \oM_{v^\perp/v}(k)$ for some isotropic $v\in\tV$.
\end{Definition}

By Lemma \ref{lem:fiberdescription}, any twistor line is isomorphic to $\bA^1$. We will eventually construct certain families of twisted K3 surfaces whose periods correspond to the twistor lines. In this section, we will record some facts about the geometry of twistor lines in $\oM_{\tV}$.

\begin{Remark}\label{rem:descriptionoffiber}
Fix a characteristic subspace $K\subset V\otimes k$. Let $U(K)$ be the group of $k$-points of the fiber of $\sU_V$ over $[K]\in \oM_V(k)$, so that
\[
  U(K)=\left\{B\in V\otimes k|B-\varphi(B)\in K+\varphi(K)\right\}/K
\]
Under the isomorphism $U(K)\iso\pi_v^{-1}([K])(k)$ of Proposition \ref{prop:isoofgroupschemes}, we obtain an explicit description of the fiber of $\pi_v$ over $\left[K\right]$, which is useful for making computations. Let $\left\{x_1,\dots,x_{\sigma_0}\right\}$ be a basis for $K$ and $B\in V\otimes k$ an element such that $B-\varphi(B)\in K+\varphi(K)$. Set
\[
  K(B)=\left\langle x_1+(x_1\innerproduct B)v,\dots,x_{\sigma_0}+(x_{\sigma_0}\innerproduct B)v,w+B+\frac{B^2}{2}w\right\rangle\subset \tV\otimes k
\]
Then $K(B)$ determines (and is determined by) $B$ modulo $K$, and $K(B)$ is the characteristic subspace corresponding to (the image of) $B$ in $U(K)$. If $K$ is strictly characteristic, then a further description of the fiber $\pi_v^{-1}([K])$ is implied by Lemma \ref{lem:bundlestructure}. Specifically, if $e$ is a generator for the line
\[
  l_K=K\cap\varphi(K)\cap\dots\cap\varphi^{\sigma_0-1}(K)
\]
then the map
\[
  B\mapsto B\innerproduct\varphi^{-\sigma_0+1}(e)
\]
gives a bijection $\sU(K)\iso k$.
\end{Remark}

Let $L\subset\oM_{\tV}$ be a twistor line that is a connected component of a fiber of $\pi_v$ over a $k$-point $[K]\in\oM_V(k)$.
\begin{Lemma}\label{lem:genericArtin}
  If $[\tK]\in L(k)$ is a $k$-point, then
  \[
    \sigma_0(\tK)=
      \begin{cases}
        \sigma_0(K)+1&\text{ if }\tK\cap\tV\subset v^\perp\\
        \sigma_0(K)&\text{ otherwise}
      \end{cases}
  \]
  In particular, the twistor line $L$ has generic Artin invariant $\sigma_0(K)+1$.
\end{Lemma}
\begin{proof}
  By assumption, $K=(\tK\cap v^\perp)/v\subset V\otimes k$. As in the proof of Lemma \ref{lem:definitionMorphism}, the quotient map $v^\perp\to v^\perp/v=V$ induces an isomorphism
  \[
    (\tK\cap\tV\cap v^\perp)/v\iso K\cap V
  \]
  It follows that $\sigma_0(\tK)=\sigma_0(K)+1$ if $\tK\cap\tV\subset v^\perp$, and $\sigma_0(\tK)=\sigma_0(K)$ otherwise.
\end{proof}

\begin{Lemma}\label{lem:incidence1}
  If $v'\in\tV$ is an isotropic vector such that $v\innerproduct v'\neq 0$, then $L$ intersects $\oM_{\tV,v'}$ in at most one point.
\end{Lemma}
\begin{proof}
  Recall that in Notation \ref{notation:vectorspace} we have fixed an orthogonal decomposition $\tV=V\oplus (U_2\otimes\bF_p)$, where $U_2\otimes\bF_p$ is generated by $v$ and $w$, which satisfy $v^2=w^2=0$ and $v\innerproduct w=-1$. This decomposition gives rise to a diagram
  \[
    \begin{tikzcd}
      \oM^{\langle v\rangle}_{\tV}\arrow{d}{\pi_v}\\
      \oM_{V}\arrow[bend left=45]{u}{\sigma_{w}}
    \end{tikzcd}
  \]
  where $\sigma_w$ is the section of $\pi_v$ defined by $K\mapsto \langle K,w\rangle$. Note that the image of $\sigma_w$ is $\oM_{\tV,w}\subset\oM_{\tV}$. From this, it is clear that $L$ intersects the locus $\oM_{\tV,w}$ in at most one point. We will transfer this result to $v'$, essentially by using the group structure given by the isomorphism of Proposition \ref{prop:isoofgroupschemes}. In terms of our orthogonal decomposition, we have $v'=w+b+(b^2/2) v$ for some $b\in V$. Consider the linear transformation $\exp(b)\colon \tV\to\tV$
  given by
  \begin{align*}
      v&\mapsto v,\\
      w&\mapsto w+b+(b^2/2)v,\hspace{.2cm}\mbox{ and}\\
      x&\mapsto x+(x\innerproduct b)v\hspace{.2cm}\mbox{ for }x\in V.
  \end{align*}
  One checks that $\exp(b)$ is in fact an isometry. We have a commutative diagram
  \[
    \begin{tikzcd}[column sep=tiny]
      \oM^{\langle v\rangle}_{\tV}\arrow{rr}{\exp(b)}[swap]{\sim}\arrow{dr}[swap]{\pi_v}&&\oM^{\langle v\rangle}_{\tV}\arrow{dl}{\pi_v}\\
      &\oM_V&
    \end{tikzcd}
  \]
  Furthermore, the image of $\exp(b)\circ\sigma_w$ is exactly $\oM_{\tV,v'}$. This gives the result.
\end{proof}

\begin{Lemma}\label{lem:incidence2}
  If $L$ has generic Artin invariant $\sigma\geq 2$, then $L$ intersects the locus of points with Artin invariant $\sigma-1$ in $p^{2\sigma-2}$ distinct points.
\end{Lemma}
\begin{proof}
  Under the isomorphism of Proposition \ref{prop:isoofgroupschemes} (see also Remark \ref{rem:descriptionoffiber}), the points in the fiber $\pi_v^{-1}([K])$ with Artin invariant $\sigma-1$ correspond to the subgroup
  \[
    \frac{V}{K\cap V}\subset\frac{V\otimes k}{K}
  \]
  The dimension of the $\bF_p$-vector space $V/(K\cap V)$ is $\sigma_0+\sigma-2$. By Lemma \ref{lem:fiberdescription}, the fiber $\pi_v^{-1}([K])$ has $p^{\sigma_0-\sigma}$ connected components. It follows that $L$ contains $p^{2\sigma-2}$ points of Artin invariant $\sigma-1$.
\end{proof}

We next study the tangent spaces of twistor lines.

\begin{Lemma}\label{lem:Ogustangentspace}
  The canonical connection on $F^*_{M_{\tV}}K_{\tV}$ induces an isomorphism
  \[
    T^1_{M_{\tV}}\iso \sHom_{M_{\tV}}\left(K_{\tV}\cap F^*_{M_{\tV}}K_{\tV},\frac{F_{M_{\tV}}^*K_{\tV}}{K_{\tV}\cap F_{M_{\tV}}^*K_{\tV}}\right)
  \]
\end{Lemma}
\begin{proof}
  See \cite{Ogus78}, Proposition 4.6.
\end{proof}

\begin{Lemma}\label{lem:relativetangentspace}
  Write $\tK=K_{\tV}|_{M^{\langle v\rangle}_{\tV}}$ and $K=K_V$. Under the isomorphism of Lemma \ref{lem:Ogustangentspace}, the sub-bundle
  \[
    T^1_{M^{\langle v\rangle}_{\tV}/M_{V}}\subset T^1_{M_{\tV}^{\langle v\rangle}}
  \]
  is identified with
  \[
    \sHom_{M^{\langle v\rangle}_{\tV}}\left(\frac{\tK\cap F^*\tK}{\tK\cap F^*\tK\cap v^\perp},\frac{F^*\tK}{\tK\cap F^*\tK}\right)\cong\frac{F^*\tK}{\tK\cap F^*\tK}
  \]
\end{Lemma}
\begin{proof}
  There are isomorphisms
  \[
    \pi_v^*T^1_{M_V}\cong\pi_v^*\sHom_{M_V}\left(K\cap F^*K,\dfrac{F^*K}{K\cap F^*K}\right)\cong\sHom_{M^{\langle v\rangle}_{\tV}}\left(\tK\cap F^*\tK\cap v^\perp,\dfrac{F^*\tK}{\tK\cap F^*\tK}\right)
  \]
  The morphism $\pi_v$ induces a commutative diagram
  \[
    \begin{tikzcd}
      T^1_{M^{\langle v\rangle}_{\tV}}\arrow{r}{\sim}\arrow{d}&\sHom_{M_{\tV}^{\langle v\rangle}}\left(\tK\cap F^*\tK,\dfrac{F^*\tK}{\tK\cap F^*\tK}\right)\arrow{d}\\
      \pi_v^*T^1_{M_V}\arrow{r}{\sim}&\sHom_{M_{\tV}^{\langle v\rangle}}\left(\tK\cap F^*\tK\cap v^\perp,\dfrac{F^*\tK}{\tK\cap F^*\tK}\right)
    \end{tikzcd}
  \]
  and the right vertical map sends a homomorphism $f\colon \tK\cap F^*\tK\to F^*\tK/\tK\cap F^*\tK$ to the composition
  \[
    \tK\cap F^*\tK\cap v^\perp\hto \tK\cap F^*\tK\too{f}\frac{F^*\tK}{\tK\cap F^*\tK}
  \]
  Thus, its kernel is
  \[
     T^1_{M^{\langle v\rangle}_{\tV}/M_{V}}=\sHom_{M^{\langle v\rangle}_{\tV}}\left(\frac{\tK\cap F^*\tK}{\tK\cap F^*\tK\cap v^\perp},\frac{F^*\tK}{\tK\cap F^*\tK}\right)
  \]
  Finally, note that pairing with $v$ gives a trivialization
  \[
    \frac{\tK\cap F^*\tK}{\tK\cap F^*\tK\cap v^\perp}\iso \cO_{M^{\langle v\rangle}_{\tV}}
  \]
\end{proof}


Consider the open subset $U_{\tV}\subset M_{\tV}$ consisting of those characteristic subspaces with maximal Artin invariant. This is the intersection of the domains of definition of the rational maps $\pi_v$ as $v$ ranges over all isotropic vectors in $\tV$. Thus, for any isotropic $v\in V$, we get a morphism $U_{\tV}\to U_{v^\perp/v}$, whose fibers over $k$-points are open subsets of twistor lines in $M_{\tV}$.

\begin{Lemma}\label{BigTangentSpace}
  Let $v_0,\dots,v_n$ be an enumeration of the isotropic vectors in $\tV$. The morphism 
  \[
    U_{\tV}\to U_{v_0^\perp/v_0}\times\dots\times U_{v_n^\perp/v_n}
  \]
  induced by the fibrations $\pi_{v_i}$ is unramified. In other words, for any $x\in U_{\tV}(k)$, the tangent vectors to the twistor lines through $x$ span $T^1_{U_{\tV},x}$
\end{Lemma}
\begin{proof}
  Take a point $x\in U_{\tV}(k)$ corresponding to a characteristic subspace $\tK\subset V\otimes k$. Fix an isomorphism $(F^*\tK/(\tK+F^*\tK))|_{x}\iso k$. If $v_i\in\tV$ is isotropic, then the tangent space to the fiber of $\pi_{v_i}$ containing $x$ is 
  \[
    \Hom\left(\frac{\tK\cap F^*\tK}{\tK\cap F^*\tK\cap v_i^\perp},k\right)\subset\Hom(\tK\cap F^*\tK,k)
  \]
  Thus, it will suffice to show that the functions $\langle \_,v_i\rangle$ span $\Hom(\tK\cap F^*\tK,k)$. Let $e\in\tK$ be a vector that spans the line $l_{\tK}=\tK\cap\varphi(\tK)\cap\dots\varphi^{\sigma_0}(\tK)$, so that $\left\{\varphi^{-\sigma_0+1}(e),\dots,e\right\}$ is a basis for $\tK\cap\varphi(\tK)$. Note that for any vector $v$ defined over $\bF_p$ we have
  \[
    (v\innerproduct e)^{p^j}=\sigma^j(v\innerproduct e)=v\innerproduct \varphi^j(e)
  \]
  for all $j$. In particular, each $v_i$ is uniquely determined by $\lambda_i=v_i\innerproduct\varphi^{-\sigma_0+1}(e)$, and the $\lambda_i$ are distinct and non-zero. We will be finished if we can show that the matrix
  \[
    A = \begin{bmatrix}
      \lambda_0 & \lambda_0^p & \dots &\lambda_0^{p^{\sigma_0-1}}\\
      \lambda_1 & \lambda_1^p & \dots &\lambda_1^{p^{\sigma_0-1}}\\
      \vdots & \vdots & \ddots & \vdots \\
      \lambda_n & \lambda_n^p & \dots &\lambda_n^{p^{\sigma_0-1}}\\
      \end{bmatrix}
  \]
  has full rank. Take a $\sigma_0\times\sigma_0$ minor of $A$, which we may assume to be $(\lambda_i^{p^{j}})_{0\leq i,j\leq \sigma_0-1}$ without loss of generality. Consider its determinant $D(\lambda_0,\dots,\lambda_{\sigma_0-1})$ as a polynomial in the $\lambda_i$. We claim that
  \[
    D(\lambda_0,\dots,\lambda_{\sigma_0-1})=c\prod_{\alpha\in\bP^{\sigma_0-1}(\bF_p)}\lambda_\alpha
  \]
  where $c\in\bF^\times$ is a non-zero constant, and for $\alpha=\left[\alpha_0:\dots:\alpha_{\sigma_0-1}\right]\in\bP^{\sigma_0-1}(\bF_p)$ we set $\lambda_\alpha=\sum_{k=0}^{\sigma_0-1}\alpha_k\lambda_k$, which is well defined up to elements of $\bF_p^\times$. To see this, note that if one of the $\lambda_\alpha$ vanishes, then $D(\lambda_0,\dots,\lambda_{\sigma_0-1})$ vanishes as well. Thus, by the nullstellensatz, the product of the $\lambda_\alpha$ divides $D(\lambda_0,\dots,\lambda_{\sigma_0-1})$. By comparing degrees we get the result. So, to show that $A$ has full rank, we need to show that we can find $\sigma_0$ vectors among the $v_i$ such that the corresponding $\lambda_i$ are linearly independent over $\bF_p$. Because the $\lambda_i$ are distinct and non-zero, we can do this if $n\geq p^{\sigma_0-1}$. By Lemma 4.12 of \cite{Ogus78}, the number of isotropic vectors in $\tV$ is $p^{2\sigma_0+1}-p^{\sigma_0+1}+p^{\sigma_0}-1$.
\end{proof}

\subsection{The relative Brauer group of a family of supersingular K3 surfaces}\label{sec:sheaves on the moduli space}

Fix a relative supersingular K3 surface $\pi\colon X\to S$. In this section we will study the relationship between the cohomology sheaves $\R^2\pi^{\fl}_*\m_p$ and $\R^1\pi^{(p)\etale}_*\nu(1)$, and relate them to the Hodge and conjugate filtrations on the de Rham cohomology of $X\to S$. In order to allow for the variation of the Artin invariant of the fibers, we will find it convenient to introduce an appropriate marking of the relative Picard group.

\begin{Definition}
  A \textit{supersingular K3 lattice} is a free abelian group $\Lambda$ of rank 22 equipped with an even symmetric bilinear form such that
  \begin{enumerate}
    \item $\disc(\Lambda\otimes\bQ)=-1$ in $\bQ^\times/\bQ^{\times 2}$,
    \item the signature of $\Lambda\otimes\bR$ is $(1,21)$, and
    \item the discriminant group $\Lambda^*/\Lambda$ is $p$-torsion.
  \end{enumerate}
\end{Definition}
For the definitions of these terms and other lattice theoretic background we refer to Chapter 14 of \cite{Huy16}. It is shown in \cite{Ogus83} that the supersingular K3 lattices are precisely those lattices that occur as the Picard group of a supersingular K3 surface. If $\Lambda$ is a supersingular K3 lattice, then there are inclusions
\[
  p\Lambda\subset p\Lambda^*\subset\Lambda\subset\Lambda^*
\]
We set
\[
  \Lambda_0=\frac{p\Lambda^*}{p\Lambda}\hspace{1cm}\Lambda_1=\frac{\Lambda}{p\Lambda^*}
\]
By condition (3), $\Lambda_0$ and $\Lambda_1$ are vector spaces over $\bF_p$. The dimension of $\Lambda_0$ over $\bF_p$ is equal to $2\sigma_0$ for some integer $1\leq\sigma_0\leq 10$, called the \textit{Artin invariant} of $\Lambda$, and the dimension of $\Lambda_1$ is $22-2\sigma_0$. The Artin invariant determines $\Lambda$ uniquely up to isometry (see Theorem 7.4 of \cite{Ogus78}). There is a natural symmetric bilinear form $\Lambda_0\otimes\Lambda_0\to\bF_p$ given by 
\[
  \overline{v}\innerproduct\overline{w}=p^{-1}\langle v,w\rangle_{\Lambda}\hspace{.1cm}\mod p
\]
for $v,w\in p\Lambda^*\subset\Lambda$. This form is non-degenerate and non-neutral, so that $V=\Lambda_0$ satisfies the assumptions of section \ref{sec:perioddomain}. We fix a supersingular K3 lattice $\Lambda$ with Artin invariant $\sigma_0$.

\begin{Definition}
  If $S$ is a $k$-scheme, we define a \textit{family of $\Lambda$-marked supersingular K3 surfaces over $S$} to be an algebraic space $X$ equipped with a smooth, proper morphism $X\to S$ whose geometric fibers are supersingular K3 surfaces, along with a morphism $m\colon \underline{\Lambda}_S\to\Pic_{X/S}$ of sheaves of groups that is compatible with the intersection forms. We let $S_\Lambda$ denote the functor whose $S$-points are isomorphism classes of families of $\Lambda$-marked supersingular K3 surfaces over $S$. 
\end{Definition}

Using Artin's representability theorems, Ogus proved the following.

\begin{Theorem}\cite[Theorem 2.7]{Ogus83}\label{thm:representable}
  The functor $S_\Lambda$ is representable by an algebraic space over $k$ that is locally of finite presentation, locally separated, and smooth of dimension $\sigma_0-1$.
\end{Theorem}

We recall the definition of Ogus's crystalline period morphism.
\begin{Definition}
  Let $\pi\colon X\to S$ be a relative supersingular K3 surface equipped with a marking $m\colon \underline{\Lambda}_S\to\Pic_{X/S}$. Composing with the Chern class map $\Pic_{X/S}\to \H^2_{dR}(X/S)$, we obtain a map $\underline{\Lambda}_S\to \H^2_{dR}(X/S)$, which induces a map $\underline{\Lambda}_S\otimes\cO_S\to \H^2_{dR}(X/S)$. Because $p\cO_S=0$, there is a natural inclusion $\uLam_{0S}\otimes\cO_S\subset\uLam_S\otimes\cO_S$. Suppose that $S$ is smooth over $k$. As described in Section 3 of \cite{Ogus83}, the kernel of the map $\uLam_S\otimes\cO_S\to \H^2_{dR}(X/S)$ is a characteristic subspace $K'\subset\uLam_{0S}\otimes\cO_S$, and that the Gauss-Manin connection induces a descent datum on $K'$ with respect to the Frobenius $F_S$. Thus we find a characteristic subspace $K\subset\uLam_{0S}\otimes\cO_S$ equipped with an isomorphism $F^*_SK\cong K'$. In particular, the moduli space $S_{\Lambda}$ is smooth, so we may apply the construction of $K$ to produce (via \'{e}tale descent) a characteristic subspace $K_{\Lambda}\subset\uLam_{0S}\otimes\cO_{S_{\Lambda}}$, which in turn gives a morphism
  \begin{equation}\label{eq:periodmorphism}
    \rho\colon S_{\Lambda}\to \oM_{\Lambda_0}
  \end{equation}
\end{Definition}

\begin{Remark}
  The morphism $\rho$ can be interpreted in terms of crystalline cohomology (see Remark \ref{rem:periodmorphism1}). We will take this viewpoint in Section \ref{sec:twistedperiodmorphism} when we define a period morphism for twisted supersingular K3 surfaces.
\end{Remark}

\begin{Remark}\label{rem:subtleremark}
  Note that this definition of the period morphism is not purely moduli-theoretic, in the sense of giving an ``intrinsic'' procedure that associates to any marked family over an arbitrary $k$-scheme $S$ a characteristic subspace. The difficulty lies in the descent through the Frobenius. If $S$ is a smooth $k$-scheme, the Gauss Manin connection induces a descent of the kernel of the map $\Lambda\otimes\cO_S\to \H^2_{dR}(X/S)$ through the Frobenius, and the resulting characteristic subspace is of formation compatible with base change along morphisms $T\to S$, where $T$ is a smooth $k$-scheme. If, however, $S$ is not smooth, then it is not clear which descent of $K'$ to choose. By defining the period morphism in terms of the universal characteristic subspace, we have in effect specified a choice of Frobenius descent of the kernel for every marked family over an arbitrary base.
\end{Remark}

\begin{Proposition}\label{prop:nuRepresentable}
  If $\pi\colon X\to S$ is a relative supersingular K3 surface, then the sheaf $\R^1\pi^{(p)\etale}_*\nu(1)$ is representable by a group algebraic space of finite presentation over $S$.
\end{Proposition}
\begin{proof}
  It will suffice to prove the result after taking an \'{e}tale cover of $S$. Thus, we may assume that there exists a marking $\uLam_S\to\Pic_{X/S}$ (see page 1522 of \cite{RS79}). Let $\rho_S\colon S\to\oM_{\Lambda_0}$ be the corresponding morphism. The sheaf $K_{\Lambda_0}$ on $\oM_{\Lambda_0}$ pulls back to a sheaf $K$ on $S$, whose associated vector bundle we will denote by $\uK$. By results of \cite{Ogus78}, we have a diagram
  \begin{equation}
    \begin{tikzcd}\label{eq:hodge}
      0\arrow{r}&F_S^*\uK\arrow{r}&\uLam_S\otimes\cO_S\arrow{r}&F^0_H\arrow{r}&F_S^*\uK\arrow{r}&0\\
      0\arrow{r}&F_S^*\uK\arrow{r}\arrow[equal]{u}&\uLam_S\otimes\cO_S\arrow{r}\arrow[equal]{u}& F^1_H \arrow{r}\arrow[hook]{u}&\uK\cap F_S^*\uK\arrow{r}\arrow[hook]{u}&0\\
      0\arrow{r}&F_S^*\uK\arrow{r}\arrow[equal]{u}&\uK+F_S^*\uK\arrow{r}\arrow[hook]{u}&F^2_H \arrow{r}\arrow[hook]{u}&0\arrow[hook]{u}&
    \end{tikzcd}
  \end{equation}
  of big \'{e}tale sheaves on $S$ with exact rows, where
  \[
    0\subset F^2_H\subset F^1_H\subset F^0_H=\H^2_{dR}(X/S)_{\etale}
  \]
  is the Hodge filtration (\ref{eq:filtrations}) on the second de Rham 
  cohomology $\H^2_{dR}(X/S)_{\etale}$. Note that the entries of this diagram 
  are representable by groups schemes on $S$. The first Chern class map 
  $\Pic_{X/S}\to \R^2\pi^{\etale}_*\Omega_{X/S}^{\etale\,\bullet}$ factors 
  through $\R^1\pi^{\etale}_*Z^1\Omega^{\etale\,\bullet}_{X/S}$. By Lemma 
  \ref{lem:closeddifferentialforms}, there is an identification 
  $\R^1\pi_*^{\etale}Z^1\Omega^{\etale\,\bullet}_{X/S}\iso F^1_H\cap F^1_C$, 
  and the inclusion
  \[
   F^1_H\cap F^1_C\hto \H^2_{dR}(X/S)_{\etale}
  \]
  together with (\ref{eq:hodge}) gives a diagram
  \begin{equation}\label{eq:cohomologyofclosedforms}
    \begin{tikzcd}
      0\arrow{r}&\uLam_{S}\otimes\cO_S/F^*_S\uK\arrow{r}\arrow{d}&F^1_H\cap F^1_C\arrow{r}\arrow{d}&\uK\cap F^*_S\uK\cap F^{*2}_S\uK\arrow{d}\arrow{r}&0\\
      0\arrow{r}&\uLam_{S}\otimes\cO_S/F^*_S\uK\arrow{r}& \H^2_{dR}(X/S)_{\etale}\arrow{r}& F^*_S\uK\arrow{r}&0
    \end{tikzcd}
  \end{equation}
  with exact rows. We emphasize that by $F^1_H\cap F^1_C$ and $\uK\cap F^*_S\uK\cap F^{*2}_S\uK$ we mean the fiber products of sheaves (or modules) on the big \'{e}tale site of $S$. The rank of these sheaves jumps on the superspecial locus, and they will not be quasi-coherent in general. Furthermore, by (\ref{eq:hodge}) there is a short exact sequence
  \[
    0\to\dfrac{\uLam_S\otimes\cO_S}{F^*_S\uK+F^{*2}_S\uK}\to F^*_S(F^1_H/F^2_H)\to F^*_S\uK\cap F^{*2}_S\uK\to 0
  \]
  We are led to the following diagram of big \'{e}tale sheaves on $S$, with exact rows and columns.
  \begin{equation}\label{eq:bigdiagram}
    \begin{tikzcd}
      &0\arrow{d}&0\arrow{d}&0\arrow{d}\\
      0\arrow{r}&\sV\arrow{r}\arrow{d}&\dfrac{\uLam_S\otimes\cO_S}{F^*_S\uK}\arrow{r}{1-F_S^*}\arrow{d}&\dfrac{\uLam_S\otimes\cO_S}{F_S^*\uK+F_S^{*2}\uK}\arrow{d}\\
      0\arrow{r}&\R^1\pi^{(p)\etale}_*\nu(1)\arrow{r}\arrow{d}&F^1_H\cap F^1_C\arrow{r}{C\circ\pi_C-F^*_S\circ\pi_H}\arrow{d}&F^*_S(F^1_H/F^2_H)\arrow{d}\\
      0\arrow{r}&\sW\arrow{d}\arrow{r}&\uK\cap F^*_S\uK\cap F^{*2}_S\uK\arrow{r}{1-F^*_S}\arrow{d}&F_S^*\uK\cap F^{*2}_S\uK\arrow{d}\\
      &0&0&0
    \end{tikzcd}
  \end{equation}
  The sheaf $\sV$ is the kernel of a map of group schemes, and so is representable. We have a map of exact sequences
  \[
    \begin{tikzcd}
      0\arrow{r}&\sW\arrow{r}\arrow[equal]{d}&\uK\cap F^*_S\uK\cap F^{*2}_S\uK\arrow[hook]{d}\arrow{r}{1-F^*_S}&F^*_S\uK\cap F^{*2}_S\uK\arrow[hook]{d}\\
      0\arrow{r}&\sW\arrow{r}&F^*_S\uK\arrow{r}{1-F^*_S}&F^*_S\uK+F^{*2}_S\uK
    \end{tikzcd}
  \]
  Thus $\sW$ is also the kernel of a map of group schemes and hence representable. By Lemma \ref{lem:sesrepresentable}, it follows that $\R^1\pi^{(p)\etale}_*\nu(1)$ is representable by an algebraic space.
  \end{proof}
  
\begin{Lemma}\label{lem:sesrepresentable}
    If $X$ is a scheme and $0\to A\to B\to C\to 0$ is a short exact sequence of sheaves of abelian groups on $X_{\Etale}$, then if $A$ and $C$ are representable by algebraic spaces, so is $B$.
\end{Lemma}
\begin{proof}
    Let us think of $B$ as a sheaf of groups on $C$ via the map $B\to C$. The map $A\to B$ gives an action of $A$ on $B$ over $C$. Because $B\to C$ is a surjection of big \'{e}tale sheaves, it admits a section \'{e}tale locally. Thus, $B$ is an $A$-torsor, and \'{e}tale locally on $C$ there is an isomorphism $B\iso A\times C$. We conclude that $B$ is an algebraic space.
\end{proof}

If $\pi\colon X\to S$ is a family of supersingular K3 surfaces, then the conditions of Proposition \ref{prop:flattoetalecohomology} are satisfied, and so we obtain a morphism
\begin{equation}\label{eq:thenotisomorphism}
  \Upsilon\colon \R^2\pi^{\fl}_*\m_p\to \R^1\pi^{(p)\etale}_*\nu(1)
\end{equation}
By Theorem \ref{thm:relativeartin} and Proposition \ref{prop:nuRepresentable}, these sheaves are represented by group algebraic spaces over $S$.

\begin{Remark}
  It follows from Proposition \ref{prop:homeomorphism} that $\Upsilon$ is a universal homeomorphism. By Remark \ref{rem:RemarkNumeroUno}, it is totally ramified relative to $S$. We will deduce in Proposition \ref{prop:aninterestingconsequence} that $\Upsilon$ is essentially the relative Frobenius over $S$.
\end{Remark}

Suppose that the family $X\to S$ is equipped with a marking $\uLam_S\to\Pic_{X/S}$. We will identify some particular subgroups of the sheaves $\R^2\pi^{\fl}_*\m_p$ and $\R^1\pi^{(p)\etale}_*\nu(1)$. We will use the notation of the proof of Proposition \ref{prop:nuRepresentable}, and in particular diagram (\ref{eq:bigdiagram}). Let $\sU_S$ be the pullback of $\sU_{\Lambda_0}$ under the induced map $\rho_S\colon S\to M_{\Lambda_0}$. Consider the exact sequence
\[
  0\to F_S^{-1}\sU_S\to \dfrac{\uLam_{0S}\otimes\cO_{S}}{F^*_S\uK}\xrightarrow{1-F^*_S} \dfrac{\uLam_{0S}\otimes\cO_{S}}{F^*_S\uK+F^{*2}_S\uK}
\]
As $F^*_S\uK\subset\uLam_{0S}\otimes\cO_{S}$, we get a commuting diagram
\[
  \begin{tikzcd}
    0\arrow{r}&F^{-1}_S\sU_{S}\arrow[hook]{d}\arrow{r}&\dfrac{\uLam_{0S}\otimes\cO_{S}}{F^*_S\uK}\arrow{r}{1-F^*_S}\arrow{d}&\dfrac{\uLam_{0S}\otimes\cO_{S}}{F^*_S\uK+F^{*2}_S\uK}\arrow{d}\\
    0\arrow{r}&\sV\arrow{r}&\dfrac{\uLam_{S}\otimes\cO_{S}}{F^*_S\uK}\arrow{r}{1-F^*_S}&\dfrac{\uLam_{S}\otimes\cO_{S}}{F^*_S\uK+F^{*2}_S\uK}
  \end{tikzcd}
\]

\begin{Definition}\label{def:definitioncirc}
  If $\pi\colon X\to S$ is a relative supersingular K3 surface and $\uLam\to\Pic_{X/S}$ is a marking, we define the subsheaf
  \[
    \R^1\pi^{(p)\etale}_*\nu(1)^o\subset \R^1\pi^{(p)\etale}_*\nu(1)
  \]
  to be the image of $F_S^{-1}\sU_S$ under the inclusions
  \[
    F_S^{-1}\sU_{S}\subset\sV\subset\R^1\pi^{(p)\etale}_*\nu(1)
  \]
  We define 
  \[
    (\R^2\pi^{\fl}_*\m_p)^o\subset \R^2\pi^{\fl}_*\m_p
  \]
  to be the preimage of $\R^1\pi^{(p)\etale}_*\nu(1)^o$ under the morphism (\ref{eq:thenotisomorphism}).
\end{Definition}

Note that these subsheaves depend on the marking, although it is suppressed from the notation.

\begin{Lemma}\label{lem:flatduality1}
  The subgroups $\R^1\pi^{(p)\etale}_*\nu(1)^o$ and $(\R^2\pi^{\fl}_*\m_p)^o$ are open and closed and are represented by group algebraic spaces. There is a short exact sequence
  \begin{equation}\label{eq:equation1}
    0\to\R^1\pi^{(p)\etale}_*\nu(1)^o\to \R^1\pi^{(p)\etale}_*\nu(1)\to\sD\to 0
  \end{equation}
  where $\sD$ is a group scheme that fits into a short exact sequence
  \[
    0\to\uLam_{1S}\to\sD\to\sW\to 0
  \]
\end{Lemma}
\begin{proof}
  Because $F_S^{-1}\sU_S$ is representable, so are $\R^1\pi^{(p)\etale}_*\nu(1)^o$ and $(\R^2\pi^{\fl}_*\m_p)^o$. We have a short exact sequence
  \[
    0\to\sV\to \R^1\pi^{(p)\etale}_*\nu(1)\to\sW\to 0
  \]
  The cokernel of the map $F_S^{-1}\sU_{S}\to\sV$ is the fixed points of $F^*_{S}$ acting on $\uLam_S\otimes\cO_S/\uLam_{0S}\otimes\cO_S$, which is just $\uLam_{1S}$. We obtain a diagram
  \begin{equation}\label{eq:flatdualitydiagram}
    \begin{tikzcd}
      &0\arrow{d}&0\arrow{d}&0\arrow{d}&\\
      0\arrow{r}&F_S^{-1}\sU_{S}\arrow{r}\arrow{d}&\sV\arrow{r}\arrow{d}&\uLam_{1S}\arrow{r}\arrow{d}&0\\
      0\arrow{r}&\R^1\pi^{(p)\etale}_*\nu(1)^o\arrow{r}\arrow{d}&\R^1\pi^{(p)\etale}_*\nu(1)\arrow{r}\arrow{d}&\sD\arrow{r}\arrow{d}&0\\
      &0\arrow{r}&\sW\arrow{r}\arrow{d}&\sW\arrow{d}\arrow{r}&0\\
      &&0&0&
    \end{tikzcd}
  \end{equation}
  with exact columns and rows. Because $\uLam_{1S}$ is separated over $S$, the morphism $\sD\to\sW$ is separated. The sheaf $\sW$ is equal to the intersection $\uK\cap\uLam_{0S}$, and in particular is separated over $S$. We conclude that $\sD$ is separated over $S$, and therefore the zero section of $\sD$ is closed. The immersion $\sW\subset\uLam_{0S}$ shows that the zero section of $\sW$ is open, and hence $\uLam_{1S}\to\sD$ is open. It follows that the zero section of $\sD$ is also open. Thus, the subgroups $\R^1\pi^{(p)\etale}_*\nu(1)^o$ and $(\R^2\pi^{\fl}_*\m_p)^o$ are open and closed.
\end{proof}


\begin{Lemma}\label{lem:smoothofdimension1}
  The group spaces $(\R^2\pi^{\fl}_*\m_p)^o$ and $\R^1\pi^{(p)\etale}_*\nu(1)^o$ are smooth over $S$ of relative dimension 1. Every geometric fiber of either has connected component isomorphic to $\bG_a$.
\end{Lemma}
\begin{proof}
  By Lemma \ref{lem:smooth}, $\R^1\pi^{(p)\etale}_*\nu(1)^o\to S$ is smooth, and by Lemma \ref{lem:fiberdescription}, every geometric fiber has connected component isomorphic to $\bG_a$. Lemma \ref{lem:flatduality1} gives that the morphism $(\R^2\pi^{\fl}_*\m_p)^o\subset\R^2\pi^{\fl}_*\m_p$ is open and closed, so by Proposition \ref{prop:relativeartin1} each geometric fiber of $(\R^2\pi^{\fl}_*\m_p)^o$ is regular of dimension 1, and has connected component isomorphic to $\bG_a$. It remains to show that $(\R^2\pi^{\fl}_*\m_p)^o\to S$ is flat. The morphism (\ref{eq:thenotisomorphism}) gives a diagram
  \[
    \begin{tikzcd}[column sep=small]
      (\R^2\pi^{\fl}_*\m_p)^o\arrow{r}{\Upsilon^o}\arrow{dr}&\R^1\pi^{(p)\etale}_*\nu(1)^o\arrow{d}\\
      &S
    \end{tikzcd}
  \]
  We have already seen that the vertical arrow is smooth. By Proposition \ref{prop:homeomorphism} the horizontal arrow is a universal homeomorphism. It follows that $(\R^2\pi^{\fl}_*\m_p)^o\to S$ is universally open. Moreover, its geometric fibers are reduced. It will suffice to prove the result in the universal case when $S=S_{\Lambda}$, so we may in addition assume that $S$ is reduced. By 15.2.3 of \cite{EGA4.3}, these conditions imply flatness.
\end{proof}


Let us now suppose that $S=\Spec k$, and that $\pi\colon X\to\Spec k$ is a supersingular K3 surface equipped with a marking $m\colon \Lambda\to\Pic(X)$. Evaluating the morphism (\ref{eq:thenotisomorphism}) on $\Spec k$, we obtain a diagram
\begin{equation}\label{eq:Lerayidentification}
  \begin{tikzcd}
    (\R^2\pi^{\fl}_*\m_p)(k)\arrow{r}{\Upsilon}[swap]{\sim}\arrow[equals]{d}&(\R^1\pi^{\etale(p)}_*\nu(1))(k)\arrow[equals]{d}\\
    \H^2(X_{\fl},\m_p)\arrow{r}{\sim}& \H^1(X_{\etale},\nu(1))
  \end{tikzcd}
\end{equation}
where the vertical arrows are the canonical identifications induced by the respective Leray spectral sequences. 

\begin{Definition}\label{def:transcendental}
  If $X$ is a K3 surface, we say that a class $\alpha\in \H^2(X,\m_p)$ is \textit{transcendental} if its image in $\H^2_{dR}(X/k)$ is orthogonal to the image of the first Chern class map
  \[
    \Pic(X)\to \H^2_{dR}(X/k)
  \]
  If $X$ is supersingular and $m\colon \Lambda\to\Pic(X)$ is a marking, then we say that $\alpha$ is \textit{transcendental with respect to} $m$ if its image in $\H^2_{dR}(X/k)$ is orthogonal to the image of
  \[
    \Lambda\to\Pic(X)\to \H^2_{dR}(X/k)
  \]
\end{Definition}

\begin{Lemma}\label{lem:transcendentalequalso}
  A class $\alpha\in\R^2\pi^{\fl}_*\m_p(k)$ is transcendental with respect to a marking $\Lambda\to\Pic(X)$ if and only if it lies in the subgroup $(\R^2\pi^{\fl}_*\m_p)^o(k)\subset \R^2\pi^{\fl}_*\m_p(k)$.
\end{Lemma}
\begin{proof}
  The subspace $\Lambda\otimes k/F^*_kK\subset \H^2_{dR}(X/k)$ has rank $22-\sigma_0$, where $\sigma_0$ is the Artin invariant of $\Lambda$. The subgroup $p\Lambda^*\subset\Lambda$ consists of those elements $x\in\Lambda$ such that $\langle x,y\rangle\equiv 0\mod{p}$ for every $y\in\Lambda$. The span of the image of $p\Lambda^*$ in $\H^2_{dR}(X/k)$ is therefore in the orthogonal complement of $\Lambda\otimes k/F^*_kK$. But the subspace $\Lambda_0\otimes k/F_k^*K\subset \H^2_{dR}(X/k)$ has rank $\sigma_0$, and the pairing on $\H^2_{dR}(X/k)$ is perfect. Therefore
  \[
    \left(\dfrac{\Lambda\otimes k}{F^*_kK}\right)^\perp=\left(\dfrac{\Lambda_0\otimes k}{F^*_kK}\right)
  \]
  as subspaces of $\H^2_{dR}(X/k)$. In particular, the subgroup of transcendental elements of $\R^2\pi^{\fl}_*\m_p(k)$ is equal to the intersection of $\R^2\pi^{\fl}_*\m_p(k)$ with $\Lambda_0\otimes k/F^*_kK$ inside of $\H^2_{dR}(X/k)$. By definition, this is equal to $(\R^2\pi^{\fl}_*\m_p)^o(k)$.
\end{proof}

Finally, let us further specialize to the case when $S=\Spec k$, and the marking $\Lambda\to\Pic(X)$ is an isomorphism.
\begin{Lemma}\label{lem:descriptionoverafield}
  If the marking $m$ is an isomorphism, then $(\R^2\pi^{\fl}_*\m_p)^o$ is the connected component of the identity.
\end{Lemma}
\begin{proof}
  Using the notation of the proof of Proposition \ref{prop:nuRepresentable}, we have a short exact sequence
  \[
    0\to\sV\to \R^1\pi^{(p)\etale}_*\nu(1)\to\sW\to 0
  \]
  Note that in this case $\sW$ is the trivial group. Thus, the identity component of $\R^1\pi^{(p)\etale}_*\nu(1)$ is identified with the image of the identity component of $F_S^{-1}(\sU_{S})$. But by Lemma \ref{lem:fiberdescription}, this group is isomorphic to $\bA^1$, and in particular connected. By Proposition \ref{prop:homeomorphism}, the morphism $\R^2\pi^{\fl}_*\m_p\to \R^1\pi^{(p)\etale}_*\nu(1)$ is a homeomorphism, so $(\R^2\pi^{\fl}_*\m_p)^o$ is the connected component of the identity.
\end{proof}

\begin{Remark}
  We will eventually show that if $\pi\colon X\to S_{\Lambda}$ is the universal marked supersingular K3 surface, then the subgroup $(\R^2\pi^{\fl}_*\m_p)^o\subset \R^2\pi^{\fl}_*\m_p$ is the connected component of the identity. More generally, the same is true for any marked family such that the marking is generically an isomorphism.
\end{Remark}

\begin{Remark}\label{rem:flatduality}
  Let us consider the $k$-points of the diagram (\ref{eq:flatdualitydiagram}). Evaluating the short exact sequence (\ref{eq:equation1}) on $k$, and applying the isomorphism $\R^2\pi^{\fl}_*\m_p(k)\iso\R^1\pi^{(p)\etale}_*\nu(1)(k)$ of diagram \ref{eq:Lerayidentification}, we get a short exact sequence
  \[
    0\to (\R^2\pi^{\fl}_*\m_p)^o(k)\to\R^2\pi^{\fl}_*\m_p(k)\to\sD(k)\to 0
  \]
  Recall that under the identification $\R^2\pi^{\fl}_*\m_p(k)=\H^2(X,\m_p)$, the subgroup $(\R^2\pi^{\fl}_*\m_p)^o(k)$ corresponds to the group $\U^2(X,\m_p)$ of (\ref{eq:artinsshortexactsequence}), and $\sD(k)$ corresponds to $\D^2(X,\m_p)$. Finally, note that as the marking is an isomorphism, the sheaf $\sW$ in (\ref{eq:flatdualitydiagram}) vanishes. Thus, the diagram (\ref{eq:flatdualitydiagram}) on $k$-points recovers (part of) diagram (\ref{eq:artinsflatduality}). We will find a strengthened form of this observation in Remark \ref{rem:strongflatduality}.
\end{Remark}


\subsection{Twisted K3 crystals}\label{sec:twistedk3crystals}

Let $W=W(k)$ be the ring of Witt vectors of $k$, $K=W[\frac{1}{p}]$
its field of fractions, and $F_W\colon W\to W$ the homomorphism
induced by the Frobenius on $k$. We begin this section by recalling
from \cite{Ogus78} the definition of a K3 crystal over $W$, and their
connection to characteristic subspaces. We then introduce a
crystalline analog of the Hodge-theoretic B-fields studied in
\cite{HS04}. Unlike in the complex case, crystalline B-fields satisfy
a non-trivial relation, given in Lemma
\ref{lem:description}. Accordingly, they seem to possess a somewhat
richer structure than those over the complex numbers. In particular,
in the supersingular case they have nontrivial moduli. Using these, we
show how to associate a K3 crystal to a pair $(X,\alpha)$, where $X$
is a K3 surface and $\alpha\in \H^2(X,\m_p)$ (we will prefer to
express the data $(X,\alpha)$ as a $\m_p$-gerbe $\sX\to X$, although
in this section this is just language). This construction is the
crystalline analog of the twisted Hodge structures defined in
\cite{HS04}. We then discuss these constructions in the relative
setting.

Suppose that $\pi\colon X\to\Spec k$ is a K3 surface over $k$. The
second crystalline cohomology group $\H^2(X/W)$ of $X$ is a free
$W$-module of rank 22, which is equipped with a $F_W$-linear
endomorphism $\Phi\colon \H^2(X/W)\to \H^2(X/W)$ induced by the
absolute Frobenius on $X$. Such a structure is called an
$F$\textit{-crystal} over $W$. It also is equipped with a perfect
pairing
\[
  \H^2(X/W)\otimes_W\H^2(X/W)\to W
\]
which satisfies a certain compatibility with $\Phi$. Following Ogus
\cite{Ogus78}, we abstract these properties in the following
definition.
\begin{Definition}\label{def:K3crystal}
  A \textit{K3 crystal of rank $n$} is a $W$-module $H$ of rank $n$
  equipped with a $F_W$-linear endomorphism $\Phi\colon H\to H$ and a
  symmetric bilinear form $H\otimes_WH\to W$ such that
  \begin{enumerate}
      \item $p^2H\subset\Phi(H)$,
      \item $\Phi\otimes k$ has rank 1,
      \item the pairing $\langle\_,\_\rangle$ is perfect, and
      \item $\langle\Phi(x),\Phi(y)\rangle=p^2F_W\langle x,y\rangle$ for all $x,y\in H$.
  \end{enumerate}
  A K3 crystal $H$ is \textit{supersingular} if $H$ is isogenous to a
  crystal such that $\Phi$ acts as multiplication by $p$. In this
  case, we use the notation $\varphi=p^{-1}\Phi$.

\end{Definition}
\begin{Remark}\label{rem:1redundant}
  Condition (1) is equivalent to the existence of a map
  $V\colon H\to H$ such that $\Phi\circ V=V\circ\Phi=p^2$ (see
  \cite[Proposition 1.6.4]{Ogus77}). If $H$ is supersingular, then
  Corollary 3.8 of \cite{Ogus78} shows that $(2)\implies(1)$.
\end{Remark}
\begin{Definition}
  The \textit{Tate module} of a K3 crystal $H$ is
  \[
    T_H=\left\{h\in H|\Phi(h)=ph\right\}
  \]
  This has a natural structure of $\bZ_p$-module, and is equipped with the restriction of the bilinear form on $H$. If $H$ is supersingular, then by Proposition 3.13 of \cite{Ogus78} the $p$-adic ordinal of the discriminant of $T_H$ is equal to $2\sigma_0$ for some integer $\sigma_0\geq 1$, called the \textit{Artin invariant} of $H$.
\end{Definition}
Supersingular K3 crystals give rise to characteristic subspaces via the following procedure. If $H$ is a supersingular K3 crystal with Tate module $T$, then we have a chain of inclusions
\[
  T\otimes W\subset H\subset T^*\otimes W
\]
\begin{Proposition}\cite[Theorem 3.20]{Ogus78}\label{prop:periodmorphism}
  If $H$ is a supersingular K3 crystal, then $T^*/T$ is a vector space over $\bF_p$ of dimension $2\sigma_0$ whose induced bilinear form is non-degenerate and non-neutral, and the image $\overline{H}$ of $H$ in $T^*\otimes W/T\otimes W=(T^*/T)\otimes k$ is a strictly characteristic subspace, as is $K_H=\varphi^{-1}(\overline{H})$. 
\end{Proposition}
In fact, Ogus shows that this procedure is reversible, so the above correspondence gives an equivalence between certain appropriately defined categories of supersingular K3 crystals and strictly characteristic subspaces.
\begin{Remark}\label{rem:periodmorphism1}
  If $H=\H^2(X/W)$, where $X$ is a supersingular K3 surface, then the characteristic subspace $\varphi(K_H)$ produced by Proposition \ref{prop:periodmorphism} is essentially the same as the characteristic subspace appearing in the definition of the period morphism. Indeed, if $H$ is any supersingular K3 crystal, then the image of $H$ under the map
  \begin{equation}\label{eq:MAP1}
    H\to\dfrac{T^*\otimes W}{T\otimes W}\xrightarrow{\cdot p}\dfrac{pT^*\otimes W}{pT\otimes W}\subset\dfrac{T\otimes W}{pT\otimes W}=T\otimes k
  \end{equation}
  is equal to the kernel of the map $T\otimes k\to H\otimes k$. Because 
  $\H^3(X/W)$ is torsion free there is a canonical identification 
  \begin{equation}\label{eq:crystallinederham}
    \H^2(X/W)\otimes k\iso \H^2_{dR}(X/k)
  \end{equation}
  of the reduction modulo $p$ of the second crystalline cohomology group of $X$ and the second de Rham cohomology 
  group of $X$ (see Summary 7.26 in \cite{BO78}). By the Tate conjecture, $\Pic(X)\otimes\bZ_p\iso T$. Thus, the image of $H$ in $(T^*/T)\otimes k$ corresponds under (\ref{eq:MAP1}) to the kernel of the Chern class map
  \[
    \Pic(X)\otimes k\to \H^2_{dR}(X/k)
  \]
\end{Remark}
 

Let $\sX\to X$ be a $\m_p$-gerbe on the K3 surface $X$. We will show how to associate to $\sX$ a K3 crystal of rank 24. This construction is essentially the isomorphism of Proposition \ref{prop:isoofgroupschemes} translated into supersingular K3 crystals. We begin by recalling the Mukai crystal associated to a K3 surface $X$, as introduced in \cite{LO15}. Let $K(1)$ denote the $F$-isocrystal with underlying vector space $K$ and Frobenius action given by multiplication by $1/p$. For any $F$-isocrystal $M$ and integer $n$, we set $M(n)=M\otimes K(1)^{\otimes n}$.

\begin{Definition}\label{def:MukaiCrystal}
The \textit{Mukai crystal} of $X$ is the $W$-module
\[
  \tH(X/W)=\H^0(X/W)(\text{-}1)\oplus \H^2(X/W)\oplus \H^4(X/W)(1)
\]
equipped with the twisted Frobenius $\tPhi:\tH(X/W)\to\tH(X/W)$. We define the \textit{Mukai pairing} on $\tH(X/W)$ by
\[
  (a,b,c).(a',b',c')=-ac'+b.b'-a'c\in\H^4(X/W)=W
\]
\end{Definition}
\begin{Notation}
  Given two classes $(a,b,c)$ and $(a',b',c')$ in $\tH(X/W)$, we now
  have two possible operations: the cup product and
  the Mukai pairing. We will reserve the notation $(a,b,c).(a',b',c')$
  for the Mukai pairing and will use the juxtaposition
  $(a,b,c)(a',b',c')$ for the cup product. For example, we will
  often translate the lattice $\tH(X/W)$ inside the rational
  cohomology $\tH(X/K)$ by taking the cup product with a class of the
  form $e^B=(1,B,B^2/2)$, and we will write this as $e^B\tH(X/W)\subset\tH(X/K)$.
\end{Notation}

Both $\H^0(X/W)$ and $\H^4(X/W)$ are canonically isomorphic (as $W$-modules) to $W$. Under these identifications, the twisted Frobenius action is given by
\[
  \tPhi(a,b,c)=(pF_W(a),\Phi(b),pF_W(c))
\]
It follows from the definitions that $\tH(X/W)$ is a K3 crystal of rank 24. Because $\H^0(\text{-}1)$ and $\H^4(1)$ have slope 1, $\tH(X/W)$ is supersingular if and only if $X$ is supersingular.

Recall the identifications of diagram (\ref{eq:Lerayidentification}). By Proposition \ref{prop:ogusexactsequence}, there is an exact sequence
\begin{equation}\label{eq:FlatCohomologytoDeRham}
  0\to \H^2(X,\m_p)\xrightarrow{d\log} F^1_H\cap F^1_C(k)\xrightarrow{C\circ\pi_C-\pi_H}F^1_H/F^2_H(k)
\end{equation}
of abelian groups. In particular, we get an injective homomorphism
\[
  \H^2(X,\m_p)\hto F^1_H\cap F^1_C(k)\hto \H^2_{dR}(X/k)
\]
and a diagram
\[
  \begin{tikzcd}[column sep=small]
    &\H^2(X/W)\arrow[two heads]{d}{\mod p}\\
    \H^2(X,\m_p)\arrow[hook]{r}{d\log}&\H^2_{dR}(X/k)
  \end{tikzcd}
\]
where the vertical arrow is induced by the canonical identification (\ref{eq:crystallinederham}) of the reduction 
modulo $p$ of the crystalline cohomology with the de Rham cohomology.

The following definition is the crystalline analog of the Hodge theoretic B-fields defined in \cite{HS04}.
\begin{Definition}
  An element $B=\frac{a}{p}\in \H^2(X/K)$ is a \textit{B-field} if $a\in \H^2(X/W)$ and the image of $a$ in $\H^2_{dR}(X/k)$ lies in the image of $\H^2(X,\m_p)$. We write $\alpha_B$ for the unique element of $\H^2(X,\m_p)$ such that $d\log(\alpha_B)\equiv a\mod p$, and we say that $B$ is a \textit{B-field lift} of $\alpha_B$. 
\end{Definition}

\begin{Remark}
  In this work we only discuss B-fields associated to $p$-torsion Brauer 
  classes. Using the de Rham-Witt theory, one can make a similar definition that works for $p^n$-torsion classes as well. As the Brauer group of a supersingular K3 surface is $p$-torsion, the mod $p$ theory presented here suffices for applications to supersingular K3 surfaces.
\end{Remark}

We will give an alternative characterization of B-fields that uses only the crystal structure on $\H^2(X/W)$. The Frobenius $\Phi\colon \H^m(X/W)\to \H^m(X/W)$ induces filtrations
\begin{align*}
  M^i\H^m(X/W)&=(p^{-i}\Phi)^{-1}(\H^m(X/W))\\
  N^{i}\H^m(X/W)&=p^{-(m-i)}\Phi(M^{m-i}\H^m(X/W))
\end{align*}
Note that there is an isomorphism $p^{-i}\Phi\colon M^i\H^m\to N^{m-i}\H^{m}$. Suppose that $\H^{m+1}(X/W)$ is torsion free, and consider the natural map
\[
  \rho\colon \H^m(X/W)\to \H^m_{dR}(X/k)
\]
given by reduction modulo $p$. The following theorem of Mazur relates the image of these filtrations under $\rho$ to the Hodge and conjugate filtrations on de Rham cohomology.

\begin{Theorem}[\cite{BO78}, Theorem 8.26 and Lemma 8.30]\label{thm:Mazur}
  Suppose that $X$ is a smooth proper variety satisfying $(*)$. If $\H^{m+1}(X/W)$ is torsion free, then
  \begin{enumerate}
    \item the image of $M^i\H^m(X/W)$ under $\rho$ is $F^i_H\H^m_{dR}(X/k)$, 
    \item the image of $N^i\H^m(X/W)$ under $\rho$ is $F^i_C\H^m_{dR}(X/k)$, and
    \item the following diagram commutes
  \[
    \begin{tikzcd}
      M^i\H^m(X/W)\arrow{d}[swap]{p^{-i}\Phi}\arrow{r}{\rho}&F^i_H\H^m_{dR}(X/k)\arrow{r}{\pi_H}&\gr^i_{F_H}\H^m_{dR}(X/k)\arrow{d}{C^{-1}}\\
      N^{m-i}\H^m(X/W)\arrow{r}{\rho}&F^{m-i}_C\H^m_{dR}(X/k)\arrow{r}{\pi_C}&\gr^{m-i}_{F_C}\H^m_{dR}(X/k)
    \end{tikzcd}
  \]
  \end{enumerate}
\end{Theorem}
\begin{Lemma}\label{lem:description}
  An element $B\in p^{-1}\H^2(X/W)$ is a B-field if and only if
  \[
    B-\varphi(B)\in \H^2(X/W)+\varphi(\H^2(X/W))
  \]
\end{Lemma}
\begin{proof}
  Suppose that $B=\frac{a}{p}$ is a B-field. By Lemma \ref{lem:closeddifferentialforms}, $\rho(a)\in F^1_H\cap F^1_C$, which implies that $a\in M^1\cap N^1$. By the exact sequence (\ref{eq:FlatCohomologytoDeRham}), we know that $C\circ\pi_C(\alpha)=\pi_H(\alpha)$. Part (3) of Theorem \ref{thm:Mazur} then shows that $\pi_C(\rho(a-\varphi(a)))=0$, which implies that $a-\varphi(a)\in N^2+pH$, so $B$ satisfies the claimed relation.
  
  Conversely, suppose that $B=\frac{a}{p}\in p^{-1}\H^2(X/W)$ is an element satisfying the relation. This implies that $\varphi(a)\in \H^2(X/W)$, so $a\in M^1$. Using that $p^2\H^2\subset\Phi(\H^2)$, we also get that $a\in\varphi(\H^2)$, so $a\in N^1$. Thus, $\rho(a)$ is contained in $F^1_H\cap F^1_C$. Part (3) of Theorem \ref{thm:Mazur} then implies that $\pi_C(\rho(a)-\rho(\varphi(a)))=0$, so $\rho(a)$ is in the image of $\H^2(X,\m_p)$, and hence $B$ is a B-field.
\end{proof}

Let $\sX\to X$ be a $\m_p$-gerbe with cohomology class $\alpha$, and let $B$ be a B-field lift of $\alpha$. Cupping with $e^B=(1,B,B^2/2)$ defines an isometry $\tH(X/K)\to\tH(X,K)$, given explicitly by
\[
  e^B(a,b,c)=\left(a,b+aB,c+b.B+a\frac{B^2}{2}\right)
\]
\begin{Definition}\label{def:twistedK3crystal}
The \textit{twisted Mukai crystal} associated to $\sX$ is the $W$-module
\[
  \tH(\sX/W)=e^B\tH(X/W)\subset \tH(X/K)
\]
\end{Definition}
Note that if $h\in \H^2(X/W)$ then $e^h\in \H^*(X/W)$, and therefore the submodule $\tH(\sX/W)$ is independent of the choice of B-field.
\begin{Notation}\label{notn:twistedcrystal-thang}
  For an integer $n$, we will write $\tH^{(n)}(\sX/W)$
  for $\tH(\sX'/W)$, where $\sX'$ is a $\m_p$ gerbe on $X$ whose
  cohomology class is $n[\sX]$. This is independent of the choice of $\sX'$.
\end{Notation}
We will show that $\tH(\sX/W)$ has a natural K3 crystal structure.

\begin{Lemma}\label{lem:frobeniusrestricts}
  The submodule $\tH(\sX/W)$ is preserved by the action of the Frobenius $\tPhi$ on $\tH(X/K)$.
\end{Lemma}
\begin{proof}
  We must show that $\tPhi(\tH(\sX/W))\subset\tH(\sX/W)$. Consider an element $e^B(a,b,c)\in\tH(\sX/W)$. A consequence of the Mukai twist is the useful relation  
  \[
    \tPhi(e^B(a,b,c))=e^{\varphi(B)}\tPhi (a,b,c)
  \]
  Thus, the lemma is equivalent to the statement that
  \[
    e^{\varphi(B)-B}\tPhi(a,b,c)\in\tH(X/W)
  \]
  for all $(a,b,c)\in\tH(X/W)$. Write $B'=\varphi(B)-B$. We have
  \[
    e^{B'}\tPhi(a,b,c)=\left(pF_W(a),\Phi(b)+pF_W(a)B',pF_W(c)+\Phi(b)\innerproduct     B'+pF_W(a)\frac{B'^2}{2}\right)
  \]
  By assumption, $B'\in \H^2(X/W)+\varphi(\H^2(X/W))$. This implies that $pB'\in \H^2(X/W)$, $\Phi(b)\innerproduct B'\in W$, and $p\frac{B'^2}{2}\in W$.
\end{proof}

\begin{Proposition}\label{prop:itsacrystal}
  The $W$-module $\tH(\sX/W)$ equipped with the endomorphism $\tPhi$ and the Mukai pairing is a K3 crystal of rank 24, which is supersingular if and only if $X$ is supersingular.
\end{Proposition}  
\begin{proof}
  Because cupping with $e^B$ is an isometry with respect to the Mukai pairing, conditions (3) and (4) are immediate. For condition (1), we must show that for all $(a,b,c)\in\tH(\sX/W)$
  \[
    e^{B-\varphi(B)}(p^2a,p^2b,p^2c)\in \tPhi(\tH(\sX/W))
  \]
  Because $\H^2(X/W)$ is a K3 crystal,
  \[
    p^2(B-\varphi(B))\in p^2\H^2(X/W)+p\varphi(\H^2(X/W))\subset\Phi(\H^2(X/W))
  \]
  and $p^2b\in\Phi(\H^2(X/W))$. To check condition (2), we must compute the image of $\tPhi\colon \tH(\sX/W)\to\tH(\sX/W)$ modulo $p\tH(\sX/W)$. This is isomorphic to the image of $e^{-B}\circ\tPhi$ modulo $p\tH(X/W)$. By assumption, $B'=\varphi(B)-B=h+\varphi(h')$ for some $h,h'\in \H^2(X/W)$. We compute
  \begin{align*}
    e^{-B}\tPhi(e^B(a,b,c))&=e^{B'}.\tPhi(a,b,c)\\
                     &\equiv e^{B'}(pF_W(a),0,0)+e^{B'}(0,\Phi(b),0)\\
                     &\equiv (0,\Phi(ah'),\Phi(ah')\innerproduct h)+(0,\Phi(b),\Phi(b)\innerproduct h)\\
                     &\equiv (0,\Phi(b+ah'),\Phi(b+ah')\innerproduct h)
  \end{align*}
  Because $\H^2(X/W)$ is a K3 crystal, the result follows.
\end{proof}

\begin{Definition}
  If $\sX\to X$ is a $\m_p$ gerbe on a supersingular K3 surface, we define the \textit{Artin invariant} $\sigma_0(\sX)$ of $\sX$ to be the Artin invariant of the supersingular K3 crystal $\tH(\sX/W)$.
\end{Definition}

\begin{Remark}
  Definition \ref{def:twistedK3crystal} makes sense in the abstract setting where $H$ is a K3 crystal and $B$ satisfies the conditions of Lemma \ref{lem:description}. Much of the rest of this section is valid in this generality as well. For the sake of exposition we have chosen to phrase our results in the geometric context.
\end{Remark}

We next introduce a twisted versions of the N\'{e}ron-Severi lattice.

\begin{Definition}\label{defn:ext-ns}
  If $\sX\to X$ is a $\m_p$-gerbe on a K3 surface, we define the
  \textit{extended N\'{e}ron-Severi group} of $X$ and $\sX$ by
  \[
    \tN(X)=\langle (1,0,0)\rangle\oplus N(X)\oplus\langle (0,0,1)\rangle
  \]
  and
  \[
    \tN(\sX)=(\tN(X)\otimes\bZ[\frac{1}{p}])\cap \tH(\sX/W)\subset\tH(X/K).
  \]
\end{Definition}
Note that $\tN(\sX)$ only depends upon the cohomology class
$[\sX]\in\H^2(X,\m_p)$. 
\begin{Notation}\label{notn:ns-thang}
  For an integer $n$, we will write $\tN^{(n)}(\sX)$
  for $\tN(\sX')$, where $\sX'$ is a $\m_p$ gerbe on $X$ whose
  cohomology class is $n[\sX]$. This is independent of the choice of $\sX'$.
\end{Notation}

In the supersingular case, we can give a very explicit presentation of $\tN(\sX)$.

\begin{Lemma}\label{lem:transcendentalgerbes}
  If $X$ is a supersingular K3 surface, then any $\alpha\in \H^2(X,\m_p)$ can be written as $\alpha=\alpha'+\beta$, where $\alpha'$ is transcendental (see Definition \ref{def:transcendental}) and $\beta$ is in the image of the boundary map $\H^1(X,\bG_m)\to \H^2(X,\m_p)$.
\end{Lemma}
\begin{proof}
  Equip $X$ with the tautological marking by $\Lambda=\Pic(X)$. Consider the short exact sequence
  \[
    0\to (\R^2\pi^{\fl}_*\m_p)^o(k)\to\R^2\pi^{\fl}_*\m_p(k)\to\sD(k)\to 0
  \]
  described in Remark \ref{rem:flatduality}. By Lemma \ref{lem:transcendentalequalso}, the subgroup $(\R^2\pi^{\fl}_*\m_p)^o(k)$ consists of exactly the transcendental classes. Because the marking is an isomorphism, $\Lambda_1\cong\sD(k)$, and so every element of $\sD(k)$ lifts to an element of $\R^2\pi^{\fl}_*\m_p(k)=\H^2(X,\m_p)$ that is in the image of the boundary map $\H^1(X,\bG_m)\to \H^2(X,\m_p)$. This gives the result.
\end{proof}

\begin{Proposition}\label{rem:nslattice}
  Suppose that $X$ is supersingular, let $\alpha\in \H^2(X,\m_p)$ be the cohomology class of $\sX\to X$, and let $\alpha_0$ be the image of $\alpha$ in $\Br(X)$. Suppose that $\alpha=\alpha'+\beta$, where $\alpha'$ is transcendental and $\beta$ has trivial Brauer class. If $\beta$ is the image of a line bundle $\sL$ under the boundary map, and $t=c_1(\sL)$, then
  \[
    \tN(\sX)=
      \begin{cases}
        \langle (1,\frac{t}{p},\frac{(t/p)^2}{2})\rangle\oplus\langle(0,D,D\innerproduct\frac{t}{p})\rangle\oplus\langle (0,0,1)\rangle&\mbox{ if }0=\alpha_0\in\Br(X),\\
        \langle (p,t,\frac{t^2}{2p})\rangle\oplus\langle(0,D,D\innerproduct\frac{t}{p})\rangle\oplus\langle (0,0,1)\rangle&\mbox{ if }0\neq\alpha_0\in\Br(X)
      \end{cases}
  \]
  In particular, if $\alpha$ is transcendental and $\alpha_0\neq 0$, then
  \[
    \tN(\sX)=\langle (p,0,0)\rangle\oplus N(X)\oplus\langle (0,0,1)\rangle
  \]
\end{Proposition}
\begin{proof}
  Suppose that $\alpha$ is transcendental, and that the class of $\alpha$ in $\Br(X)$ is non-zero. Let $B$ be a B-field lift of $\alpha$. An element $(a,b,c)\in\tN(X)\otimes\bQ$ lies in $\tN(\sX)$ if and only if
  \[
    e^{-B}(a,b,c)=(a,b-aB,c-b\innerproduct B+a\frac{B^2}{2})
  \]
  is in $\tH(X/W)$. This implies $a\in\bZ$. Note that $b=\frac{l}{p}$ for some $l\in N(X)$. If $b-aB=h$ for some $h\in \H^2(X/W)$, then $aB=\frac{l}{p}-h$. Because the class of $\alpha$ in $\Br(X)$ is non-zero, this implies that $a$ is divisible by $p$. Hence, $b=aB+h\in N(X)$. Finally, as $\alpha$ is transcendental, $b\innerproduct B\in W$ and $p\frac{B^2}{2}\in W$. Therefore, $c\in\bZ$. We have shown that in this case
  \[
    \tN(\sX)=\langle (p,0,0)\rangle\oplus N(X)\oplus\langle (0,0,1)\rangle
  \]
  Next, let $\alpha\in \H^2(X,\m_p)$ be arbitrary. By Lemma \ref{lem:transcendentalgerbes}, we can write $\alpha=\alpha'+\beta$, where $\alpha'$ is transcendental and $\beta$ is the image of a line bundle $\sL$ under the boundary map $\H^1(X,\bG_m)\to \H^2(X,\m_p)$. Let $B'$ be a B-field lift of $\alpha'$ and $t=c_1(\sL)$. Then $B=B'+\frac{t}{p}$ is a B-field lift of $\alpha$, and we have an isomorphism
  \[
    e^{\frac{t}{p}}\colon e^{B'}\tH(X/W)\iso e^B\tH(X/W)
  \]
  of K3 crystals. This induces an isomorphism on Tate modules, and one checks that it also induces an isomorphism on extended N\'{e}ron-Severi groups. The result therefore follows by the previous case.
\end{proof}

 Write $T(X), \tT(X)$, and $\tT(\sX)$ for the Tate modules of $\H^2(X/W), \tH(X/W)$, and $\tH(\sX/W)$. 
 
\begin{Proposition}\label{prop:twistedNeronseveri}
  The inclusion $\tN(\sX)\hto \tH(\sX/W)$ factors through $\tT(\sX)$. If $X$ is supersingular, the induced map
  \[
    \tN(\sX)\otimes\bZ_p\to\tT(\sX)
  \]
  is an isomorphism.
\end{Proposition}
\begin{proof}
  The first claim follows from the Tate twists in the definition of $\tH(X/W)$. Suppose that $X$ is supersingular. It follows from the definitions that
  \[
    \tT(\sX)=(\tT(X)\otimes\bQ)\cap\tH(\sX/W)
  \]
  By the Tate conjecture (or by assumption), the map $N(X)\otimes\bZ_p\iso T(X)$ is an isomorphism. The calculations of Proposition \ref{prop:twistedNeronseveri} then apply as written to show that the $\bZ_p$-span of the given vectors is equal to $\tT(\sX)$, which gives the result.
\end{proof}

\begin{Corollary}
  If $\sX\to X$ is a $\m_p$-gerbe over a supersingular K3 surface, then $\sigma_0(\sX)=\sigma_0(X)+1$ if the Brauer class of $\sX$ is non-zero, and $\sigma_0(\sX)=\sigma_0(X)$ otherwise.
\end{Corollary}

We will next discuss these constructions in the relative setting. Given a relative K3 surface $\pi\colon X\to S$ and a $\m_p$-gerbe $\sX\to X$, we would like to define a relative twisted Mukai crystal $\tH(\sX/S)$. The correct notion of a crystal over a non-perfect base carries significant technicalities. We will restrict our attention to the following situation.
\begin{Situation}\label{situation:liftedsituation}
  Suppose that $S=\Spec A$ is affine, where $A$ is a smooth $k$-algebra. Fix a smooth, $p$-adically complete lift $S'=\Spec A'$ of $S$ over $W$, together with a lift $F_{S'}$ of the absolute Frobenius of $S$ (that is, a morphism $F_{S'}\colon S'\to S'$ that reduces to $F_S$ modulo $p$).
\end{Situation}
If $S$ is perfect, then there is a unique such choice of $A'$, given by the Witt vectors $A'=W(A)$. In general, there will be many choices of $S'$. Note however that our assumption that $S$ is smooth ensures that they are at least locally isomorphic.
\begin{Definition}\label{def:crystal}
  Suppose that we are in Situation \ref{situation:liftedsituation}. A \textit{crystal on $S/W$ (with respect to $S'$)} is a pair $(M,\nabla)$, such that
  \begin{enumerate}
      \item $M$ is a finitely generated $p$-adically complete projective $A'$-module, and
      \item $\nabla\colon M\to M\,\widehat{\otimes}\,\widehat{\Omega}^1_{S'/W}$ is a connection such that
      \item $\nabla$ is integrable and topologically quasi-nilpotent.
  \end{enumerate}
  An \textit{F-crystal on $S/W$ (with respect to $(S',F_{S'})$)} is a triple $(M,\nabla,\Phi)$ where $(M,\nabla)$ is a crystal on $S/W$ and 
  \[
    \Phi\colon F_{S'}^*(M,\nabla)\to (M,\nabla)
  \]
  is a horizontal\footnote{Recall that if $(\sE,\nabla)$ and $(\sF,\nabla')$ are modules with connection, a map $f\colon \sE\to\sF$ is \textit{horizontal} if $\nabla'\circ f=(f\otimes\id)\circ\nabla$.} map of $A'$-modules that becomes an isomorphism after inverting $p$.
\end{Definition}

For the definitions of these terms, we refer the reader to \cite[07GI]{stacks-project}.

\begin{Remark}
  This is essentially the same as Definition 1.1.3 of \cite{Ogus77}, although we have chosen to work with complete objects, while Ogus works with formal objects on formal schemes. In addition, Ogus omits the nilpotence condition, as he works only with $F$-crystals, where it is automatic by \cite[Corollary 1.7]{Ogus77}.
\end{Remark}

\begin{Remark}
  As explained in Remark 1.8 of \cite{Ogus77}, the definition of an $F$-crystal 
  on $S/W$ is in a certain sense independent of  the choice of $S'$ and 
  $F_{S'}$. Suppose that $(M,\nabla)$ is a crystal on $S/W$ with respect to the 
  lifting $S'$. If $T'$ is another smooth, complete lift of $S$ over $W$, then 
  locally $S'$ and $T'$ are isomorphic over $S$. The connection $\nabla$ 
  induces an isomorphism between the pullbacks of $M$ along any two such local 
  isomorphisms, and the integrability assumption implies that these isomorphisms 
  satisfy the cocycle condition. Thus, $(M,\nabla)$ induces in a canonical way 
  a crystal on $S/W$ with respect to $T'$.
  
  If $(M,\nabla,\Phi)$ is an $F$-crystal on $S/W$ with respect to $(S',F_{S'})$, and $G_{S'}$ is another choice of lifting of the Frobenius, then the connection induces a horizontal isomorphism $\varepsilon\colon G_{S'}^*(M,\nabla)\to F_{S'}^*(M,\nabla)$. The triple $(M,\nabla,\Phi\circ\varepsilon)$ is then an $F$-crystal on $S/W$ with respect to $(S',G_{S'})$.
  
  This independence is explained by the site-theoretic approach to crystals and crystalline cohomology, as developed in generality in \cite{MR0384804}. See \cite[Proposition 6.8]{BO78} and \cite[07JH]{stacks-project} for the equivalence of these approaches.
\end{Remark}

\begin{Definition}\cite[Section 5]{Ogus78}\label{def:relativecrystal}
  A \textit{K3 crystal on }$S/W$ (with respect to $(S',F_{S'})$)  of rank $n$ is an F-crystal $(H,\nabla,\Phi)$ on $S/W$ (with respect to $(S',F_{S'})$) where $H$ is an $A'$-module of rank $n$, endowed with a horizontal symmetric pairing $H\otimes H\to\cO_{S'}$ such that
  \begin{enumerate}
    \item there exists a horizontal map $V\colon (H,\nabla)\to F_{S'}^*(H,\nabla)$ 
    satisfying $\Phi\circ V=V\circ\Phi=p^2$,
    \item the $\gr^\bullet_FH$ are locally free $\cO_{S}$-modules and $\gr^1_FH$ has rank one,\footnote{For the definition of the Hodge filtration on an abstract $F$-crystal we refer to \cite{Ogus77}.}
    \item the pairing is perfect, and
    \item $\langle \Phi(x),\Phi(y)\rangle=p^2F_{S'}^*\langle x,y\rangle$ for any two sections $x,y$ of $F^*_{S'}H$.
  \end{enumerate}
  We say that $H$ is \textit{supersingular} if for all geometric points $s\to 
  S$ the restricted crystal $H(s)$ is a supersingular K3 crystal in the sense of Definition \ref{def:K3crystal}.
\end{Definition}


\begin{Definition}
  Suppose that we are in Situation \ref{situation:liftedsituation}, and that $X\to S$ is a relative K3 surface. We let
    \[
      (\H^2(X/S'),\nabla_{S'},\Phi_{S'})
    \]
  be the $F$-crystal on $S/W$ relative to $(S',F_{S'})$ corresponding to the second crystalline cohomology of $X\to S$ with its canonical $F$-structure induced by the Frobenius. When equipped with the pairing induced by the cup product, this is a K3 crystal of rank 22, which is supersingular if and only if $X\to S$ is supersingular. We let 
  \[
    (\tH(X/S'),\tnabla_{S'},\tPhi_{S'})
  \]
  be the $F$-crystal on $S/W$ relative to $(S',F_{S'})$ corresponding to the total crystalline cohomology of $X\to S$ with the twisted $F$-structure (as in Definition \ref{def:MukaiCrystal}). When equipped with the Muaki pairing, this is a K3 crystal of rank 24, which is supersingular if and only if $X\to S$ is supersingular.
\end{Definition}

\begin{Remark}
  As explained in \cite{Ogus78}, the relative analogs of Proposition \ref{prop:periodmorphism} and Remark \ref{rem:periodmorphism1} hold as well.
\end{Remark}

Set $\H^2(X/S'_K)=\H^2(X/S')\otimes_{W}K$ and $\tH(X/S'_K)=\tH(X/S')\otimes_WK$. We will extend Definition \ref{def:twistedK3crystal} to the relative setting. We have maps
\[
  \R^2\pi^{\fl}_*\m_p\too{\Upsilon}\R^1\pi^{(p)\etale}_*\nu(1)\xrightarrow{d\log}\R^2\pi^{\etale}_*\Omega^{\bullet}_{X/S}
\]
There is a canonical isomorphism
\[
  \H^2(X/S')\otimes_{A'}A\iso \H^2_{dR}(X/S)
\]
As in the case when $S=\Spec k$, we find a diagram
\begin{equation}\label{eq:diagram1b}
  \begin{tikzcd}[column sep=small]
    &\Gamma(S',\H^2(X/S'))\arrow[two heads]{d}\\
    \Gamma(S,\R^2\pi^{\fl}_*\m_p)\arrow{r}&\Gamma(S,\H^2_{dR}(X/S))
  \end{tikzcd}
\end{equation}
although the horizontal map is no longer necessarily injective.

\begin{Definition}
  Let $X\to S$ be a relative K3 surface and $\sX\to X$ a $\m_p$-gerbe with cohomology class $\alpha\in\Gamma(S,\R^2\pi^{\fl}_*\m_p)$. Suppose that we are in Situation \ref{situation:liftedsituation}. Let $\sB\in\Gamma(S',\H^2(X/S'_K))$ be a section such that $p\sB\in\Gamma(S',\H^2(X/S'))$ and the image of $p\sB$ in $\Gamma(S,\H^2_{dR}(X/S))$ is equal to the image of $\alpha$. Consider the composition
  \[
    \tH(X/S')\xrightarrow{h\mapsto e^{\sB}\otimes h}\tH(X/S'_K)\otimes\tH(X/S')\to\tH(X/S'_K)
  \]
  where the second map is given by the cup product. We define
  \[
    \tH(\sX/S')=e^\sB \tH(X/S')\subset \H^*(X/S'_K)
  \]
  to be its image. As before, note that this does not depend on our choice of $\sB$. This definition makes sense more generally for a cohomology class $\alpha\in \Gamma(S,\R^2\pi^{\fl}_*\m_p)$ that is not in the image of $\H^2(X,\m_p)$, and hence may only be represented by a gerbe flat locally on $S$.
\end{Definition}

\begin{Proposition}\label{prop:compatiblewithbasechange2}
  The $A'$-module $\tH(\sX/S')$ is of formation compatible with base change, in 
  the following sense. Suppose $B$ is a smooth $k$-algebra and $B'$ is a 
  $p$-adically complete lift of $B$ to $W$. Write $T=\Spec B$ and $T'=\Spec 
  B'$. Given a commutative diagram 
  \[
    \begin{tikzcd}
      T'\arrow{r}\arrow{d}&S'\arrow{d}\\
      T\arrow{r}&S
    \end{tikzcd}
  \]
  of $W$-schemes, the natural map
  \[
    \tH(\sX_T/T')\to\tH(\sX/S')\otimes_{A'}B'
  \]
  is an isomorphism of $B'$-modules.
\end{Proposition}
\begin{proof}
  The crystalline cohomology of a relative K3 surface is of formation 
  compatible with base change, in the sense that the natural map
  \[
    \tH(X/S')\otimes_{A'}B'\to\tH(X_T/T')
  \]
  is an isomorphism. Let $\sB_T=\sB\otimes 1\in\tH(X/S')\otimes_{A'}B'$ be the pullback of $\sB$. The cup product is compatible with base change as well, so
  \[
    \tH(\sX/S')\otimes_{A'}B'=e^{\sB}\tH(X/S')\otimes_{A'}B'\iso e^{\sB_T}\tH(X_T/T')
  \]
  Let $\alpha_T$ be the cohomology class of the $\m_p$-gerbe $\sX_T\to X_T$. It follows that $\sB_T$ is a B-field lift of $\alpha_T$, and therefore that $e^{\sB_T}\tH(X_T/T')=\tH(\sX_T/T')$.
\end{proof}

To endow $\tH(\sX/S')$ with the structure of a K3 crystal, we need to give it a 
connection. Under our assumption that $S$ is smooth, the second de Rham cohomology of $X\to S$ is equipped with the 
Gauss-Manin connection
\[
  \nabla_0\colon \H^2_{dR}(X/S)\to \H^2_{dR}(X/S)\otimes\Omega^1_{S/k}
\]
For any $D\in\Gamma(S,T^1_{S/k})$, composing $\nabla_0$ with $D$ gives a map $\nabla_0(D)\colon \H^2_{dR}(X/S)\to \H^2_{dR}(X/S)$.
Morevoer, via the isomorphism
\[
  \H^2(X/S')\otimes_W k\iso \H^2_{dR}(X/S),
\]
the connection 
$\nabla_{S'}$ reduces to $\nabla_0$.

\begin{Lemma}\label{lem:alphaflat}
  If $\alpha\in \Gamma(S,\R^2\pi^{\fl}_*\m_p)$ is a cohomology class, then $\nabla_0(\beta)=0$, where $\beta=d\log\circ\Upsilon(\alpha)$ is the image of $\alpha$ in $\Gamma(S,\H^2_{dR}(X/S))$ under the horizontal map of (\ref{eq:diagram1b}).
\end{Lemma}
\begin{proof}
We recall from Proposition \ref{prop:flattoetalecohomology} that $\Upsilon$ is defined to be the composite of morphisms
\[
  \R^2\pi^{\fl}_*\m_p\to\R^1\pi^{\etale}_*(\cO_X^\times/\cO_X^{\times p})\to\R^1\pi^{(p)\etale}_*\nu(1)
\]
We consider the induced maps on global sections
\[
  \Gamma(S,\R^2\pi^{\fl}_*\m_p)\to\Gamma(S,\R^1\pi^{\etale}_*(\cO_X^\times/\cO_X^{\times p}))\to\Gamma(S,\R^1\pi^{(p)\etale}_*\nu(1))
\]
We will show that the image of any global section of $\R^1\pi^{\etale}_*(\cO_X^\times/\cO_X^{\times p})$ is horizontal (in fact, this is equivalent to the result, because under the assumption that $S$ is smooth the first map is an isomorphism on global sections). We will do this using the explicit description of the Gauss-Manin connection in terms of cocycles due to Katz. We will follow the presentation in Section 3 of \cite{MR0291177}, and also refer to the slightly different formulation in \cite{MR0237510}.

Let $r$ be the relative dimension of $X\to S$ (in our case $r=2$). Choose a finite covering $\left\{U_{\alpha}\right\}$ of $X$ by open affines such that each $U_{\alpha}$ is \'{e}tale over $\bA^r_S$, and such that on each $U_{\alpha}$ the sheaf $\Omega^1_{X/S}$ is a free $\cO_X$-module with basis $\left\{dx_1^{\alpha},\dots,dx_r^{\alpha}\right\}$. We consider the double complex $\sC^{\bullet,\bullet}$ with terms
\[
  \sC^{p,q}=\sC^p(\left\{U_{\alpha}\right\},\Omega^q_{X/S})
\]
where $\sC^p(\left\{U_{\alpha}\right\},\Omega^q_{X/S})$ is the set of alternating \v{C}ech cochains on the covering $\left\{U_{\alpha}\right\}$. This is equipped with a vertical differential $d\colon \sC^{p,q}\to\sC^{p,q+1}$ and a horizontal differential $\delta\colon \sC^{p,q}\to\sC^{p+1,q}$, defined in \cite{MR0237510}. We consider the associated total complex $\sK^{\bullet}$ with terms
\[
  \sK^{r}=\sum_{p+q=n}\sC^{p,q}
\]
and differential $d+\delta$. As explained in \cite{MR0237510}, the hypercohomology of $\sK^{\bullet}$ computes the de Rham cohomology of $X\to S$.

Let $D\in\Der_k(\cO_S,\cO_S)$ be a $k$-derivation of $\cO_S$. For each index $\alpha$, let $D_{\alpha}\in\Der_k(U_{\alpha},U_{\alpha})$ be the unique extension of $D$ which kills $dx_1^{\alpha},\dots,dx_r^{\alpha}$. In \cite{MR0291177}, Katz defines maps $\widetilde{D}\colon \sC^{p,q}\to\sC^{p,q}$ and $\lambda(D)\colon \sC^{p,q}\to\sC^{p+1,q-1}$, and shows that the sum
\[
  \widetilde{D}+\lambda(D)\colon \sK^{r}\to\sK^{r}
\]
upon passage to cohomology computes the map $\nabla_0(D)$.

We now make our computation. Our cohomology class in $\Gamma(S,\R^1\pi^{\etale}_*(\cO_X^\times/\cO_X^{\times p})$ may be represented by a cocycle $f_{i,j}\in\sC^1(\left\{U_{\alpha}\right\},\cO_X^\times/\cO_X^{\times p})$. After possibly shrinking our cover, we may find lifts $\zeta_{i,j}\in\Gamma(U_i\cap U_j,\cO_X^\times)$ of each of the $f_{i,j}$. Because we are in characteristic $p$, the quotients $d\log(\zeta_{i,j})$ are independent of our choice of lifts, and hence give rise to a cocycle
\[
  d\log(\zeta_{i,j})\in\sC^{1,1}
\]
that represents the image of our cohomology class in the second de Rham cohomology. Fix a $k$-derivation $D$ of $\cO_S$. We wish to show that the cocycle
\[
  (\widetilde{D}+\lambda(D))(d\log(\zeta_{i,j}))\in\sC^{1,1}\oplus\sC^{2,0}
\]
is in the image of $d+\delta$. Consider the cocycle
\[
  D_{i}\log(\zeta_{i,j})=\dfrac{D_i(\zeta_{i,j})}{\zeta_{i,j}}\in\sC^{1,0}
\]
Note that, because $D_i$ is a derivation, this is also independent of the choice of $\zeta_{i,j}$. One computes that
\[
  (d+\delta)(D_{i}\log(\zeta_{i,j}))=(\widetilde{D}+\lambda(D))(d\log(\zeta_{i,j}))
\]
which gives the result.
\end{proof}

\begin{Proposition}
  The submodule $\tH(\sX/S')\subset\tH(X/S'_K)$ is horizontal with respect to the connection $\tnabla_{S'}\otimes K$.
\end{Proposition}
\begin{proof}
  The assertion is equivalent to the statement that for all $(a,b,c)\in\tH(X/S')$,
  \[
    (e^\sB\otimes\id)\circ\tnabla_{S'}(e^{\text{-}\sB}(a,b,c))\in\tH(X/S')\,\widehat{\otimes}\,\widehat{\Omega}^1_{S'/W}
  \]
  By the compatibility of the connection with the cup product, we have
  \begin{align*}
    (e^\sB\otimes\id)\circ\tnabla_{S'}(e^{\text{-}\sB}(a,b,c))&=\tnabla_{S'}(a,b,c)+(e^\sB\otimes\id)((a,b,c)\tnabla_{S'}(e^{\text{-}\sB}))\\
                                                 &=\tnabla_{S'}(a,b,c)+(e^\sB\otimes\id)((a,b,c)(0,\nabla_{S'}(\text{-}\sB),\sB.\nabla_{S'}(\sB)))\\
                                                 &=\tnabla_{S'}(a,b,c)+(e^\sB\otimes\id)(0,a\nabla_{S'}(\text{-}\sB),a\sB.\nabla_{S'}(\sB)+b.\nabla_{S'}(\text{-}\sB))\\
                                                 &=\tnabla_{S'}(a,b,c)+(0,a\nabla_{S'}(\text{-}\sB),b.\nabla_{S'}(\text{-}\sB))
  \end{align*}
  By Lemma \ref{lem:alphaflat}, we have $\nabla_{S'}(\sB)\in \H^2(X/S')\,\widehat{\otimes}\,\widehat{\Omega}^1_{S'/W}$, and the result follows.
\end{proof}

Thus, we obtain a connection
\[
  \tnabla_{S'}\colon \tH(\sX/S')\to \tH(\sX/S')\,\widehat{\otimes}\,\widehat{\Omega}^1_{S'/W}
\]
on $\tH(\sX/S')$. It is immediate that it satisfies condition (3) of Definition \ref{def:crystal}. The remaining properties can be shown by the same methods as in the punctual case, using the relative version of Theorem \ref{thm:Mazur}. We will instead deduce them by reduction to the punctual case. Recall that given a closed point $s\in S$, there is a unique map $s'\colon (\Spec W,F_{W})\to (S',F_{S'})$ called the \textit{Teichm\"{u}ller lifting} of $s$. If $s'$ is a Teichm\"{u}ller lifting of a closed point $s$, then the restricted $F$-crystal
\[
  (s'^*\tH(X/S'),s'^*\tnabla_{S'}, s'^*\tPhi_{S'})
\]
is identified with $(\tH(X_{s}/W),\tPhi)$. 

\begin{Lemma}\label{lem:relativefrobeniusrestricts}
  The submodule $\tH(\sX/S')$ is preserved by the action of the Frobenius $\tPhi_{S'}$ on $\tH(X/S'_K)$.
\end{Lemma}
\begin{proof}
  We will show that the quotient
  \[
    M=(\tPhi_{S'}(F_{S'}^*\tH(\sX/S'))+\tH(\sX/S'))/\tH(\sX/S')
  \]
  of $A'$-modules is zero. Consider a closed point $s\in S$, with 
  Teichm\"{u}ller lift 
  \[
    s'\colon (\Spec W(s),F_{W(s)})\to (S',F_{S'})
  \]
  The restriction of $M$ to such a point vanishes by Lemma \ref{lem:frobeniusrestricts}. It follows that $M\otimes W/pW=0$. This means that multiplication by $p$ on $M$ is an isomorphism. Because $M$ is $p$-adically complete, we conclude that $M=0$.
\end{proof}

\begin{Proposition}
  The triple $(\tH(\sX/S'),\tnabla_{S'},\tPhi_{S'})$ is an $F$-crystal. Equipped with the Mukai pairing, it is a K3 crystal on $S/W$ with respect to $(S',F_{S'})$ of rank 24, which is supersingular if and only if $X\to S$ is supersingular.
\end{Proposition}
\begin{proof}
  That $(\tH(\sX/S'),\tnabla_{S'},\tPhi_{S'})$ is an $F$-crystal follows immediately from the fact that \newline$(\tH(X/S'),\tnabla_{S'},\tPhi_{S'})$ is an $F$-crystal.
  Next, we will check the conditions of Definition \ref{def:relativecrystal}. 
  Because $\tH(X/S')$ is a K3 crystal, it admits a map $\tV_{S'}\colon \tH(X/S')\to 
  F_{S'}^*\tH(X/S')$ satisfying condition (1). We will show that $\tV_{S'}$ 
  restricts to such a map on $\tH(\sX/S')$. Indeed, by Proposition 
  \ref{prop:itsacrystal}, this is true for the restriction of $\tH(\sX/S')$ to 
  any closed point of $S$. By the same argument as Lemma 
  \ref{lem:relativefrobeniusrestricts}, this implies the result. Because 
  $\gr^\bullet_FH$ is of formation compatible with base change, condition (2) 
  also follows from Proposition \ref{prop:itsacrystal}. Conditions (3) and (4) 
  hold because they are true for $\tH(X/S')$, and cupping with $e^\sB$ is an 
  isometry with respect to the Mukai pairing.
  
  Finally, supersingularity is by definition a condition on the fibers, so the final claim follows from Proposition \ref{prop:itsacrystal}.
\end{proof}


\begin{Remark}
  Note that we have only defined the crystal $\tH(\sX/S')$ in Situation \ref{situation:liftedsituation}. If $S$ is smooth but not necessarily affine, the compatibility of our construction with respect to base change shows that our locally defined crystals glue to a crystal on $S/W$, in the sense of \cite[07IS]{stacks-project}.
\end{Remark}

\subsection{Moduli of twisted supersingular K3 surfaces and the twisted period morphism}\label{sec:twistedperiodmorphism}

In this section we discuss the moduli space of marked twisted supersingular K3 surfaces. We then use the twisted K3 crystals of the previous section to define a twisted crystalline period morphism. We compute its differential, and show that it is \'{e}tale.

Let $\Lambda$ be a supersingular K3 lattice and let $\pi\colon X\to S_{\Lambda}$ denote the universal marked supersingular K3 surface. For psychological reasons, we adopt the following notation.
\begin{Definition}
    \begin{align*}
        \sS_{\Lambda}&=\R^2\pi^{\fl}_*\m_p\\
        \sS_{\Lambda}^o&=(\R^2\pi^{\fl}_*\m_p)^o
    \end{align*}
\end{Definition}

Recall from Definition \ref{def:definitioncirc} that $(\R^2\pi^{\fl}_*\m_p)^o\subset\R^2\pi^{\fl}_*\m_p$ is the subgroup whose fiber over a marked K3 surface is the group of transcendental classes.
\begin{Theorem}
  The functor $\sS_{\Lambda}^o$ is representable by an algebraic space that is locally of finite presentation, locally separated, and smooth over $\Spec k$. The morphism $\sS_{\Lambda}^o\to S_{\Lambda}$ is a smooth group space of relative dimension 1, and the connected component of any of its geometric fibers is isomorphic to $\bG_a$.
\end{Theorem}
\begin{proof}
  This follows from Theorem \ref{thm:representable} and Lemma \ref{lem:smoothofdimension1}.
\end{proof}

The $k$-points of the moduli space $\sS_{\Lambda}^o$ are given by pairs $(\sX,m)$ where $\sX$ is a $\m_p$-gerbe over a supersingular K3 surface $X$, and $m$ is a marking of $X$ such that $\sX$ is transcendental with respect to $m$. Using the assumption that $\sX$ is transcendental, we can give an alternate description of such a pair as a $\m_p$-gerbe $\sX$ equipped with a certain marking of the extended N\'{e}ron-Severi group of $\sX$.

\begin{Definition}
  An \textit{extended supersingular K3 lattice} is a free abelian group $\tLam$ of rank 24 equipped with an even symmetric bilinear form such that
  \begin{enumerate}
    \item $\disc(\tLam\otimes\bQ)=1$ in $\bQ^\times/\bQ^{\times 2}$,
    \item the signature of $\tLam\otimes\bR$ is $(2,22)$, and
    \item the discriminant group $\tLam^*/\tLam$ is $p$-torsion.
  \end{enumerate}
\end{Definition}
These lattices behave similarly to supersingular K3 lattices. In particular, if $\tLam$ is an extended supersingular K3 lattice, then $\tLam_0=p\tLam^*/p\tLam$ is a vector space over $\bF_p$ of dimension $2\sigma_0$ for some integer $1\leq\sigma_0\leq 11$, called the \textit{Artin invariant} of $\tLam$. This vector space has a natural bilinear form, which is non-degenerate and non-neutral.

\begin{Notation}\label{def:hyperbolicplane2}$ $\par\nobreak\ignorespaces
  \begin{itemize}
      \item We let $U_2$ denote the rank 2 lattice which is generated by two elements $e,f$ satisfying $e^2=f^2=0$ and $e\innerproduct f=-1$.
      \item Given a lattice $L$ and an integer $n$, we write $L(n)$ for the lattice with underlying abelian group $L$ but with the form multiplied by $n$. Thus, $U_2(p)$ denotes the lattice generated by two elements $e,f$ satisfying $e^2=f^2=0$ and $e\innerproduct f=-p$.
      \item We set $\tLam=\Lambda\oplus U_2(p)$, where $\Lambda$ is our fixed supersingular K3 lattice.
  \end{itemize}
\end{Notation}

\begin{Lemma}\label{lem:supersingularlattices}
  If $\Lambda$ is a supersingular K3 lattice of Artin invariant $\sigma_0$, then $\Lambda\oplus U_2$ and $\Lambda\oplus U_2(p)$ are extended supersingular K3 lattices of Artin invariants $\sigma_0$ and $\sigma_0+1$. Moreover, every extended supersingular K3 lattice is of this form for some $\Lambda$.
\end{Lemma}
\begin{proof}
  In \cite{Ogus78} it is shown that the Artin invariant determines a supersingular K3 lattice up to isometry. The proof implies the same for extended supersingular K3 lattices. It is then easy to compute that $\Lambda\oplus U_2$ and $\Lambda\oplus U_2(p)$ have the required properties.
\end{proof}
By the calculations of Proposition \ref{rem:nslattice}, we see that extended supersingular K3 lattices are exactly those lattices that occur as the extended N\'{e}ron-Severi group of a twisted supersingular K3 surface.

\begin{Lemma}\label{lem:inducedmarking}
  Let $X$ be a supersingular K3 surface with a marking $m\colon \Lambda\to\Pic(X)$. If $\sX\to X$ is a $\m_p$-gerbe that is transcendental with respect to $m$, then the map
  \[
    \tLam=\Lambda\oplus U_2(p)\to\tH(\sX/W)
  \]
  given by $e\mapsto (0,0,1)$ and $f\mapsto (p,0,0)$ factors through $\tN(\sX)$.
\end{Lemma}
\begin{proof}
  We must check that if $l\in\Lambda$ and $s\in\bZ$ then $(p,l,s)\in\tN(\sX)$. If $B$ is a B-field lift of $\alpha$, we compute
  \[
    e^{-B}(p,l,s)=(p,l-pB,s-l.B+p\frac{B^2}{2})
  \]
  As $\alpha$ is transcendental relative to the marking $m$, $l.B\in W$ and $p\frac{B^2}{2}\in W$, which gives the result.
\end{proof}

We will next define our twisted period morphism. We first briefly recall some material from \cite{Ogus78}. Let $T$ be a $\bZ_p$-lattice that is isometric to the Tate module of a supersingular K3 crystal of rank $n$ (in \cite{Ogus78} these are called ``K3 lattices''). For instance, we may take $\Lambda\otimes\bZ_p$ where $\Lambda$ is a supersingular K3 lattice (and $n=22$), or $\tLam\otimes\bZ_p$, where $\tLam$ is an extended supersingular K3 lattice (and $n=24$). Suppose that we are in Situation \ref{situation:liftedsituation}, and that $H_{S'}$ is a supersingular K3 crystal of rank $n$ on $S/W$ with Tate module $T_H$. A $T$\textit{-structure} on $H_{S'}$ is an isometry $T\otimes\cO_{S'}\to T_H$. The image of $H$ in $T^*\otimes\cO_{S'}/T\otimes\cO_{S'}\cong T_0\otimes\cO_S$ descends uniquely through the absolute Frobenius of $S$, giving rise to a characteristic subspace $K_H\subset T_0\otimes\cO_S$.

\begin{Definition}\label{def:twistedperiodmorphism}
  Let $\tpi\colon \tX\to\sS_{\Lambda}^o$ be the pullback of $\pi\colon X\to S_{\Lambda}$ to $\sS_{\Lambda}^o$, and let $\alpha\in (\R^2\tpi^{\fl}_*\m_p)(\sS_{\Lambda}^o)$ be the restriction of the universal cohomology class. Given an affine \'{e}tale open $S$ of $\sS^o_{\Lambda}$ and data as in Situation \ref{situation:liftedsituation}, we consider the twisted K3 crystal $\tH(\sX/S')$ corresponding to $(\tX,\alpha)|_S$.
  Let $T=\tLam\otimes\bZ_p$. By Lemma \ref{lem:inducedmarking} the induced map $\tuLam\otimes\cO_{S'}\to\tH(\sX/S')$ factors through $T\otimes\cO_{S'}$, and thus gives rise to a $T$-structure. As described in \cite{Ogus78}, the image of $\tH(\sX/S')$ in 
  \[
    \dfrac{T^*\otimes\cO_{S'}}{T\otimes\cO_{S'}}=\dfrac{\tLam^*}{\tLam}\otimes\cO_S\cong\tuLam_0\otimes\cO_S
  \]
  descends through the Frobenius, and we obtain a characteristic subspace $K\subset\tuLam_0\otimes\cO_S$. By Proposition \ref{prop:compatiblewithbasechange2}, this is independent of our choices of lifts, and in particular glues to a global characteristic subspace
  \[
    \tK_{\Lambda}\subset\tuLam_{0}\otimes\cO_{\sS_{\Lambda}^o}
  \]
  Note that for every geometric point $s\in\sS_{\Lambda}^o$ we have $e\notin(\tK_{\Lambda})_s$. We therefore obtain a morphism
  \begin{equation}\label{eq:twistedperiodmorphism}
    \trho\colon \sS_{\Lambda}^o\to \oM^{\langle e\rangle}_{\tLam_0}
  \end{equation}
  which we call the \textit{twisted period morphism}.
\end{Definition}
  
We refer to Remark \ref{rem:subtleremark} for some of the subtleties of this definition. We remark that the non-twisted period morphism (\ref{eq:periodmorphism}) has a similar interpretation in terms of crystalline cohomology (see Remark \ref{rem:periodmorphism1}). The twisted period morphism and the usual one are related by a commutative diagram
\begin{equation}\label{eq:commutingsquare}
  \begin{tikzcd}
    \sS_{\Lambda}^o\arrow{r}{p}\arrow{d}[swap]{\trho}&S_\Lambda\arrow{d}{\rho}\\
    \oM_{\tLam_0}^{\langle e\rangle}\arrow{r}{\pi_{e}}&\oM_{\Lambda_0}
  \end{tikzcd}
\end{equation}

  Let $\sX\to X$ be a $\m_p$-gerbe over a supersingular K3 surface and let $m\colon\Lambda\to \Pic(X)$ be a marking. By Lemma \ref{lem:inducedmarking}, we obtain a map $\tLam\to\tH(\sX/W)$, and hence a $T=\tLam\otimes\bZ_p$-structure $\tLam\otimes\bZ_p\to\tH(\sX/W)$ on the supersingular K3 crystal $\tH(\sX/W)$. We will let $K(\sX)\subset\tLam_0\otimes k$ be the corresponding characteristic subspace. In other words, this gives the image of the point corresponding to $\sX$ under $\trho$.

\begin{Proposition}\label{prop:etale}
  The twisted period morphism $\trho$ is \'{e}tale.
\end{Proposition}
\begin{proof}
   Consider the diagram
  \begin{equation}\label{eq:diagramwithasquare}
    \begin{tikzcd}
      \sS_{\Lambda}^o\arrow{dr}[swap]{\trho}\arrow{r}{(\trho,p)}\arrow[bend left=30]{rr}{p}&\oM_{\tLam_0}^{\langle e\rangle}\times_{\oM_{\Lambda_0}}S_\Lambda\arrow{r}\arrow{d}&S_\Lambda\arrow{d}{\rho}\\
      &\oM_{\tLam_0}^{\langle e\rangle}\arrow{r}{\pi_{e}}&\oM_{\Lambda_0}
    \end{tikzcd}
  \end{equation}
  with Cartesian square induced by the commuting square (\ref{eq:commutingsquare}). Ogus has shown in \cite{Ogus78} that $\rho$ is \'{e}tale. Thus, it will suffice to show that $(\trho,p)$ is \'{e}tale. The base change of $(\trho,p)$ by the absolute Frobenius $F_{S_{\Lambda}}\colon S_{\Lambda}\to S_{\Lambda}$ gives a map
  \[
    (\trho,p)^{(p/S_{\Lambda})}\colon (\sS^o_{\Lambda})^{(p/S_{\Lambda})}\to\left(\overline{M}_{\tLam_0}^{\langle e\rangle}\times_{\overline{M}_{\Lambda_0}}S_{\Lambda}\right)^{(p/S_{\Lambda})}
  \]
  which is a \'{e}tale if and only if $(\trho,p)$ is \'{e}tale. Let $\sU_{S_{\Lambda}}$ be the pullback of $\sU_{\Lambda_0}$ to $S_{\Lambda}$. By Proposition \ref{prop:isoofgroupschemes}, we have an isomorphism
  \[
    \R^1\pi^{(p)\etale}_*\nu(1)^o=F_{S_{\Lambda}}^{-1}\sU_{S_{\Lambda}}\iso\left(\overline{M}_{\tLam_0}^{\langle e\rangle}\times_{\overline{M}_{\Lambda_0}}S_{\Lambda}\right)^{(p)}
  \]
  We will identify the two for the remainder of this proof. By the definition of the twisted period morphism, we have a commuting diagram
  \begin{equation}\label{eq:bythedefinition}
    \begin{tikzcd}
      \sS^o_{\Lambda}\arrow{r}{F_{\sS^o_{\Lambda}/S_{\Lambda}}}\arrow[bend right=25]{rr}{\Upsilon^o}&(\sS^o_{\Lambda})^{(p)}\arrow{r}{F_{S_{\Lambda}}^{-1}(\trho,p)}&(\overline{M}_{\tLam_0}^{\langle e\rangle}\times_{\overline{M}_{\Lambda_0}}S_{\Lambda})^{(p)}
    \end{tikzcd}
  \end{equation}
  Because its source and target are smooth, $(\trho,p)$ is \'{e}tale if and only if it is \'{e}tale on the fiber over every geometric point $s$ of $S_{\Lambda}$. Thus, we may replace $S_{\Lambda}$ with the spectrum $\Spec k'$ of an algebraically closed field, and the universal K3 surface $\pi\colon X\to S_{\Lambda}$ with its fiber $\pi_s\colon X_s\to\Spec k'$. As in Lemma \ref{lem:LemmaNumeroUno}, we consider the diagram
  \begin{equation}\label{eq:5sideddiagram}
    \begin{tikzcd}
      (\R^2\pi_{s*}^{\fl}\m_p)^o\arrow{r}{F_{\sS^o_s/k'}}\arrow{d}&(\R^2\pi_{s*}^{(p)\fl}\m_p)^o\arrow{d}{\iota}\arrow{r}{F_{k'}^{-1}(\trho,p)}&\R^1\pi_{s*}^{(p)\etale}\nu(1)^o\arrow{dl}{\delta}\\
      \R^2\pi_{s*}^{\etale}\cO_{X_s}^\times\arrow{r}&\R^2\pi_{s*}^{(p)\etale}\cO_{X^{(p)}_s}^\times&
    \end{tikzcd}
  \end{equation}
  We do not know yet that this diagram commutes, merely that the outer compositions agree (by Lemma \ref{lem:LemmaNumeroUno}) and that the square commutes. This implies that
  \[
    \delta\circ F^{-1}_{k'}(\trho,p)\circ F_{\sS^o_s/k'}=\iota\circ F_{\sS^o_s/k'}
  \]
  By Theorem \ref{thm:grothendieck} and the Leray spectral sequence, there is an isomorphism 
  \[
    \R^2\pi_{s*}^{\etale}\cO_{X_s}^\times\cong\varepsilon_{s*}\R^2\pi^{\fl}\cO_{X_s}^\times
  \]
  Thus, the functors $\R^2\pi_{s*}^{\etale}\cO_{X_s}^\times$ and $\R^2\pi_{s*}^{(p)\etale}\cO_{X^{(p)}_s}^\times$ on the category of schemes over $\Spec k'$ are sheaves in the flat topology. Being representable, the functors in the top row of diagram (\ref{eq:5sideddiagram}) are also sheaves in the flat topology. The relative Frobenius $F_{\sS^o_s/k'}$ is faithfully flat, and hence an epimorphism in the category of sheaves on the big flat site. Thus, we have
  \[
    \delta\circ F^{-1}_{k'}(\trho,p)=\iota
  \]
  and so in fact the diagram commutes. After passing to completions at the identity section, both $\iota$ and $\delta$ become isomorphisms (see Remark \ref{rem:RemarkNumeroUno}). Therefore, $F^{-1}_{k'}(\trho,p)$ induces an isomorphism on completions at the identity section, and so is \'{e}tale.
\end{proof}

\begin{Proposition}\label{prop:commutativetriangle}
  The map
  \[
    (\trho,p)\colon \sS_{\Lambda}^o\to\oM^{\langle e\rangle}_{\tLam_0}\times_{\oM_{\Lambda_0}}S_\Lambda
  \]
  induced by the commuting square (\ref{eq:commutingsquare}) is an isomorphism. In other words, the square is Cartesian.
\end{Proposition}
\begin{proof}
  For readability, let us we write $M^+=\oM^{\langle e\rangle}_{\tLam_0}$ and $M=\oM_{\Lambda_0}$. With the same identifications as in the proof of Proposition \ref{prop:etale}, we then have a commutative diagram
  \begin{equation}\label{eq:Cartesian}
    \begin{tikzcd}
      \sS_{\Lambda}^o\arrow{r}[swap]{(\trho,p)}\arrow{dr}[swap]{\trho}\arrow[bend left=15]{rr}{\Upsilon^o}\arrow[bend left=28]{rrr}{p}&M^+\hspace{-.1cm}\times_{M}\hspace{-.1cm}S_\Lambda\arrow{d}\arrow{r}&[-14pt]M^+\hspace{-.1cm}\times_{M}\hspace{-.1cm}S_\Lambda\arrow{r}\arrow{d}&S_{\Lambda}\arrow{d}{\rho}\\
      &M^+\arrow{r}{F_{M^+/M}}\arrow{dr}[swap]{F_{M^+}}&M^{+(p/M)}\arrow{r}{\pi_e^{(p/M)}}\arrow{d}{W_{M^+/M}}&M\arrow{d}{F_M}\\
      &&M^+\arrow{r}{\pi_e}&M
    \end{tikzcd}
  \end{equation}
  where the three squares are Cartesian. Because $F_{M^+/M}$ and $\Upsilon^o$ are universal homeomorphisms, it follows that $(\trho,p)$ is a universal homeomorphism as well. By Proposition \ref{prop:etale}, $(\trho,p)$ is \'{e}tale. By Zariski's Main Theorem (for algebraic spaces, see \cite[082K]{stacks-project}), an \'{e}tale universal homeomorphism is an isomorphism.
\end{proof}

Proposition \ref{prop:commutativetriangle} could be viewed as a relative Torelli theorem for the twisted period morphism over the non-twisted period morphism. Let us describe a few consequences for an individual supersingular K3 surface $\pi\colon X\to\Spec k$. We define groups schemes $\sU$ and $\sH$ over $k$ by the exact sequences
\[
  0\to\sU\to\dfrac{\Pic(X)_0\otimes k}{K(X)}\xrightarrow{1-F^*}\dfrac{\Pic(X)_0\otimes k}{K(X)+F^*K(X)}
\]\[
  0\to\sH\to\dfrac{\Pic(X)\otimes k}{K(X)}\xrightarrow{1-F^*}\dfrac{\Pic(X)\otimes k}{K(X)+F^*K(X)}
\]
where $K(X)$ is the characteristic subspace corresponding to $X$ (equipped with the trivial marking $\id:\Pic(X)\to\Pic(X)$).
\begin{Proposition}\label{prop:aninterestingconsequence}
  There are natural isomorphisms $\R^2\pi^{\fl}_*\m_p\iso\sH$ and $\R^1\pi^{(p)\etale}_*\nu(1)\iso \sH^{(p/\Spec k)}$ of group schemes fitting into a commutative diagram
  \[
    \begin{tikzcd}
      \R^2\pi^{\fl}_*\m_p\isor{d}{}\arrow{r}{\Upsilon}&\R^1\pi^{(p)\etale}_*\nu(1)\isor{d}{}\\
      \sH\arrow{r}{F_{\sH/k}}&\sH^{(p/\Spec k)}
    \end{tikzcd}
  \]
\end{Proposition}
\begin{proof}
  Let $\Lambda=\Pic(X)$, and equip $X$ with the trivial marking $\id:\Lambda\to\Pic(X)$. We obtain the right hand vertical isomorphism from diagram (\ref{eq:flatdualitydiagram}) (as the marking is an isomorphism, the sheaf $\sW$ is trivial). By Proposition \ref{prop:commutativetriangle} (and the isomorphism of Proposition \ref{prop:isoofgroupschemes}) the twisted period morphism induces an isomorphism
  \[
    (\R^2\pi^{\fl}_*\m_p)^o\iso\sU
  \]
  which by the definition of the twisted period morphism fits into a commuting diagram
  \[
    \begin{tikzcd}
      (\R^2\pi^{\fl}_*\m_p)^o\arrow[hook]{r}\isor{d}{}&\R^2\pi^{\fl}_*\m_p\arrow{r}{\Upsilon}\arrow[dashed]{d}&\R^1\pi^{(p)\etale}_*\nu(1)\isor{d}{}\\
      \sU\arrow[hook]{r}&\sH\arrow{r}{F_{\sH/k}}&\sH^{(p/\Spec k)}
    \end{tikzcd}
  \]
  As $\Upsilon$ and $F_{\sH/k}$ are universal homeomorphisms, they induce isomorphisms on the respective groups of connected components. We may therefore produce by translation an isomorphism filling in the dashed arrow.
\end{proof}
\begin{Remark}
  The exact sequence (of group schemes!)
  \[
    0\to\R^2\pi^{\fl}_*\m_p\to\dfrac{\Pic(X)\otimes k}{K(X)}\xrightarrow{1-F^*}\dfrac{\Pic(X)\otimes k}{K(X)+F^*K(X)}
  \]
  produced by Proposition \ref{prop:aninterestingconsequence} is a strengthened form of the equality mentioned in Remark 3.26 of \cite{Ogus78} (we warn that the referenced formula appears to contain some typos).
\end{Remark}

\begin{Remark}\label{rem:strongflatduality}
  The natural inclusion $\sU\subset\sH$ gives rise to a short exact sequence
  \[
    0\to\sU\to\sH\to\sD\to 0
  \]
  where $\sD$ is the quotient (this is the same as the sheaf $\sD$ in the diagram (\ref{eq:flatdualitydiagram})). Write $\Lambda=\Pic(X)$. We have a commutative diagram
  \[
    \begin{tikzcd}
      0\arrow{r}&\dfrac{p\Lambda^*}{p\Lambda}\arrow{r}\arrow{d}&\dfrac{\Lambda}{p\Lambda}\arrow{r}\arrow{d}&\dfrac{\Lambda}{p\Lambda^*}\arrow{r}\isor{d}{}&0\\
      0\arrow{r}&\sU\arrow{r}&\sH\arrow{r}&\sD\arrow{r}&0
    \end{tikzcd}
  \]
  with exact rows. Evaluating on $k$-points and applying Proposition \ref{prop:aninterestingconsequence}, we recover part of diagram (\ref{eq:artinsflatduality}) (see also Remark \ref{rem:flatduality}).
\end{Remark}

A basic fact about the Gauss-Manin connection is that it may be computed in terms of the Kodaira-Spencer map (see \cite{MR0337959}). A similar result is true in our twisted setting. As explained in \cite{Ogus78}, the connections on the crystals $\tH(\sX/S')$ defined on an affine cover of $\sS_{\Lambda}^o$ induce a map
\begin{equation}\label{eq:kodairaspencer}
    T^1_{\sS^o_{\Lambda}/k}\to\sHom_{\sS^o_{\Lambda}}\left(\tK_{\Lambda}\cap F^*\tK_{\Lambda},\dfrac{F^*\tK_{\Lambda}}{\tK_{\Lambda}\cap F^*\tK_{\Lambda}}\right)
\end{equation}
which, under the identification in Lemma \ref{lem:Ogustangentspace} of the tangent space to the period domain, is exactly the differential
\[
  d\trho\colon T^1_{\sS^o_{\Lambda}/k}\to\trho^*T^1_{\oM_{\tLam_0}^{\langle e\rangle}/k}
\]
of $\trho$. The diagram (\ref{eq:diagramwithasquare}) gives rise to a diagram
    \begin{equation}\label{eq:tangentspacesdiagram}
      \begin{tikzcd}
        0\arrow{r}&T^1_{\sS_{\Lambda}^o/S_{\Lambda}}\arrow{r}\arrow{d}{d(\trho,p)}&T^1_{\sS_{\Lambda}^o/k}\arrow{r}\arrow{d}{d\trho}&p^*T^1_{S_{\Lambda}/k}\arrow{r}\arrow{d}{p^*(d\rho)}&0\\
        0\arrow{r}&\trho^*T^1_{\oM^{\langle e\rangle}_{\tLam_0}/\oM_{\Lambda_0}}\arrow{r}&\trho^*T^1_{\oM^{\langle e\rangle}_{\tLam_0}/k}\arrow{r}&\trho^*\pi_e^*T^1_{\oM_{\Lambda_0}/k}\arrow{r}&0
      \end{tikzcd}
    \end{equation}
  of tangent sheaves with exact rows. Pulling back the identification of Lemma \ref{lem:relativetangentspace}, we find identifications
\[
  \begin{tikzcd}
    T^1_{\sS_{\Lambda}^o/S_{\Lambda}}\arrow{r}\arrow[hook]{d}&\sHom_{\sS^o_{\Lambda}}\left(\dfrac{\tK_{\Lambda}\cap F^*\tK_{\Lambda}}{\tK_{\Lambda}\cap F^*\tK_{\Lambda}\cap e^\perp},\dfrac{F^*\tK_{\Lambda}}{\tK_{\Lambda}\cap F^*\tK_{\Lambda}}\right)\arrow[hook]{d}\\
    T^1_{\sS^o_{\Lambda}/k}\arrow{r}&\sHom_{\sS^o_{\Lambda}}\left(\tK_{\Lambda}\cap F^*\tK_{\Lambda},\dfrac{F^*\tK_{\Lambda}}{\tK_{\Lambda}\cap F^*\tK_{\Lambda}}\right)
  \end{tikzcd}
\]
As described in the proof of Lemma \ref{lem:relativetangentspace}, pairing with $e$ gives rise to isomorphisms
\[
  \sHom_{\sS^o_{\Lambda}}\left(\dfrac{\tK_{\Lambda}\cap F^*\tK_{\Lambda}}{\tK_{\Lambda}\cap F^*\tK_{\Lambda}\cap e^\perp},\dfrac{F^*\tK_{\Lambda}}{\tK_{\Lambda}\cap F^*\tK_{\Lambda}}\right)\cong\dfrac{F^*\tK_{\Lambda}}{\tK_{\Lambda}\cap F^*\tK_{\Lambda}}\cong p^*\dfrac{F^*K_{\Lambda}}{K_{\Lambda}\cap F^*K_{\Lambda}}\cong p^*\R^2\pi_*\cO_{X}
\]
Thus, $d(\trho,p)$ induces a map
\begin{equation}\label{eq:isomorphism1}
    T^1_{\sS_{\Lambda}^o/S_{\Lambda}}\to p^*\R^2\pi_*\cO_X
\end{equation}
which we have proved to be an isomorphism.

\begin{Proposition}
  The isomorphism (\ref{eq:isomorphism1}) agrees with the ``Kodaira-Spencer'' isomorphism constructed in Lemma \ref{lem:tangentspacetomodulispace}.
\end{Proposition}
\begin{proof}
  Omitted.
\end{proof}


\section{Supersingular Twistor Space}\label{sec:moduli}

\subsection{Moduli of twisted sheaves}\label{sec:derived}

Let $X$ be a K3 surface over our fixed algebraically closed field $k$ of characteristic $p\geq 3$, and let $\sX\to X$ be a $\m_p$-gerbe. In this section we will study the derived category of twisted sheaves on $\sX$ with the aim of extending various results that are well known over the complex numbers to our setting. In particular, we define Chern characters for twisted sheaves on $\sX$, consider the action of derived equivalences on the twisted Mukai crystals, and establish the existence of moduli spaces of twisted sheaves.

In establishing these results, we will need to reduce certain statements for twisted K3 surfaces to the non-twisted case. A considerable simplification occurs here if we assume that $X$ is supersingular. This is because any $\m_p$-gerbe over a supersingular K3 surface sits in a canonical flat family with a $\m_p$-gerbe over the same surface that has trivial Brauer class. Indeed, consider the connected component $\bA^1\subset\R^2\pi^{\fl}_*\m_p$ of the group of $\m_p$-gerbes on $X$ such that $\sX$ is a fiber of the tautological family $\tsX\to\bA^1$ (see Section \ref{sec:BrauerGroup}, where these are discussed as the basic examples of twistor families). We say that $\sX$ \textit{deforms the trivial gerbe} if this component is the identity component. In any case, the map $\bA^1(k)\to\Br(X)$ induced by the Kummer sequence is surjective (for instance, by diagram (\ref{eq:flatdualitydiagram})) so some fiber of the family $\tsX\to\bA^1$ has trivial Brauer class.

The corresponding reductions in the finite height case are significantly more involved, and require lifting to the complex numbers and comparison with the Hodge-theoretic constructions of \cite{HS04}. Thus, although the main results of this section are true for with no assumptions on the height of $X$, we will for the most part restrict our attention to the supersingular case.


Let $p\colon \sX\to X$ be a $\m_p$-gerbe on a smooth projective variety over $k$. Suppose that $\sX$ has the resolution property. We have in mind the case when $\sX$ is a $\m_p$-gerbe over a K3 surface, or the exterior sum of two $\m_p$-gerbes on the product of two K3 surfaces.

\begin{Definition}\label{def:twistedcherncharacter}
  If $\sE$ is a locally free sheaf of positive rank on $\sX$, the \textit{twisted Chern character} of $\sE$ is
  \[
    \ch_{\sX}(\sE)=\sqrt[\leftroot{-2}\uproot{2}p]{\ch(p_*(\sE^{\otimes p}))}
  \]
  where by convention we choose the $p$-th root so that $\rk(\ch_{\sX}(\sE))=\rk(\sE)$. We define the twisted Mukai vector of $\sE$ by
  \[
    v_{\sX}(\sE)=\ch_{\sX}(\sE).\sqrt{\Td(X)}
  \]
\end{Definition}
If $\sE$ is 0-twisted, then its twisted Chern character is the same as the usual Chern character of its pushforward to $X$:
\[
  \ch_{\sX}(\sE)=\ch(p_*\sE)
\]
In most cases we will consider, $\sE$ will be a twisted sheaf. The twisted Chern character determines an additive map 
\[
  \ch_{\sX}\colon K^{(1)}(\sX)\to A^*(X)_{\bQ}
\]
where $K^{(1)}(\sX)$ is the Grothendieck group of locally free twisted sheaves on $\sX$ and $A^*(X)_{\bQ}$ is the numerical Chow theory of $X$ tensored with $\bQ$. We record a few straightforward lemmas extending properties of the usual Chern characters to our twisted Chern characters.

\begin{Lemma}\label{lem:riemannrochEulerVersion}
  If $\sE$ and $\sF$ are locally free twisted sheaves on $\sX$, then
  \[
    \chi(\sE,\sF)=\deg(\ch_{\sX}(\sE^\vee\otimes\sF).\Td(X))
  \]
\end{Lemma}
\begin{proof}
  The sheaves $\sHom(\sE,\sF)$ and $\sE^\vee\otimes\sF$ are $0$-twisted. Thus,
  \[
    \ch_{\sX}(\sE^\vee\otimes\sF)=\ch(p_*(\sE^\vee\otimes\sF))
  \]
  and the result follows from the Grothendieck-Riemann-Roch theorem as in the nontwisted case.
\end{proof}

If $f\colon\sX\to\sY$ is a map of $\m_p$-gerbes, we obtain a homomorphism
\[
  f^*\colon K^{(1)}(\sY)\to K^{(1)}(\sX)
\]
determined by $\left[\sE\right]\mapsto\left[\mathbf{L}f^*\sE\right]$.

\begin{Lemma}\label{lem:thisisalemma}
  Let $\pi\colon T\to X$ be a morphism of smooth projective $k$-schemes, and consider the Cartesian diagram
  \[
    \begin{tikzcd}
      \sX_T\arrow{r}{\pi_{\sX}}\arrow{d}&\sX\arrow{d}\\
      T\arrow{r}{\pi}&X
    \end{tikzcd}
  \]
  If $\alpha\in K^{(1)}(\sX)$, then $\pi^*\ch_{\sX}(\alpha)=\ch_{\sX_T}(\pi_{\sX}^*\alpha)$.
\end{Lemma}
\begin{proof}
  If $\sE$ is a locally free twisted sheaf on $\sX$, then
  \[
    \pi^*\ch_{\sX}(\sE)=\pi^*\sqrt[\leftroot{-2}\uproot{2}p]{\ch(p_*(\sE^{\otimes p}))}=\sqrt[\leftroot{-2}\uproot{2}p]{\ch(\mathbf{L}\pi^*p_*(\sE^{\otimes p}))}=\sqrt[\leftroot{-2}\uproot{2}p]{\ch(p_{T*}\mathbf{L}\pi_{\sX}^*(\sE^{\otimes p}))}=\ch_{\sX_T}(\mathbf{L}\pi_{\sX}^*\sE)
  \]
  This implies the result.
\end{proof}

Let $\sX\to X$ and $\sY\to Y$ be $\m_p$-gerbes over smooth projective varieties satisfying the resolution property. In the following lemma we formulate a version of the Grothendieck-Riemann-Roch formula for the projection $\sX\times\sY\to\sX$.
\begin{Lemma}\label{lem:riemannroch}
  For any $\alpha\in K^{(1,0)}(\sX\times\sY)$ we have
  \[
    \ch_{\sX}(\pi_{\sX*}\alpha)=\pi_{X*}(\ch_{\sX\times\sY}(\alpha).\Td(\pi_{X}))
  \]
  where $\pi_{\sX}\colon \sX\times\sY\to\sX$ is the projection and $\pi_X\colon X\times Y\to X$ is the induced map on coarse spaces.
\end{Lemma}
\begin{proof}
  We observe that by the usual Grothendieck-Riemann-Roch formula, the result holds when $\sX$ is trivial. In the general case, we choose a finite flat cover $f_X\colon U\to X$ by a smooth projective $k$-scheme $U$ such that the gerbe $\sX_U=\sX\times_XU$ is trivial (for instance, we may take $f_X$ to be the absolute Frobenius). Consider the Cartesian diagrams
  \[
    \begin{tikzcd}
      \sX_U\times\sY\arrow{r}{f_{\sX\times\sY}}\arrow{d}[swap]{\pi_{\sX_U}}&\sX\times\sY\arrow{d}{\pi_{\sX}}\\
      \sX_U\arrow{r}{f_{\sX}}&\sX
    \end{tikzcd}\hspace{2cm}
    \begin{tikzcd}
      X_U\times Y\arrow{r}{f_{X\times Y}}\arrow{d}[swap]{\pi_U}&X\times Y\arrow{d}{\pi_X}\\
      U\arrow{r}{f_X}&X
    \end{tikzcd}
  \]
  Using Lemma \ref{lem:thisisalemma}, we compute
  \begin{align*}
    f_X^*(\pi_{X*}(\ch_{\sX\times\sY}(\alpha).\Td(\pi_{X})))&=\pi_{U*}f_{X\times Y}^*(\ch_{\sX\times\sY}(\alpha).\Td(\pi_{X}))\\
                                           &=\pi_{U*}(\ch_{\sX_U\times\sY}(f_{\sX\times\sY}^*\alpha).\Td(\pi_U))\\
                                           &=\ch_{\sX_U}(\pi_{\sX_U*}(f_{\sX\times\sY}^*\alpha))\\
                                           &=\ch_{\sX_U}(f_{\sX}^*\pi_{\sX*}\alpha)\\
                                           &=f_X^*\ch_{\sX}(\pi_{\sX*}\alpha)
  \end{align*}
  The map $f_{X*}f_X^*\colon A^*(X)_{\bQ}\to A^*(X)_{\bQ}$ is given by multiplication by the degree of $f_X$, so this computation implies the result.
\end{proof}
\begin{Remark}
  The proof gives the same formula in various other situations. We have not attempted to give a general formulation. However, we note that one must be somewhat careful in applying the formula of Lemma \ref{lem:riemannroch}, as it is easily seen to fail, for instance, for a $(1,1)$ twisted sheaf on $\sX\times\sY$, or for a twisted sheaf on $\sX$ and the coarse space morphism $\sX\to X$.
\end{Remark}

We next discuss the relationship between the derived category and cohomological equivalences. Generally speaking, we find the same behaviors as in the untwisted case (see Chapter 5 of \cite{Huy06}) and in the twisted case over the complex numbers (see \cite{HS04}). We denote by $D^{(n)}(\sX)$ the bounded derived category associated to $\Coh^{(n)}(\sX)$. We consider a perfect complex $\sP^\bullet\in D^{(1,1)}(\sX\times\sY)$ of twisted sheaves. Using $\sP^\bullet$ as a kernel, we get a Fourier-Mukai transform 
\[
  \Phi_{\sP^\bullet}\colon D^{(\text{-}1)}(\sX)\to D^{(1)}(\sY)
\]
and an induced transform
\[
  \Phi_{[\sP^\bullet]}\colon K^{(\text{-}1)}(\sX)\to K^{(1)}(\sY)
\]
on K theory. The twisted Mukai vector $v_{\sX\times\sY}(\sP^\bullet)\in A^*(X\times Y)\otimes\bQ$ gives a cohomological transform 
\[
  \Phi^{\cris}_{v_{\sX\times\sY}(\sP^\bullet)}\colon \H^*(X/K)\to \H^*(Y/K)
\]
\begin{Lemma}\label{lem:actiononcohomologycompatible}
  If $\sP^{\bullet}\in D^{(1,1)}(\sX\times\sY)$ is a perfect complex of twisted sheaves, then the diagram
  \[
    \begin{tikzcd}[column sep=large]
      K^{(\text{-}1)}(\sX)\arrow{r}{\Phi_{[\sP^\bullet]}}\arrow{d}[swap]{v_{\sX}}&K^{(1)}(\sY)\arrow{d}{v_{\sY}}\\
      \H^*(X/K)\arrow{r}{\Phi^{\cris}_{v_{\sX\times\sY}(\sP^{\bullet})}}&\H^*(Y/K)
    \end{tikzcd}
  \]
  commutes.
\end{Lemma}
\begin{proof}
  In the untwisted case, this follows by applying the Grothedieck-Riemann-Roch formula to the projections $X\times Y\to X$ and $X\times Y\to Y$ (see Corollary 5.29 of \cite{Huy06}). Using the twisted Grothendieck-Riemann-Roch formula of Lemma \ref{lem:riemannroch}, the same proof immediately gives the result in the twisted case as well.
\end{proof}

\begin{Proposition}\label{prop:cohomologicaltransform}
  Suppose that $\sX\to X$ and $\sY\to Y$ are $\m_p$-gerbes over K3 surfaces, and $\sP^\bullet\in D^{(1,1)}(\sX\times\sY)$ is a perfect complex of twisted sheaves inducing an equivalence of categories $\Phi_{\sP^\bullet}:D^{(\text{-}1)}(\sX)\to D^{(1)}(\sY)$. The cohomological transform
  \[
    \Phi^{\cris}_{v_{\sX\times\sY}(\sP^\bullet)}\colon\tH(X/K)\to\tH(Y/K)
  \]
  is an isomorphism of $K$-vector spaces, an isometry with respect to the Mukai pairing, and commutes with the respective Frobenius operators.
\end{Proposition}
\begin{proof}
  Using Lemma \ref{lem:actiononcohomologycompatible} and \ref{lem:riemannrochEulerVersion}, one shows that the cohomological transform is an isomorphism of vector spaces and an isometry exactly as in the untwisted case (see Proposition 5.33 and Proposition 5.44 of \cite{Huy06}). To see that it is compatible with the crystal structure, recall that the $i$-th component $\ch(\sE)_i$ of the Chern character of a sheaf $\sE$ on a smooth projective variety satisfies
  \[
    \Phi(\ch(\sE)_i)=p^i\ch(\sE)_i
  \]
  The result follows upon expanding $v_{\sX\times\sY}(\sP^\bullet)$ under the Kunneth decomposition.
\end{proof}

Let us record a first example of an equivalence of derived categories and its action on cohomology. Suppose that $\sX\to X$ is a $\m_p$-gerbe over a K3 surface, and let $\sL\in\Pic(X)$ be a line bundle. In Definition \ref{def:relativeroot}, we defined a $\m_p$-gerbe $\sX\{\sL^{1/p}\}\to X$ over $X$, which is equipped with a universal line bundle $\sM\in\Pic^{(\text{-}1,1)}(\sX\times_X\sX\{\sL^{1/p}\})$. Consider the natural map 
\[
  \iota\colon\sX\times_X\sX\{\sL^{1/p}\}\to\sX\times\sX\{\sL^{1/p}\}.
\]
\begin{Proposition}\label{prop:twistedlinebundleequivalence}
  The Fourier-Mukai transform induced by the sheaf $\iota_*\sM\in\Coh^{(\text{-}1,1)}(\sX\times\sX\{\sL^{1/p}\})$ is an equivalence of categories
  \[
    \Phi_{\iota_*\sM}\colon D^{(1)}(\sX)\to D^{(1)}(\sX\{\sL^{1/p}\})
  \]
  The induced map on cohomology is given by
  \[
    e^{l/p}\colon\tH(\sX/W)\to\tH(\sX\{\sL^{1/p}\}/W)
  \]
  where $l\in\Pic(X)$ is the first Chern class of $\sL$.
\end{Proposition}
\begin{proof}
  To see that $\Phi_{\iota_*\sM}$ is an equivalence, note that it is inverse to $\Phi_{\iota_*\sM^\vee}$. To determine the action on cohomology, we use that the induced map on cohomology 
  \[
    \Phi^{\cris}_{\iota_*\sM}\colon\tH(X/K)\to\tH(X/K)
  \]
  satisfies
  \[
    \left(\Phi^{\cris}_{\iota_*\sM}\right)^p=\Phi^{\cris}_{\iota_*\sL}=e^l
  \]
  The result follows upon noting that $\Phi^{\cris}_{\iota_*\sM}$ sends sheaves of positive rank to sheaves of positive rank.
\end{proof}
In fact, in this example $\Phi_{\iota_*\sM}$ induces also an equivalence of abelian categories. For the remainder of this section, we will further specialize to the case when $X$ and $Y$ are supersingular. As noted at the beginning of this section, many of these results hold without this assumption, but the proofs are more involved.

\begin{Proposition}\label{prop:twistedcherncharacter}
  If $p\colon\sX\to X$ is a $\m_p$-gerbe on a supersingular K3 surface over $k$, then for any twisted sheaf $\sE\in\Coh^{(1)}(\sX)$, the twisted Chern character $\ch_{\sX}(\sE)$ lies in the twisted extended N\'{e}ron-Severi lattice $\tN(\sX)$.
\end{Proposition}
\begin{proof}
  Let $\alpha\in \H^2(X,\m_p)$ be the class of $\sX$. We assume first that $\sX$ is essentially trivial, so that there exists a line bundle $\sL$ such that $B=\frac{c_1(\sL)}{p}$ is a B-field lift of $\alpha$. The boundary map $\H^1(X,\bG_m)\to \H^2(X,\m_p)$ takes a line bundle to its gerbe of $p$-th roots, so on $\sX$ there is a universal line bundle $\sM$ equipped with an isomorphism $\sM^{\otimes p}\iso p^*\sL$. Using $\sM$, we compute
  \[
    \sL^\vee\otimes p_*(\sE^{\otimes p})\cong p_*(\sM^\vee\otimes \sE)^{\otimes p}
  \]
  which implies that
  \[
    \ch_{\sX}(\sE)=e^{\frac{c_1(\sL)}{p}}\ch(p_*(\sM^\vee\otimes \sE))
  \]
  which gives the result. 
  
  Next, suppose that $\alpha$ is not essentially trivial, so that $\left[\alpha\right]\in\Br(X)$ has order $p$. We make some reductions. It will suffice to show the result for $\sE$ a torsion free twisted sheaf of positive rank. Let $\eta$ denote the generic point of $X$. There exists a $\sX_{\eta}$-twisted sheaf $\sF$ of rank $p$ and a surjection $\sE|_{\eta}\to \sF$.
  Thus we get a map $\sE\to\eta_*\sF$. We find a surjection $\sE\to \sF'$ for some coherent $\sX$-twisted sub-sheaf $\sF'$ of rank $p$, and by modding out by torsion we may assume $\sF'$ is torsion free. As $\sE$ is torsion free, so is the kernel of this map. So, by induction we may reduce to the case when $\sE$ is a torsion free twisted sheaf of rank $p$. We claim that such an $\sE$ is simple, that is, that $\Hom(\sE,\sE)\cong k$. Indeed, as $\sE$ is torsion free, the natural map $\Hom(\sE,\sE)\to\Hom(\sE^{\vee\vee},\sE^{\vee\vee})$ is an injection, and because $\left[\alpha\right]\in\Br(X)$ has order $p$, any locally free sheaf of rank $p$ is simple. There is a deformation theory for torsion free twisted sheaves with unobstructed determinant with obstruction space
  \[\ker(\Ext_{X}^2(\sE,\sE)\too{\Tr}\H^2(X,\cO_{X}))\]
  Under Serre duality, $\Tr$ is dual to the natural map $\H^0(X,\cO_X)\to\Hom(\sE,\sE)$, which is an isomorphism because $\sE$ is simple. Consider the connected component $\bA^1$ of the group of $\m_p$-gerbes on $X$ that contains $\sX$ as the fiber over some $t\in\bA^1(k)$. The coarse space of the corresponding twistor family $\widetilde{\sX}\to\bA^1$ is the trivial family $X\times\bA^1\to\bA^1$, so in particular the determinant of $\sE$ is unobstructed. Hence, $\sE$ is unobstructed. By the Grothendieck existence theorem for twisted sheaves, we find a twisted sheaf $\sE'$ on $\widetilde{\sX}_{k\left[\left[x-t\right]\right]}$ that is flat over $\Spec(k\left[\left[x-t\right]\right])$ and whose restriction to the closed fiber is isomorphic to $\sE$. By Proposition 2.3.1.1 of \cite{Lie04}, the stack of coherent twisted sheaves on the morphism $\widetilde{\sX}\to\bA^1$ is in particular limit preserving. Thus we may apply Artin approximation (see \cite[07XB]{stacks-project}) to produce an \'{e}tale morphism
  \[
    f:(U,u)\to (\bA^1,t)
  \]
  and a coherent $\widetilde{\sX}\times_{\bA^1}U$-twisted sheaf $\sE''$ that is $U$-flat such that the restrictions of $\sE''$ and $f^*\sE$ to the closed fiber $(\widetilde{\sX}\times_{\bA^1}U)\times_Uu\cong\widetilde{X}_t$ are isomorphic. Taking the normalization of $\bA^1$ in $U$, we find a factorization $U\to C\to\bA^1$ such that $U\to C$ is an open immersion and $C\to\bA^1$ is finite and flat. Thus, every connected component of $C$ maps surjectively onto $\bA^1$. Consider a flat extension of $\sE''$ to $\sX_C$. The twisted Chern class is constant in a flat family, and the origin $0\in\bA^1$ corresponds to an essentially trivial $\m_p$-gerbe. We therefore obtain by the previous case that 
  \[
    \ch_{\sX}(\sE)\in\tN(\widetilde{\sX_0})\subset\tN(X)\otimes\bQ
  \]
  But $\sE$ has rank $p$, so by the explicit presentations in Proposition \ref{rem:nslattice} we see that $\ch_{\sX}(\sE)\in\tN(\sX)$. This completes the proof.
  
\end{proof}

\begin{Proposition}\label{prop:TwistedTate}
  If $\sX\to X$ is a $\m_p$-gerbe on a supersingular K3 surface, then the map
  \[
    \ch_{\sX}\colon K^{(1)}(\sX)\to\tN(\sX)
  \]
  is surjective.
\end{Proposition}
\begin{proof}
  If the Brauer class of $\sX$ is trivial, then $\sX$ is the gerbe of $p$-th roots of some line bundle $\sL$ on $X$. This means that there is an invertible $\sX$-twisted sheaf $\sM$ equipped with an isomorphism $\sM^{\otimes p}\iso p^*\sL$. We have a commutative diagram
  \[
    \begin{tikzcd}
      \Coh(X)\arrow{r}{p^*(\_)\otimes\sM}\arrow{d}[swap]{\ch}&\Coh^{(1)}(\sX)\arrow{d}{\ch_{\sX}}\\
      \tN(X)\arrow{r}{\cdot\frac{c_1(\sL)}{p}}&\tN(\sX)
    \end{tikzcd}
  \]
  where the top horizontal arrow is an equivalence of categories, and the lower horizontal arrow is an isometry, and in particular and isomorphism. Because the left vertical arrow is surjective, so is the right vertical arrow.
  
  Next, suppose that the Brauer class of $\sX$ is non-trivial. We refer to Proposition \ref{rem:nslattice} for an explicit description of $\tN(\sX)$. The structure sheaf of a closed point gives a $\sX$-twisted sheaf with twisted Chern class $(0,0,1)$. If $D\subset X$ is a closed integral subscheme of dimension 1, then by Tsen's Theorem the Brauer class of the gerbe $\sX_D=\sX\times_XD\to D$ is trivial, and hence there is an invertible $\sX_D$-twisted sheaf. The pushforward of such a sheaf under the map $\sX_D\hto\sX$ gives a $\sX$-twisted sheaf whose twisted Chern character is of the form $(0,D,s)$. Finally, by a theorem of Grothendieck, there exists a locally free $\sX$-twisted sheaf $\sE$ of rank $p$.
\end{proof}
\begin{Remark}
Combining Proposition \ref{prop:TwistedTate} with Proposition \ref{prop:twistedNeronseveri} shows that the natural map
  \[
    \im\left(\ch_{\sX}\colon K^{(1)}(\sX)\to\tN(\sX)\right)\otimes\bZ_p\to\tT(\sX)
  \]
is an isomorphism, where $\tT(\sX)$ is the Tate module of $\tH(\sX/W)$. Thus, the analog of the Tate conjecture holds for twisted supersingular K3 surfaces.
\end{Remark}

We next discuss stability conditions.
Fix a polarization $H$ of $X$. Recall the following definition (see Definition 2.2.7.6 of \cite{Lie04}).
\begin{Definition}
  If $\sE$ is a $\sX$-twisted sheaf, then the \textit{geometric Hilbert polynomial} of $\sE$ is the function
  \[
    P_{\sE}(m)=\deg(\ch_{\sX}(\sE(m)).\Td(X))
  \]
  where $\sE(m)=\sE\otimes p^*H^{\otimes m}$. Let $\alpha_{d}$ be the leading coefficient of $P_{\sE}$, and write
  \[
    p_{\sE}(m)=\frac{1}{\alpha_d}P_{\sE}(m)
  \]
\end{Definition}

As explained in Section 2.2.7 of \cite{Lie04}, $P_{\sE}$ is a numerical polynomial with the usual properties. In particular, we can use it to define stability and semistable of sheaves.

\begin{Definition}
  A $\sX$-twisted sheaf $\sE$ is \textit{stable} (resp. \textit{semistable}) if it is pure and for all proper non-trivial subsheaves $\sF\subset\sE$
  \[
    p_{\sF}(m)<p_{\sE}(m)
  \]
   (resp. $\leq$) for $m$ sufficiently large.
\end{Definition}

We have the following result.

\begin{Proposition}\label{prop:genericpolarization}
  If $v=(r,l,s)\in\tN(\sX)$ is a primitive Mukai vector, then there exists a locally finite collection of hyperplanes $W\subset N(X)\otimes\bR$ such that if $H$ does not lie on any of these hyperplanes, then any $H$-semistable $\sX$-twisted sheaf $\sE$ with $v_{\sX}(\sE)=v$ is $H$-stable.
\end{Proposition}
\begin{proof}
  We need some form of the Bogomolov inequality for $\m_H$-semistable twisted sheaves on $\m_p$-gerbes. For our purposes, little more than the existence of such a lower bound will suffice. By Lemma 3.2.3.13 of \cite{Lie04}, if $\sX\to X$ is a $\m_p$-gerbe on a smooth proper surface, then there exists a constant $C\geq 0$ such that for any polarization $H$ of $X$ and any $H$-semistable $\sX$-twisted sheaf $\sE$ of rank $r$,
  \begin{equation}\label{eq:langersinequality}
    \Delta(\sE)\geq -Cr^4
  \end{equation}
  The discriminant $\Delta(\sE)$ of a twisted sheaf is defined in by Definition 3.2.1.1 of \cite{Lie04}. For the remainder of the proof we will closely follow the proof of Theorem 4.C.3 of \cite{HL10}. We say that a class $\xi\in N(X)$ is of type $v$ if
  \[
    -r^2(\Delta+2Cr^4)\leq \xi^2<0
  \]
  The wall determined by $\xi$ is the set
  \[
    W_{\xi}=\left\{H\in\sC|\xi.H=0\right\}\subset\sC
  \]
  We say that $W_{\xi}$ is of type $v$ if $\xi$ is of type $v$. The proof of Theorem 4.C.2 of \cite{HL10} shows that the set of walls of type $v$ is locally finite in the ample cone $\sC$. Suppose that $\sE$ is a $\sX$-twisted $H$-semistable sheaf. If $\sE$ fails to be $H$-stable, then there exists a subsheaf $\sE'\subset\sE$ with Mukai vector $v'=(r',l',s')$ such that $r'<r$ and $p_{\sE'}(m)\equiv p_{\sE}(m)$. Suppose that $r>0$. We have
  \[
    p_{\sE}(m)=\frac{m^2}{2}+m\frac{l.H}{rH^2}+\frac{s}{rH^2}+\frac{1}{H^2}
  \]
  so this condition is equivalent to
  \[
    \xi_{\sE',\sE}.H=0\hspace{1cm}\mbox{ and }\hspace{1cm}r's=rs'
  \]
  where
  \[
    \xi=\xi_{\sE',\sE}=r'l-rl'
  \]
  We may assume that $\sE'$ is saturated, so that $\sE''=\sE/\sE'$ is torsion-free of rank $r''=r'-r$. By the Hodge index theorem, either $\xi=0$ or $\xi^2<0$. If $\xi=0$, then we compute
  \[
    \frac{r}{r'}v'=v
  \]
  But $r'<r$, so this contradicts our assumption that $v$ was primitive. Therefore, $\xi^2<0$. We have the identity
  \[
    \Delta(\sE)-\frac{r}{r'}\Delta(\sE')-\frac{r}{r''}\Delta(\sE'')=-\frac{\xi^2}{r'r''}
  \]
  Applying (\ref{eq:langersinequality}) to $\sE'$ and $\sE''$, we get
  \[
    -\Delta(\sE)r^2-2Cr^6\leq \xi^2<0
  \]
  We conclude that if such an $\sE'$ exists, then $H$ lies on a wall of type $v$. Therefore, if $H$ is not on any wall of type $v$, then $\sE$ is stable.
  
  The $r=0$ case is proved exactly as in Theorem 10.2.5 of \cite{Huy16}.
\end{proof}

\begin{Definition}
  A polarization $H$ is \textit{v-generic} if any $H$-semistable twisted sheaf with Mukai vector $v$ is $H$-stable.
\end{Definition}

For future use, we record the following observation.

\begin{Lemma}\label{lem:genericpolarizationinfamily}
    Suppose that $X$ is supersingular, and let $\tsX\to\A^1$ be a universal family of $\m_p$-gerbes over a connected component of $\R^2\pi^{\fl}_*\m_p$. Suppose that $v=(r,l,s)\in\tN(X)$ is a primitive isotropic vector such that, for each closed point $t\in\bA^1$, $v$ is contained in the image of the inclusion $\tN(\sX)\subset\tN(X)\otimes\bQ$. There exists a polarization $H$ of $X$ such that in each fiber $H$ is $v$-generic.
\end{Lemma}
\begin{proof}
    The extended N\'{e}ron-Severi groups of the fibers are given by Proposition \ref{rem:nslattice}. In particular, we see that, among the fibers of the universal family over each connected component of $\R^2\pi^{\fl}_*\m_p$, there are only finitely many possibilities for the twisted extended N\'{e}ron-Severi group (viewed as a subgroup of $\tN(X)\otimes\bQ$). For each, we find by Proposition \ref{prop:genericpolarization} a locally finite union of hyperplanes in $N(X)\otimes\bR$. But the union of all of these is again a locally finite union of hyperplanes, and we may therefore find an $H$ with the desired property.
\end{proof}

If $X\to S$ is a relative K3 surface, we define the \textit{relative extended N\'{e}ron-Severi lattice}
\[
  \tuN_{X/S}=\underline{\bZ}_S\oplus\Pic_{X/S}\oplus\underline{\bZ}_S
\]
which we equip with the Mukai pairing. We make the following definition.

\begin{Definition}\label{def:modulistack0}
  Let $X\to S$ be a relative K3 surface, $\sX\to X$ a $\m_p$-gerbe, $v$ a global section of $\tuN_{X/S}\otimes\bQ$, and $H$ a relative polarization of $X/S$. The \textit{moduli space of} $\sX$-\textit{twisted stable sheaves with twisted Mukai vector} $v$ is the stack $\sM_{\sX/S}(v)$ on $S$ whose objects over an $S$-scheme $T$ are $T$-flat $\sX_T$-twisted sheaves $\sE$ locally of finite presentation such that for each geometric point $t\in T$ the fiber $\sE_t$ is $H_t$-stable and has twisted Mukai vector $v_t$.
\end{Definition}

In the special case when $S=\Spec k$, we will write $\sM_{\sX}(v)$ for $\sM_{\sX/\Spec k}(v)$.

\begin{Proposition}\label{prop:rel}
  Let $\sX\to X$ be a $\m_p$-gerbe on a proper smooth family of K3 surfaces over a Henselian DVR $R$, and let $H$ and $v$ be as in Definition \ref{def:modulistack0}. Suppose that $v$ restricts to a primitive element of the twisted N\'{e}ron-Severi group of each geometric fiber. If $H$ is $v$-generic in each geometric fiber of $\sX/R$, then the moduli space $\sM_{\sX}(v)$ of $H$-stable twisted sheaves on $\sX$ with twisted Mukai vector $v$ is either empty or a $\G_m$-gerbe over a proper smooth scheme over $R$.
  
  In particular, every fiber of $\sM_{\sX}(v)\to\Spec R$ is a $\G_m$-gerbe over a K3 surface if and only if one geometric fiber is a $\G_m$-gerbe over a K3 surface.
\end{Proposition}
\begin{proof}
  Any stable twisted sheaf on a fiber is a smooth point of the morphism  $\sM_{\sX}(v)\to\Spec R$. Moreover, Langton's theorem (and the genericity of $H$) tells us that $\sM_{\sX}(v)$ is a $\G_m$-gerbe over its sheafification $M_{\sX}(v)$ and that $M_{\sX}(v)\to\Spec R$ is proper. It follows that $\sM_{\sX}(v)$ is a $\G_m$-gerbe over a proper smooth $R$-scheme.
  
  Since $\sM_{\sX}(v)\to\Spec R$ is smooth, its Stein factorization is finite étale. It follows that all geometric fibers have the same number of connected components. This concludes the proof.
\end{proof}

\begin{Theorem}\label{thm:moduli0}
  Let $\sX\to X$ be a $\m_p$-gerbe on a supersingular K3 surface and $v=(r,l,s)\in\tN(\sX)$ a primitive Mukai vector with $v^2=0$. If $H$ is sufficiently a generic, then the moduli space $\sM_{\sX}(v)$ of $H$-stable twisted sheaves on $\sX$ with twisted Mukai vector $v$ is either empty or satisfies
  \begin{enumerate}
    \item $\sM_{\sX}(v)$ is a $\G_m$-gerbe over a supersingular K3 surface $M_{\sX}(v)$,
    \item the universal sheaf $\sP$ on $\sM_{\sX}(v)\times\sX$ induces a Fourier-Mukai equivalence
    \[
      \Phi_{\sP}\colon D^{(\text{-}1)}(\sM_{\sX}(v))\to D^{(1)}(\sX)
    \]
    and
    \item the $\G_m$-gerbe $\sM_{\sX}(v)\to M_{\sX}(v)$ is trivial if and only if there exists a vector $w\in\tN(\sX)$ such that $v.w$ is coprime to $p$.
  \end{enumerate}
\end{Theorem}
\begin{proof}
  We begin by showing that $\sM_{\sX}(v)$ is either empty or a $\G_m$-gerbe over a K3 surface. We first reduce to the case when $\sX$ has trivial Brauer class. Consider the connected component of the universal family of $\m_p$-gerbes on $X$ that contains $\sX$ as a fiber. If $\sX$ has nontrivial Brauer class, then $v$ extends to a section of $\tuN_{X/\bA^1}\otimes\bQ$ that in each geometric fiber is isotropic and contained in the extended twisted N\'{e}ron-Severi lattice. Note that $v$ may become non-primitive in some fiber of the universal family, necessarily corresponding to a gerbe with trivial Brauer class. However, it follows however from the calculations of Proposition \ref{rem:nslattice} that there exists some fiber with trivial Brauer class where $v$ remains primitive. As in Lemma \ref{lem:genericpolarizationinfamily}, we may find a relative polarization that is $v$-generic in each fiber where $v$ remains primitive. We now apply Proposition \ref{prop:rel} twice to compare the original fiber to to the geometric generic fiber, and the geometric generic fiber to our chosen fiber with trivial Brauer class. It therefore suffices to prove the result when $\sX$ has trivial Brauer class. This follows from results of Mukai \cite{Muk84} (see also Corollary 10.2.4 and Proposition 10.2.5 of \cite{Huy16}). We conclude that $\sM_{\sX}(v)$ is either empty or a $\G_m$-gerbe over a K3 surface.

  Condition (2) then follows from a criterion of Bridgeland (Theorem 2.3 and Theorem 3.3 of \cite{MR1651025}). 
  By Theorem \ref{thm:moduli0} and Proposition \ref{prop:cohomologicaltransform} we know therefore that the rational cohomological transform is an isomorphism of $K$-vector spaces that is compatible with the Mukai pairing and the Frobenius operators. This shows that $M_{\sX}(v)$ is supersingular, completing the proof of (1). 
  
  Finally, we show (3). Let $\sM$ be a $\m_p$-gerbe over $M_{\sX}(v)$ whose associated $\G_m$-gerbe is $\sM_{\sX}(v)$. There is then a canonical equivalence of categories $\Coh^{(1)}(\sM)\cong\Coh^{(1)}(\sM_{\sX}(v))$. Under this equivalence, the universal sheaf of the moduli problem gives rise to a sheaf
  \[
    \sP\in \Coh^{(1,1)}(\sM\times\sX)
  \]
  The induced cohomological transform gives an isometry
  \[
    g\colon \tN(\sM)\to\tN(\sX)
  \]
  sending $(0,0,1)$ to $v$. Suppose that the gerbe $\sM_{\sX}(v)\to M_{\sX}(v)$ is trivial. This is equivalent to the existence of an invertible twisted sheaf, say $\sL$, on $\sM$. We have $\ch_{\sM}(\sL).(0,0,1)=1$, so $g(\ch_{\sM}(\sL)).v=1$. Conversely, suppose that there exists a $w\in\tN(\sX)$ such that $(v.w,p)=1$. By Proposition \ref{prop:TwistedTate}, there exists a $\sX$-twisted sheaf $\sE$ with $v_{\sX}(\sE)=w$. Consider the perfect complex
  \[
    \R q_*(p^*(\sE)\ltensor\sP^\vee)\in D^{(\text{-}1)}(\sM)
  \]
  The rank of this complex over a geometric point $x\in M_{\sX}(v)$ is
  \[
    \chi(\sE,\sP_x)=-v.w
  \]
  The existence of a such a complex implies that
  \[
    \left[\sM_{\sX}(v)\right]\in\Br(M_{\sX}(v))[v.w]
  \]
  The result follows.
\end{proof}

We treat the question of when our moduli spaces are non-empty. Our strategy is to lift to characteristic 0 and appeal to results of Yoshioka on the existence of semi-stable sheaves with prescribed invariants. Yoshioka's results rely on the global structure of the moduli space of K3 surfaces over the complex numbers, and hence are analytic in nature. This represents the unique points in this paper where we (implicitly) use analytic techniques. It would be interesting to try to remove this dependence.
\begin{Proposition}\label{prop:nonempty}
    Let $\sX\to X$ be a $\m_p$-gerbe on a supersingular K3 surface. If $v=(r,l,s)\in\tN(\sX)$ is a primitive Mukai vector with $v^2=0$, then the moduli space $\sM_{\sX}(v)$ with respect to a $v$-generic polarization is non-empty if one of the following holds:
    \begin{enumerate}
      \item $r>0$.
      \item $r=0$ and $l$ is effective.
      \item $r=l=0$ and $s>0$.
    \end{enumerate}
\end{Proposition}
\begin{proof}
  By Proposition \ref{prop:rel}, we reduce as in Theorem \ref{thm:moduli0} to the case when $\sX$ has trivial Brauer class. Thus, $\sX=\left\{\sL^{1/p}\right\}$ for some line bundle $\sL$.
  
  After tensoring with an appropriate line bundle, we may assume that $l$ is ample. By Theorem A.1 of \cite{LO15}, we find a complete DVR $R$ with fraction field $K$ of characteristic 0 and residue field $k$, a relative K3 surface $X_R\to\Spec R$, a line bundle $\sL_R$ on $X_R$, a vector $v_R\in A^*(X_R)\otimes\bQ$, and a class $H_R\in\Pic(X_R)$, together with an isomorphism $X\cong X_{R}\otimes k$ under which $\sL_R$ restricts to $\sL$, $v_R$ restricts to $v$, and $H_R$ restricts to $H$.
  
  
  Set $\sX_R=\left\{\sL_R^{1/p}\right\}$, so that $\sX_R\to\Spec R$ is a relative twisted K3 surface whose special fiber is isomorphic to $\sX$. Note that in each geometric fiber $v_R$ restricts to an element of the twisted N\'{e}ron-Severi group that is primitive and isotropic. Moreover, $H_R$ is $v_R$-generic in each geometric fiber. By Theorem 3.16 of \cite{Yos06}, the geometric generic fiber is non-empty. By Proposition \ref{prop:rel} this gives the result.
  
  
\end{proof}

Finally, we record the following lemma, which appears to be well known.

\begin{Lemma}\label{lem:locallyfree}
  If $v=(r,l,s)$ is primitive and $v^2=0$ and $H$ is $v$-generic, then any object $\sE\in\sM_{\sX}(v)(T)$ is locally free.
\end{Lemma}
\begin{proof}
  By the local criterion for flatness, it will suffice to prove this when $T=\Spec k$. If $\sE$ is stable, then $\sE^{\vee\vee}$ is also stable. The Mukai vector of $\sE^{\vee\vee}$ is $(r,l,s')$ for some $s'\leq s$, with equality if and only if $\sE$ is locally free. The sheaf $\sE^{\vee\vee}$ gives a point in the moduli space of stable sheaves on $\sX$ with Mukai vector $(r,l,s')$. These moduli spaces are always smooth, and are either empty or of the expected dimension $l^2-2rs'$. If $s'<s$, this number is negative, a contradiction.
\end{proof}


\subsection{Twistor families of positive rank}\label{sec:constructionoftwistorlines}

We maintain Notation \ref{def:hyperbolicplane2}, so that $\Lambda$ is a fixed supersingular K3 lattice and $\tLam=\Lambda\oplus U_2(p)$. In Section \ref{sec:twistorlines}, we defined a twistor line $f_v\colon \bA^1\to \oM_{\tLam_0}$ to be a connected component of a fiber of $\pi_v$ for some isotropic vector $v\in\tLam_0$, where $\pi_v$ fits into the diagram
\[
  \begin{tikzcd}[column sep=tiny]
    \oM_{\tLam_0}\arrow[hookleftarrow]{rr}\arrow[dashed]{dr}& &\oM^{\langle v\rangle}_{\tLam_0}\arrow{dl}{\pi_v}\\
       & \oM_{v^\perp/v} &
  \end{tikzcd}
\]
We will make the following (somewhat preliminary) definition.
\begin{Definition}\label{def:twistorfamily}
  A  family of twisted supersingular K3 surfaces over an open subset $U\subset\bA^1$ is a \textit{twistor family} if it admits a marking such that the induced map from $U$ to the period domain is an isomorphism onto an open subset of a twistor line.
\end{Definition}

In Theorem \ref{thm:modularinterpretation}, we will give another characterization of twistor families as relative moduli spaces of twisted sheaves on universal families of $\m_p$-gerbes over an open subset of a connected component of the group of $\m_p$-gerbes on a supersingular K3 surface. For reasons that will become clear shortly, we will make the following distinction.

\begin{Definition}\label{def:kindsoftwistorlines}
  A twistor line (or twistor family) is \textit{of positive rank} if $v.e\neq 0$, and is \textit{Artin-Tate} if $v.e=0$.
\end{Definition}

We will restrict our attention in this section to the positive rank case, and postpone our study of the Artin-Tate case to Section \ref{sec:artintate}.
The main result of this section is Theorem \ref{thm:constructionoftwistorlines}, which says that, under certain restrictions, twistor lines of positive rank lift to families of twisted K3 surfaces, in a strong sense. We will apply Theorem \ref{thm:constructionoftwistorlines} in Section \ref{sec:torellitheorem} to deduce consequences for moduli spaces of marked twisted supersingular K3 surfaces. 

We have already constructed moduli spaces of twisted sheaves in Theorem \ref{thm:moduli0}, which for a generic choice of polarization are $\G_m$-gerbes over supersingular K3 surfaces. In order to examine these stacks in terms of our period morphism, we will identify an associated $\m_p$-gerbe.

\begin{Definition}\label{def:modulistack}
  Let $\sX\to X$ be a $\m_p$-gerbe over a supersingular K3 surface, $v=(r,l,s)\in\tN(\sX)$ a Mukai vector with $r>0$, and $H$ a polarization on $X$. Let $\sL$ be a line bundle on $\sX$ with first Chern class $l$. We define a stack $\sM^{\det}_{\sX}(v)$ on $\Spec k$ whose objects over a $k$-scheme $T$ are pairs $(\sE,\phi)$ where $\sE$ is a $T$-flat $\sX_T$-twisted sheaf that is locally of finite presentation such that for each geometric point $t\in T$ the fiber $\sE_t$ is $H$-stable and has Mukai vector $v$, and $\phi\colon \det\sE\iso\sL_T$ is an isomorphism of invertible sheaves on $\sX_T$.
\end{Definition}

Note that this definition is restricted to the case of positive rank.

\begin{Proposition}\label{prop:integralcorrespondence}
   Let $\sX\to X$ be a $\m_p$-gerbe on a supersingular K3 surface and $v=(p,l,s)\in\tN(\sX)$ a primitive Mukai vector of rank $p$ satisfying $v^2=0$. The cohomological correspondence
    \[
      \Phi^{\cris}_{v_{\sM\times\sX}(\sP)}\colon \tH^{(\text{-}1)}(\sM^{\det}_{\sX}(v)/W)\to\tH(\sX/W)
    \]
    (see Notation \ref{notn:twistedcrystal-thang}) induced by the universal twisted sheaf
   \[
     \sP\in\Coh^{(1,1)}(\sM^{\det}_{\sX}(v)\times\sX)
   \]
   is an isomorphism of K3 crystals that is compatible with the inclusions of the respective twisted N\'{e}ron-Severi lattices.
\end{Proposition}
\begin{proof}
  Define the subgroup
  \[
    \tN(\sM^{\det}_{\sX}(v)\times\sX)\defeq\tN(\sM^{\det}_{\sX}(v))\uplus\tN(\sX)\subset A^*(M_{\sX}(v)\times X)\otimes\bQ
  \]
  By the results of Theorem \ref{thm:moduli0} and Proposition \ref{prop:cohomologicaltransform}, it remains only to show that
  \[
    v_{\sM^{\det}_{\sX}(v)\times\sX}(\sP)\in\tN(\sM^{\det}_{\sX}(v)\times\sX)
  \]
  Let us first suppose that the Brauer classes of the gerbes $\sX$ and $\sM^{\det}_{\sX}(v)$ are both trivial. By tensoring with twisted invertible sheaves, we may reduce to the case that the gerbes $\sX$ and $\sM^{\det}_{\sX}(v)$ are themselves trivial. The result then follows by applying the Grothendieck-Riemann-Roch theorem to the respective projections as in the untwisted case (see Lemma 10.6 of \cite{Huy06}).
  
  We now prove the result in general. Consider the connected component $\tsX\to\bA^1$ of the universal family of $\m_p$-gerbes on $X$ containing $\sX$ as a fiber. Let $U\subset\bA^1$ be the open locus where $v$ remains primitive in the twisted N\'{e}ron-Severi group of the fiber. Expand $v$ in the basis considered in Proposition \ref{rem:nslattice}, and let $s\in\bZ$ be the coefficient of $(0,0,1)$. By tensoring with a line bundle of non-zero degree, we may assume without loss of generality that $s$ is not divisible by $p$. By Lemma \ref{lem:genericpolarizationinfamily}, we may choose $H$ to be $v$-generic in each fiber over a point in $U$. Consider the relative moduli space $\sM_{\tsX_U}(v)\to U$. Let
  \[
    \tsP\in\Coh^{(1,1)}(\sM^{\det}_{\tsX/\bA^1}(v)\times_{\A^1}\tsX)
  \]
  be the relative universal twisted sheaf.  Consider a point $t\in U$ such that the gerbe $\tsX_{t}$ has trivial Brauer class. There exists then a vector $w\in\tN(\tsX_t)$ with rank 1. As $p$ does not divide $s$, part (3) of Theorem \ref{thm:moduli0} gives that the fiber $(\sM_{\tsX}(v))_t$ of the moduli space also has trivial Brauer class. Thus, by the previous case the result holds for $\tsP_t$. By the constancy of the twisted Chern character in flat families, the result also holds for $\sP$.
\end{proof}
\begin{Remark}
  The previous result presumably holds for the cohomological correspondence induced by any Fourier-Mukai equivalence. The above special case will suffice for our purposes in this section.
\end{Remark}

We make the following observation regarding markings of tautological families of $\m_p$-gerbes.

\begin{Lemma}\label{lem:universaltwistorlines}
  Let $\pi\colon X\to\Spec k$ be a supersingular K3 surface and $\sX_0\to X$ a $\m_p$ gerbe. Let $\bA^1\subset \R^2\pi^{\fl}_*\m_p$ be the connected component containing $\sX_0$ as a fiber, and $\sX\to\bA^1$ the corresponding tautological family. If the Brauer class of $\sX_0$ is non-trivial, then any marking $\tLam\to\tN(\sX_0)$ extends to a map $\tuLam_{\bA^1}\to\tuN_{X\times\bA^1/\bA^1}\otimes\bQ$ which on each geometric fiber $t\in\bA^1$ restricts to a map $\tLam\to\tN(\sX_t)$.
\end{Lemma}
\begin{proof}
  This follows immediately from the calculations of Proposition \ref{rem:nslattice}.
\end{proof}

We record some purely lattice-theoretic facts. Consider an isometric inclusion $\tLam\subset\tN$ of extended supersingular K3 lattices (or of supersingular K3 lattices). We have a chain of inclusions
\[
  \tLam\subset\tN\subset\tN^*\subset\tLam^*
\]
\begin{Lemma}\label{lem:latticefact1}
  The $\bF_p$-vector spaces $\tN/\tLam$ and $\tLam^*/\tN^*$ are naturally dual. In particular, they have the same dimension.
\end{Lemma}
\begin{proof}
  Applying $\Hom_{\bZ}(\_,\bZ)$ to the short exact sequence
  \[
    0\to\tLam\to\tN\to\tN/\tLam\to 0
  \]
  we get a short exact sequence
  \[
    0\to\tN^*\to\tLam^*\to\Ext_{\bZ}^1(\tN/\tLam,\bZ)\to 0
  \]
  Applying $\Hom_{\bZ}(\tN/\tLam,\_)$ to the short exact sequence $0\to\bZ\to\bQ\to\bQ/\bZ\to 0$, we find a natural isomorphism
  \[
    \Hom_{\bF_p}(\tN/\tLam,\bF_p)\iso\Ext_{\bZ}^1(\tN/\tLam,\bZ)
  \]
  This gives the result.
\end{proof}
Consider the containments
\[
  \dfrac{p\tN}{p\tLam}\subset\dfrac{p\tN^*}{p\tLam}\subset\dfrac{p\tLam^*}{p\tLam}=\tLam_0
\]
of $\bF_p$-vector spaces.
\begin{Lemma}\label{lem:latticefact2}
  As subspaces of $\tLam_0$, we have
  \[
    \left(\dfrac{p\tN}{p\tLam}\right)^\perp=\dfrac{p\tN^*}{p\tLam}
  \]
\end{Lemma}
\begin{proof}
    It is immediate that the right hand side is contained in the left hand side. We will show that they have the same dimensions.
  Let $\sigma_0$ be the Artin invariant of $\tLam$ and $\sigma$ the Artin invariant of $\tN$, so that $p\tLam^*/p\tLam$ has dimension $2\sigma_0$, and $p\tN^*/p\tN$ has dimension $2\sigma$. Consider the short exact sequences
  \[
    0\to \dfrac{p\tN^*}{p\tLam}\to\dfrac{p\tLam^*}{p\tLam}\to\dfrac{p\tLam^*}{p\tN^*}\to 0\hspace{2cm}
    0\to \dfrac{p\tN}{p\tLam}\to\dfrac{p\tN^*}{p\tLam}\to\dfrac{p\tN^*}{p\tN}\to 0
  \]
  Combined with Lemma \ref{lem:latticefact1}, we find that
  \[
    \dim_{\bF_p}\left(\dfrac{p\tN}{p\tLam}\right)=\dim_{\bF_p}\left(\dfrac{p\tLam^*}{p\tN^*}\right)=\sigma_0-\sigma,\mbox{ and}\hspace{.5cm}\dim_{\bF_p}\left(\dfrac{p\tN^*}{p\tLam}\right)=\sigma_0+\sigma
  \]
  This gives the result.
\end{proof}


\begin{Lemma}\label{lem:isotropicvectors}
  If $N$ is a supersingular K3 lattice or an extended supersingular K3 lattice, and $v\in N_0$ is an isotropic vector, then there exists a primitive isotropic vector $l\in pN^*\subset N$ whose image in $N_0$ is equal to $v$.
\end{Lemma}
\begin{proof}
  The existence of an isotropic vector in $N_0$ implies that $\sigma_0(N)\geq 2$. By the explicit presentation of supersingular K3 lattices in \cite{RS79} and Lemma \ref{lem:supersingularlattices}, we see that there is an orthogonal decomposition $N=N'\oplus U_2(p)$. Let $\left\{e,f\right\}$ be the standard generators for $U_2(p)$, and let $w$ be the image of $e$ in $N_0$. By Witt's Lemma (Theorem \ref{thm:Witt}), there exists an isometry $g\in O(N_0)$ taking $w$ to $v$. By a result of Nikulin (see Theorem 14.2.4 of \cite{Huy16}), the map
  \[
    O(N)\to O(N_0)
  \]
  is surjective, so we find an orthogonal transformation $h\in O(N)$ inducing $g$. The vector $h(e)\in N$ is primitive and isotropic, and its image in $N_0$ is equal to $V$, as desired.
\end{proof}

We are now ready to prove the main result of this section. This should be viewed as a (partial) supersingular analog of Proposition 3.9, Chapter 7 of \cite{Huy16}, which describes those curves in the Hodge-theoretic period domain that lift to twistor spaces of complex analytic K3 surfaces.

\begin{Theorem}\label{thm:constructionoftwistorlines}
  Let $x\in\sS_{\Lambda}^o(k)$ be a $k$-point. Let $L\subset\oM_{\tLam_0}$ be a twistor line corresponding to an isotropic vector $v\in\tLam_0$, and let $f_v\colon L\to\oM_{\tLam_0}$ be the inclusion. Suppose that $\trho(x)\in L$ and that the Artin invariant of $\trho(x)$ is equal to the generic Artin invariant of $L$. Set $U=L\cap\oM_{\tLam_0}^{\langle e\rangle}$. If $v.e\neq 0$, then there exists a lift $f\colon U\to \sS_{\Lambda}^o$ such that the diagram
  \[
    \begin{tikzcd}
      &\sS_{\Lambda}^o\arrow{d}{\trho}\\
      U\arrow{ur}{f}\arrow{r}{f_v|_U}&\oM_{\tLam_0}^{\langle e\rangle}
    \end{tikzcd}
  \]
  commutes and $x\in f(U)$.
\end{Theorem}
\begin{Remark}
  Let us explain the idea behind the proof. Let $\sX\to X$ be a $\m_p$-gerbe corresponding to $x$. Using $v$, we form an appropriate moduli space of sheaves, say $\sY\to Y$, on $\sX$ (in fact, we may need to first replace $\sX$ by another $\m_p$-gerbe with the same Brauer class). Let $\tsY\to\bA^1$ be the corresponding universal twistor family containing $\sY$ as a fiber. By our relative twisted Torelli theorem (Proposition \ref{prop:commutativetriangle}), the induced map from the base $\bA^1$ to the period domain identifies $\bA^1$ with a twistor line. We show that $\sX$ is also a moduli space of twisted sheaves on $\sY$, an instance of Mukai duality. By taking an appropriate relative moduli space of sheaves on the family $\tsY$, we obtain a family of twisted surfaces containing the original surface $\sX$ as a fiber. Our assumptions on $v$ will allow us to ensure that moduli space $\sY$ parametrizes sheaves of rank $p$, and that the gerbe $\sY$ has Brauer class of order $p$. Together, these conditions enable us to ignore stability conditions at a key moment. It is possible that this could be eliminated with enough knowledge of the stability of ``wrong-way'' slices of universal sheaves. After proving our twisted crystalline Torelli theorem (a consequence of this result!) we will explain in Theorem \ref{thm:modularinterpretation} how to remove certain of our assumptions on $v$.
\end{Remark}
\begin{proof}
  The $k$-point $x$ of $\sS_\Lambda^o$ corresponds to a $\m_p$-gerbe
  $\sX\to X$ over a supersingular K3 surface $X$ along with a marking
  $m\colon\Lambda\to\Pic(X)$. Let us identify $\Lambda$ with its image
  in $\Pic(X)$, so that $m$ is just the canonical inclusion. By Lemma \ref{lem:inducedmarking}, we obtain an induced inclusion $\tLam\subset\tN(\sX)$ identifying $e$ with $(0,0,1)$ and $f$ with $(p,0,0)$.
  
  \begin{Claim}
    There exists an element $x=(p,l,l^2/2p)\in p\tLam^*$ such that the image of $x$ in $\tLam_0$ is a non-zero scalar multiple of $v$.
  \end{Claim}
  By Lemma \ref{lem:isotropicvectors}, we may find a primitive
  isotropic vector $(r,l,s)\in p\tLam^*\subset\tLam$ whose image in
  $\tLam_0$ is $v$. Note that $r$ is necessarily of the form $pa$ for
  some integer $a$, and $a$ must be invertible modulo $p$ by our
  assumption that $v.e\neq 0$. Let $b$ be an integer such that
  $ab\equiv 1\mod{p}$. Consider the vector
  $x=(p,bl,ab^2s)\in\tLam$. It is immediate that this vector is
  isotropic. As $(r,l,s)$ is primitive, so is $(p,bl,ab^2s)$. Finally,
  note that $ab^2-b$ is divisible by $p$, and that
  \[
    x=(p,bl,ab^2s)=b(pa,l,s)+p\left(1-ab,0,s\frac{ab^2-b}{p}\right)
  \]
  Because $1-ab$ is divisible by $p$,
  $(1-ab,0,s(ab^2-b)/p)\in\tLam$. It follows that $x\in p\tLam^*$ and
  that the image of $x$ in 
  $\tLam_0$ is a non-zero scalar multiple of $v$. This completes the proof of the claim.
  
  We will henceforth let $x=(p,l,l^2/2p)\in p\tLam^*\subset\tLam$ denote a fixed vector with these properties.
  \begin{Claim}
    The image of $x$ in $\tN(\sX)$ is primitive and isotropic, and is contained in $p\tN(\sX)^*$.
  \end{Claim}
  Let $K$ be the characteristic subspace corresponding to $\trho(x)$. The twistor line $L$ is by definition contained in the locus of characteristic subspaces not containing $v$, so $v\notin K$. But as subspaces of $\tLam_0$ we have
  \[
    K\cap\tLam_0=\dfrac{p\tN(\sX)}{p\tLam}
  \]
  so the image of $x$ in $\tN(\sX)$ remains primitive. By Lemma \ref{lem:latticefact1}, we have
  \[
    \left(\dfrac{p\tN(\sX)}{p\tLam}\right)^\perp=\dfrac{p\tN(\sX)^*}{p\tLam}
  \]
  By Lemma \ref{lem:genericArtin}, our assumption on the Artin invariant of $x$ implies that $v\in (K\cap\tLam_0)^\perp$, and therefore the image of $x$ in $\tN(\sX)$ is contained in $p\tN(\sX)^*$. This completes the proof of the claim.
  
  Using this lift, we will first construct a particular family of twisted supersingular K3 surfaces over an open subset of $\bA^1$ that contains $\sX$ as a fiber. This family will come with a natural marking by $\tLam$, and we will show that it satisfies the conclusions of the theorem. 
  
  Fix a line bundle $\sL$ on $X$ with first Chern class $l$, and consider the stack $\sX'=\sX\{\sL^{\vee 1/p}\}$ (see Definition \ref{def:relativeroot}). By Proposition \ref{prop:twistedlinebundleequivalence}, the universal bundle induces a derived equivalence, and the corresponding map on cohomology is given by
  \[
    e^{-l/p}:\tH(\sX/W)\to\tH(\sX'/W)
  \]
  In particular, note that $e^{-l/p}(p,l,l^2/2p)=(p,0,0)$. By Theorem
  \ref{thm:moduli0}, the moduli space $\sY=\sM^{\det}_{\sX'}(p,0,0)$
  of stable twisted sheaves on $\sX'$ with respect to a generic
  polarization is a $\m_p$-gerbe over a supersingular K3 surface $Y$,
  and there is a universal object
  \[
    \sP\in \Coh^{(1,1)}(\sY\times\sX')
  \]
  equipped with an isomorphism
  \[
    \det\sP\iso\cO_{\sY\times\sX'}
  \]
  The isometry
  \[
    g\colon\tN^{(\text{-}1)}(\sY)\to\tN(\sX')
  \]
  (see Notation \ref{notn:ns-thang})
  induced by the kernel $v(\sP)$ satisfies $(0,0,1)\mapsto (p,0,0)$
  and $(p,0,0)\mapsto (0,0,1)$. We have $(p,0,0)\in p\tN(\sX')^*$, so
  by Theorem \ref{thm:moduli0} the Brauer class of the gerbe
  $\sY\to Y$ is non-trivial. Let $\pi\colon Y\to\Spec k$ be the
  structure map, and consider the connected component
  $\bA^1\subset \R^2\pi^{\fl}_*\m_p$ that contains $\sY$ as a
  fiber. Let $\tsY\to\bA^1$ be the universal family. Let
  $V\subset\bA^1$ be the locus of points $t\in\bA^1(k)$ where the image of the class
  $(p,0,0)=g^{-1}(0,0,1)$ in $\tN^{(\text{-}1)}(\sY_t)$ remains
  primitive. Note that if $\sY$ deforms the trivial gerbe then $V$ is the complement of the origin, and otherwise $V=\bA^1$.
  We form the relative moduli space
  \[
    \tsZ'=\sM^{\det}_{\tsY|_V/V}(p,0,0)\to V
  \]
  of twisted sheaves on the fibers of $\tsY|_V\to V$ that are stable with respect to a sufficiently generic polarization, which comes with a universal object
  \[
    \tsQ\in \Coh^{(1,1)}(\tsZ'\times_{V}\tsY|_V)
  \]
  Let $\sZ'$ be the fiber corresponding to $\sY$. We will construct an isomorphism $\sX'\iso\sZ'$. By Lemma \ref{lem:locallyfree}, the universal twisted sheaf $\sP$ is locally free, and in particular flat over $\sY$. Because the Brauer class of $\sY\to Y$ is non-trivial, stability conditions for rank $p$ twisted sheaves on $\sY$ are vacuous. Thus, $\sP$ gives (by descent) an object of $\sM^{\det}_{\sY}(p,0,0)(\sX')$, and hence a morphism
  \[
    \Theta=[\sP]\colon \sX'\to\sZ'
  \]
  We claim that this is an isomorphism. The restriction $\sQ$ of $\tsQ$ to $\sZ'\times\sY$ induces a Fourier-Mukai equivalence $D^{(\text{-}1)}(\sZ')\to D^{(1)}(\sY)$. The sheaf $\sQ^\vee\in\Coh^{(\text{-}1,\text{-}1)}(\sZ'\times\sY)$ also induces a Fourier-Mukai equivalence $D^{(1)}(\sZ')\to D^{(\text{-}1)}(\sY)$ (see for instance Theorem 1.6.15 of \cite{Navas10}). Consider the ``wrong-way'' Fourier-Mukai transforms
  \begin{align*}
    \Phi_{\sQ}^o\colon D^{(\text{-}1)}(\sY)\to D^{(1)}(\sZ')&\hspace{1cm}&\Phi_{\sQ^\vee}^o\colon D^{(1)}(\sZ')\to D^{(\text{-}1)}(\sY)
  \end{align*}
  One shows as in Theorem 1.6.15 of \cite{Navas10} that these are both equivalences, and that there is a natural isomorphism of functors $\Phi_{\sQ^\vee}^o\circ\Phi_{\sQ}^o\cong [-2]$ (see Remark 5.8 of Chapter 5 in \cite{Huy06}). By construction, the map
  \[
    \Theta\times\id\colon \sX'\times\sY\to\sZ'\times\sY
  \]
  satisfies $(\Theta\times\id)^*\sQ\cong\sP$. Using this and the projection formula, we obtain a natural isomorphism
  \[
    \mathbf{L}\Theta^*\circ\Phi_{\sQ}^o\cong\Phi_{\sP}
  \]
  where
  \[
    \Phi_{\sP}\colon D^{(\text{-}1)}(\sY)\to D^{(1)}(\sX')
  \]
  is the natural Fourier-Mukai transform associated to $\sP$. We therefore find that
  \[
    \mathbf{L}\Theta^*\cong\Phi_{\sP}\circ\Phi^o_{\sQ^\vee}[2]
  \]
  In particular, we conclude that $\mathbf{L}\Theta^*$ is an
  equivalence of categories. By looking at structure sheaves of closed
  points, this shows that the map $\Theta$ induces an isomorphism
  $X\iso Z$. As $\sP$ is a $\sX'$-twisted sheaf, its
  $\cO_{\sX'}^\times$-action is identified with its canonical action by the
  inertia group of $\sX'$. It follows that $\Theta$ is in fact a
  morphism of $\m_p$-gerbes, and hence is an isomorphism.
  
  We have realized $\sX'$ as a fiber of a family of twisted K3 surfaces. We wish to do the same for $\sX$. Consider the isometries
  \[
    \tN(\sZ')\xrightarrow{h}\tN^{(\text{-}1)}(\sY)\xrightarrow{g}\tN(\sX')
  \]
  induced by $\Phi_{\sQ}$ and $\Phi_{\sP}$. Let $\tZ'$ be the coarse space of $\tsZ'$ and $Z'$ the coarse space of $\sZ'$. Let $\sL'=\Theta_*(\sL)\in\Pic(Z')$. We claim that $\sL'$ is the restriction of a line bundle on $\tZ'$. Indeed, consider the class $(p,-l,l^2/2p)\in\tN(\sX')$. We have that
  \[
    h^{-1}\circ g^{-1}(p,-l,l^2/2p)=(p,-l',l'^2/2p)
  \]
  where $l'$ is the first Chern class of $\sL'$. Any class in
  $\tN^{(\text{-}1)}(\sY)$ extends to a section of
  $\tuN_{Y_{\bA^1}/\bA^1}^{(\text{-}1)}\otimes\bQ$ which restricts
  over each geometric point $t\in\bA^1$ to a class in
  $\tN^{(\text{-1})}(\tsY_t)$. Thus, $(p,-l',l'^2/2p)$ extends to a
  section of $\tuN_{\tZ'/V}\otimes\bQ$ which restricts over each
  geometric point $t\in V$ to a class in $\tN(\tsZ'_t)$. It follows
  that $l'$ extends to a class in the relative Picard group of $\tZ'$,
  and so $\sL'$ extends to a line bundle, say $\tsL'$, on
  $\tZ'$. Define
  \[
    \tsZ=\tsZ'\{\tsL'^{1/p}\}
  \]
  This is a $\m_p$-gerbe over $\tZ'$, and is equipped with a morphism
  $\tsZ\to V$ (see Definition \ref{def:relativeroot}). If $\sZ$ is the
  fiber corresponding to $\sZ'$, then 
  $\Theta$ induces an isomorphism
  \[
    \sX=\sX'\{\sL^{1/p}\}\iso\sZ'\{\sL'^{1/p}\}=\sZ
  \]
  Thus, $\sX$ is isomorphic to a fiber of the family $\tsZ\to V$.
  
  We will now complete the proof. The derived equivalences we have constructed induce isometries
  \begin{equation}\label{eq:isometries}
    \begin{tikzcd}
      \tN(\sX)\arrow{r}{e^{-l/p}}&\tN(\sX')\arrow{r}{g^{-1}}\arrow[bend right=20]{rr}[swap]{\Theta_*}&\tN^{(\text{-}1)}(\sY)\arrow{r}{h^{-1}}&\tN(\sZ')\arrow{r}{e^{l'/p}}&\tN(\sZ)
    \end{tikzcd}
  \end{equation}
  By our construction, the family $\tsZ$ carries a natural marking by $\Lambda$. Indeed, the marking $\tLam=\Lambda\oplus U_2(p)\to\tN(\sX)$ of the extended N\'{e}ron-Severi group of $\sX$ induces a marking $\tLam\to\tN(\sX')$, and hence a map $\tLam\to\tN^{(\text{-}1)}(\sY)$, which by Lemma \ref{lem:universaltwistorlines} extends to a map $\tuLam_{\bA^1}\to\tuN_{Y_{\bA^1}/\bA^1}^{(\text{-}1)}\otimes\bQ$ restricting over each geometric point $t\in\bA^1$ to a map $\tLam\to\tN^{(\text{-1})}(\tsY_t)$. We thus find the same for the families $\tsZ'$ and $\tsZ$. The composition $\tN(\sX)\to\tN(\sZ)$ satisfies $(0,0,1)\mapsto (0,0,1)$ and $(p,0,0)\mapsto (p,0,0)$, so there is an induced marking $\Lambda\to\Pic_{\tsZ/V}$. The resulting marked family gives a morphism
  \[
    f_0\colon V\to\sS^o_{\Lambda}
  \]
  The isomorphism $\sX\iso\sZ$ is compatible with the respective markings by $\Lambda$, so the image of $f_0$ contains $x$. Consider the composition
  \[
    \trho\circ f_0\colon V\to \oM^{\langle e\rangle}_{\tLam_0}
  \]
  It follows directly from our constructions that the image of $V$ under $\trho\circ f_0$ is $U$, and moreover that the map $V\to U$ is a bijection on closed points. To obtain the stronger result that we have claimed, we need to show that in fact $\trho\circ f_0$ maps $V$ isomorphically onto $U$. To see this, consider the induced isometry $\tLam\to\tN(\sY)$, which satisfies $(p,l,l^2/2p)\mapsto (0,0,1)$ and $e=(0,0,1)\mapsto (p,0,0)$. Set $\Lambda'=(p,l,s)^{\perp}/(p,l,s)$ (a subquotient of $\tLam$). Set $\tLam'=\Lambda'\oplus U_2(p)$, and for clarity let us denote the standard basis of this copy of $U_2(p)$ by $e',f'$. The family $\tsY\to\bA^1$ carries a natural marking by $\Lambda'$, and by Proposition \ref{prop:commutativetriangle} this family induces a diagram
  \[
    \begin{tikzcd}
      &\sS_{\Lambda'}^o\arrow{d}{\trho}\\
      \bA^1\arrow{ur}{f_1}\arrow[hook]{r}&\oM_{\tLam'_0}^{\langle e'\rangle}
    \end{tikzcd}
  \]
  where the lower horizontal arrow is the inclusion of the twistor line (up to our identification of this line with $\bA^1$). Let $K(\tsY)$ be the sub-bundle of $\tLam_0\otimes\cO_{\bA^1}$ corresponding to the relative Mukai crystal of $\tsY$ with its $T=\tLam\otimes\bZ_p$-structure, and let $K(\tsZ)$ be the sub-bundle of $\tLam_0\otimes\cO_V$ corresponding to the relative Mukai crystal of $\tsZ$ with its $\tLam\otimes\bZ_p$-structure. Applying Proposition \ref{prop:cohomologicaltransform} (and Proposition \ref{prop:integralcorrespondence}), we see that $K(\tsY)|_V$ and $K(\tsZ)$ are equal as sub-bundles of $\tLam_0\otimes\cO_{V}$. Furthermore, the map $\tLam'\to\tLam$ induced by the derived equivalence between $\sY$ and $\sZ$ and satisfies $e'\mapsto (p,l,l^2/2p)$ and $f'\mapsto (0,0,1)$. It follows that the map $f_0\colon V\to\sS^o_{\Lambda}$ induced by $\tsZ\to V$ maps $V$ isomorphically onto $U$. Composing with the inverse of this isomorphism, we find a lift $f$ of $f_v|_U$ and a diagram 
    \[
    \begin{tikzcd}
      &\sS_{\Lambda}^o\arrow{d}{\trho}\\
      U\arrow{ur}{f}\arrow{r}{f_v|_U}&\oM_{\tLam_0}^{\langle e\rangle}
    \end{tikzcd}
  \]
  As previously noted, the image of $f$ contains $x$, so this completes the proof.

\end{proof}

To apply this theorem, we will use the following technical lemma. Let us adopt Notation \ref{notation:vectorspace}, except that we shall write $e,f$ for the generators $v,w$ of $U_2$. Let $\tK\subset\tV\otimes k$ be a characteristic subspace such that $e\notin\tK$, and let $K=\pi_e(K)$. Let $\sigma$ be the Artin invariant of $K$.
\begin{Lemma}\label{lem:technical}
  If $\sigma\geq 2$, then there exists an isotropic vector $v\in\tV$ such that $v.e\neq0$, $v\notin\tK$, and $v\in (\tK\cap \tV)^\perp$.
\end{Lemma}
\begin{proof}
  As discussed in Remark \ref{rem:descriptionoffiber}, there exists an element $B\in V\otimes k$, uniquely determined up to elements of $K$, such that
  \[
    \tK=\left\langle x_1+(x_1.B)e,\dots,x_{\sigma_0}+(x_{\sigma_0}.B)e,f+B+\dfrac{B^2}{2}e\right\rangle
  \]
  Consider the vector space $(K\cap V)^\perp/(K\cap V)$, which has dimension $2\sigma$. The natural form on this space is non-degenerate, and as $\sigma\geq2$, it has dimension greater than or equal to 4. Thus, it contains a nonzero isotropic vector (by, for instance, Proposition \ref{prop:classification}). We may therefore find a nonzero isotropic vector $x\in (K\cap V)^\perp$ such that $x\notin K\cap V$.
  
  Let $y_1,\dots,y_{\sigma_0-\sigma}$ be a basis for $K\cap V$. We consider first the case that the image of $B$ in $V\otimes k/K$ is not contained in the subgroup 
  \[
    \dfrac{V}{K\cap V}=\dfrac{V+K}{K}\subset\dfrac{V\otimes k}{K}
  \]
  We have that
  \[
    \tK\cap\tV=\left\langle y_1+(y_1.B)e,\dots,y_{\sigma_0-\sigma}+(y_{\sigma_0-\sigma}.B)e\right\rangle
  \]
  It is immediate that the vector $v=f+x$ has the desired properties.
  
  Next, suppose that the image of $B$ is contained in $V/K\cap V$. Without loss of generality, we may then assume that $B\in V$, and therefore
  \[
    \tK\cap\tV=\left\langle y_1+(y_1.B)e,\dots,y_{\sigma_0-\sigma}+(y_{\sigma_0-\sigma}.B)e,f+B+\dfrac{B^2}{2}e\right\rangle
  \]
  One checks that the vector $v=f+x+B+(x.B+B^2/2)e$ has the desired properties.
  
  
  
\end{proof}

\subsection{Artin-Tate twistor families}\label{sec:artintate}

In this section we discuss Artin-Tate twistor families. We begin by interpreting moduli spaces of rank 0 twisted sheaves in terms of elliptic fibrations. We then interpret the fibers of Artin-Tate families in terms of torsors, giving a geometric manifestation of the isomorphism $\Br\cong\Sha$.

\begin{Definition}
  An \textit{elliptic fibration} on a surface $X$ is a proper flat morphism $f\colon X\to\P^1$ whose generic fiber is a smooth genus $1$ curve. An elliptic fibration is \textit{Jacobian} if it admits a section $\P^1\to X$, and is \textit{non-Jacobian} otherwise. A \textit{multisection of degree} $n$ is an integral subscheme $\Sigma\subset X$ such that $f|_{\Sigma}\colon \Sigma\to\P^1$ is finite flat of degree $n$.
\end{Definition}

 
Let $f\colon X\to\P^1$ be an elliptic fibration on a supersingular K3 surface, and $p\colon\sX\to X$ a $\m_p$-gerbe. Let $E\subset X$ be a smooth fiber of $f$. Fix an integer $s$, and consider the Mukai vector
\[
  v=(0,E,s)
\]
Note that $v$ is primitive and $v^2=0$. Fix a generic polarization $H$ on $X$. By Theorem \ref{thm:moduli0} the moduli stack $\sM_{\sX}(v)$ is a $\G_m$-gerbe over a supersingular K3 surface. We wish to identify a natural choice of a corresponding $\m_p$-gerbe, as in Definition \ref{def:modulistack} in the positive rank case. Fix a multisection $\Sigma\subset X$ whose degree is equal to the index of the fibration $X\to\P^1$. In our case, $X$ is supersingular, so the index is equal to $1$ or $p$. Let $\widetilde{\Sigma}=\sX\times_X\Sigma$ be the restriction of $\sX$ to $\Sigma\subset X$. By Tsen's theorem, the Brauer class of the gerbe $\widetilde{\Sigma}$ is trivial. We fix also an invertible twisted sheaf $\sL$ on $\widetilde{\Sigma}$. We define the following relative stacks over $\P^1$.

\begin{Definition}\label{sec:modul-interpr-isom-10}
\noindent
\begin{enumerate}
\item Let $\ms R^{(\Sigma,\sL)}_{\ms X}(s)\to\P^1$ be the stack whose objects over
  $T\to\P^1$ are pairs $(\sE,\phi)$, where $\sE$ is a $T$-flat $\ms X\times_{\P^1}T$-twisted
  quasi-coherent sheaf of finite presentation such that for each geometric
  point $t\in T$, the pushforward of the fiber $\sE_t$ along the natural
  closed immersion
  \[
    \ms X\times_{\P^1}t\to\ms X\times t
  \]
  is $H$-stable with twisted Mukai vector $v=(0,E,s)$, and $\phi$ is an isomorphism
  \[
    \phi\colon \det(q_{\widetilde{\Sigma}_T*}(\sE|_{\widetilde{\Sigma}_T}\otimes\sL^\vee))\iso \det(q_{\Sigma_T*}\cO_{\Sigma_T})
  \]
   where $\widetilde{\Sigma}_T=\widetilde{\Sigma}\times_{\P^1}T$, $\Sigma_T=\Sigma\times_{\P^1}T$, and $q_{\widetilde{\Sigma}_T}\colon\widetilde{\Sigma}_T\to T$ and $q_{\Sigma_T}\colon\Sigma_T\to T$ are the projections.
 \item Let $\ms R_{\ms X}(s)\to\P^1$ be the same but omitting the isomorphisms $\phi$.
 \end{enumerate}
\end{Definition}

\begin{Remark}
 The reader will note that we define the stability condition in terms of the pushforward of the family to $\ms X$, rather than in the usual classical way, in terms of a relative polarization on $X$ over $\P^1$. This is done in order to avoid dealing with Hilbert polynomials on gerbes -- which are not purely cohomological in nature -- in the case of a singular variety (such as a singular fiber of the pencil). 
\end{Remark}

\begin{Example}\label{ex:trivial-example}
  Consider the case when $\ms X\to X$ is the trivial gerbe $X\times\B\m_p$. There is then an invertible
  $\ms X$-twisted sheaf $\sL'$ such that $\ms L'^{\tensor
    p}\cong\ms O_{\ms X}$. Tensoring with $\ms L'^\vee$ and pushing
  forward to $X$ defines an isomorphism between $\ms M_{\ms X}(v)$ and
  the stack of coherent pure $1$-dimensional sheaves on $X$ with
  determinant $\ms O(E)$ and second Chern class $s$. As shown in
  Section 4 of \cite{MR1629929}, this stack is isomorphic to the relative moduli stack of stable sheaves on the fibers of $f\colon X\to\P^1$ of rank $1$ and degree $s$.
\end{Example}

\begin{Lemma}
  \label{sec:modul-interpr-isom-8}
  Pushforward defines an isomorphism
  \[
    \varphi\colon \ms R_{\ms X}(s)\to\ms M_{\ms X}(v)\hspace{1cm}
  \]
  of $k$-stacks.
\end{Lemma}
\begin{proof}
  First we define the morphism. Fix a $k$-scheme $T$. A point of $\ms R_{\ms X}=\ms R_{\ms X}(s)$ is given by a lift $T\to\P^1$ and a $T$-point as in Definition \ref{sec:modul-interpr-isom-10}. But the stability condition is preserved under the pushforward
  \[
    \ms X\times_{\P^1}T\to\ms X\times T
  \]
  by definition of $\ms R_{\ms X}$. Hence, pushing forward along this
  morphism gives an object of $\ms M_{\ms X}=\ms M_{\ms X}(v)$, giving the morphism $\varphi$.

  To show that $\varphi$ is an isomorphism of stacks, we will show that it is a
  proper monomorphism (hence a closed immersion) that is surjective
  on $k$-points (hence an isomorphism, as $\ms M_{\ms X}$ is smooth). We first
  make the following claim.

  \begin{Claim}
    Given a morphism $a\colon T\to\P^1$ and a family $F$ in $\ms R_{\ms
      X}(T)$ with pushforward $\iota_\ast F$ on $\ms X\times T$, we
    can recover the graph of $a$ as the Stein factorization of the morphism
    $$\Supp(\iota_\ast F)\to \P^1\times T,$$
    where $\Supp(\iota_\ast F)$ denotes the scheme-theoretic support
    of $\iota_\ast F$.
  \end{Claim}
  \begin{proof}
    Since $X\to\P^1$ is cohomologically flat in dimension $0$, the
    claim follows if the natural map
    $$\ms O_{\ms X\times_{\P^1}T}\to\End(F)$$
    is injective (as it is automatically compatible with base change
    on $T$). By the assumption about the determinant of the fibers of
    $F$, for each geometric point $t\to T$ we know that $\ms O_{\ms
      X_t}\to\End(F_t)$ is injective (as $F_t$, supported on one fiber
    of $X\to\P^1$, must have full support for the determinant on $\ms
    X$ to be correct). The result now follows from Lemma 3.2.3 of
    \cite{MR2579390}.
  \end{proof}

  Suppose $$t_1,t_2\colon T\to\P^1$$
  are two morphisms and $F_i$ is an object of $\ms R_{\ms X}(t_i)$ for
  $i=1,2$. Write $$\iota_i\colon \ms X\times_{\P^1,t_i}T\to\ms X\times T$$
  for the two closed immersions. If
  $\varphi(F_1)\cong\varphi(F_2)$ then their scheme-theoretic supports
  agree, whence their Stein factorizations agree. By the Claim, the
  two maps $t_1$ and $t_2$ must be equal. But then $F_1$ and $F_2$
  must be isomorphic because their pushforwards are isomorphic.

  Now let us show that $\varphi$ is proper. Fix a complete dvr $R$ over
  $k$ with fraction field $K$, a morphism $\Spec K\to\P^1$, and an
  object $F_K\in\ms R_{\ms X}(K)$. Let $\Spec R\to\P^1$ be the unique
  extension ensured by the properness of $\P^1$ and let $$\iota\colon \ms
  X\tensor_{\P^1}R\to\ms X\tensor_k R$$ be the natural closed
  immersion. We wish to show that the unique stable
  limit $F$ of $\iota_\ast F_K$ has the form $\iota_\ast F$ for an
  $R$-flat family of coherent $\ms X_R$-twisted sheaves (as the
  stability condition then follows by definition).

  Let $\ms I$ be the ideal of the image of $\iota$. By assumption, the
  map of sheaves
  $$\nu\colon \ms O_{\ms X\tensor_k R}\to\send(F)$$
  kills $\ms I$ in the generic fiber over $R$. Since $F$ is $R$-flat,
  so is $\send(F)$ (as $R$ is a dvr). Thus, the image of $\ms I$ in
  $\send(F)$ is $R$-flat. But this image has trivial generic fiber,
  hence must be trivial. It follows that $\nu$ kills $\ms I$, whence
  $F$ has a natural structure of pushforward along $\iota$, as desired.
  
\end{proof}

\begin{Corollary}\label{cor:gerbesoversheaves}
    If $X\to\P^1$ has index $p$, then the stack
    \[
      \ms R^{(\Sigma,\sL)}_{\sX}(s)\to R_{\sX}(s)
    \]
    is a $\m_p$-gerbe over a supersingular K3 surface, with associated $\G_m$-gerbe $\ms R_{\sX}(s)\to R_{\sX}(s)$.
\end{Corollary}
\begin{proof}
    By Theorem \ref{thm:moduli0}, $\sM_{\sX}(v)\to M_{\sX}(v)$ is a $\G_m$-gerbe over a supersingular K3 surface. By Lemma \ref{sec:modul-interpr-isom-8}, $\ms R_{\ms X}(s)\to R_{\ms X}(s)$ is as well. There is a natural forgetful map
    \[
      \ms R^{(\Sigma,\sL)}_{\sX}(s)\to\ms R_{\ms X}(s)
    \]
    Consider a morphism $T\to\P^1$ and an object $(\sE,\phi)$ over $T$. The sheaf $\sE|_{\widetilde{\Sigma}_T}\otimes\sL^\vee$ on $\widetilde{\Sigma}_T$ has rank 1 and is 0-twisted. Because $\Sigma\to\P^1$ is finite flat of degree $p$, the pushforward
    \[
      q_{\widetilde{\Sigma}_T*}(\sE|_{\widetilde{\Sigma}_T}\otimes\sL^\vee)
    \]
    is of rank $p$. An automorphism $f\in\Aut(\sE)\iso\G_m$ acts on its determinant via the $p$-th power map. It follows that $\ms R^{(\Sigma,\sL)}_{\sX}(s)$ is a $\m_p$-gerbe over its sheafification, and that $\ms R_{\ms X}(s)$ is its associated $\G_m$-gerbe.
\end{proof}

Our interpretation of moduli spaces of rank 0 sheaves in terms of elliptic fibrations means that Artin-Tate twistor families come with additional geometric structure. Recall the following special case of a theorem of Artin and Tate (Theorem 3.1 of \cite{MR1610977}).

\begin{Theorem}[Artin-Tate]\label{sec:modul-interpr-isom}
  The edge map in the $E^2$ term of the Leray spectral sequence for $\G_m$ on $f\colon X\to\P^1$ yields an isomorphism
  \[
    \Br(X)\to\Sha(k(t),\Pic_{X_{k(t)}/k(t)})
  \]
  resulting in a natural surjection
  \[
    \Sha(k(t),\Jac(X_{k(t)}))\surj\Br(X)
  \]
  with kernel isomorphic to $\Z/i\Z$, where $i$ is the index of the generic fiber $X_{k(t)}$ over $k(t)$.
\end{Theorem}
In particular, if $X\to\P^1$ has a section, the latter arrow is an isomorphism.
\begin{proof}
This is Proposition 4.5 (and ``cas particulier (4.6)'') of \cite{MR0244271}.
\end{proof}

It follows by descent theory that any element of $\Sha(k(t),\Jac(X_{k(t)}))$ corresponds to an \'etale form $X'$ of $X$, and $X'$ is also a K3 surface. We will describe this isomorphism geometrically using the twistor families. This description will allow us to take a tautological family of $\m_p$-gerbes on $X$ and produce a family of K3 surfaces that are each forms of a given elliptic fibration on $X$.


The following are very well known for \'etale cohomology with coefficients of order prime to
the residue characteristic. They are also true for fppf cohomology with coefficients in $\m_p$, as we
record here.

\begin{Lemma}
  \label{sec:few-remarks-cohom-2}
  Suppose $C$ is a proper smooth curve over an algebraically closed
  field $k$. The Kummer sequence induces a canonical isomorphism
  \[
    \Pic(C)/p\Pic(C)=\Z/p\Z\simto\H^2(C,\m_p)
  \]
\end{Lemma}
\begin{proof}
  By Tsen's theorem $\H^2(C,\G_m)=0$. Since the $p$-th power Kummer
  sequence is exact on the fppf site, the result follows.
\end{proof}

\begin{Proposition}
  \label{P:curve-coho}
  Suppose $Z\to G$ is a proper smooth morphism of finite presentation of relative dimension $1$ 
  with $G$ connected and $\alpha\in\H^2(Z,\m_p)$. There exists a unique
  element $a\in\Z/p\Z$ such that for every
  geometric point $g\to G$, the restriction $\alpha_{Z_g}\in\Z/p\Z$ is
  equal to $a$ via the isomorphism of Lemma \ref{sec:few-remarks-cohom-2}.
\end{Proposition}
\begin{proof}
  It suffices to prove this under the assumption that $G$ is the
  spectrum of a complete dvr.
  \begin{Lemma}
    When $G$ is the spectrum of a complete dvr, we have that $\H^2(Z,\G_m)=0$.
  \end{Lemma}
  \begin{proof}
    First, since $Z$ is regular the group is torsion. Thus, any
    $\G_m$-gerbe is induced by a $\m_n$-gerbe for some $n$. Fix a
    $\m_n$-gerbe $\ms Z\to Z$. By Tsen's theorem and Lemma 3.1.1.8 of \cite{MR2388554}, there is
    an invertible $\ms Z_g$-twisted sheaf $L$, where $g$ is the closed
    point of $g$. The obstruction to deforming such a sheaf lies in
    $$\H^2(Z,\ms O)=0,$$
    so that $L$ has a formal deformation over the completion of $\ms
    Z$. By the Grothendieck Existence Theorem for proper Artin stacks,
    Theorem 11.1 of \cite{MR2312554} (or, in this case, the classical Grothendieck Existence Theorem for
    coherent modules over an Azumaya algebra representing $\alpha$),
    this formal deformation algebraizes, trivializing the class of
    $\ms Z$ in $\H^2(Z,\G_m)$, as desired.
  \end{proof}
  Applying the lemma and the Kummer sequence, we see that
  $$\Pic(Z)/p\Pic(Z)\simto\H^2(Z,\m_p).$$
  But restricting to a fiber defines a canonical isomorphism
  $$\Pic(Z)/p\Pic(Z)\simto\Pic(Z_g)/p\Pic(Z_g)\simto\Z/p\Z$$
  independent of the point $g$. The result follows.
\end{proof}

\begin{Corollary}
  \label{C:coho-class-indep}
  Suppose $E\subset Z\to T$ is a family of smooth genus $1$ fibers in
  a proper flat family of elliptic surfaces of finite presentation
  over a connected base. Given a class
  $$\alpha\in\H^2(Z,\m_p),$$
  there is an element $a\in\Z/p\Z$ such that for every geometric point
  $t\in T$, the restriction of $\alpha$ to $E_t$ equals $a$ via the
  isomorphism of Lemma \ref{sec:few-remarks-cohom-2}.
\end{Corollary}


\begin{Proposition}
  \label{sec:modul-interpr-isom-13}
  If $\ms X$ is a $\m_p$-gerbe that deforms the trivial gerbe, then for any $s$ the morphism $\ms R_{\ms X}(s)\to\P^1$ is an \'etale form of the morphism $\ms R_X(s)\to\P^1$.
\end{Proposition}

\begin{proof}[Proof of Proposition \ref{sec:modul-interpr-isom-13}]
   This turns out to be surprisingly subtle, and uses our assumption that $\ms X$ deforms the trivial gerbe in an essential way. We begin with a lemma.

  \begin{Lemma}
    \label{sec:modul-interpr-isom-15}
    Given a fiber $D\subset X$ of $f$, there is an invertible $\ms
    X\times_X D$-twisted sheaf $\Lambda$ of rank $1$ such that for
    each smooth curve $C\to D$ the restriction of $\Lambda$ to $C$ has
    degree $0$.
  \end{Lemma}
  \begin{proof}
    We may replace $D$ with its induced reduced structure, so we will
    assume that $D$ is a reduced curve supported on a fiber of
    $f$. Write $\pi_D\colon D\to\spec k$ for the structure morphism.
    Let $\sY\to\A^1$ be the universal family over the identity component of the group of $\m_p$-gerbes on $X$, so that $\sY_0$ is the trivial gerbe and $\sY_{t_0}\cong\sX$ for some $t_0\in\A^1$.
    Let $\ms D\to D\times\A^1$ be the restriction of $\sY$ under the map $D\times\A^1\to X\times\A^1$. The
    gerbe $\ms D$ gives rise to a morphism of fppf sheaves
    \[
      \A^1\to \R^2\pi_{D*}\m_p.
    \]
    The Kummer sequence shows that there is an isomorphism of sheaves
    \[
      \Pic_{D/k}/p\Pic_{D/k}\simto\R^2\pi_{D\ast}\m_p.
    \]
    Thus, the gerbe $\ms D$ gives rise to a morphism
    \[
      h\colon \A^1\to\Pic_{D/k}/p\Pic_{D/k}
    \]
    under which $0$ maps to $0$.

    On the other hand, there is a multidegree morphism of $k$-spaces
    \[
      \Pic_{D/k}\surj\prod_{i=1}^m\Z,
    \]
    where $m$ is the number of irreducible components of $D$. (This
    map comes from taking the degree of invertible sheaves pulled back
    to normalizations of components, and surjectivity is a basic
    consequence of the ``complete gluing'' techniques of \cite{MR1432058}.) This
    gives rise to a morphism
    \[
      \deg_p\colon \Pic_{D/k}/p\Pic_{D/k}\to\prod_{i=1}^m\Z/p\Z
    \]
    of sheaves. 
      Composing with $h$, it follows from the connectedness of $\A^1$ that $t_0\in\A^1$ must map into
      the kernel of $\deg_p$.

      By Tsen's theorem, there is an
      invertible $\ms D$-twisted sheaf $\sL$, and the above
      calculation shows that $\sL^{\tensor p}$ is the pullback of
      an invertible sheaf $L$ on $D$ such that for each irreducible
      component $D_i\subset D$, the pullback of $L$ to the
      normalization of $D_i$ has degree $p\delta_i$. Let
      \[
        \lambda_i\in\Pic(D_i)
      \]
      be an invertible sheaf whose pullback
      to the normalization has degree $-\delta_i$. A simple gluing argument
      shows that there is an invertible sheaf 
      \[
        \lambda\in\Pic(D)
      \]
      such that
      \[
        \lambda|_{D_i}\cong\lambda_i
      \]
      for each $i$. Replacing $\sL$ by $\sL\tensor\lambda$
      yields an invertible $\ms D$-twisted sheaf whose restriction to
      each $D_i$ has degree $0$, yielding the desired result (as any
      non-constant $C\to D$ factors through a $D_i$).
  \end{proof}

  To show that $\ms R_{\ms X}(s)$ is an \'etale form of $\ms R_{X}(s)$, we may
  base-change to the Henselization $U$ of $\P^1$ at a closed point. By
  Lemma \ref{sec:modul-interpr-isom-15}, there is an invertible
  twisted sheaf $\sL_u$ on the closed fiber 
  \[
    X_u\subset X_U
  \]
  whose restriction to each irreducible component has degree
  $0$. Since the obstruction to deforming such a sheaf lies in 
  \[
    \H^2(X_u,\ms O)=0,
  \]
  we know that $\sL_u$ deforms to an invertible $\ms X_U$-twisted sheaf $\sL$ whose restriction to
  any smooth curve in any fiber of $X_U$ over $U$ has degree $0$.

  Tensoring with $\sL$ gives an isomorphism of stacks
  \[
    \Sh_{\ms X_U/U}(0,E,s)\simto\Sh_{X_U/U}(0,E,s).
  \]
  We claim that this isomorphism preserves $H$-stability. Since
  stability is determined by Hilbert polynomials, it suffices to prove
  the following.

  \begin{Claim}
    \label{sec:modul-interpr-isom-14}
    For any geometric point $u\to U$ and any coherent $\ms X$-twisted
    sheaf $G$, the geometric Hilbert polynomial of $\iota_\ast G$
    equals the Hilbert polynomial of $\sL^{\vee}\tensor G$. In
    particular, $G$ is stable if and only if $\sL^{\vee}\tensor G$
    is stable.
  \end{Claim}
  \begin{proof}
    Since $G$ is filtered by subquotients supported on the reduced
    structure of a single irreducible component of $D$, it suffices to
    prove the result for such a sheaf. Let $\nu\colon C\to D$ be the
    normalization of an irreducible component. The sheaves
    $\nu_\ast\nu^\ast G$ and $G$ differ by a sheaf of finite
    length. Thus, it suffices to prove the result for twisted sheaves
    on $C$ and twisted sheaves of finite length. In either case, we
    are reduced to showing the following:
    given a finite morphism $q\colon S\to X_u$ from a smooth
    $\kappa(u)$-variety, let $\ms S\to S$ be the pullback of $\ms
    X_u\to X_u$. Then for any coherent $\ms S$-twisted sheaf $G$, the
    geometric Hilbert polynomial of $q_\ast G$ equals the Hilbert
    polynomial of $\sL^{\vee}\tensor q_\ast G$.

    Using the Riemann-Roch theorem for geometric Hilbert polynomials (Lemma \ref{lem:riemannrochEulerVersion}),
    the classical Riemann-Roch theorem, and the projection formula, we see
    that it is enough to prove that the geometric Hilbert polynomial
    of $G$ (with respect to the pullback of $H$ to $S$) equals the
    usual Hilbert polynomial of $\sL_{\ms S}^{\vee}\tensor G$,
    under the assumption that $\sL_{\ms S}$ has degree $0$.

    Let $L\in\Pic^0(S)$ be the sheaf whose pullback to $\ms S$
    isomorphic to $\sL^{\tensor p}$. Using the isomorphism between
    $K(\ms S)\tensor\Q$ and $K(X)\tensor\Q$, the geometric Hilbert
    polynomial of $G$ is identified with the usual Hilbert polynomial
    of the class
    \[
      (\sL^{\vee}\tensor G)\tensor\frac{1}{p}L^\vee.
    \]
    But $L\in\Pic^0(S)$, so this Hilbert polynomial is the same as the
    Hilbert polynomial of $\sL^{\vee}\tensor G$, as claimed.
  \end{proof}

  We conclude that $\ms R_{\ms X}(s)$ is an \'etale form of the moduli space $\ms R_X(s)$, as claimed.
\end{proof}

Recall that the subgroup $\U^2(X,\m_p)\subset \H^2(X,\m_p)$ is the $k$-points of the connected component of the identity of $\R^2\pi^{\fl}_*\m_p$, and therefore classifies those $\m_p$-gerbes on $X$ that deform the trivial gerbe.

\begin{Proposition}\label{prop:Br-iso-Sha}
  Let $X\to\P^1$ be an elliptic supersingular K3 surface with Jacobian $J(X)\to\P^1$. If $\ms X$ is a $\m_p$-gerbe that deforms the trivial gerbe, then for any $s$ the morphism $R_{\sX}(s)\to\P^1$ is an \'{e}tale form of the morphism $J(X)\to\P^1$. The resulting map $U^2(X,\m_p)\to\Sha(k(t),\Jac(X_{k(t)}))$ defined by
  \[
    \ms X\mapsto[R_{\ms X}(1)]
  \]
  fits into a commutative diagram
  \[
    \begin{tikzcd}
      \U^2(X,\m_p)\arrow{d}\arrow[two heads]{dr}&\\
      \Sha(k(t),\Jac(X_{k(t)}))\arrow[two heads]{r}&\Br(X)
    \end{tikzcd}
  \]
  where the horizontal arrow is the Artin-Tate map of Proposition \ref{sec:modul-interpr-isom}.
\end{Proposition}
\begin{proof}
  By Proposition \ref{sec:modul-interpr-isom-13}, $\ms R_{\sX}(s)\to\P^1$ is an \'{e}tale form of $\ms R_{X}(s)\to\P^1$. The same is therefore true for $R_{\sX}(s)\to\P^1$ and $R_X(s)\to\P^1$. The latter is isomorphic as an elliptic surface to the Jacobian $J^s(X)\to\P^1$, which is an \'{e}tale form of $J(X)\to\P^1$ (and also of $X\to\P^1$) (see Chapter 11, Remark 4.4 of \cite{Huy16}). This gives the first claim.

  For the second claim, we consider the corresponding genus 1 curves $X_{\eta}$ and $R_\eta=(R_{\ms X}(1))_{\eta}$ over the generic point $\eta=\Spec k(t)$ of $\P^1$. The Leray spectral sequence and Tsen's theorem show that the edge map gives an isomorphism
  \[
    \Br(X_\eta)\simto\H^1(\eta,\Jac(X_\eta))
  \]
  which we can describe concretely as follows. Over $\widebar{k(t)}$ the gerbe $\ms X_\eta\to X_\eta$ has trivial Brauer class, hence carries an invertible twisted sheaf $\sL$ such that $\sL^{\tensor p}$ has degree $pa(\ms X)$, where
  \[
    a(\ms X)\in\frac{1}{p}\Z/\Z
  \]
  is the unique element that corresponds to the cohomology class of the restriction of $\ms X$ to any smooth fiber $E\subset X$ of $f$ under the natural isomorphism
  \[
    \frac{1}{p}\Z/\Z\simto\Z/p\Z.
  \]
  By our assumption that $\ms X$ deforms the trivial gerbe, $a(\sX)=0$. Given an element $\sigma$ of the Galois group of $\widebar{k(t)}$ over $k(t)$, there is an invertible sheaf $L_\sigma\in\Pic(X_{\widebar{k(t)}})$ such that $\sigma^\ast\sL\tensor\sL^{\vee}\cong L_{\sigma}|_{\ms X_\eta}$. This defines a $1$-cocycle in the sheaf $\Pic_{X_\eta/\eta}$, and its cohomology class is the image of a unique class in $\H^1(\eta,\Jac(X_\eta))$, as desired.

  On the other hand, tensoring with $\sL^{\vee}$ gives an isomorphism between the stack of invertible $\ms X_\eta$-twisted sheaves of degree $1$ and the stack of invertible sheaves on $X_\eta$ of degree $1$. The latter stack is a gerbe over $X_\eta$, and the Galois group induces the cocycle given by the translation action of $\Jac(X_\eta)$ on $X_\eta$. But this gives the edge map in the Leray spectral sequence. 
  
  Finally, the surjectivity of the diagonal arrow is implied by Theorem \ref{thm:Artinsflatduality}.
\end{proof}


\section{Applications}\label{sec:applications}

\subsection{A crystalline Torelli theorem for twisted supersingular K3 surfaces}\label{sec:torellitheorem}

In this section we use twistor lines to prove a crystalline Torelli theorem for twisted supersingular K3 surfaces over the algebraically closed field $k$ of characteristic $p\geq 3$. Our approach is inspired by Verbitsky's proof of a Torelli theorem for hyperk\"{a}hler manifolds over the complex numbers using classical twistor space theory.\footnote{See \cite{MR3161308} for the proof, and \cite{MR3051203} for a skillful exposition of these ideas. Also see \cite{Huy16}, where the same strategy is applied in a simplified form to deduce the Torelli theorem for complex K3 surfaces.} 

As we will explain, this gives in particular an alternative proof of Ogus's crystalline Torelli theorem (Theorem III of \cite{Ogus83}). Let us briefly discuss the differences between our proof and that of Ogus. A key input in \cite{Ogus83} is the result of Rudakov, Zink, and Shafarevich that supersingular K3 surfaces in characteristic $p\geq 5$ do not degenerate (see Theorem 3 of \cite{RTS83} and Theorem \ref{thm:rudakovshafaerevich}). This is applied to show that the moduli space $S_{\Lambda}$ is almost proper over $k$ (meaning it satisfies the surjectivity part of the valuative criterion with respect to DVRs). Using that the period domain is proper, it is then deduced that the period morphism is almost proper, and eventually (after adding ample cones) an isomorphism. Because of this, the main result of \cite{Ogus83} is restricted to characteristic $p\geq 5$, although the surrounding theory works in characteristic $p\geq 3$. In our proof, we will show using twistor lines that the period morphism is an isomorphism without using the almost properness of the moduli space. Thus, we do not need to restrict to characteristic $p\geq 5$, and we are able to extend Ogus's Torelli theorem to characteristic 3 (characteristic 2 carries its own difficulties, which we will not discuss here). In fact, as we observe in Section \ref{sec:misc-applications}, we can reverse the flow of information in Ogus's proof to obtain in particular a proof of the non-degeneration result of \cite{RTS83} in characteristic $p\geq 3$. This result was previously known only for $p\geq 5$, with partial results for $p=2,3$ via different methods (see \cite{RS83}).

The period domain $M_{\Lambda_0}$ is a projective variety, and in particular separated. However, the moduli space $S_\Lambda$ is not separated. To correct this, we will follow Ogus \cite{Ogus83} and equip characteristic subspaces with ample cones. If $\Lambda$ is a supersingular K3 lattice, we set
\[
  \Delta_\Lambda=\left\{\delta\in \Lambda|\delta^2=-2\right\}
\]
Associated to an element $\delta\in\Delta_{\Lambda}$ is the reflection
\[
  s_{\delta}(w)=w+(\delta.w)\delta
\]
The \textit{Weyl group} of $\Lambda$ is the subgroup $W_{\Lambda}\subset O(\Lambda)$ generated by the $s_{\delta}$. We set
\[
  V_\Lambda=\left\{x\in \Lambda\otimes\bR|x^2>0\mbox{ and }x.\delta\neq 0\mbox{ for all }\delta\in\Delta_\Lambda\right\},
\]
and let $C_\Lambda$ be the set of connected components of $V_\Lambda$.
\begin{Proposition}[\cite{Ogus83}, Proposition 1.10]\label{prop:transitiveaction}
  The group $\pm W_{\Lambda}$ acts simply transitively on $C_{\Lambda}$.
\end{Proposition}
If $\Lambda=\Pic(X)$ for some K3 surface $X$, then it is shown in \cite{Ogus83} that the ample cone of $X$ corresponds to a connected component of $V_\Lambda$, and hence to a certain element of $C_\Lambda$. Let $K\in M_{\Lambda_0}(S)$ be a characteristic subspace. For each (possibly non-closed) point $s\in S$, we set
\[
  \Lambda(s)=\left\{x\in\Lambda\otimes\bQ|px\in\Lambda\mbox{ and }\overline{px}\in K(s)\right\}
\]
We caution that this notation is not completely compatible with that of \cite{Ogus83}.
\begin{Definition}
  If $K\in M_{\Lambda_0}(S)$ is a characteristic subspace, then an \textit{ample cone} for $K$ is a choice of elements $\alpha_s\in C_{\Lambda(s)}$ for each $s\in S$, such that if $s$ specializes to $s_0$ then $\alpha_{s_0}\subset\alpha_{s}$. We let $P_{\Lambda}$ denote the functor on schemes over $\oM_{\Lambda_0}$ whose value on a scheme is the set of ample cones of the corresponding characteristic subspace.
\end{Definition}

\begin{Proposition}[\cite{Ogus83}, Proposition 1.16]\label{prop:ampleconespace}
  The functor $P_{\Lambda}$ is representable and locally of finite type over $k$, and the natural map $P_{\Lambda}\to \oM_{\Lambda_0}$ is \'{e}tale and surjective.
\end{Proposition}

Ogus also shows that $\rho\colon S_\Lambda\to\oM_{\Lambda_0}$ factors through $P_{\Lambda}$, resulting in a map
\[
  \rho^a\colon S_\Lambda\to P_{\Lambda}
\]
This map is shown to be \'{e}tale and separated. Ogus's crystalline Torelli theorem asserts that this map is an isomorphism. We seek an extension of these ideas to the twisted setting. We define a functor $\sP_{\tLam}$ by the Cartesian diagram
\[
  \begin{tikzcd}
    \sP_{\tLam}\arrow{r}{\pi_{e}^a}\arrow{d}&P_{\Lambda}\arrow{d}\\
    \oM^{\langle e\rangle}_{\tLam_0}\arrow{r}{\pi_e}&\oM_{\Lambda_0}
  \end{tikzcd}
\]
It is immediate from Proposition \ref{prop:ampleconespace} that the functor $\sP_{\tLam}$ is representable and locally of finite type over $k$, and that the natural map $\sP_{\tLam}\to\oM^{\langle e\rangle}_{\tLam_0}$ is \'{e}tale and surjective.

\begin{Proposition}\label{prop:thesquareistotallyCartesian}
  The morphism $\trho\colon \sS_{\Lambda}^o\to \oM^{\langle e\rangle}_{\tLam_0}$ factors through $\sP_{\tLam}$. The resulting diagram
  \begin{equation}\label{eq:periodsquare}
    \begin{tikzcd}
        \sS_{\Lambda}^o\arrow{r}{p}\arrow{d}[swap]{\trho^a}&S_\Lambda\arrow{d}{\rho^a}\\
        \sP_{\tLam}\arrow{r}{\pi_{e}^a}&P_{\Lambda}
      \end{tikzcd}
  \end{equation}
  is Cartesian, and the map $\trho^a\colon \sS_{\Lambda}^o\to \sP_{\tLam}$ is separated and \'{e}tale.
\end{Proposition}
\begin{proof}
  By the universal property of the fiber product and Proposition \ref{prop:commutativetriangle}, we get a diagram
    \[
      \begin{tikzcd}
        \sS_{\Lambda}^o\arrow{r}{p}\arrow{d}{\trho^a}\arrow[bend right=35]{dd}[swap]{\trho}&S_\Lambda\arrow{d}[swap]{\rho^a}\arrow[bend left=35]{dd}{\rho}\\
        \sP_{\tLam}\arrow{r}{\pi_{e}^a}\arrow{d}&P_{\Lambda}\arrow{d}\\
        \oM^{\langle e\rangle}_{\tLam_0}\arrow{r}{\pi_{e}}&\oM_{\Lambda_0}
      \end{tikzcd}
    \]
  where the lower and outer squares are Cartesian. It follows that the upper square is Cartesian. Because $\rho^a$ is separated and \'{e}tale, the same is true for $\trho^a$.
\end{proof}

Our crystalline Torelli theorem for twisted supersingular K3 surfaces asserts that the twisted period morphism $\trho^a$ is an isomorphism. We prove this in several steps.

Let us define the Artin invariant of a $k$-point of $\sP_{\tLam}$ (respectively, $P_{\Lambda}$) to be the Artin invariant of its image in $\oM^{\langle e\rangle}_{\tLam_0}$ (respectively $\oM_{\Lambda_0}$). Using Kummer surfaces, Ogus showed in \cite{Ogus78} that the crystalline Torelli theorem is true at all points of Artin invariant $\sigma_0\leq 2$ (this part of his argument needs only that $p\geq 3$). Using Proposition \ref{prop:thesquareistotallyCartesian}, we will extend this to the twisted setting.

\begin{Proposition}\label{prop:torelliatonepoint}
  The fiber of $\trho^a$ over any $k$-point of Artin invariant $\sigma_0\leq 2$ is a singleton. Each connected component of $\sP_{\tLam}$ contains a point of Artin invariant $\leq 2$.
\end{Proposition}
\begin{proof}
  It is shown in \cite{Ogus83} that the fiber of $\rho^a$ over any $k$-point of Artin invariant $\sigma_0\leq 2$ is a singleton (see Step 4 of the proof of Theorem III'). Consider a point $x\in \sP_{\tLam}$ of Artin invariant $\leq 2$. The image of $x$ in $P_{\Lambda}$ also has Artin invariant $\leq 2$, so by Proposition \ref{prop:thesquareistotallyCartesian} the fiber $(\trho^{a})^{-1}(x)$ is a singleton. 
  
  Following Ogus \cite{Ogus83}, let us say that a morphism of schemes is \textit{almost proper} if it satisfies the surjectivity part of the valuative criterion with respect to DVR's. In the proof of Proposition 1.16 of \cite{Ogus83}, Ogus shows that $P_{\Lambda}\to\oM_{\Lambda_0}$ is almost proper. This property is clearly preserved under base change, so the morphism
  \[
    \sP_{\tLam}\to\oM^{\langle e\rangle}_{\tLam_0}
  \]
  is also almost proper. Let $\sP\subset\sP_{\tLam}$ be a connected component. By Proposition \ref{prop:ampleconespace}, its image in $\oM^{\langle e\rangle}_{\tLam_0}$ is open. The scheme $\oM^{\langle e\rangle}_{\tLam_0}$ has two connected components, each of which are irreducible and contain a (unique) point of Artin invariant $1$ (see the discussion following Definition \ref{def:charsub}). We may therefore find a DVR $R$ with residue field $k$ and fraction field $K$ and a diagram
  \[
    \begin{tikzcd}
      \Spec K\arrow{r}\arrow{d}&\sP\arrow{d}\\
      \Spec R\arrow{r}\arrow[dashed]{ur}&\oM^{\langle e\rangle}_{\tLam_0}
    \end{tikzcd}
  \]
  such that the closed point of $\Spec R$ is mapped to a point of Artin invariant $\leq 2$. As the right vertical arrow is almost proper, we may find the dashed arrow, and this gives the result.
\end{proof}

We next show that $\trho^a$ is surjective.

\begin{Proposition}\label{prop:surjective}
  The morphism $\trho^a\colon \sS_{\Lambda}^o\to \sP_{\tLam}$ is surjective.
\end{Proposition}
\begin{proof}
  Clearly, $\trho^a$ is surjective if and only if $\rho^a$ is surjective. By Proposition \ref{prop:transitiveaction}, the group $\pm W_\Lambda$ acts transitively on the set of ample cones of $\Lambda$. Moreover, this group preserves any characteristic subspace. Therefore, $\rho^a$ is surjective if and only if $\rho$ is surjective (see \cite{Ogus83}, proof of Theorem III' step 4 for a similar argument). Because $\rho$ is surjective if and only if $\trho$ is surjective, we are reduced to showing that
  \[
    \trho\colon \sS_{\Lambda}^o\to\oM^{\langle e\rangle}_{\tLam_0}
  \]
  is surjective. We will induct on the Artin invariant. For each $1\leq\sigma_0\leq 10$ fix a supersingular K3 lattice $\Lambda^{\sigma_0}$, and define $\tLam^{\sigma_0+1}=\Lambda^{\sigma_0}\oplus U_2(p)$ (this is an extended supersingular K3 lattice of Artin invariant $\sigma_0+1$). Let
  \[
    \rho_{\sigma_0}\colon S_{\Lambda^{\sigma_0}}\to \oM_{\Lambda^{\sigma_0}_0}\hspace{1cm}\mbox{and}\hspace{1cm}\trho_{\sigma_0+1}\colon \sS_{\Lambda^{\sigma_0}}^o\to \oM^{\langle e\rangle}_{\tLam^{\sigma_0+1}_0}
  \]
  be the period morphisms. As described in Step 4 of the proof of Theorem III$'$ of \cite{Ogus83}, $\rho_1$ is surjective. We will show that
  \begin{enumerate}
    \item  if $\rho_{\sigma_0}$ is surjective, then $\trho_{\sigma_0+1}$ is surjective, and
    \item  if $\trho_{\sigma_0}$ is surjective, then $\rho_{\sigma_0}$ is surjective.
  \end{enumerate}
  Claim (1) follows immediately from Proposition \ref{prop:commutativetriangle}. To show (2), suppose that $\trho_{\sigma_0}$ is surjective, and take a point $p\in \oM_{\Lambda^{\sigma_0}}(k)$ corresponding to a characteristic subspace $K\subset\Lambda^{\sigma_0}_0\otimes k$. The orthogonal sum $\Lambda^{\sigma_0}\oplus\, U_2$ is an extended supersingular K3 lattice of Artin invariant $\sigma_0$. Hence, it is isometric to $\tLam^{\sigma_0}$, and the vector spaces $\Lambda^{\sigma_0}_0$ and $\tLam^{\sigma_0}_0$ are isometric. Choose an isotropic vector $w\in\Lambda^{\sigma_0}_0$ such that $w\notin K$. By Witt's Lemma (Theorem \ref{thm:Witt}) there exists an isometry $g_0\colon\Lambda^{\sigma_0}_0\iso\tLam^{\sigma_0}_0$ taking $w$ to $e$. By a result of Nikulin (see Theorem 14.2.4 of \cite{Huy16}), this map is induced by an isometry
  \[
    g\colon \Lambda^{\sigma_0}\oplus U_2\iso\tLam^{\sigma_0}
  \]
  Consider the corresponding isomorphism
  \[
    g_0^*\colon \oM_{\tLam^{\sigma_0}_0}\iso \oM_{\Lambda^{\sigma_0}_0}
  \]
  Suppose $g_0^*(q)=p$. As $w\notin K$, it follows that $q\in \oM^{\langle e\rangle}_{\tLam^{\sigma_0}_0}$. We are assuming that $\trho_{\sigma_0}$ is surjective, so we may find a point $x\in\sS^o_{\Lambda^{\sigma_0-1}}(k)$, corresponding to a marked twisted supersingular K3 surface $\sX\to X$, such that $\trho_{\sigma_0}(x)=q$.
  Let $a,b\in U_2$ be the standard basis satisfying $a^2=b^2=0$ and $a.b=-1$. Set $v=g(a)$. The image of $v$ under the marking $\tLam^{\sigma_0}\to\tN(\sX)$ remains primitive and isotropic. After possibly replacing $v$ with $-v$, we may arrange so that the image of $v$ in $\tN(\sX)$ satisfies one of the conditions of Proposition \ref{prop:nonempty}. By Theorem \ref{thm:moduli0}, the moduli space $\sM_{\sX}(v)$ is a trivial $\G_m$-gerbe over the supersingular K3 surface $M_{\sX}(v)$. Consider a quasi-universal sheaf
  \[
    \sP\in\Coh^{(0,1)}(M_{\sX}(v)\times\sX)
  \]
  As in Proposition \ref{prop:integralcorrespondence}, the induced cohomological transform induces an isomorphism
  \[
    \Phi^{\cris}_{v_{M\times\sX}(\sP)}\colon\tH(M_{\sX}(v)/W)\iso\tH(\sX/W)
  \]
  of K3 crystals, which sends $\tN(M_{\sX}(v))$ to $\tN(\sX)$. (This does not follow directly from Proposition \ref{prop:integralcorrespondence} as stated. However, the proof of this result is very similar to the proof of Proposition \ref{prop:integralcorrespondence}.) The restriction $h\colon\tN(M_{\sX}(v))\to\tN(\sX)$ of $\Phi^{\cris}_{v_{M\times\sX}(\sP)}$ sends $(0,0,1)$ to $v$. The composition
  \[
    \Lambda^{\sigma_0}\oplus U_2\too{g}\tLam^{\sigma_0}\to\tN(\sX)\too{h^{-1}}\tN(M_{\sX}(v))
  \]
  sends $a$ to $(0,0,1)$, and hence induces a map
  \[
    \Lambda^{\sigma_0}=a^\perp/a\to (0,0,1)^{\perp}/(0,0,1)=N(M_{\sX}(v))
  \]
  The K3 surface $M_{\sX}(v)$ with this marking gives the desired point in $S_{\Lambda^{\sigma_0}}$ whose image under $\rho$ is $p$.
\end{proof}


We now know that $\trho^a\colon \sS_{\Lambda}^o\to\sP_{\tLam}$ is separated, \'{e}tale, and surjective, and that there exists a point in each connected component of $\sP_{\tLam}$ whose fiber is a singleton. Combining this with our construction of twistor families in Theorem \ref{thm:constructionoftwistorlines}, we will show that $\trho^a$ is an isomorphism.

\begin{Theorem}\label{thm:twistedtorelli}
  If $p\geq3$, then the twisted period morphism $\trho^a\colon \sS_{\Lambda}^o\to \sP_{\tLam}$ is an isomorphism.
\end{Theorem}
\begin{proof}
  Let $\sP$ be a connected component of $\sP_{\tLam}$, $\sS$ its preimage in $\sS_{\Lambda}^o$, and $\trho'\colon \sS\to\sP$ the restriction of $\trho^a$. Consider the diagonal 
  \[
    \Delta_{\trho'}\colon\sS\to \sS\times_{\sP}\sS
  \]
  Because $\trho'$ is separated, $\Delta_{\trho'}$ is a closed immersion, and because $\trho'$ is \'{e}tale, it is also an open immersion. Thus, $\Delta_{\trho'}$ is an isomorphism onto its image $\Delta=\Delta_{\trho'}(\sP)$, which is open and closed. We will show that $\Delta=\sS\times_{\sP}\sS$. Take a $k$-point $(x_0,x_1)\in \sS\times_{\sP}\sS$. We will construct a connected subvariety $C\subset \sS\times_{\sP}\sS$ that contains $(x_0,x_1)$ and a $k$-point whose image in $\sP$ has Artin invariant $\leq 2$.
  
  Let $\sigma$ be the Artin invariant of $x_0$ and $x_1$. Suppose that $\sigma>2$. By induction, it will suffice to construct a connected subvariety containing $(x_0,x_1)$ and a point of Artin invariant $\sigma-1$.
  Let $\tK$ be the characteristic subspace corresponding to $\trho(x_0)=\trho(x_1)$. By Lemma \ref{lem:technical}, we may find an isotropic vector $v$ such that $v.e\neq0$, $v\notin\tK$, and $v\in (\tK\cap\tLam_0)^{\perp}$. That is, $v$ satisfies the assumptions of Theorem \ref{thm:constructionoftwistorlines} (see Lemma \ref{lem:genericArtin}). Applying Theorem \ref{thm:constructionoftwistorlines} to both $x_0$ and $x_1$, we find an open subscheme $U\subset L$ and lifts 
  \[
    \varphi_0,\varphi_1\colon U\to\sS
  \]
  of $f$, which agree upon composing with $\trho$. Consider the diagram
  \[
    \begin{tikzcd}
      U\arrow{dr}{\varphi_i}\arrow[bend left=25]{drr}{\varphi'_i}\arrow[bend right=25]{ddr}&&\\
      &\sS_{\Lambda}^o\arrow{r}{p}\arrow{d}[swap]{\trho^a}&S_{\Lambda}\arrow{d}{\rho}\\
      &\sP_{\tLam}\arrow{r}&P_{\Lambda}
    \end{tikzcd}
  \]
  where $i=0,1$. We will show that by performing an appropriate sequence of elementary modifications, we can modify the $\varphi_i$ so that they agree upon composing with $\trho^a$. By definition of $\sP_{\tLam}$, this is the same as $\rho\circ\varphi_0'=\rho\circ\varphi_1'$.
  
  By Proposition \ref{prop:transitiveaction}, the Weyl group of $\Lambda$ acts transitively on the set of ample cones $C_{\Lambda}$. By precomposing the marking of one of our families by an appropriate reflection, we may therefore ensure that $\rho\circ\varphi_0'=\rho\circ\varphi_1'$ at the generic point $\eta$ of $U$. The same is then true also on the open locus of points in $U$ with Artin invariant $\sigma$. Suppose that $s\in U$ is one of the finitely many closed points with Artin invariant $\sigma-1$. There is an inclusion $\Lambda\subset\Lambda(s)$, which induces an inclusion $V_{\Lambda}\supset V_{\Lambda(s)}$. The morphisms $\rho\circ\varphi_i'$ applied to $\eta$ give an ample cone $\alpha\subset V_{\Lambda}$, and applied to $s$ give two possibly different ample cones $\alpha_0,\alpha_1\subset V_{\Lambda(s)}$, which have the property that $\alpha_0,\alpha_1\subset\alpha$. By Proposition \ref{prop:transitiveaction}, there exists a sequence of elements $\delta_i\in\Lambda(s)$ such that the composite of the associated reflection sends $\alpha_0$ to $\alpha_1$. We may assume that the $\delta_i$ are not in the image of the Picard group of the generic fiber. By Remark 8.2.4 of \cite{Huy16}, the Weyl group is generated by irreducible effective curves, so we may furthermore assume that each $\delta_i$ is the class of an irreducible effective curve. By Proposition 2.8 of \cite{Ogus83}, taking the elementary modification of the family $X\to U$ corresponding to $\varphi_1'$ with respect to the $\delta_i$ produces a new family $Y\to U$, along with isomorphisms
  \[
    \Theta_{U-s}\colon Y_{U-s}\iso X_{U-s}\hspace{.5cm}\mbox{ and }\hspace{.5cm}\Theta_{s}\colon Y_{s}\iso X_{s}
  \]
  Moreover, these isomorphisms act on the Picard groups in such a way so that $Y\to U$ inherits a marking by $\Lambda$, and the induced morphism $U\to S_{\Lambda}\to P_{\Lambda}$ now agrees with $\varphi_0$ at $s$. Replacing $\varphi_1$ with the induced $U\to\sS_{\Lambda}^o$, we have modified our lift $\varphi_1$ so that $\rho\circ\varphi_0'=\rho\circ\varphi_1'$ at every point of Artin invariant $\sigma$, and at $s$. Furthermore, note that $x_1$ is still in the image of $\varphi_1$. Repeating this procedure at each of the finitely many points of $U$ with Artin invariant $\sigma-1$, we find the desired lifts $\varphi_0,\varphi_1$ satisfying $\trho^a\circ\varphi_0=\trho^a\circ\varphi_1$.
  
  We have constructed a morphism
  \[
    (\varphi_0,\varphi_1)\colon U\to\sS\times_{\sP}\sS
  \]
  whose image contains $(x_0,x_1)$. By Lemma \ref{lem:incidence1}, $U$ is equal to $L\cong\bA^1$ minus at most one point. Therefore, by Lemma \ref{lem:incidence2}, the image of $U$ contains a point, say $(x'_0,x'_1)$, where $x'_0$ and $x'_1$ have Artin invariant $\sigma-1$. Continuing in this manner, we find a chain of (open subsets of) twistor lines $C$ and a map $C\to\sS\times_{\sP}\sS$ whose image contains $(x_0,x_1)$, such that the image of the composition $C\to\sS\times_{\sP}\sS\to\sP$ contains a point of Artin invariant $\leq 2$. 
  
  By Proposition \ref{prop:torelliatonepoint}, the preimage of this point under $\trho'$ is a singleton. It follows that $C$ intersects $\Delta$, and hence $C$ is contained in $\Delta$. As $(x_0,x_1)$ was arbitrary, this shows that $\Delta=\sS\times_{\sP}\sS$. The diagonal morphism $\Delta_{\trho'}$ is therefore an isomorphism, so $\trho'$ is a monomorphism. Because $\trho'$ is \'{e}tale, it is an open immersion. By Proposition \ref{prop:surjective}, $\trho'$ is surjective. It follows that $\trho'$ is an isomorphism.
\end{proof}

Note that the pullback of the section $\sigma_f\colon \oM_{\Lambda_0}\subset\oM^{\langle e\rangle}_{\tLam_0}$ along $\rho$ gives the locus $S_{\Lambda}\subset\sS_{\Lambda}^o$ of marked surfaces with trivial $\m_p$-gerbe. Thus, Theorem \ref{thm:twistedtorelli} yields an alternative proof of Ogus's crystalline Torelli theorem (Theorem III$'$ of \cite{Ogus83}, as well as its consequences Theorems I,II,II$'$,II$''$, and III), and extends these results to characteristic $p\geq 3$.
\begin{Corollary}\label{cor:ogust}
    If $p\geq 3$, then the period morphism $\rho^a\colon S_{\Lambda}\to P_{\Lambda}$ is an isomorphism.
\end{Corollary}

Of course, because the diagram (\ref{eq:periodsquare}) is Cartesian, the opposite implication holds as well. Let us translate our result from characteristic subspaces to K3 crystals. Under the inclusion $\tH(\sX/W)\subset\tH(X/K)$, the codimension filtration on $\tH(X/K)$ induces a filtration on $\tH(\sX/W)$.
Theorem \ref{thm:twistedtorelli} implies the following pointwise statement (compare to Theorem II of \cite{Ogus83}).

\begin{Theorem}
  Let $k$ be an algebraically closed field of characteristic $p\geq 3$ and $\sX$ and $\sY$ be $\m_p$-gerbes on supersingular K3 surfaces over $k$. If $\Theta\colon \tH(\sX/W)\to\tH(\sY/W)$ is an isomorphism of $W$-modules such that
  \begin{enumerate}
    \item $\Theta$ is compatible with the bilinear forms and commutes with the respective Frobenius operators,
    \item $\Theta$ preserves the extended N\'{e}ron-Severi groups, so that there is a commutative diagram
      \[
        \begin{tikzcd}
          \tN(\sX)\arrow{r}\arrow[hook]{d}&\tN(\sY)\arrow[hook]{d}\\
          \tH(\sX/W)\arrow{r}{\Theta}&\tH(\sY/W)
        \end{tikzcd}
      \]
    \item $\Theta$ preserves the codimension filtration, and
    \item $\Theta$ maps an ample class to an ample class,
  \end{enumerate}
  then $\Theta$ is induced by a unique isomorphism $\sX_{\bG_m}\iso\sY_{\bG_m}$ of associated $\bG_m$-gerbes. Conversely, any such isomorphism induces a map $\Theta$ satisfying these conditions.
\end{Theorem}

Using the same methods as in \cite{Ogus78}, one can deduce the following (compare to Theorem I of \cite{Ogus83}).
\begin{Theorem}
  If $p\geq3$ and $\sX$ and $\sY$ are $\m_p$-gerbes on supersingular K3 surfaces, then $\sX$ and $\sY$ are isomorphic if and only if there exists a filtered isomorphism
  \[
    \tH(\sX/W)\iso\tH(\sY/W)
  \]
  of $W$-modules that is compatible with the bilinear forms and with the Frobenius operators.
\end{Theorem}

\subsection{Some consequences for moduli of supersingular K3 surfaces}
\label{sec:misc-applications}

In this section we record some further consequences of our methods for
families of supersingular K3 surfaces. As discussed in the
introduction to Section \ref{sec:torellitheorem}, we obtain an
alternate proof of the following theorem of Rudakov, Zink, and
Shafarevich (Theorem 3 of \cite{RTS83}). This result is new in
characteristic 3.

\begin{Theorem}\label{thm:rudakovshafaerevich}
  Let $k$ be an algebraically closed field of characteristic $p\geq 3$,
  and let $R=k[[t]]$ and $K=k((t))$. If $X$ is a K3 surface over $K$
  whose geometric fiber has Picard rank 22, then there exists a finite
  extension $R'/R$ and a smooth surface $X'$ over $R'$ such that
  $X'_{K'}\cong X_{K'}$.
\end{Theorem}
\begin{proof}
  After taking a finite extension $R'$ of $R$, we may
  arrange so that the family $X\to\Spec K$ admits a marking by some
  supersingular K3 lattice $\Lambda$ (see page 1522 of \cite{RS79}),
  and hence gives a morphism $K'\to S_{\Lambda}$. As explained in
  Proposition 1.16 of \cite{Ogus83}, the forgetful map
  $P_{\Lambda_0}\to\oM_{\Lambda_0}$ satisfies the existence part of
  the valuative criterion with respect to DVR's. Because the period
  domain $\oM_{\Lambda_0}$ is proper over $\Spec k$, we conclude by
  Theorem \ref{thm:twistedtorelli} that $S_{\Lambda}\to\Spec k$
  satisfies the existence part of the valuative criterion with respect
  to DVR's. This gives the result.
\end{proof}

By Theorem \ref{thm:twistedtorelli}, and our identification of twistor
lines in the period domain, we see that the moduli of supersingular K3
surfaces of Artin invariant $\leq\sigma_0+1$ is rationally fibered
over the moduli of supersingular K3 surfaces of Artin invariant
$\leq\sigma_0$, and the fibers of this map are twistor families. We
record a few observations along these lines.

Let $1\leq\sigma_0\leq 8$ be an integer, fix a supersingular K3
lattice $\Lambda^{\sigma_0+1}$ of Artin invariant $\sigma_0+1$ and set
$\tLam^{\sigma_0+1}=\Lambda^{\sigma_0+1}\oplus U_2$. Pick a primitive
isotropic vector $v=(p,l,s)\in p\tLam^{\sigma_0+1}$, let
$\Lambda^{\sigma_0}=v^\perp/v$, and set
$\tLam^{\sigma_0+1}=\Lambda^{\sigma_0}\oplus U_2(p)$. To a
$\Lambda^{\sigma_0+1}$-marked supersingular K3 surface $X$ we may
associate the moduli space $\sM_X(v)$, whose extended N\'{e}ron-Severi
group has an induced marking by $\tLam^{\sigma_0+1}$. This is
reflected at the level of period domains by a diagram
\begin{equation}\label{eq:diagram22}
  \begin{tikzcd}[column sep=tiny]
    \oM_{\Lambda^{\sigma_0+1}_0}^{\langle v\rangle}\arrow{rr}{\sim}\arrow{dr}[swap]{\pi_{v}}&& \oM_{\tLam^{\sigma_0+1}_0}^{\langle e\rangle}\arrow{dl}{\pi_{e}}\\
    &\oM_{\Lambda^{\sigma_0}_0}&
    \end{tikzcd}
\end{equation}
This shows that the crystalline period domain for marked supersingular
K3 surfaces of Artin invariant $\leq\sigma_0+1$ is covered by open
subsets that are isomorphic to the sheaf of groups
$(\R^2\pi^{\fl}_*\m_p)^o$ over the crystalline period domain for
supersingular K3 surfaces of Artin invariant $\leq\sigma_0$. With a
little extra work, we may obtain a statement at the level of moduli
spaces.
\begin{Definition}
  If $W\subset V$ is a subspace, we let $U_V^W\subset \oM_V$ be the
  locus of characteristic subspaces $K\subset V\otimes k$ such that
  $K\cap V\subset W$. In particular, $U_V^0=U_V$ is the open subset of
  strictly characteristic subspaces.
\end{Definition}
\begin{Lemma}\label{lem:rationalsections}
  For any totally isotropic subspace $W\subset\Lambda_0$, the morphism
  $P_{\Lambda}\to \oM_{\Lambda_0}$ admits a section over the open
  subset $U_{\Lambda_0}^W\subset \oM_{\Lambda_0}$.
\end{Lemma}
\begin{proof}
  A section over an open subset $U\subset \oM_{\Lambda_0}$ corresponds
  to a choice of ample cone for the restriction of the universal
  characteristic subspace $K_{\Lambda_0}$ to $U$. Let
  $W\subset\Lambda_0$ be a totally isotropic subspace. Define
  \[
    \Lambda_W=\left\{x\in\Lambda\otimes\bQ|px\in\Lambda\mbox{ and
      }\overline{px}\in W\right\}
  \]
  This is a supersingular K3 lattice, and we have
  $\Lambda\subset\Lambda_W$. For any (possibly non-closed) point
  $s\in U_{\Lambda_0}^W$, we have that $K(s)\cap \Lambda_0\subset W$,
  so there are inclusions
  \[
    \Lambda\subset\Lambda(s)\subset\Lambda_W
  \]
  It follows that
  \[
    V_{\Lambda_W}\subset V_{\Lambda(s)}\subset V_{\Lambda}
  \]
  Pick an ample cone $\alpha\in C_{\Lambda_W}$. For each $s$, let
  $\alpha_s\in C_{\Lambda(s)}$ be the unique connected component
  containing the image of $\alpha$ under the above inclusion. These
  ample cones are compatible, and hence we get an ample cone for the
  restriction of $K_{\Lambda_0}$ to $U_{\Lambda_0}^W$. In other words,
  we obtain for each $\alpha\in C_{\Lambda_W}$ a diagram
  \[
    \begin{tikzcd}
      &P_{\Lambda}\arrow{d}\\
      U_{\Lambda_0}^W\arrow{ur}{[\alpha]}\arrow[hook]{r}&\oM_{\Lambda_0}
    \end{tikzcd}
  \]
  which in particular implies the result.
\end{proof}

Combining this with our observation about the period domains, we
obtain the following statement.

\begin{Proposition}
  The moduli spaces $S_{\Lambda^{\sigma_0+1}}$ and
  $\sS_{\Lambda^{\sigma_0}}$ are birational.
\end{Proposition}
\begin{proof}
  By Lemma \ref{lem:rationalsections}, the morphisms
  $\sP_{\tLam^{\sigma_0+1}}\to \oM^{\langle
    e\rangle}_{\tLam^{\sigma_0+1}_0}$ and
  $P_{\Lambda^{\sigma_0+1}}\to \oM_{\Lambda^{\sigma_0+1}}$ admit
  sections defined on open subsets of the target. A section of an
  \'{e}tale morphism is an open immersion, so the isomorphism
  $\oM_{\Lambda^{\sigma_0+1}_0}\iso \oM_{\tLam^{\sigma_0+1}_0}$
  (\ref{eq:diagram22}) induces the desired birational correspondence.
\end{proof}

This iterated bundle structure has the following pointwise consequence.

\begin{Proposition}
  Every supersingular K3 surface can be obtained from the unique
  supersingular K3 surface $X$ of Artin invariant $\sigma_0=1$ by
  iteratively taking moduli spaces of twisted sheaves. That is, for
  any supersingular K3 surface $Y$, there exists a sequence
  $v_1,\dots,v_n$ of Mukai vectors and $\alpha_1,\dots,\alpha_n$ of
  cohomology classes such that $v_i\in\tN(X_i)$ and
  $\alpha_i\in \H^2(X_i,\m_p)$, where $X_1=X$ and
  \[
    X_{i+1}=M_{(X_i,\alpha_i)}(v_i)
  \]
  such that $X_n=Y$ (in fact, we may take $n=\sigma_0(Y)-1$).
\end{Proposition}

We also obtain consequences for other moduli spaces of supersingular
K3 surfaces. Given a positive integer $d$, let $M_{2d}$ denote the
moduli space of polarized K3 surfaces of degree $2d$, and let
$S_{2d}\subset M_{2d}$ denote the supersingular locus.

\begin{Theorem}
  For any positive integer $d$, every irreducible component of
  $S_{2d}$ is unirational and every general point has tangent space
  spanned by the tangent spaces of twistor lines.
\end{Theorem}
\begin{proof}
  Fix a component $S$ of $S_{2d}$.  Let $i$ be the maximal Artin invariant that
  appears in $S$, so that there is a dense open subscheme $S'\subset S$
  parametrizing polarized supersingular K3 surfaces of degree $2d$ with Artin
  invariant $i$. By Lemma \ref{lem:rationalsections} and Corollary
  \ref{cor:ogust}, we have there is an open subscheme
  $U$ of a crystalline period domain $\overline
  M_{\Lambda_0}$ and a dominant unramified morphism $U\to
  S'$. By Proposition \ref{prop:periodspaceisunirational},
  $U$ is unirational, and by Lemma \ref{BigTangentSpace}, the tangent space 
  to any point of $U$ is generated by twistor lines. This gives
  the desired conclusion.
\end{proof}

\begin{Question}
  We are not sure if the supersingular locus of $M_{2d}$ is irreducible.
\end{Question}

Finally, we will show that every twistor family is a moduli space of
twisted sheaves on a universal twistor family. In particular, this
will remove the assumptions on $v$ in Theorem
\ref{thm:constructionoftwistorlines}. To treat the Artin-Tate case, we
will use the following lemma.
\begin{Lemma}\label{lem:lotsoffibrations}
  Suppose that $p\geq 5$ and that $v\in\tN(\sX)_0$ is an isotropic
  vector with $v.e=0$, where $e=(0,0,1)$. If $v$ is not a multiple of
  $e$, then there exists an elliptic fibration $X\to\P^1$ with smooth
  fiber $E$ and an integer $s$ such that the image of the vector
  $(0,E,s)\in p\tN(\sX)^*\subset\tN(\sX)$ in $\tN(\sX)_0$ is a scalar
  multiple of $v$.
\end{Lemma}
\begin{proof}
  Consider any lift of $v$ to an isotropic vector
  $(0,l,s)\in p\tN(\sX)^*\subset\tN(\sX)$. We will modify this vector
  while ensuring that its image in $\tN(\sX)_0$ remains a non-zero
  scalar multiple of $v$. We claim that we may assume that $l$ is
  primitive. Suppose that $l=nl'$ for some primitive $l'$ and some
  $n$, which is necessarily invertible modulo $p$. We pick an integer
  $s'$ such that $ns'\equiv s$ modulo $p$. Replacing $(0,l,s)$ with
  $(0,l',s')$ achieves our goal. Next, after possibly multiplying by
  $-1$, we may arrange so that $l$ is in the closure of the positive
  cone, and hence is effective. By Proposition 3.3.10 of \cite{Huy16},
  we find a sequence of $(-2)$-curves $C_i$ such that
  \[
    s_{C_n}\circ\dots\circ s_{C_1}(l)
  \]
  is linearly equivalent to a smooth fiber $E\subset X$ of an elliptic
  fibration $f\colon X\to\P^1$. Because the reflections $s_{C_i}$ act
  trivially on $\tN(X)_0$, the Mukai vector $(0,E,s)$ also lifts $v$.
\end{proof}

\begin{Theorem}\label{thm:modularinterpretation}
  If $p\geq 5$, then every twistor family is a moduli space of twisted
  sheaves on (an open subset of) a connected component of a universal
  family of $\m_p$-gerbes on some supersingular K3 surface. If
  $p\geq 3$, the same is true if $v.e\neq0$ (with the notation of
  Theorem \ref{thm:constructionoftwistorlines}).
\end{Theorem}
\begin{proof}
  It is likely possible to give a direct proof for this result along
  the lines of Theorem \ref{thm:constructionoftwistorlines}. Allowing
  ourselves to use the twisted crystalline Torelli theorem, the proof
  is greatly simplified. Choose a lift $x\in\tLam$ as in Theorem
  \ref{thm:constructionoftwistorlines} (in the positive rank case), or
  with the properties of Lemma \ref{lem:lotsoffibrations} (in the
  Artin-Tate case). The remainder of the proof follows Theorem
  \ref{thm:constructionoftwistorlines}, except that we deduce the
  existence of the desired isomorphism from the Torelli theorem.
\end{proof}

\bibliographystyle{plain}
\bibliography{biblio}

\begin{thebibliography}{10}

\bibitem{MR2786662}
Dan Abramovich, Martin Olsson, and Angelo Vistoli.
\newblock Twisted stable maps to tame {A}rtin stacks.
\newblock {\em J. Algebraic Geom.}, 20(3):399--477, 2011.

\bibitem{Artin74}
M.~Artin.
\newblock Supersingular {$K3$} surfaces.
\newblock {\em Ann. Sci. \'Ecole Norm. Sup. (4)}, 7:543--567 (1975), 1974.

\bibitem{MR0419450}
M.~Artin and J.~S. Milne.
\newblock Duality in the flat cohomology of curves.
\newblock {\em Invent. Math.}, 35:111--129, 1976.

\bibitem{MR0384804}
Pierre Berthelot.
\newblock {\em Cohomologie cristalline des sch\'emas de caract\'eristique
  {$p>0$}}.
\newblock Lecture Notes in Mathematics, Vol. 407. Springer-Verlag, Berlin-New
  York, 1974.

\bibitem{BO78}
Pierre Berthelot and Arthur Ogus.
\newblock {\em Notes on crystalline cohomology}.
\newblock Princeton University Press, Princeton, N.J.; University of Tokyo
  Press, Tokyo, 1978.

\bibitem{Ill96}
Jos\'e Bertin, Jean-Pierre Demailly, Luc Illusie, and Chris Peters.
\newblock {\em Introduction \`a la th\'eorie de {H}odge}, volume~3 of {\em
  Panoramas et Synth\`eses [Panoramas and Syntheses]}.
\newblock Soci\'et\'e Math\'ematique de France, Paris, 1996.

\bibitem{BS15}
Bhargav Bhatt and Peter Scholze.
\newblock Projectivity of the {W}itt vector affine {G}rassmannian.
\newblock {\em Invent. Math.}, 209(2):329--423, 2017.

\bibitem{MR1629929}
Tom Bridgeland.
\newblock Fourier-{M}ukai transforms for elliptic surfaces.
\newblock {\em J. Reine Angew. Math.}, 498:115--133, 1998.

\bibitem{MR1651025}
Tom Bridgeland.
\newblock Equivalences of triangulated categories and {F}ourier-{M}ukai
  transforms.
\newblock {\em Bull. London Math. Soc.}, 31(1):25--34, 1999.

\bibitem{Charles13}
Fran\c{c}ois Charles.
\newblock The {T}ate conjecture for {K}3 surfaces over finite fields.
\newblock {\em Invent. Math.}, 194(1):119--145, 2013.

\bibitem{Charles14}
Fran\c{c}ois Charles.
\newblock Birational boundedness for holomorphic symplectic varieties,
  {Z}arhin's trick for {$K3$} surfaces, and the {T}ate conjecture.
\newblock {\em Ann. of Math. (2)}, 184(2):487--526, 2016.

\bibitem{dejong-gabber}
A.~J. de~Jong.
\newblock A result of {G}abber.
\newblock 2003.

\bibitem{MR638598}
P.~Deligne.
\newblock Rel\`evement des surfaces {$K3$}\ en caract\'eristique nulle.
\newblock In {\em Algebraic surfaces ({O}rsay, 1976--78)}, volume 868 of {\em
  Lecture Notes in Math.}, pages 58--79. Springer, Berlin-New York, 1981.
\newblock Prepared for publication by Luc Illusie.

\bibitem{EGA4.3}
J.~Dieudonn\'e and A.~Grothendieck.
\newblock {\'E}l\'ements de g\'eom\'etrie alg\'ebrique. {IV}. \'etude locale
  des sch\'emas et des morphismes de sch\'emas. {III}.
\newblock {\em Inst. Hautes \'Etudes Sci. Publ. Math.}, (28):255, 1966.

\bibitem{Ekedahl}
T.~{Ekedahl}.
\newblock {On non-liftable Calabi-Yau threefolds}.
\newblock {\em ArXiv Mathematics e-prints}, June 2003.

\bibitem{MR0344253}
Jean Giraud.
\newblock {\em Cohomologie non ab\'elienne}.
\newblock Springer-Verlag, Berlin-New York, 1971.
\newblock Die Grundlehren der mathematischen Wissenschaften, Band 179.

\bibitem{MR0244271}
Alexander Grothendieck.
\newblock Le groupe de {B}rauer. {III}. {E}xemples et compl\'ements.
\newblock In {\em Dix expos\'es sur la cohomologie des sch\'emas}, volume~3 of
  {\em Adv. Stud. Pure Math.}, pages 88--188. North-Holland, Amsterdam, 1968.

\bibitem{MR1859189}
Larry~C. Grove.
\newblock {\em Classical groups and geometric algebra}, volume~39 of {\em
  Graduate Studies in Mathematics}.
\newblock American Mathematical Society, Providence, RI, 2002.

\bibitem{Huy06}
D.~Huybrechts.
\newblock {\em Fourier-{M}ukai transforms in algebraic geometry}.
\newblock Oxford Mathematical Monographs. The Clarendon Press, Oxford
  University Press, Oxford, 2006.

\bibitem{MR3051203}
Daniel Huybrechts.
\newblock A global {T}orelli theorem for hyperk\"ahler manifolds [after {M}.
  {V}erbitsky].
\newblock {\em Ast\'erisque}, (348):Exp. No. 1040, x, 375--403, 2012.
\newblock S\'eminaire Bourbaki: Vol. 2010/2011. Expos\'es 1027--1042.

\bibitem{Huy16}
Daniel Huybrechts.
\newblock {\em Lectures on {K}3 surfaces}, volume 158 of {\em Cambridge Studies
  in Advanced Mathematics}.
\newblock Cambridge University Press, Cambridge, 2016.

\bibitem{HL10}
Daniel Huybrechts and Manfred Lehn.
\newblock {\em The geometry of moduli spaces of sheaves}.
\newblock Cambridge Mathematical Library. Cambridge University Press,
  Cambridge, second edition, 2010.

\bibitem{HS04}
Daniel Huybrechts and Paolo Stellari.
\newblock Equivalences of twisted {$K3$} surfaces.
\newblock {\em Math. Ann.}, 332(4):901--936, 2005.

\bibitem{MR0491704}
Luc Illusie, editor.
\newblock {\em Cohomologie {$l$}-adique et fonctions {$L$}}.
\newblock Lecture Notes in Mathematics, Vol. 589. Springer-Verlag, Berlin-New
  York, 1977.
\newblock S\'eminaire de G\'eometrie Alg\'ebrique du Bois-Marie 1965--1966 (SGA
  5), Edit\'e par Luc Illusie.

\bibitem{Ill79}
Luc Illusie.
\newblock Complexe de de\thinspace {R}ham-{W}itt et cohomologie cristalline.
\newblock {\em Ann. Sci. \'Ecole Norm. Sup. (4)}, 12(4):501--661, 1979.

\bibitem{MR0291177}
Nicholas~M. Katz.
\newblock Nilpotent connections and the monodromy theorem: {A}pplications of a
  result of {T}urrittin.
\newblock {\em Inst. Hautes \'Etudes Sci. Publ. Math.}, (39):175--232, 1970.

\bibitem{MR0337959}
Nicholas~M. Katz.
\newblock Algebraic solutions of differential equations ({$p$}-curvature and
  the {H}odge filtration).
\newblock {\em Invent. Math.}, 18:1--118, 1972.

\bibitem{MR0237510}
Nicholas~M. Katz and Tadao Oda.
\newblock On the differentiation of de {R}ham cohomology classes with respect
  to parameters.
\newblock {\em J. Math. Kyoto Univ.}, 8:199--213, 1968.

\bibitem{KMP16}
Wansu Kim and Keerthi Madapusi~Pera.
\newblock 2-adic integral canonical models.
\newblock {\em Forum Math. Sigma}, 4:e28, 34, 2016.

\bibitem{Lie15}
M.~{Lieblich}.
\newblock {Rational curves in the moduli of supersingular K3 surfaces}.
\newblock {\em ArXiv e-prints}, July 2015.

\bibitem{Lie14}
Max Lieblich.
\newblock On the unirationality of supersingular {K}3 surfaces.

\bibitem{Lie04}
Max Lieblich.
\newblock Moduli of twisted sheaves.
\newblock {\em Duke Math. J.}, 138(1):23--118, 2007.

\bibitem{MR2388554}
Max Lieblich.
\newblock Twisted sheaves and the period-index problem.
\newblock {\em Compos. Math.}, 144(1):1--31, 2008.

\bibitem{MR2579390}
Max Lieblich.
\newblock Compactified moduli of projective bundles.
\newblock {\em Algebra Number Theory}, 3(6):653--695, 2009.

\bibitem{LO15}
Max Lieblich and Martin Olsson.
\newblock Fourier-{M}ukai partners of {K}3 surfaces in positive characteristic.
\newblock {\em Ann. Sci. \'Ec. Norm. Sup\'er. (4)}, 48(5):1001--1033, 2015.

\bibitem{Liedtke15}
Christian Liedtke.
\newblock Supersingular {K}3 surfaces are unirational.
\newblock {\em Invent. Math.}, 200(3):979--1014, 2015.

\bibitem{MP15}
Keerthi Madapusi~Pera.
\newblock The {T}ate conjecture for {K}3 surfaces in odd characteristic.
\newblock {\em Invent. Math.}, 201(2):625--668, 2015.

\bibitem{Mau14}
Davesh Maulik.
\newblock Supersingular {K}3 surfaces for large primes.
\newblock {\em Duke Math. J.}, 163(13):2357--2425, 2014.
\newblock With an appendix by Andrew Snowden.

\bibitem{Milne76}
J.~S. Milne.
\newblock Duality in the flat cohomology of a surface.
\newblock {\em Ann. Sci. \'Ecole Norm. Sup. (4)}, 9(2):171--201, 1976.

\bibitem{MR1432058}
Laurent Moret-Bailly.
\newblock Un probl\`eme de descente.
\newblock {\em Bull. Soc. Math. France}, 124(4):559--585, 1996.

\bibitem{Muk84}
Shigeru Mukai.
\newblock Symplectic structure of the moduli space of sheaves on an abelian or
  {$K3$}\ surface.
\newblock {\em Invent. Math.}, 77(1):101--116, 1984.

\bibitem{MR2514037}
David Mumford.
\newblock {\em Abelian varieties}, volume~5 of {\em Tata Institute of
  Fundamental Research Studies in Mathematics}.
\newblock Published for the Tata Institute of Fundamental Research, Bombay; by
  Hindustan Book Agency, New Delhi, 2008.
\newblock With appendices by C. P. Ramanujam and Yuri Manin, Corrected reprint
  of the second (1974) edition.

\bibitem{Navas10}
Hermes Jackson~Martinez Navas.
\newblock Fourier–mukai transform for twisted sheaves.
\newblock {\em (Dissertation, Bonn)}, 2010.

\bibitem{Ogus77}
Arthur Ogus.
\newblock {$F$}-crystals and {G}riffiths transversality.
\newblock In {\em Proceedings of the {I}nternational {S}ymposium on {A}lgebraic
  {G}eometry ({K}yoto {U}niv., {K}yoto, 1977)}, pages 15--44. Kinokuniya Book
  Store, Tokyo, 1978.

\bibitem{Ogus78}
Arthur Ogus.
\newblock Supersingular {$K3$}\ crystals.
\newblock In {\em Journ\'ees de {G}\'eom\'etrie {A}lg\'ebrique de {R}ennes
  ({R}ennes, 1978), {V}ol. {II}}, volume~64 of {\em Ast\'erisque}, pages 3--86.
  Soc. Math. France, Paris, 1979.

\bibitem{Ogus83}
Arthur Ogus.
\newblock A crystalline {T}orelli theorem for supersingular {$K3$}\ surfaces.
\newblock In {\em Arithmetic and geometry, {V}ol. {II}}, volume~36 of {\em
  Progr. Math.}, pages 361--394. Birkh\"auser Boston, Boston, MA, 1983.

\bibitem{MR2312554}
Martin Olsson.
\newblock Sheaves on {A}rtin stacks.
\newblock {\em J. Reine Angew. Math.}, 603:55--112, 2007.

\bibitem{RS79}
A.~N. Rudakov and I.~R. Shafarevich.
\newblock Supersingular {$K3$}\ surfaces over fields of characteristic {$2$}.
\newblock {\em Izv. Akad. Nauk SSSR Ser. Mat.}, 42(4):848--869, 1978.

\bibitem{RS83}
A.~N. Rudakov and I.~R. Shafarevich.
\newblock Surfaces of type {$K3$}\ over fields of finite characteristic.
\newblock In {\em Current problems in mathematics, {V}ol. 18}, pages 115--207.
  Akad. Nauk SSSR, Vsesoyuz. Inst. Nauchn. i Tekhn. Informatsii, Moscow, 1981.

\bibitem{RTS83}
A.~N. Rudakov, T.~Zink, and I.~R. Shafarevich.
\newblock The influence of height on degenerations of algebraic surfaces of
  type {K}3.
\newblock {\em Mathematics of the USSR-Izvestiya}, 20(1):119, 1983.

\bibitem{MR0429918}
Tetsuji Shioda.
\newblock On elliptic modular surfaces.
\newblock {\em J. Math. Soc. Japan}, 24:20--59, 1972.

\bibitem{MR0435084}
Tetsuji Shioda.
\newblock Algebraic cycles on certain {$K3$} surfaces in characteristic{$p$}.
\newblock pages 357--364, 1975.

\bibitem{Shioda77}
Tetsuji Shioda.
\newblock Some results on unirationality of algebraic surfaces.
\newblock {\em Mathematische Annalen}, 230:153--168, 1977.

\bibitem{stacks-project}
The {Stacks Project Authors}.
\newblock Stacks {P}roject.
\newblock \url{http://stacks.math.columbia.edu}, 2017.

\bibitem{MR1610977}
John Tate.
\newblock On the conjectures of {B}irch and {S}winnerton-{D}yer and a geometric
  analog.
\newblock In {\em S\'eminaire {B}ourbaki, {V}ol.\ 9}, pages Exp.\ No.\ 306,
  415--440. Soc. Math. France, Paris, 1995.

\bibitem{MR0225778}
John~T. Tate.
\newblock Algebraic cycles and poles of zeta functions.
\newblock In {\em Arithmetical {A}lgebraic {G}eometry ({P}roc. {C}onf. {P}urdue
  {U}niv., 1963)}, pages 93--110. Harper \& Row, New York, 1965.

\bibitem{MR3161308}
Misha Verbitsky.
\newblock Mapping class group and a global {T}orelli theorem for hyperk\"ahler
  manifolds.
\newblock {\em Duke Math. J.}, 162(15):2929--2986, 2013.
\newblock Appendix A by Eyal Markman.

\bibitem{Yos06}
K\=ota Yoshioka.
\newblock Moduli spaces of twisted sheaves on a projective variety.
\newblock In {\em Moduli spaces and arithmetic geometry}, volume~45 of {\em
  Adv. Stud. Pure Math.}, pages 1--30. Math. Soc. Japan, Tokyo, 2006.

\end{thebibliography}
\end{document}